\documentclass[11pt,paper=a4]{article}
 \usepackage{indentfirst, latexsym,bm}
 \usepackage{amsmath}
 \usepackage{pifont}
 \usepackage{amsfonts}
 \usepackage{mathrsfs}
 \usepackage{array}
 \usepackage{multirow} 
 \usepackage{graphicx}
 \usepackage{subfigure}
 \usepackage{picinpar}
 \usepackage{setspace}
 \usepackage[textsize=footnotesize]{todonotes}
 \RequirePackage[colorlinks,citecolor=blue,urlcolor=blue,linkcolor=blue]{hyperref}
 \usepackage[top=2cm,bottom=2cm, outer=2cm, inner=2cm]{geometry}

 \RequirePackage[numbers]{natbib}
 \usepackage{enumitem} 
 
 \usepackage{booktabs}
 \usepackage{authblk}
 
 \usepackage{threeparttable}
 \usepackage{mathrsfs,amsfonts,amsmath}
 
 \usepackage{color}
 
 \makeatletter
 \def\namedlabel#1#2{\begingroup
 	#2%
 	\def\@currentlabel{#2}%
 	\phantomsection\label{#1}\endgroup
 }
 \makeatother

 \numberwithin{figure}{section}
 
 \newcommand\email[1]{\rm\href{mailto:#1}{ \nolinkurl{#1}}}

 \newcommand*{\doublebar}[1]{\bar{\bar{#1}}}

 \renewcommand{\theequation}{\arabic{section}.\arabic{equation}}
 
 \newtheorem{theorem}{Theorem}[section]
 \newtheorem{definition}[theorem]{Definition}
 \newtheorem{lemma}[theorem]{Lemma}
 \newtheorem{corollary}[theorem]{Corollary}
 \newtheorem{proposition}[theorem]{Proposition}
 \newtheorem{remark}[theorem]{Remark}
 \newtheorem{condition}[theorem]{Condition}
 \newtheorem{example}{Example}[section]

 \def\blemma{\begin{lemma}}\def\elemma{\end{lemma}}
 \def\bproposition{\begin{proposition}}\def\eproposition{\end{proposition}}
 \def\ttheorem{\begin{theorem}}\def\etheorem{\end{theorem}}
 \def\bcorollary{\begin{corollary}}\def\ecorollary{\end{corollary}}
 \def\bremark{\begin{remark}}\def\eremark{\end{remark}}
 \def\bcondition{\begin{condition}}\def\econdition{\end{condition}}

 \def\benumerate{\begin{enumerate}}\def\eenumerate{\end{enumerate}}
 \def\bitemize{\begin{itemize}}\def\eitemize{\end{itemize}}

 \def\beqlb{\begin{eqnarray}}\def\eeqlb{\end{eqnarray}}
 \def\beqnn{\begin{eqnarray*}}\def\eeqnn{\end{eqnarray*}}
 \def\ar{\!\!\!&}

 \def\proof{\noindent{\it Proof.~~}}\def\qed{\hfill$\Box$\medskip}

\begin{document}
  
 \title{\bf\Large Stochastic Volterra Equations for Local Times of Spectrally Positive L\'evy Processes with Gaussian Components}
 	
 \author{Wei Xu\footnote{School of Mathematics and Statistics, Beijing Institute of Technology, China. Email: xuwei.math@gmail.com}  
 }

 \maketitle

 \begin{abstract}
 Following our previous work \cite{Xu2024b}, this paper continues to investigate the evolution dynamics of local times of spectrally positive L\'evy processes with Gaussian components in the spatial direction. 
 We prove that conditioned on the finiteness of the first time at which the local time at zero exceeds a given value, local times at positive line are equal in law to the unique solution of a stochastic Volterra equation driven by a Gaussian white noise and two Poisson random measures with convolution kernel given in terms of the scale function. 
 Also, we obtain several equivalent stochastic equations by using the potential theoretic techniques and prove the strong existence and uniqueness by using the generalized Yamada-Watanabe theorems. 
 
  \medskip
   
 Armed with the stochastic Volterra representation, we then establish a comparison principle for the local times of spectrally positive L\'evy processes with various drifts or stopped when local times at zero exceed different given values, which proposes a stochastic flow enjoying the branching property. 
 And also, we explore some novel properties of local times in the spatial direction including uniform moment estimates, $(1/2-\varepsilon)$-H\"older continuity and maximal inequality.   
 By using the method of duality, we provide an exponential-affine representation of the Laplace functional in terms of the unique non-negative solution of a path-dependent nonlinear Volterra equation associated with the Laplace exponent of L\'evy process. 
 This gives another perspective on the evolution dynamics of local times in the spatial direction. 
 
 \bigskip
  
  
  \medskip
  

 \end{abstract}
 

  \section{Introduction and main results}
 \label{Sec.Introduction}
 \setcounter{equation}{0}
 
 
 As an important research branch in mathematics and probability, local times of L\'evy processes have been deeply studied and also widely applied in various fields including random trees, queueing systems, ruin theory and so on; see \cite{BarndorffMikoschResnick2001,Kyprianou2014} for a review. 
 More specifically, Aldous first introduced the celebrated Brownian continuum random tree as scaling limit of discrete random trees in \cite{Aldous1991} and reconstructed it in \cite{Aldous1993} from Brownian excursions whose local times are identical to the tree-width. 
 Further, Le Gall and Le Jan \cite{LeGallLeJan1998b,LeGallLeJan1998a} coded the genealogy of general continuous-state branching processes via the exploration process defined by the local times at zero of the reflected processes of the time-reversed processes associated to a spectrally positive L\'evy process with Laplace exponent being identically the branching mechanism. 
 Later,  the well-known L\'evy continuum random tree was constructed by Duquesne and Le Gall \cite{DuquesneLeGall2002} from the spectrally positive L\'evy process. 
 Furthermore, Duquesne proposed in a pioneer work \cite{Duquesne2006} that each spectrally positive L\'evy process stopped at hitting zero can be seen as the contour process of a compact real tree.  
 Therefore, more in-depth study of the local times of spectrally positive L\'evy processes would contribute to understand the corresponding real trees.

 Since the well-known  \textsl{Ray-Knight theorems} were proved independently by Ray \cite{Ray1963} and Knight \cite{Knight1963} to connect Brownian local times with various Bessel processes, understanding the distributions and the inner structure of local times of symmetric Markov processes have motivated an abundance of amazing and wonderful work during the past decades, including Dynkin's  and Eisenbaum's isomorphism theorems (see \cite{Dynkin1984,EisenbaumKaspi1995}), the alternate Ray-Knight theorems (see \cite{EisenbaumKaspi1994,EisenbaumKaspiMarcusRosenShi2000,SabotTarres2016}), sufficient and necessary conditions for the Markov property (see \cite{EisenbaumKaspi1993}) and the joint continuity (see \cite{MarcusRosen1992a,MarcusRosen1992b}). We refer to the monograph \cite{MarcusRosen2006} for a survey on the local times of symmetric Markov processes.  
 Recently, the isomorphism theorems also have been established for non-symmetric Markov processes in \cite{FitzsimmonsRosen2014,LeJanMarcusRosen2015} by using Markovian loop soups and permanental processes; see also \cite{LeJanMarcusRosen2017} for details. 
 In particular, more explicit results have been obtained for the local times of L\'evy processes with the symmetry being unnecessary, e.g., the joint continuity  \cite{BarlowHawkes1985}, Hilbert transform \cite{FitzsimmonsGetoor1992}, H\"older regularity \cite{Barlow1988}, laws of the iterated logarithm \cite{Bertoin1995} and so on.
 
 In recent years, the inner structure and evolution dynamic of local times of spectrally positive L\'evy processes in the spatial direction have attracted considerable attention because of their close connection to random trees and random maps. 
 Besides the Brownian case, the discontinuous trajectories of L\'evy processes result in the fail of Markovianity and the intractability to their local times; see \cite{EisenbaumKaspi1993}. 
 Lambert \cite{Lambert2010} first used an excursion of compound Poisson process with unit negative drift and positive jumps to code a binary, splitting tree, and identified that local times of the excursion equal in distribution to a homogeneous, binary  \textsl{Crump-Mode-Jagers branching process}. 
 Later, Lambert and Uribe Bravo \cite{LambertUribeBravo2018} considered the totally ordered measured trees and coded them by the spectrally positive L\'evy processes. Their results also illustrate the link between local times of spectrally positive L\'evy processes and the corresponding tree-width processes. 
 Meanwhile, as the preceding compound Poisson processes converges weakly to a spectrally positive L\'evy process after a suitable time-spatial scaling,  Lambert and Simatos \cite{LambertSimatos2015} proved  the finite-dimensional convergence of their local times. 
 By marking each jumps of a driftless spectrally positive stable process by a random path, Forman et al. \cite{FormanPalRizzoloWinkel2018} established a locally uniform approximation for their local times and proved the finiteness of all moments of the H\"older coefficient. 
 
 In this work, we are mainly interested in the macro-evolution mechanism and sample path properties of local times of spectrally positive L\'evy processes in the spatial direction, also known as \textsl{the second Ray-Knight theorems}. 
 In our previous work \cite{Xu2024b}, local times of spectrally positive stable processes stopped when the local time at zero exceeds a given value were intuitively described by the unique solution of a stochastic Volterra equation driven by Poisson random measures, which also gave a detailed account of the perturbations caused by each jump on local times. 
 In a recent work, Rivero and Contreras \cite{RiveroContreras2024} also generalized the first and second Ray-Knight theorems for spectrally negative L\'evy processes and establish several Poisson representations of their local times on the whole real line by using excursion theory. 
 As a continuation of \cite{Xu2024b}, in this paper we first establish a stochastic Volterra representation as well as a comparison principle for the local times of spectrally positive L\'evy processes with Gaussian components, and then investigate their distribution properties as well as  trajectory regularity.

 \subsection{Main results}
 
 Consider a spectrally positive L\'evy process $\xi:=\{ \xi(t):t\geq 0 \}$ defined on a filtrated probability space $(\Omega,\mathscr{F},\mathscr{F}_t,\mathbf{P})$ with a non-positive drift and a Gaussian component.
 It can be fully characterized by the \textsl{Laplace exponent} $\varPhi$ that is of the form
 \beqlb\label{LaplaceExponent}
 \varPhi (\lambda):= \log \mathbf{E}\big[\exp\{ -\lambda \xi(1) \} \big] =  b\cdot\lambda + c\cdot \lambda^2  +   \int_0^\infty \big(e^{-\lambda y}-1+\lambda y \big) \, \nu (dy) ,\quad \lambda \geq 0 ,
 \eeqlb  
 for some \textsl{L\'evy triplet} $(b,c,\nu)$ with $b\geq 0$, $c>0$ and $\nu(dy)$ being the L\'evy measure on $(0,\infty)$ satisfying
 \beqlb\label{LevyTriplet} 
 \int_0^\infty \big(y\wedge y^2 \big) \, \nu(dy)= \int_0^\infty \big(1 \wedge (2y) \big) \,\bar\nu(y)\,dy <\infty , 
 \eeqlb 
 where $\bar\nu(y):= \nu([y,\infty))$ is the \textsl{tail-function} of $\nu(dy)$. 
 The function $\varPhi$ is zero at zero and increases strictly to infinity at infinity.
 It is strictly convex and infinitely differentiable on $(0,\infty)$ with 
 $ \varPhi^{\,\prime} (0)=-\mathbf{E}\big[\xi(1)\big] = b$.  
 The process $\xi$ drifts to $-\infty$ or is recurrent according as $b>0$ or $=0$.
  
 For every $t\geq 0$, let $\mu_{\xi,t}(dy)$ be the \textsl{occupation measure} of $\xi$ on the time interval $[0,t]$ given for every non-negative and measurable function $f$ on $\mathbb{R}$ by
 \beqlb\label{OccupationMeasure}
 \int_0^t f\big(\xi(s)\big) \, ds \overset{\rm a.s.}= \int_\mathbb{R} f(x)\, \mu_{\xi,t}(dx).
 \eeqlb
 Theorem 1 in \cite[p.126]{Bertoin1996} induces that the random measure $\mu_{\xi,t}(dx)$ is absolutely continuous almost surely with respect to the Lebesgue measure and the density, denoted by $\{L_\xi(x,t):x\in \mathbb{R} \}$, is square integrable.
 The identity \eqref{OccupationMeasure} turns into the well-known  \textsl{occupation density formula}
 \beqlb\label{OccupationDensityF}
 \int_0^t f \big( \xi(r) \big)\, dr \overset{\rm a.s.}= \int_\mathbb{R} f(x) L_\xi(x,t) \, dx, \quad t\geq 0.
 \eeqlb
 
 The existence of Gaussian component yields that $\varPhi(\lambda)=O(\lambda^2) $ as $\lambda \to\infty$ and 
 $ \int_1^\infty 1/\varPhi(\lambda)\, d\lambda <\infty$, 
 which along with Lemma 2.2 in \cite{LambertSimatos2015} induces that the two-parameter process  
 \medskip\smallskip\\ \medskip  \smallskip 
 \centerline{\bf $L_\xi:=\{L_\xi(x,t):x\in \mathbb{R}, t\geq 0 \}$ is jointly continuous almost surely.} 
 \noindent In particular, the process $\{L_\xi(0,t):t\geq 0\}$ is continuous and non-decreasing with $L_\xi(0,\infty) \in (0,\infty]$ a.s.
 This allows us to define the \textsl{inverse local time} $\tau_\xi^L:=\{\tau_\xi^L(\zeta):\zeta\geq 0\}$ at level $0$ by
 $\tau_\xi^L(\zeta)=\infty$ if $\zeta> L_\xi(0,\infty)$ and
 \beqlb \label{eqn.1009}
 \tau_\xi^L(\zeta):= \inf\big\{ s\geq 0:  L_\xi(0,s)\geq \zeta \big\}, \quad \mbox{if } \zeta\in \big[0,L_\xi(0,\infty)\big].
 \eeqlb
 For any $\zeta > 0$, the transience and recurrence of $\xi$ induce that $L_\xi(0,\infty) =\infty$ and $ \tau_\xi^L(\zeta) <\infty$ a.s.  when $b=0$ or $L_\xi(0,\infty) <\infty$ a.s. and $\mathbf{P}(\tau_\xi^L(\zeta)=\infty)>0$  when $b>0$.
 In this work, we are mainly interested in the evolution dynamic of local times in the spatial direction conditioned on the finiteness of the first time at which the local time at zero exceeds a given value, i.e., 
 \medskip\smallskip\\ \medskip   
 \centerline{\bf $L^\xi_\zeta$ is the process $\big\{L_\xi(x,\tau_\xi^L(\zeta)):x\geq 0  \big\}$ conditioned on $\tau_\xi^L(\zeta)<\infty$. }\smallskip


 Our first main result establishes a stochastic Volterra equation driven by a Gaussian white noise and two Poisson random measures for $L^\xi_\zeta$ with integrands being in terms of the \textsl{scale function} $W:=\{W(x):x\in \mathbb{R} \}$ associated to  $\varPhi$. 
 The function $W$ is a non-negative function that is identically zero on $  (-\infty,0)$ and characterized on $[0,\infty)$ as a continuous, strictly increasing function with Laplace transform 
 \beqlb\label{ScaleFunction}
 \int_0^\infty e^{-\lambda x}W(x)\, dx= \frac{1}{\mathit{\Phi}(\lambda )}, 
 \quad \lambda >0;
 \eeqlb 
 see \cite{ChanKyprianouSavov2011,KuznetsovKyprianouRivero2013,Kyprianou2014} for details.
 It is continuous on $\mathbb{R}$ and twice continuously differentiable on $(0,\infty)$ with
 \beqlb \label{eqn.201}
 W^{\prime}(0)=\frac{1}{c} ,
 \qquad
 \sup_{x\in \mathbb{R}} W'(x)\leq \frac{1}{c}
 \qquad \mbox{and}\qquad
  W(\infty)=\begin{cases} 
  	\infty, & \mbox{if $b=0$}; \vspace{3pt}\\
  1/b, & \mbox{if $b>0$}.
  \end{cases}   
 \eeqlb 
 This uniform upper bound comes from Corollary~\ref{Coro.UpperBound}.  
  The next theorem will be proved in Section~\ref{Sec.SVETightness}.

 \begin{theorem}\label{Main.Thm01}
 	For each $\zeta\geq 0$, the process $L^\xi_\zeta$ equals in distribution to the unique non-negative continuous solution of the stochastic Volterra equation
 	\beqlb\label{MainThm.SVE}
 	X_\zeta(t) \ar=\ar \zeta\cdot c \cdot W'(t)  
 	+  \int_0^\zeta\int_0^\infty \big( W(t)-W(t-y) \big) \, N_0(dz, dy)   + \int_0^t \int_0^{X_\zeta(s)}   W'(t-s) \, B_c(ds,dz) \cr
 	\ar\ar + \int_0^t \int_0^{X_\zeta(s)} \int_0^\infty  \big( W(t-s)-W(t-s-y) \big) \, \widetilde{N}_\nu(ds,dz,dy), \quad t\geq 0,
 	\eeqlb
 	where $N_0(dz, dy)$ is a  Poisson random measure on $(0,\infty)^2$ with intensity $\bar\nu(y)\,dz\,dy$, $B_c(ds,dz)$ is a Gaussian white noise on $(0,\infty)^2$ with intensity $2c\cdot ds\, dz$ and $\widetilde{N}_\nu(ds,dz,dy)$ is a compensated Poisson random measure on $(0,\infty)^3$ with intensity $ds\, dz\,\nu(dy)$. 
 \end{theorem}

 Definitions of $W$ and $(N_0, B_c, \widetilde{N}_\nu)$ tell that the stochastic Volterra equation \eqref{MainThm.SVE} is fully determined by the \textsl{characteristic vector}  $ (\zeta\, ;\, b,c,\nu)$.
 The two upper bounds in \eqref{LevyTriplet} and \eqref{eqn.201} induce that
 \beqlb  \label{eqn.100}
 \int_0^\infty \big(W(t)-W(t-y) \big)\,\bar{\nu}(y)\,dy + \int_0^\infty  \big( W(t)-W(t-y) \big)^2\, \nu(dy) 
 \ar\leq\ar C \int_0^\infty  (t\wedge y)^2 \, \nu(dy) <\infty,
 \eeqlb
 uniformly in $t $ on compacts.  
 This along with \eqref{eqn.201} allows us to consider the three stochastic integrals in \eqref{MainThm.SVE} as It\^o's integrals that have been deeply explored in \cite{IkedaWatanabe1989,Walsh1986}. 
 Hence it is natural to formulate the definition of solutions to \eqref{MainThm.SVE} by extending those of stochastic differential equations in \cite[Chapter IV.1]{IkedaWatanabe1989} and \cite{Kurtz2014}. 
 More specifically, by a \textsl{solution} of \eqref{MainThm.SVE}, we mean a  process $X_\zeta \in C(\mathbb{R}_+;\mathbb{R}_+)$ defined on a filtrated probability space on which three mutually independent driving noises $(N_0, B_c, \widetilde{N}_\nu)$ as in Theorem~\ref{Main.Thm01} are defined such that \eqref{MainThm.SVE} holds almost surely.  
 The \textsl{uniqueness} is said to hold for \eqref{MainThm.SVE} if any two solutions equal in distribution.
 Also, we recall the quadruple $(X_{\zeta},N_0, B_c, \widetilde{N}_\nu)$ a solution of \eqref{MainThm.SVE} to emphasize the particular role of driving noises. 
 
 \begin{remark}\label{Remark.CIR}
 In the Brownian case, i.e., $\nu(\mathbb{R}_+)=0$, 
 we have $W'(x)= \frac{1}{c} \cdot \exp \{ -\frac{b}{c}\cdot x  \}$. 
 By the equality $e^{ -\frac{b}{c}\cdot (t-s)  }=1-\int_s^t e^{ -\frac{b}{c}\cdot (r -s) }\,dr$ and the stochastic Fubini's theorem; see \eqref{eqn.SFT01}, we can write \eqref{MainThm.SVE}   as
 \beqlb\label{eqn.CIR}
  X_\zeta(t) 
  \ar=\ar \zeta -\int_0^t \frac{b}{c} \cdot  X_\zeta(s)\, ds + \int_0^t \int_0^{X_\zeta(s)}  \frac{1}{c}  \, B_c(ds,dz), \quad t\geq 0,
 \eeqlb
 The strong non-negative continuous solution uniquely exists and is a branching diffusion starting with $\zeta$ ancestors. In particular, when $b=0$ and $c=1/2$, it is a square of $0$-dimensional Bessel process and \eqref{eqn.CIR} is identical to (1.4) in \cite{AidekonHuShi2024} with $\mu=1$. 
 \end{remark}
 
 Usually, the last two stochastic Volterra integrals in \eqref{MainThm.SVE} are not local martingales but have mean zero. 
 Taking expectations on both sides of \eqref{MainThm.SVE} and then using Fubini's theorem, 
 \beqlb\label{eqn.103}
 \mathbf{E}\big[  L^\xi_\zeta(x)   \big]
 = \mathbf{E}\big[ X_\zeta(x) \big] 
 \ar=\ar \zeta\cdot c \cdot W'(x) +  \zeta \cdot \int_0^\infty \big( W(x)-W(x-y) \big) \bar\nu(y)\,dy\cr
 \ar=\ar \zeta\cdot c \cdot W'(x)  +  \zeta \cdot  \int_0^\infty W'(x-y)   \doublebar\nu(y) \,dy,
 \eeqlb
 where $\doublebar{\nu}(x):=\int_{[x,\infty)} \bar{\nu}(y)\,dy$ is the integrated tail-function of $\bar\nu$.
 Furthermore, applying the identity in  Lemma~\ref{Lemma.Identity01} to \eqref{eqn.103} induces that 
 \beqlb\label{eqn.10301}
 \mathbf{E}\big[  L^\xi_\zeta(x) \big]
 = \mathbf{E}\big[ X_\zeta(x) \big] 
 \ar=\ar \zeta\cdot \big(1- b\cdot W(x)\big) \in [0,\zeta],
 \eeqlb
 which is identically $\zeta$ when $b=0$ or decreases strictly to $0$ as $x\to\infty$ when $b>0$. 
 The finite height of the process $\xi$ stopped at $\tau^\xi_\zeta$ induces that $L^\xi_\zeta(x) \to 0$ a.s. as $x\to\infty$ and the point $0$ is a cemetery state. 
 Analogue to the criticality for branching processes, we say $L^\xi_\zeta$ is \textsl{critical} if $b=0$ or \textsl{subcritical} if $b>0$.

 \begin{remark}\label{Remark.EquivalentSVE}
 With the help of \eqref{eqn.103} and \eqref{eqn.10301}, we merge the first term on the right side of \eqref{MainThm.SVE} with the compensator of the second term to obtain the following alternate representation 
  \beqlb\label{MainThm.SVENew}
  X_\zeta(t) \ar=\ar  \zeta\cdot \big(1- b \cdot W(x)\big) + \int_0^\zeta\int_0^\infty \big( W(t)-W(t-y) \big) \, \widetilde{N}_0(dz, dy)  + \int_0^t \int_0^{X_\zeta(s)}   W'(t-s) \, B_c(ds,dz) \cr
  \ar\ar + \int_0^t \int_0^{X_\zeta(s)} \int_0^\infty  \big( W(t-s)-W(t-s-y) \big) \, \widetilde{N}_\nu(ds,dz,dy), \quad t\geq 0,
  \eeqlb 
 with $\widetilde{N}_0(dz, dy):=N_0(dz, dy)- \bar\nu(y)\,dz\,dy$. 
 Additionally, 
 by Theorem~7.1 in \cite[p.84]{IkedaWatanabe1989}  one can find a Brownian motion $B$ on an extension of the original probability space such that almost surely
 	\beqnn
 	\int_0^t \int_0^{X_\zeta(s)}   W'(t-s) \, B_c(ds,dz) = \int_0^t  W'(t-s) \sqrt{2c\cdot X_\zeta(s)} \, dB(s) , \quad t\geq 0. 
 	\eeqnn
 \end{remark}
 
 
 \begin{figure} \label{Figure.01}
 	\centering
 	\includegraphics[scale=0.4]{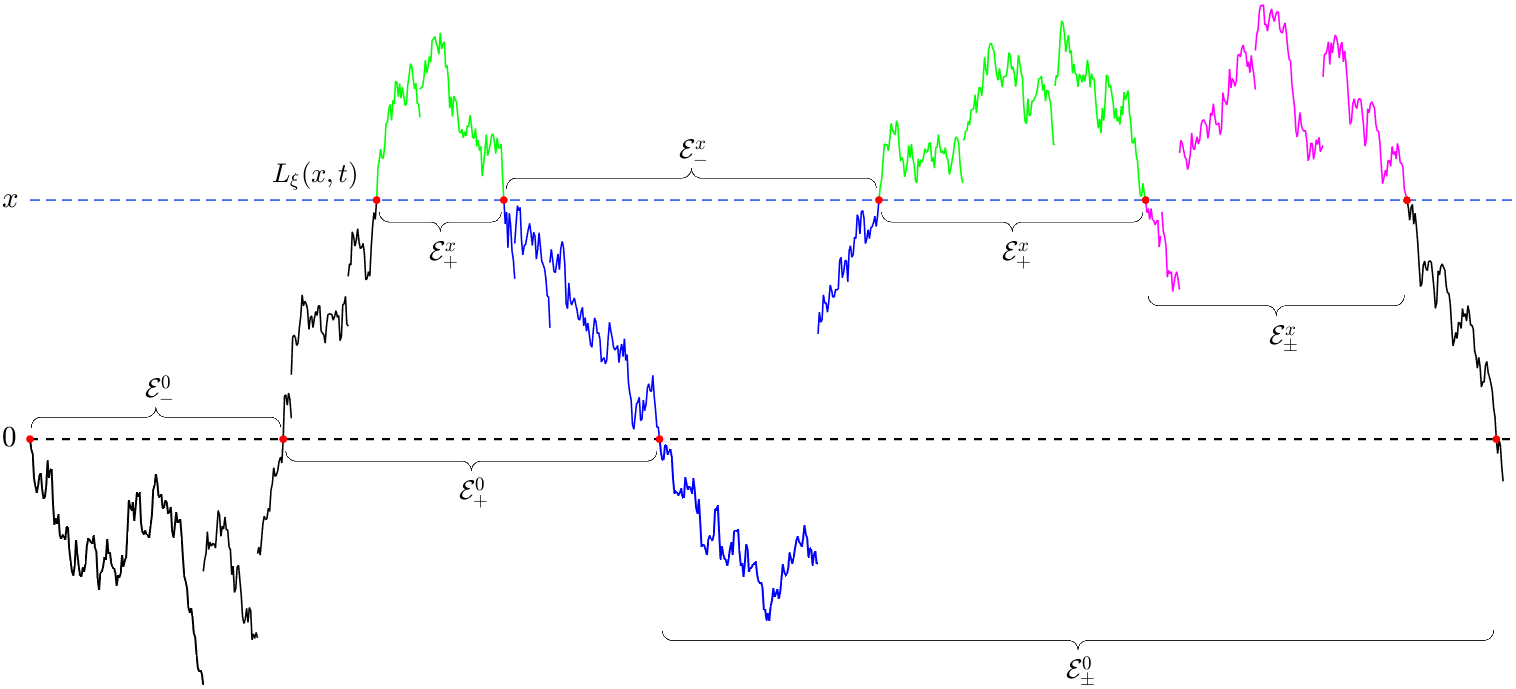}
 	\caption{A sample path of spectrally positive L\'evy processes containing three typical kinds of excursions away from level $x$: (i) $\mathcal{E}_+^x$  excursions that are completely above $x$ (the green trajectories); (ii) $\mathcal{E}_-^x$  excursions that are completely blow $x$ (the blue trajectory); (iii) $\mathcal{E}_\pm^x$ excursions that move below $x$ until jumping into $(x,\infty)$ and then stay above $x$ up to hitting $x$ (the magenta trajectory).}
 \end{figure}
 By using results in Chapter IV in \cite{Bertoin1996} and the spatial homogeneity, the local times of $\xi$ at any level $x$ (up to a constant multiplicative factor) also can be defined by approximations involving the numbers of excursions away from $x$ with intervals of certain types. 
 Because of the existence of Gaussian component and the lake of negative jumps, excursions away from $x$ can classified into three distinguishable types: (i) $\mathcal{E}_+^x$ consists of all excursions that are completely above $x$; (ii) $\mathcal{E}_-^x$ consists of all excursions that are completely blow $x$; (iii) $\mathcal{E}_{\pm}^x$ consists of all excursions that move below $x$ until jumping into $(x,\infty)$ and then stay above $x$ up to hitting $x$; see Figure~\ref{Figure.01}. 
 We refer to \cite{RiveroContreras2024} for more details. 
 Roughly speaking, the first two types come from the Gaussian component and the third type results from the positive jumps. 
 Our proofs further elaborate on their various contributions to the local times at different levels and  the equation \eqref{MainThm.SVE}. 
 More precisely, excursions $\mathcal{E}_-^0$ contribute only to the local times of negative levels and hence have no connection to \eqref{MainThm.SVE}. 
 In contrast, excursions $\mathcal{E}_+^0$ contribute only to the local times of positive levels as well as three of the four terms on the right side of \eqref{MainThm.SVE} in different ways, i.e., the first term on the right side of  \eqref{MainThm.SVE} describe their average contribution to the local times at different levels; meanwhile, perturbations of their sub-excursions in $\mathcal{E}_-^x\cup\mathcal{E}_+^x$  and $\mathcal{E}_\pm^x$ respectively compose part of the last two stochastic Volterra integrals in \eqref{MainThm.SVE}. 
 Different from the previous two types, excursions $\mathcal{E}_\pm^0$ contribute to the local times of both positive and negative levels.  
 Their average contribution at positive levels along with the randomness of their overshoots composes the first stochastic Volterra integral in \eqref{MainThm.SVE}. 
 Also, perturbations of their sub-excursions in $\mathcal{E}_-^x\cup\mathcal{E}_+^x$  and $\mathcal{E}_\pm^x$ compose the rest part of the last two stochastic Volterra integrals in \eqref{MainThm.SVE}, respectively.

 \begin{remark}\label{Remark.104}
 In view of \eqref{eqn.201} and \eqref{eqn.100} the first two terms on the right side of \eqref{MainThm.SVE} converge respectively to $\zeta$ and $0$ as $t\to0+$, which indicates that excursions  $\mathcal{E}_\pm^0$ make no contribution to the local times at level $0$. This along with the spatial homogeneity asserts that local times at any level $x$ are fully determined by excursions  $\mathcal{E}_-^x\cup\mathcal{E}_+^x$. 
 \end{remark}

 The monotonicity of $L_\xi$ in the time variable induces that when $b=0$, both $\tau_\xi^L(\zeta)$ and $L^\xi_\zeta$ increase almost surely in $\zeta$. 
 However, this property fails when $b>0$, since the definition of $L^\xi_\zeta$ replies on the conditional law $\mathbf{P}(\cdot\,|\, \tau_\xi^L(\zeta)<\infty)$. 
 On the other hand, the large negative drift will pull $\xi$ back into negative line quickly. 
 Hence, it is natural to conjecture that local times $L_\xi$ would decrease as $b$ increases. 
 As the second main result, these observations motivate us to establish in the next theorem a comparison principle for local times of spectrally positive L\'evy processes with different drifts or stopped when local times at zero exceed different given values. 
 The proof can be found in Section~\ref{Sec.SVETightness}.
  
  
 \begin{theorem} \label{Thm.Comparison}
 For  $b_1,b_2,\zeta_1, \zeta_2\geq 0$ and $i=1,2$, let $\xi_i$ be a spectrally positive L\'evy process with L\'evy triplet $(b_i,c,\nu)$. 
 There exists a filtrated probability space endowed with two processes $X_{\zeta_1}^{b_1},X_{\zeta_2}^{b_2} \in C(\mathbb{R}_+;\mathbb{R}_+)$ and three driving noises $(N_0,B_c,\widetilde{N}_\nu)$ defined as in Theorem~\ref{Main.Thm01} such that the following two claims hold.  
 \begin{enumerate}
 \item[(1)] For $i=1,2$, the process $L^{\xi_i}_{\zeta_i}$ equals in distribution to $ X_{\zeta_i}^{b_i}$ and the quadruple $(X_{\zeta_i}^{b_i},N_0,B_c, \widetilde{N}_\nu)$ is the unique solution of \eqref{MainThm.SVE} with characteristic vector $(\zeta_i\,;\, b_i,c,\nu)$.
 		
 \item[(2)] If $b_1\geq b_2$ and $\zeta_1\leq \zeta_2$, we have  $\mathbf{P}\big( X_{\zeta_1}^{b_1}(t)\leq X_{\zeta_2}^{b_2}(t) ,t\geq 0 \big)=1$. 
 \end{enumerate} 
 \end{theorem}
 
 Comparison principles have been widely established for classic deterministic/stochastic differential equations (see \cite{Hartman2002,IkedaWatanabe1989,Situ2005}) and ordinary Volterra equations (see \cite{Brunner2004,GripenbergLondenStaffans1990}). 
 However, by comparison, much less is known for stochastic Volterra equations, e.g., a comparison principle was first obtained in \cite{Tudor1989} with only Volterra drift and then generalized to the general case in \cite{CanadasFriesen2024,FerreyraSundar2000} under some regular conditions on the kernel and coefficients. 
 To the best of our knowledge, all results in the existing literature do not apply to our stochastic equation \eqref{MainThm.SVE}. 
 
 Compared to the comparison principles in the aforementioned references for differential equations  with explicit, manifest initial states and drifts, the equation \eqref{MainThm.SVE} not only has a random and time-dependent initial state, it also seems to have no drift. 
 This makes our comparison principle puzzling and incomprehensible. 
 In order to make it more understandable in intuition and integrate it into the classic theory, in the next theorem we provide an equivalent representation of \eqref{MainThm.SVE} by separating the impact of drift $b$ on the local times from that of diffusion and jumps. 
 Consider the function  
 \beqlb\label{Phi0}
 \varPhi_0(\lambda):=\varPhi(\lambda)-b\cdot \lambda,
 \quad \lambda\geq 0,
 \eeqlb
 which is the Laplace exponent of a driftless spectrally positive L\'evy process with triplet $(0,c,\nu)$. 
 Denote by $W_0$ and $W'_0$ the scale function associated to $\varPhi_0$ and its right-derivative respectively. 
 
  \begin{theorem} \label{MainThm.EquavilentRep}
 	The stochastic Volterra equation \eqref{MainThm.SVE} is equivalent to 
 	\beqlb \label{Main.SVE02} 
 	X_\zeta(t)\ar=\ar \zeta \cdot c\cdot W'_0(t)  +  \int_0^\zeta\int_0^\infty  \big( W_0(t)-W_0(t-y) \big) \, N_0(dz, dy) \cr
 	\ar\ar - \int_0^t b\cdot W'_0(t-s) X_\zeta(s)ds  + \int_0^t \int_0^{X_{\zeta}(s)}   W'_0(t-s) \, B_c(ds,dz) \cr
 	\ar\ar   
 	 + \int_0^t \int_0^{X_\zeta(s)} \int_0^\infty \big( W_0(t-s)-W_0(t-s-y) \big) \, \widetilde{N}_\nu(ds,dz,dy) .\quad  t\geq 0.
 	\eeqlb
 \end{theorem}
 
 Similarly as in Remark~\ref{Remark.EquivalentSVE}, one can use the identity in Lemma~\ref{Lemma.Identity01} with $b=0$ to combine the first term on right side of \eqref{Main.SVE02} together with the compensator of the second term and then write \eqref{Main.SVE02} as 
 \beqlb \label{Main.SVE0201} 
  X_\zeta(t)\ar=\ar  \zeta   +  \int_0^\zeta\int_0^\infty \big( W_0(t)-W_0(t-y) \big) \, \widetilde{N}_0(dz, dy)    \cr
  \ar\ar -\int_0^t 	b \cdot W'_0(t-s) X_\zeta(s)ds + \int_0^t \int_0^{X_\zeta(s)}   W'_0(t-s) \, B_c(ds,dz) \cr
  \ar\ar + \int_0^t \int_0^{X_\zeta(s)} \int_0^\infty  \big( W_0(t-s)-W_0(t-s-y) \big)\, \widetilde{N}_\nu(ds,dz,dy),  \quad t\geq 0.
 \eeqlb
 Recall that the scale function $W_0$ is independent of the drift parameter $b$. 
 The time-independent initial state and non-positive drift make it natural to expect comparison principles for \eqref{Main.SVE0201} and then \eqref{MainThm.SVE}.

 Analogous to the stochastic dominance between  random variables, we say the process $Y_1\in D(\mathbb{R}_+;\mathbb{R})$ is \textsl{absolutely dominant over} $Y_2\in D(\mathbb{R}_+;\mathbb{R})$, denoted as $Y_1\geq Y_2$, if $\mathbf{P}(Y_1(t)\geq Y_2(t), t\geq 0)=1$. 
 Moreover, $Y_1$ is said to be \textsl{first-order stochastically dominant over}  $Y_2 $, denoted as $Y_1  \succeq Y_2$,  if
 \beqnn
 \mathbf{P}\Big(  \int_0^\infty Y_1(s)\mu(ds)\geq x\Big)  \geq \mathbf{P}\Big(  \int_0^\infty Y_2(s)\mu(ds)\geq x\Big) ,\quad x\in \mathbb{R},
 \eeqnn
 for any measure $\mu(ds)$ on $\mathbb{R}_+$. 
 A family of processes $\{Y_\zeta\}_{\zeta\geq 0}$ is said to be \textsl{absolutely increasing} or \textsl{stochastically increasing} [resp. \textsl{absolutely decreasing} or \textsl{stochastically decreasing}] if $Y_{\zeta_2} \geq  Y_{\zeta_1}$ or $Y_{\zeta_2} \succeq Y_{\zeta_1}$ [resp. $Y_{\zeta_1} \geq  Y_{\zeta_2}$ or $Y_{\zeta_1} \succeq Y_{\zeta_2}$] when $\zeta_2\geq \zeta_1$.
 
 Our comparison principle claims that the local time processes $\{ L^\xi_\zeta \}_{\zeta\geq 0}$ are absolutely increasing if $b=0$ or stochastically increasing if $b>0$. 
 This indicates that a class of stochastic flows are hidden in the local times $L_\xi$. 
 Recently, Contreras and Xu \cite{ContrerasXu2024} reconstructed them and investigated their (conditionally) asymptotic properties by using excursion theory and fluctuation theory of L\'evy processes. 
 In the past decades, various stochastic flows have been derived from L\'evy processes and their related models. For instance, Ray-Knight theorems for (reflected) Brownian motions introduce the flows of squared Bessel processes, Their pathwise construction was recently given by A\"id\'ekon et al \cite{AidekonHuShi2025}, in which a flow called \textsl{Jacobi flow} was also defined by strong solutions of stochastic differential equations that were constructed via a perturbed reflecting Brownian motion. 
 For a general spectrally positive L\'evy process, 
 Ray-Knight theorems in \cite{DuquesneLeGall2002} states that local times of the reflected processes of its time-reversed processes at $0$ equal in law to the flows of continuous-state branching processes, which were reconstructed in \cite{DawsonLi2012} as strong solutions of stochastic differential equations driven by a Gaussian white noise and a Poisson random measure.

 Without the existence and uniqueness of strong solution, one cannot construct the stochastic flows in $L_\xi$ by using \eqref{MainThm.SVE} as in the aforementioned literature. 
 Here the \textsl{strong solution} of \eqref{MainThm.SVE} is defined as a solution $(X_\zeta, N_0, B_c, \widetilde{N}_\nu)$ with $X_\zeta$ being equal almost surely to a Borel measurable function of $(N_0, B_c, \widetilde{N}_\nu)$. 
 Although the scale function $W$ is twice continuously differentiable
 on $(0,\infty)$, the second derivative $W''$ may be singular around $0$, i.e.  Theorem~3.1 in \cite{BehmeOechslerSchilling2022} shows that 
 \beqnn
 W''(0+)= \frac{b-\doublebar{\nu}(0+)}{c^2}, 
 \eeqnn
 which is finite if and only if $\doublebar{\nu}(0+)<\infty$ that is equivalent to $\int_0^\infty y \, \nu(dy)<\infty$. 
 Moreover, note that all coefficients in \eqref{MainThm.SVE} are $1/2$-H\"older continuous. 
 Consequently, we cannot obtain the strong solutions by using the standard Picard iteration and the Euler method as usual. 
 Fortunately, in the case of $\doublebar{\nu}(0+)<\infty$, the function $W''$ is continuous on $\mathbb{R}_+$ and also differentiable on $(0,\infty)$ with derivative $W'''$ enjoying the regularity of $\bar{\nu}$; see Theorem 2 in \cite{ChanKyprianouSavov2011}. 
 This allows us to improve the method developed in \cite{AlfonsiSzulda2024,PromelScheffels2023} to prove the pathwise uniqueness for \eqref{MainThm.SVE} and then establish the strong existence and uniqueness by using the generalized Yamada-Watanabe theorems in \cite{YamadaWatanabe1971}; see the next theorem whose the proof is given in Section~\ref{Sec.PathUniqueness}. 

 \begin{theorem} \label{Thm.StrongUniqueness}
 	If $\doublebar{\nu}(0+)<\infty$, the stochastic Volterra equation \eqref{MainThm.SVE} has a unique strong solution. Moreover, the solution is a semimartingale with the following representation
 	\beqlb\label{eqn.1004}
 	\left\{ \begin{split}
 	X_\zeta(t) 
 	&= \zeta \cdot c\cdot W'(t)  
 	+  \int_0^\zeta\int_0^\infty \big( W(t) -W(t-y)\big) \, N_0(dz, dy) \cr
 	&\quad + \int_0^t \Big( W''(0+) \cdot M(s) +  \int_0^s W'''(s-r)M(r)\,dr \Big) \, ds + \frac{1}{c}\cdot M(t),\cr
 	M(t) &= \int_0^t   \int_0^r \int_0^{X_\zeta(s)} \int_{r-s}^\infty    \widetilde{N}_\nu(ds,dz,dy)\,dr + \int_0^t \int_0^{X_\zeta(s)} \, B_c(ds,dz).
 	\end{split}
 	\right.
 	\eeqlb  
 \end{theorem}

 Armed with the stochastic Volterra equation \eqref{MainThm.SVE}, we are allowed to use the instruments and methods provided by modern probability theory and stochastic analysis to revisit or refine the well-known results for local times of L\'evy processes and also investigate their unknown properties. 
 The main contribution in the second part of this work
 is to illustrate the strength of  \eqref{MainThm.SVE} by using them to study the distribution properties and sample path regularity of local times in the spatial direction. 
 As the first result, in the next theorem we provide a uniform upper bound for all moments of local times of any levels. 
 
 \begin{theorem}\label{Thm.Moment}
  For each $p\geq 1$, there exists a constant $C>0$ such that for any $\zeta,x\geq 0$,
  \beqlb\label{eqn.Moment} 
  \mathbf{E}\Big[\big| L^\xi_\zeta(x) \big|^p  \Big] 
  \leq C\cdot (\zeta \vee \zeta^p) \cdot \big(1+W(x)\big)^{2p-2} . 
  \eeqlb
 \end{theorem}
 
 In the critical case ($b=0$), the uniform upper bound in \eqref{eqn.Moment} can be changed into $ C\cdot (\zeta \vee \zeta^p)\cdot (1+x)^{2p-2}$ since $W(x) \leq x/c$; see \eqref{eqn.201}.
 Differently, when $b>0$, the third claim in \eqref{eqn.201} induces that it can be bounded by $C\cdot (\zeta \vee \zeta^p)$ uniformly  in $x\geq 0$.  
 Moreover, it should be stressed that our uniform upper bound is far from being optimal and can be further improved if some $L^p$-estimates are provided for $W''$, e.g., for a drifted Brownian, in view of Remark~\ref{Remark.CIR} we can replace the upper bound by $c_1 \cdot (\zeta \vee \zeta^p)\cdot \exp\{-c_2\cdot x\}$   for some constants $c_1,c_2>0$.

  The second property we are interested in is the H\"older regularity of the trajectories of $L^\xi_\zeta$. Finding condition for the joint continuity of the local times motivated a large number of magnificent work in the last century; see \cite[Chapter V.5]{Bertoin1996} for a survey. 
  For instance, the joint H\"older continuity of local times of $\alpha$-stable L\'evy processes with $\alpha\in (1,2]$ was first obtained in  \cite{Boylan1964} by discretizing the occupation measures. The optimal H\"older exponent was given in \cite{Barlow1988} and the finiteness of all moments of the H\"older coefficient in the driftless case was proved in the recent work \cite{FormanPalRizzoloWinkel2018} . 
  In the general case, local times in the spatial direction were proved to be H\"older continuous in \cite{BarlowHawkes1985,BassKhoshnevisan2006} if the reciprocal of characteristic function has real part decreasing at infinity as a power function.  
  
  In the next theorem, we prove the H\"older continuity of $L^\xi_\zeta$ by using the inequalities of stochastic integral to \eqref{MainThm.SVE} and then the Kolmogorov continuity theorem. 
  Moreover, we also provide a uniform upper bound for all moments of the H\"older coefficients and the local maximal of $L^\xi_\zeta$ by using the Garsia-Rodemich-Rumsey inequality.
  For $\kappa\in (0,1]$ and $x> 0$,  the \textsl{$\kappa$-H\"older coefficient} of a H\"older continuous function $f$ on $[0,x]$ is defined by
  \beqnn
  \big\|f\big\|_{C^{0,\kappa}_x} :=\sup_{0\leq y<z\leq x} \frac{|f(y)-f(z)|}{|y-z|^\kappa}.
  \eeqnn
  
  \begin{theorem}\label{Thm.Regularity}
  	The process $L^\xi_\zeta$ is locally H\"older continuous with index strictly less than $1/2$. 
  	Moreover,  for each $\kappa \in (0,1/2)$ and $p\geq 1$, there exists a constant $C>0$ such that for any $\zeta, x\geq 0$, 
  	\beqlb\label{Regularity.01}
  	\mathbf{E}\Big[\big\| L^\xi_\zeta \big\|_{C^{0,\kappa}_x}^{p}\Big] \leq  C \cdot (\zeta \vee \zeta^{p}) \cdot (1+x)^{p(3-\kappa) }
  	\quad \mbox{and}\quad
  	\mathbf{E}\Big[ \sup_{z\in[0,x]} \big|L^\xi_\zeta(z)\big|^p \Big] \leq C \cdot (\zeta \vee \zeta^{p}) \cdot (1+x)^{3p} . 
  	\eeqlb
  \end{theorem}
  
  Compared to the $(\frac{\alpha-1}{2} -\varepsilon)$-H\"older continuity in the $\alpha$-stable case; see Theorem~2.9 in \cite{Xu2024b}, the Gaussian component significantly improves the regularity of $L^\xi_\zeta$  even though more perturbations are added. 
  Actually, this is not surprising and inexplicable,  since the real part the reciprocal of characteristic function is $O(|z|^{-2})$ at infinity and hence  $L_\xi$ is $(\frac{1}{2}-\varepsilon)$-H\"older continuous in the spatial variable; see \cite{BarlowHawkes1985,BassKhoshnevisan2006} and \cite[Exercise 4, p.151]{Bertoin1996}. 
  As we mentioned in Remark~\ref{Remark.104}, similarly to Brownian local times, local times at level $x$ are only contributed by excursions $\mathcal{E}^x_+\cup \mathcal{E}^x_-$ and hence should enjoy the same regularity. 
  
  The last property of $L^\xi_\zeta$ we explore in this work is its Laplace functional.  
  In contrast to Brownian local times whose the Laplace transform can be given as an exponential affine function of the initial state and the solution of a Riccati equation by using their branching property, 
  the evolution dynamic of $L^\xi_\zeta$ is much intractable because of the lack of Markovanity.
  Fortunately, an intuitive comparison between the stochastic integrals in \eqref{MainThm.SVE} and \eqref{eqn.CIR} reveals that  $L^\xi_\zeta$ should enjoy the analogous ``affine property''  of Brownian local times. 
  By generalizing the duality method in \cite{Xu2024b},  we establish in the next theorem an explicit representation of Laplace functionals 
  $\mathbf{E}\big[ \exp\{- \int_{[0,x]} L^\xi_\zeta(x-s)\,\mu(ds)\} \big]$ in term of the unique solution of the \textsl{path-dependent nonlinear Volterra equation} 
  \beqlb\label{NonlinearVolterra}
  V_\mu (t) = \int_{[0,t]} W'(t-s) \,\mu(ds)-\int_0^t \boldsymbol{\mathcal{R}}\circ V_\mu (s) W'(t-s) \, ds,\quad t\geq 0 , 
  \eeqlb
  where $\mu(ds)$ is a $\sigma$-finite measure  on $\mathbb{R}_+$ and $\boldsymbol{\mathcal{R}}$ is a nonlinear operator acting on a locally bounded
  function $f$ by
  \beqlb\label{OperatorR}
  \boldsymbol{\mathcal{R}}\circ f(t):= c \cdot \big(f(t)\big )^2 + \int_0^\infty
  \bigg(\exp\Big\{- \int_{(t-y)^+}^{t} f(r)\, dr\Big\}-1  +\int_{(t-y)^+}^{t} f(r)\, dr\bigg) \nu ( dy ) . 
  \eeqlb
  
  \begin{theorem} \label{Thm.LaplaceF}
  	For each $\sigma$-finite measure $\mu(ds)$ on $\mathbb{R}_+$, 
  	the nonlinear Volterra equation \eqref{NonlinearVolterra} has a unique global solution $V_\mu \in D(\mathbb{R}_+;\mathbb{R}_+) $ and the Laplace functional of $L^\xi_\zeta$ admits the representation
  	\beqlb\label{LaplaceFun01}
  	\mathbf{E}\bigg[ \exp\Big\{- \int_{[0,x]} L^\xi_\zeta(x-s)\,\mu(ds) \Big\} \bigg]
  	= \exp\big\{  -\zeta \cdot \boldsymbol{\mathcal{F}}\circ V_\mu(x)   \big\}, 
  	\quad x\geq 0,
  	\eeqlb
  	where $\boldsymbol{\mathcal{F}}$ is a nonlinear operator acting on a locally bounded
  	function $f$ by
  	\beqlb\label{LaplaceFun02}
  	\boldsymbol{\mathcal{F}}\circ f(t)
  	\ar=\ar c\cdot f(t) + \int_0^\infty \bigg(1- \exp\Big\{-\int_{(t-y)^+}^{t} f(r)\,dr\Big\}   \bigg)\, \bar\nu(y)\,dy ,\quad t\geq 0 . 
  	\eeqlb

  \end{theorem}
  
  It is worth mentioning that Rivero and Contreras \cite{RiveroContreras2024} used the excursion theory to provide a different representation for the Laplace functional of $L^\xi_\zeta$ as follows
  \beqnn
  \lefteqn{	\mathbf{E}\bigg[ \exp\Big\{- \int_{[0,x]} L^\xi_\zeta(x-s)\,\mu(ds) \Big\} \bigg]
  }\ar\ar\cr
  \ar=\ar  \exp\bigg\{-\zeta\cdot \int_{D_0}\Big( 1-\exp\Big\{-\int_{[0,x]} \ell_{\mathbf{e}}(x-s)\,\mu(ds) \Big\} \Big)\boldsymbol{n}_0(d\mathbf{e})   \cr
  \ar\ar \quad -\zeta \cdot \int_0^\infty \Big( 1-\exp\Big\{-\int_0^y dr\int_{D_0}\Big( 1-\exp\Big\{-\int_{[r,x]} \ell_{\mathbf{e}}(x-s)\,\mu(ds) \Big\}   \Big) \underline{\boldsymbol{n}}(d\mathbf{e}) \Big\}\Big) \bar{\nu}(y)dy \bigg\},
  \eeqnn
  where $\ell_{\mathbf{e}}(s)$ is the total local times of excursion $\mathbf{e}$ at level $s$, $\boldsymbol{n}_0$ and $\underline{\boldsymbol{n}}$ are the measures of excursions  away from $0$ for $\xi$ and its reflected process at infimum respectively. 
  Intuitively, the two terms in the Laplace exponent should be identical to the corresponding terms in \eqref{LaplaceFun01}-\eqref{LaplaceFun02}. 
  Due to our limited acknowledge of excursion theory, the equivalence between these two representations remains to be addressed.

  Finally, we end this work with some discussion of the total local times in the subcritical case ($b>0$). 
  Proposition~4 in \cite[p.130]{Bertoin1996} tells that the process $\tau_\xi^L$ is a subordinator killed at an independent exponential time and its Laplace transform is of the form
  \beqlb\label{eqn.500}
  \mathbf{E}\big[\exp\big\{-\lambda\cdot \tau_\xi^L(\zeta)\big\}\big] = \exp\big\{-\zeta/u^\lambda(0)\big\},\quad \lambda  ,\, \zeta\geq 0,
  \eeqlb
  where $u^\lambda:=\{u^\lambda(y):y\in\mathbb{R}\}$ is the density of the \textit{$\lambda$-resolvent kernel} of $\xi$. 
  This induces that $$\mathbf{P}\big(L_\xi(0,\infty)\geq \zeta\big)=1-\mathbf{P}\big(\tau_\xi^L(\zeta) =\infty\big)= \exp\big\{-\zeta/u^0(0)\big\}$$
  and hence $L_\xi(0,\infty)$ is exponentially distributed with mean $u^0(0)$. 
  Define the process
  \beqlb \label{eqn.Linf}
  L^\xi_\infty:= \big\{L_\xi(x,\infty):x\geq 0  \big\}\overset{\rm a.s.}=
  \big\{L_\xi \big(x,\tau_\xi^L(L_\xi(0,\infty)) \big):x\geq 0 \big\} . 
  \eeqlb
  All preceding results and their proofs can be generalized to $L^\xi_\infty$ with the constant $\zeta$ replaced by an exponentially distributed random variable $\varrho$ with mean $u^0(0)$; see the next corollary as one      example.  
  
  \begin{corollary}\label{Main.Thm.b}
  	When $b>0$, the process $L_\infty^\xi$ equals in distribution to the unique continuous solution of (\ref{MainThm.SVE}) with $\zeta = \varrho$.
  	For any $\sigma$-finite measure $\mu(ds)$ on $\mathbb{R}_+$, it has Laplace functional 
  	\beqlb\label{LaplaceFun011}
  	\mathbf{E}\bigg[ \exp\Big\{- \int_{[0,x]} L_\infty^\xi(x-s)\,\mu(ds) \Big\} \bigg]
  	= \Big( 1+u^0(0)\cdot \boldsymbol{\mathcal{F}}\circ V_\mu(x) \Big)^{-1} , 
  	\quad x\geq 0 . 
  	\eeqlb

  \end{corollary}



 \smallskip
 {\it \textbf{Organization of this paper.}} \ \ 
 In Section~\ref{Sec.PreMainResult}, we provide some auxiliary results for scale functions and stochastic Volterra integrals. 
 In Section~\ref{Sec.CompoundPoisson}, we recall two stochastic Volterra representations for the local times of compound Poisson processes and then show that after a suitable scaling, they behave like those of the desired L\'evy process. 
 In Section~\ref{Sec.SVETightness}, we first establish the stochastic Volterra representation for the process $L^\xi_\zeta$ by proving that the stochastic Volterra equation solved by the local times of compound Poisson processes converges weakly to \eqref{MainThm.SVE}. 
 Thereafter, we prove the comparison principle with the help of our weak convergence arguments. 
 The existence and uniqueness of strong solution are proved in Section~\ref{Sec.PathUniqueness}. 
 Section~\ref{Sec.MomentEstimate} is devoted to establish the moment estimates given in Theorem~\ref{Thm.Moment} and the general equivalent representation of  \eqref{MainThm.SVE} that will induces Theorem~\ref{MainThm.EquavilentRep} as a by-product. 
 Theorem~\ref{Thm.Regularity} is proved in Section~\ref{Sec.EquivalentSVE}. 
 The well-posedness of nonliear Volterra equation \eqref{NonlinearVolterra} and the affine representation of the Laplace functional of $L^\xi_\zeta$ are established in Section~\ref{Sec.NonlinearVE}, and also, the uniqueness in law of solution to \eqref{MainThm.SVE} is proved as a corollary. 
 In Appendix~\ref{AppendixHMartinagle}, we recall some basic theory of stochastic integrals driven by semimartingales indexed by a Banach space.

  \smallskip
  
 {\it \textbf{Notation.}} \ \ 
 Let $\mathbb{N}$ be the space of all natural numbers including $0$. 
 For any $x\in\mathbb{R}$, let $x^+:=x\vee 0$, $x^-:=x\wedge 0$ and $[x]$ be the integer part of $x$. 
 Let $\overset{\rm f.d.d.}\longrightarrow$, $\overset{\rm u.c.}\to$,  $\overset{\rm a.s.}\to$, $\overset{\rm d}\to$  and  $\overset{\rm p}\to$ be the convergence in the sense of finite dimensional distributions, the uniform convergence on compacts, almost sure convergence,  convergence in distribution and convergence in probability respectively.
 We also use $\overset{\rm a.s.}=$, $\overset{\rm d}=$ and $\overset{\rm p}=$ to denote almost sure equality, equality in distribution and equality in probability respectively.

 For a Banach space $\mathbb{V}$ with a norm $\|\cdot\|_\mathbb{V}$, let $D([0,\infty),\mathbb{V})$ be the space of all c\`adl\`ag $\mathbb{V}$-valued functions endowed with the Skorokhod topology and $C([0,\infty),\mathbb{V})$ the space of all continuous $\mathbb{V}$-valued functions endowed with the uniform topology.
 For any $\mathcal{T}\subset[0,\infty)$ and $p\in(0,\infty]$, let $L^p(\mathcal{T}; \mathbb{V})$ be the space of all $\mathbb{V}$-valued measurable functions $f$ on $\mathcal{T}$ satisfying   
 $$
 \big\|f\big\|_{L^p_\mathcal{T}}^p:=  \int_\mathcal{T} \big\|f(x)\big\|_{\mathbb{V}}^p dx  <\infty.
 $$
 We also write $\|f\|_{L^p_T}$ for $\|f\|_{L^p_{[0,T]}}$ and $\|f\|_{L^{^p}}$ for $\|f\|_{L^p_\infty}$. 
 Let $L^p_{\rm loc}(\mathbb{R}_+; \mathbb{V}) :=  \cap_{T\geq 0}L^p([0,T];\mathbb{V})$. 
 For two functions $f,g$ and a $\sigma$-finite measure $\mu$ on $\mathbb{R}_+$,
 the two convolutions $f*g$ and $f*d\mu$ are defined by
 \beqnn
 f*g(x):= \int_0^x f(x-y)g(y)\, dy
 \quad \mbox{and}\quad 
  f*d\mu(x):= \int_{[0,x]} f(x-y)\,\mu(dy) 
 \quad x\geq 0.
 \eeqnn


 
 We use $C$ to denote a positive constant whose value might change from line to line.

  \section{Preliminaries}
 \label{Sec.PreMainResult}
 \setcounter{equation}{0}
  
  In this section, we first provide some auxiliary results for scale functions and then give some moment estimates and stochastic Fubini's theorems for stochastic Volterra integrals with respect to Gaussian white noise and Poisson random measure.

 \subsection{Scale functions}
  
%
 
 Associated to the functions $\bar{\nu}$ and $\doublebar{\nu}$ we recall two alternate representations of the Laplace exponent $\varPhi $
 \beqlb\label{EquivalentLE}
 \varPhi(\lambda)\ar=\ar  b\cdot\lambda + c\cdot \lambda^2  +  \lambda \cdot \int_0^\infty \big(1-e^{-\lambda x}  \big)  \bar\nu (x) \, dx \cr
 \ar=\ar b\cdot\lambda + c\cdot \lambda^2  +  \lambda^2 \cdot \int_0^\infty  e^{-\lambda x}   \, \doublebar{\nu}(x)\, dx ,
 \eeqlb
 which appear frequently in \cite{Bertoin1996,Kyprianou2014} and can be obtained by using Fubini's theorem to \eqref{LaplaceExponent}. 
 Using integration by parts to \eqref{ScaleFunction} gives that
 \beqlb\label{eqn.202}
  \int_0^\infty e^{-\lambda x}W'(x)\,dx 
  = \int_0^\infty \lambda e^{-\lambda x}W(x)\,dx =\frac{\lambda}{\varPhi(\lambda)} , 
  \quad \lambda >0,
 \eeqlb
 which also holds for all 
 $\lambda \in \mathbb{C}_+:= \{ x+\mathtt{i}z\,:\, x\geq 0 \mbox{ and } z\in\mathbb{R}\}$ 
 if $b>0$ and all 
 $\lambda \in \mathbb{C}_+ \setminus \{ 0 \}$ if $b=0$; see the proof of Theorem~2.1 in \cite{KuznetsovKyprianouRivero2013}. 
 We first provide an identity that has been used to  derive \eqref{eqn.10301}. 

 \begin{lemma}\label{Lemma.Identity01}
 For all $x\geq 0$, we have $  b\cdot W(x) + c\cdot W'(x)  + W'*\doublebar{\nu}(x) =1$.
 \end{lemma}
 \proof 
 To prove the identity,  we define a $\sigma$-finite measure $\mu_b(dy) := c\cdot \delta_0(dy) + \big(b+ \doublebar{\nu}(y)\big)\,dy$ on $\mathbb{R}_+$ 
 with $\delta_0(dy)$ being the Dirac measure at point $0$. 
 It has Laplace transform 
 \beqnn
  \int_{\mathbb{R}_+} e^{-\lambda y}\, \mu_b(dy) = \frac{b}{\lambda } + c + \int_0^\infty e^{-\lambda x}\doublebar{\nu}(x)\, dx = \frac{\varPhi(\lambda)}{\lambda^2}.
 \eeqnn
 Here the second equality follows from \eqref{EquivalentLE}.  
 This along with \eqref{eqn.202} induces that for $\lambda >0$, 
 \beqnn
  \int_{\mathbb{R}_+} e^{-\lambda x}W'*d\mu_b(x)\,dx= \int_{\mathbb{R}_+} e^{-\lambda x}W'(x)\,dx\cdot \int_{\mathbb{R}_+} e^{-\lambda y}\, \mu_b(dy) = \frac{1}{\lambda}, 
 \eeqnn
 which identifies that $W'*d\mu_b(x)=1$ for almost every $x\geq 0$. 
 Additionally, note that
 \beqnn
 W'*d\mu_b(x) = b\cdot W(x) + c\cdot W'(x)  + W'*\doublebar{\nu}(x),
 \eeqnn
 which tells that $W'*d\mu_b$ is continuous on $\mathbb{R}_+$ and hence $W'*d\mu_b(x)=1$ for all $x\geq 0$. 
 \qed

   \begin{lemma}\label{Lemma.201}
 	For any $\lambda\in \mathbb{R}_+\cap (-b,\infty)$, there exists a constant $C>0$ such that for any $\delta>0$,
 	\beqnn
 	\int_\mathbb{R} \big|e^{-\lambda(x+\delta)} W'(x+\delta)-e^{-\lambda x} W'(x)  \big|^2\,dx 
 	\leq C\cdot \delta. 
 	\eeqnn
 \end{lemma}
 \proof 
 Note that $\varPhi(\mathtt{i}z+\lambda) \sim b\lambda$ as $|z|\to 0$ and is $O(|z|^2)$ as $|z|\to \infty$, there exists a constant $C>0$ such that 
 \beqlb \label{eqn.2002}
 \bigg|\int_\mathbb{R}  e^{-(\mathrm{i}z+\lambda)x} W'(x)\, dx  \bigg|
 = \bigg| \frac{\mathtt{i}z+\lambda}{\varPhi(\mathtt{i}z+\lambda)} \bigg|  \leq C\cdot \Big(  1\wedge \frac{1}{|z|}\Big),\quad z\in\mathbb{R}.
 \eeqlb
 By using the Fourier isometry along with the square integrability of  $e^{-\lambda x} W'(x)$ and then the change of variables, 
 \beqnn
 \lefteqn{\int_\mathbb{R} \Big|e^{-\lambda(x+\delta)} W'(x+\delta)-e^{-\lambda x} W'(x)  \Big|^2\,dx }\quad\ar\ar\cr
 \ar=\ar \int_\mathbb{R} \bigg|\int_\mathbb{R} e^{-\mathtt{i}z x -\lambda(x+\delta) } W'(x+\delta) \,dx  - \int_\mathbb{R} e^{-\mathtt{i}z x -\lambda x  } W'(x) \,dx\bigg|^2\, dz\cr
 \ar=\ar \int_\mathbb{R} \Big|\big(e^{\mathrm{i}z\delta}-1\big)\int_\mathbb{R}  e^{-(\mathrm{i}z+\lambda)x} W'(x)\, dx   \bigg|^2\, dz\cr
 \ar\leq\ar C\int_\mathbb{R}\big((z\delta)^2\wedge 1\big) \Big(  1\wedge \frac{1}{|z|^2}\Big)\, dz\cr
 \ar\leq\ar C\int_{|z|\leq 1 } (z\delta)^2  \, dz
 + C\int_{1<|z|\leq 1/\delta} \delta^2  \, dz +C\int_{|z|> 1/\delta}   \frac{1}{|z|^2} \, dz,
 \eeqnn
 which can be uniformly bounded by $C\cdot\delta$. 
 \qed

 \begin{corollary}\label{Coro.51101}
 	For each $T\geq 0$,  there exist two constants $C>0$ and $n_0\geq 1$ such that for any  $\delta \in (0,1)$,
 	\beqnn
 	 \int_{-\infty}^T \big| W'(x+\delta) -  W'(x)\big|^2\, dx 
 	\leq C\cdot \delta.
 	\eeqnn
 \end{corollary}
 \proof For any $\beta>0$, by the fact that $W'(x)=0$ for all $x<0$ we have 
 \beqnn
  \int_{-\infty}^T \big| W'(x+\delta) -  W'(x)\big|^2\, dx
  \ar=\ar  \int_{-1}^T \big| W'(x+\delta) -  W'(x)\big|^2\, dx\cr
  \ar\leq\ar  e^{2\beta (T+1)} \int_{-1}^T e^{-2\beta (x+\delta)}\big|   W'(x+\delta) -    W'(x)\big|^2\, dx\cr
  \ar\leq\ar e^{2\beta (T+1)} \int_{-1}^T  \big|  e^{-\beta (x+\delta)}  W'(x+\delta) -   e^{-\beta x}  W'(x)  \big|^2\, dx \cr
  \ar\ar + e^{2\beta (T+1)} \int_0^T \big| \big( e^{-\beta \delta}-1\big)e^{-\beta x} W'(x) \big|^2\, dx . 
 \eeqnn
 By Lemma~\ref{Lemma.201}, the first term on the right side can be bounded by $C\cdot \delta $ uniformly in $\delta\in(0,1)$. 
 By \eqref{eqn.201} and the inequality $|e^{-z}-1|\leq z$ for any $z\geq 0$, the second term can be bounded by 
 \beqnn
 e^{2\beta (T+1)}\int_0^T \beta^2 e^{-2\beta x}   \, dx \cdot   \frac{\delta^2}{c^2} 
 \leq C\cdot \delta^2 , 
 \eeqnn
 uniformly in $\delta\in(0,1)$ and hence the desired upper bound holds. 
 \qed  
 
 Recall the function $\varPhi_0$ defined in \eqref{Phi0}.
 Let us consider the function 
 \beqnn
 \varPhi_\beta(\lambda):=\varPhi_0(\lambda)+ \beta\cdot \lambda= \varPhi(\lambda)+ (\beta-b)\cdot \lambda,\quad \lambda \geq 0,
 \eeqnn
 which is the Laplace exponent of a spectrally positive L\'evy process with triplet $(\beta,c,\nu)$. 
 Denote by $W_\beta$ and $W'_\beta$ the scale function associated to $\varPhi_\beta$ and its right-derivative respectively. 
 
 \begin{lemma}\label{Lemma.Identity02}
 	For all $x\geq 0$, we have 
 	$ W'(x) = W'_\beta(x) + (\beta-b)\cdot W'_\beta * W'(x)$. 
 \end{lemma}
 \proof It is obvious $W' = W'_\beta$ if $\beta=b$. When $\beta\neq b$,  for any $\lambda\in \mathbb{R}$ such that $\varPhi_\beta(\lambda),\varPhi(\lambda)>0$, we have
 \beqnn
 \frac{\lambda}{\varPhi_\beta(\lambda)} +(\beta-b)\cdot \frac{ \lambda}{\varPhi_\beta(\lambda)} \cdot \frac{\lambda}{\varPhi(\lambda)} 
 =   \frac{\lambda \cdot \big( \varPhi(\lambda)+(\beta-b)\cdot \lambda\big)}{\varPhi_\beta(\lambda)\varPhi(\lambda)} 
 = \frac{\lambda}{\varPhi(\lambda)} ,
 \eeqnn
 which along with \eqref{eqn.202} and the continuity of $W^{\prime}_\beta $, $ W'$ on $\mathbb{R}_+$ yields the desired identity.
 \qed

 \subsection{Stochastic integrals driven by martingale measures} 
 \label{AppendixFubini} 
 
 
 Let $B(ds,dz)$ be a $(\mathscr{F}_t)$-Gaussian white noise on $(0,\infty)^2$ with intensity $ds\,dz$ and $\widetilde{N}(ds,dy,dz)$ be a compensated $(\mathscr{F}_t)$-Poisson random measure  on $(0,\infty)^3$ with intensity $ds\,\mu(dy)\,dz$ and $\mu(dy)$ being a $\sigma$-finite measure on $(0,\infty)$. 
 For a $(\mathscr{F}_t)$-predictable non-negative process $X$, two measurable functions $f$ on $\mathbb{R}_+$ and $g$ on $\mathbb{R}_+^2$ satisfying that for any $t\geq 0$,
 \beqnn
 \int_0^t |f(s)|^2 ds + \int_0^t \int_0^\infty |g(s,y)|^2 ds\, \mu(dy) <\infty , 
 \eeqnn
 we consider the following two stochastic Volterra integrals
 \beqnn
  \begin{split}
 B(f,t)
 &:=  \int_0^t \int_0^{X(s)} f(t-s)B(ds,dz),\cr   
 \widetilde{N}(g,t)
 &:=  \int_0^t  \int_0^{X(s)}\int_0^\infty g(t-s,y)\widetilde{N}(ds,dz,dy),
 \end{split} 
 \qquad t\geq 0.
 \eeqnn

 \begin{proposition} \label{Prop.BDG}
 	For two constants $p\geq 1$ and $T>0$, assume that
 	\beqnn
 	\sup_{t\in[0,T]}\mathbf{E}\Big[\big|X(t)\big|^{p}\Big]<\infty
 	\quad \mbox{and}\quad 
 	\int_0^T ds \int_0^\infty \big|g(s,y)\big|^{2p}\mu(dy) <\infty.
 	\eeqnn
 	Then there exists a constant $C>0$ depending only on $p$ such that
 	\beqlb \label{BDG0}
 	\sup_{t\in[0,T]}\mathbf{E}\Big[ \big| B(f,t) \big|^{2p}  \Big] 
 	\leq C\cdot \sup_{t\in[0,T]}\mathbf{E} \Big[\big|X(t)\big|^p \Big] \cdot \bigg| \int_0^T   \big|f(s)\big|^{2} ds  \bigg|^{p}
 	\eeqlb
 	and 
 	\beqlb \label{BDG}
 	\sup_{t\in[0,T]}\mathbf{E}\Big[ \big| \widetilde{N}(g,t) \big|^{2p}  \Big] 
 	\ar\leq\ar C  \cdot \sup_{t\in[0,T]}\mathbf{E} \Big[\big|X(t)\big|^p \Big] \cdot \bigg| \int_0^T ds \int_0^\infty\big|g(s,y)\big|^{2} \mu(dy)  \bigg|^{p} \cr
 	\ar\ar + C \cdot \sup_{t\in[0,T]}\mathbf{E}\Big[\big|X(t)\big|\Big]\cdot \int_0^T ds \int_0^\infty\big|g(s,y)\big|^{2p} \mu(dy) .
 	\eeqlb 
 \end{proposition}
 \proof For each $t_0\in[0,T]$, we consider two auxiliary martingales $B_{t_0}(f,\cdot)$ and $\widetilde{N}_{t_0}(g,\cdot)$ defined by
 \beqnn
 \begin{split}
 	B_{t_0}(f,t)
 	&:=  \int_0^t \int_0^{X(s)} f(t_0-s)B(ds,dz),\cr   
 	\widetilde{N}_{t_0}(g,t)
 	&:=  \int_0^t  \int_0^{X(s)}\int_0^\infty g(t_0-s,y)\widetilde{N}(ds,dz,dy),
 \end{split} 
 \qquad t\geq 0.
 \eeqnn
 It is obvious that  $ B (f,t_0)\overset{\rm a.s.}=B_{t_0}(f,t_0)$ and $ \widetilde{N} (g,t_0)\overset{\rm a.s.}=\widetilde{N}_{t_0}(g,t_0)$.
 Hence 
 \beqnn
 \sup_{t\in[0,T]}\mathbf{E}\Big[ \big| B(f,t) \big|^{2p}  \Big]  \ar=\ar
 \sup_{t_0\in[0,T]}\mathbf{E}\Big[ \big|B_{t_0} (f,t_0) \big|^{2p}  \Big] ,\cr
 \sup_{t\in[0,T]}\mathbf{E}\Big[ \big| \widetilde{N}(g,t) \big|^{2p}  \Big]  
 \ar=\ar
 \sup_{t_0\in[0,T]}\mathbf{E}\Big[ \big|\widetilde{N}_{t_0} (g,t_0) \big|^{2p}  \Big] .
 \eeqnn
 By using the Burkholder-Davis-Gundy inequality and then H\"older's inequality as well as the change of variables, there exists a constant $C>0$ that depends only on $p$ such that 
 \beqnn
 \mathbf{E}\Big[ \big|B_{t_0} (f,t_0) \big|^{2p}  \Big]
 \ar\leq\ar   \mathbf{E}\Big[ \sup_{t\in[0,t_0]} \big|B_{t_0} (f,t) \big|^{2p}  \Big] \cr
 \ar\leq\ar C\cdot\mathbf{E}\bigg[\Big|\int_0^{t_0} X(s)\cdot \big|f(t_0-s)\big|^2 ds\Big|^p\bigg]\cr
 \ar\leq\ar \mathbf{E}\bigg[ \Big|\int_0^{t_0}  \big|f(t_0-s)\big|^2 ds\Big|^{p-1}\cdot  \int_0^{t_0} \big|X(s)\big|^p \cdot \big|f(t_0-s)\big|^2 ds \bigg]\cr
 \ar\leq\ar C\cdot \sup_{t\in [0,t_0]} \mathbf{E}\Big[ \big|X(t)\big|^p \Big]  \cdot \Big|\int_0^{t_0}  \big|f(s)\big|^2 ds\Big|^p,
 \eeqnn
 which immediately induces that 
 \beqnn
 \sup_{t\in[0,T]}\mathbf{E}\Big[ \big| B(f,t) \big|^{2p}  \Big]  =\sup_{t_0\in[0,T]}\mathbf{E}\Big[ \big|B_{t_0} (f,t_0) \big|^{2p}  \Big]
 \ar\leq\ar   C\cdot \sup_{t\in [0,T]} \mathbf{E}\Big[ \big|X(t)\big|^p \Big]  \cdot \Big|\int_0^{T}  \big|f(s)\big|^2 ds\Big|^p. 
 \eeqnn
 We now start to prove \eqref{BDG}. 
 By using Theorem~D.1 in \cite{Xu2024b} and then the change of variables to the following last inequality, we have 
 \beqnn
 \mathbf{E}\Big[ \big|\widetilde{N}_{t_0} (g,t_0) \big|^{2p}  \Big] 
 \ar\leq\ar 
 \mathbf{E}\Big[ \sup_{t\in[0,t_0]} \big|\widetilde{N}_{t_0} (g,t) \big|^{2p}  \Big] \cr
 \ar\leq\ar 
 C  \cdot \sup_{t\in[0,t_0]}\mathbf{E} \Big[\big|X(t)\big|^p \Big] \cdot \bigg| \int_0^{t_0}   ds \int_0^\infty\big|g(t_0-s,y)\big|^{2} \mu(dy)\bigg|^{p} \cr
 \ar\ar + C \cdot \sup_{t\in[0,t_0]}\mathbf{E}\Big[\big|X(t)\big|\Big]\cdot \int_0^{t_0}  ds  \int_0^\infty\big|g(t_0-s,y)\big|^{2p} \mu(dy) \cr
 \ar\leq\ar C  \cdot \sup_{t\in[0,T]}\mathbf{E} \Big[\big|X(t)\big|^p \Big] \cdot \bigg| \int_0^T  ds \int_0^\infty\big|g(s,y)\big|^{2} \mu(dy) \bigg|^{p} \cr
 \ar\ar + C \cdot \sup_{t\in[0,T]}\mathbf{E}\Big[\big|X(t)\big|\Big]\cdot \int_0^T  ds \int_0^\infty\big|g(s,y)\big|^{2p} \mu(dy)  ,
 \eeqnn
 for some constant $C$ depending only on $p$. 
 Hence the desired inequality \eqref{BDG} holds. 
 \qed

 The next two stochastic Fubini's theorems follow directly from Theorem~2.6 in \cite[p.296]{Walsh1986} and Theorem~D.2 in \cite{Xu2024b} respectively. 
 
 \begin{proposition}\label{StoFubiniThm}
 Let $T\geq 0$, $m(dt)$ be a $\sigma$-finite measure on $\mathbb{R}_+$ and $h$ be a measurable function on $\mathbb{R}_+$ such that
 \beqnn
  \int_0^T \big|h(T-t)\big|  \cdot \bigg( \int_0^t |f(s)|^2\, ds +\int_0^t  \int_0^\infty \big| g(s,y) \big|^2\, \mu(dy)\,ds \bigg)\,m(dt)  <\infty .
 \eeqnn
 If $\sup_{t\in[0,T]}  X(t)  <\infty$ a.s., we have
 \beqlb\label{eqn.SFT01}
  \int_0^T h(T-t)  B(f,t)\, m(dt) 
  \ar=\ar  \int_0^T \int_0^{X(s)} \bigg( \int_s^{T} h(T-t)f(t-s)\, m(dt) \bigg)B(ds,dz) 
 \eeqlb
 and 
 \beqlb\label{eqn.SFT02}
  \int_0^T h(T-t)   \widetilde{N}(g,t)\, m(dt)  
  \ar=\ar \int_0^T  \int_0^{X(s)}\int_0^\infty \bigg(\int_s^{T}h(T-t) g(t-s,y)\, m(dt) \bigg) \widetilde{N}(ds,dz,dy).
 \eeqlb
 	
 \end{proposition}

  \section{Local times of compound Poisson processes}
 \label{Sec.CompoundPoisson}
 \setcounter{equation}{0}
 
 In this section, we first recall some properties of compound Poisson processes with unit negative drift and positive jumps and then two stochastic Volterra representations established in \cite{Xu2024b} for their local times. 
 Thereafter, we show that after a suitable scaling, their local times behave like those of  $\xi$.
 
 \subsection{Compound Poisson processes}
 
 Let $\varPi (dy)$ be a probability law on $(0,\infty)$ with tail-distribution  $\overline{\varPi}(y):=\varPi\big([y,\infty)\big)$ and  finite mean 
 \beqnn
 \big\|\overline{\varPi} \big\|_{L^1} := \int_0^\infty \overline{\varPi}(y) \, dy
 = \int_0^\infty y \, \varPi(dy) <\infty.
 \eeqnn 
 Consider a compound Poisson process $Y:=\{Y(t):t\geq 0 \}$ defined on  $(\Omega, \mathscr{F}, \mathscr{F}_t, \mathbf{P})$ with a drift $-1$, arrival rate $\gamma>0$ and jump-size distribution $\varPi(dy)$.
 It is a spectrally positive L\'evy process with bounded variation and Laplace exponent $\varphi (\lambda):= \log \mathbf{E}[\exp\{ -\lambda Y(1) \}]$ being of the form
 \beqnn
 \varphi (\lambda)=  \lambda +  \gamma \cdot \int_0^\infty \big( e^{-\lambda y}-1 \big) \, \varPi (dy) = \lambda \cdot \Big(1 -\gamma \cdot \int_0^\infty e^{-\lambda y}\, \overline{\varPi} (y)\, dy \Big) ,\quad \lambda \geq 0.
 \eeqnn
 The function $\varphi $ is infinitely differentiable and strictly convex on $(0,\infty)$ with $\varphi(0)=0$ and $\varphi(\infty)=\infty$. 
 In particular, it is strictly increasing on $[0,\infty)$ if and only if 
 \beqnn
 \varphi^\prime(0)= -\mathbf{E} \big[ Y(1) \big]= 1 - \gamma\cdot \big\|\overline{\varPi} \big\|_{L^1} \geq 0. 
 \eeqnn 
 The process $Y$ drifts to $-\infty$, $\infty$ or is recurrent according as $\varphi^\prime(0)>0$, $<0$ or $=0$.  
 
 Denote by $\tau^+_{Y}:=\inf\big\{t> 0:Y(t)\geq 0\big\}$ the first passage time of $Y$ into $[0,\infty)$. 
 Actually, the process $Y$ always gets into the positive half line by jumping, i.e., $Y(\tau^+_{Y}-)<0$ and $Y(\tau^+_{Y})>0$ a.s.
 By Theorem~17(ii) in \cite[p.204]{Bertoin1996}, 
 \beqnn
 \mathbf{P}\big( Y(\tau^+_{Y})\in dy \,\big|\,\tau^+_{Y}<\infty \big) =\varPi^*(dy),  
 \eeqnn
 where $\varPi^*(dy) $ is the \textit{size-biased distribution} of $\varPi(dy)$ with probability density function 
 \beqnn
 \pi^*(y) := \mathbf{1}_{\{ y> 0\}} \cdot \frac{ \overline{\varPi}(y)}{ \big\|\overline{\varPi} \big\|_{L^1} }. 
 \eeqnn
 
 Denote by $L_{Y}:=\{L_{Y}(x,t):x\in \mathbb{R},\, t\geq 0 \}$ the local times of $Y$ satisfying the occupation density formula (\ref{OccupationDensityF}). 
 It also can be defined as the number of hitting times of $Y$ at each level. 
 More precisely, the following identity holds 
 \beqnn
 L_{Y} \overset{\rm a.s.}= \big\{\# \{s\in(0, t]: Y(s-)=x \}: x\in\mathbb{R},\, t\geq 0 \big\}
 \eeqnn
 with $\# A $ being the number of elements in the set $A$; see  \cite[p.128]{Bertoin1996}. 
 The right-continuity of $Y$ induces that the two-parameter process $L_{Y}$ is jointly c\`adl\`ag almost surely; see \cite{FitzsimmonsPort1989}. 
 The total local time $L_{Y}(x,\infty)$ at level $x$ is almost surely infinite for some and hence all $x\in\mathbb{R}$ if and only if $Y$ is recurrent ($\varphi'(0)=0$).  
 The right-inverse local time $\tau^L_{Y}:=\{\tau^L_{Y}(\zeta): \zeta\geq 0\}$  at level $0$ can be defined similarly as in \eqref{eqn.1009}. 
 For $a,\theta >0$, let $a\cdot Y(\theta\cdot):= \{aY(\theta t):t\geq 0\}$. 
 The identity \eqref{OccupationDensityF} along with the change of variables and the right-continuity of local times implies the following two equivalences 
 \beqlb\label{TimeSpatialChange01}
 L_{a\cdot Y(\theta\cdot)}  \overset{\rm a.s.}= \big\{(a\theta)^{-1} \cdot L_{Y}(x/a,\theta t):x\in\mathbb{R},t\geq 0 \big\}
 \quad \mbox{and}\quad 
 \tau^L_{a\cdot Y(\theta\cdot)}(\zeta) \overset{\rm a.s.}= \theta^{-1}\cdot  \tau^L_{Y}(a\theta\zeta) .
 \eeqlb 
 It is obvious that $\tau^L_Y$ only jumps at positive integer points.  
 For $k\in \mathbb{N}$, we are interested in an $\mathbb{N}$-valued c\`adl\`ag process $ L^{Y}_{k}:=\{L^{Y}_{k}(x):x\geq 0  \}$ defined as   
 \medskip\smallskip\\ \medskip   
 \centerline{\bf $L^{Y}_{k}$ is the process $  \{ L_{Y}(x,\tau^L_{Y}(k)):  x\geq 0\}$  under $\mathbf{P}(\,\cdot\,|\tau^L_{Y}(k)<\infty)$. }\smallskip
 
 \subsection{Stochastic Volterra representation of $L^{Y}_{k}$} 

 Lambert \cite{Lambert2010} proved that the process $L^Y_k$ equals in distribution to a homogeneous, binary Crump-Mode-Jagers branching process starting from $k$ ancestors, in which the residual life of ancestors is distributed as $\varPi^*$, offsprings have a common life-length distribution $\varPi$ and each individual gives birth to its children   according to a Poisson process with rate $\gamma$. 
 Further, Xu \cite{Xu2024b} established the following Hawkes representations for the process $L^Y_k$ by linking the branching process to a marked Hawkes process. 
  
 \begin{lemma}\label{Lemma.HR}
 Assume that $\varphi^\prime(0)\geq 0$, the process $L^Y_k$ equals in distribution to the unique strong solution to the stochastic Volterra equation
 \beqlb\label{HawkesRep}
 Z_k(t)\ar=\ar \sum_{j=1}^k \mathbf{1}_{\{\ell_j> t\}} + \int_0^t \int_0^{Z_k (s-)} \int_0^\infty  \mathbf{1}_{\{y>t-s\}} \, N_\varPi (ds,dz,dy), \quad t\geq 0 , 
 \eeqlb
 where $\{ \ell_j \}_{j\geq 1}$ is a sequence of i.i.d $\mathscr{F}_0$-measurable random variables distributed as $\varPi^*(dy)$ and $N_\varPi(ds,dz,dy)$ is a $(\mathscr{F}_t)$-Poisson random measure on $(0,\infty)^3$ with intensity $ \gamma \cdot ds\, dz\,\varPi(dy)$.   
 \end{lemma}
     
 Associated to the function $\gamma \cdot \overline{\varPi}$ we define the \textit{resolvent of the second kind} $R_\varPi$ by  the unique global solution in $D(\mathbb{R}_+;\mathbb{R}_+)$ to the linear Volterra equation
 \beqlb\label{Resolvent.01}
 R_\varPi(t) =  \gamma  \cdot \overline{\varPi}(t) +   \gamma \cdot \overline{\varPi} * R_\varPi (t),\quad t\geq 0.
 \eeqlb
 This equation is also known as \textit{resolvent equation} or \textit{renewal equation}.
 Moreover, its unique solution admits the following Neumann series expansion 
 \beqlb\label{Resolvent.0101}
 R_\varPi(t) = \sum_{k=1}^\infty \big( \gamma  \cdot \overline{\varPi} \big)^{(*k)}(t),\quad t\geq 0,
 \eeqlb
 where $f^{(*k)}$ denotes the $k$-th convolution of $f$.  
 For convention, we extend $ R_\varPi$ to the whole line by setting $ R_\varPi(t)=0$ for $t<0$.
 The function $R_\varPi$ is integrable if and only if $\gamma  \cdot \big\|\overline{\varPi} \big\|_{L^1} <1$; equivalently, $\varphi^\prime(0)> 0$, in which case, 
 \beqnn
  \big\| R_\varPi\big\|_{L^1} =\frac{\gamma  \cdot \big\|\overline{\varPi} \big\|_{L^1}}{1-\gamma  \cdot \big\|\overline{\varPi} \big\|_{L^1}}= \frac{\gamma  \cdot \big\|\overline{\varPi} \big\|_{L^1}}{\varphi^\prime(0)}<\infty. 
 \eeqnn
 In addition, associated to $R_\varPi$ we also introduce a two-parameter function $R$ on $\mathbb{R}^2$ defined by $R(t,y)=0$ if $t<0$ or $y<0$ and
 \beqlb\label{Resolvent.02}
 R(t,y)\ar=\ar \mathbf{1}_{\{y>t\}} + \int_0^t R_\varPi (t-s)\cdot \mathbf{1}_{\{y>s\}}\, ds
 =\mathbf{1}_{\{y>t\}} + \int_{t-y}^t R_\varPi (s) \, ds,\quad t,\, y\geq 0. 
 \eeqlb
 We refer to \cite{HorstXu2023,Xu2024b} for more detailed and intuitive explanations of $ R_\varPi $ and $R$.
 Armed with these two functions, we now recall the second stochastic Volterra equation established in \cite{Xu2024b} for $ L^{Y}_{k}$ .

 \begin{lemma}\label{Lemma.SVR}
  The stochastic Volterra equation \eqref{HawkesRep} is equivalent to 
  \beqlb\label{HawkesRepresentation}
   Z_k(t)
   \ar=\ar  
   \sum_{j=1}^{k} R(t,\ell_{j})  + \int_0^t \int_0^{Z_k (s-)} \int_0^\infty  R (t-s,y) \, \widetilde{N}_\varPi (ds,dz,dy), \quad t\geq 0,
  \eeqlb
 	where  $\widetilde{N}_\varPi(ds,dz,dy):= N_\varPi(ds,dz,dy)- \gamma \cdot ds\, dz\,\varPi(dy)$. 
 \end{lemma}
 
  \subsection{Poisson approximation}\label{Sec.CPApproximation}
  
 We now show that the local times of L\'evy process $\xi$ can be well-approximated by those of compound Poisson processes after a suitable scaling.   
 For $n\geq 1$, assume that the $n$-th compound Poisson process $Y_n:=\{ Y_n(t) :t\geq 0 \}$ has a drift $-1$, arrival rate $\gamma_n>0$ and jump-size distribution $\varPi_n(dy)$ on $(0,\infty)$. 
 Its Laplace exponent $\varphi_n$, local times $L_{Y_n}$ and right-inverse local time $\tau^L_{Y_n}$ at level $0$ and the processes $L^{Y_n}_k$ are defined as before. 
 
 Associated to a sequence of positive scaling parameters $\{ \kappa_n \}_{n\geq 1}$ such that  $\kappa_n\to\infty$ as $n\to\infty$, we define a sequence of rescaled processes $\{ \xi^{(n)} \}_{n\geq 1}$ by
 \beqlb\label{RescaledCP}
  \xi^{(n)}(t):= \frac{1}{n} \cdot Y_n(\kappa_n\cdot t),\quad t\geq 0,\, n\geq 1. 
 \eeqlb  
 Denote by $\overline\varPi_n(y)$ and $\varPi_n^*(dy)$ the tail-distribution and size-biased distribution of $\varPi_n(dy)$ respectively. 
 The process $\xi^{(n)}$ is a spectrally positive L\'evy process with Laplace exponent 
 \beqnn
  \varPhi^{(n)}(\lambda) 
  := \kappa_n\cdot \varphi_n\Big(\frac{\lambda}{n}\Big) 
  \ar=\ar \frac{\kappa_n}{n} \big( 1-\gamma_n \cdot \big\|\overline\varPi_n\big\|_{L^1} \big)\cdot \lambda + \kappa_n \gamma_n  \int_0^\infty \Big( e^{-\frac{\lambda}{n} y}-1+ \frac{\lambda}{n}  y \Big) \, \varPi_n (dy). \qquad 
 \eeqnn
 By Corollary~4.3 in \cite[p.440]{JacodShiryaev2003}, the rescaled process
 $ \xi^{(n)}$ converges weakly to $\xi$ in $D(\mathbb{R}_+;\mathbb{R})$ as $n\to\infty$ under the following necessary and sufficient condition.
 
 \begin{condition}\label{Main.Condition}
 The limit $\varPhi^{(n)}(\lambda)  \to \varPhi (\lambda)$ holds for all $\lambda \geq 0$ as $n\to\infty$.
 \end{condition}
  
 Define the local times $L_{\xi^{(n)}}$,  the right-inverse local time $\tau^L_{\xi^{(n)}}$ at level $0$ and the process $L^{\xi^{(n)}}_\zeta$ for $\zeta\geq 0$ as before. 
 By the occupation density formula \eqref{OccupationDensityF}, 
 \beqnn
 \int_\mathbb{R} f(x) L_{\xi^{(n)}}(x,t) dx \overset{\rm a.s.}=\int_0^t f\big(\xi^{(n)}(r)\big)dr 
 \overset{\rm d}\to \int_0^t f\big(\xi(r)\big)dr \overset{\rm a.s.}= \int_\mathbb{R} f(x) L_\xi(x,t) dx , 
 \eeqnn
 as $n\to\infty$ for all bounded and continuous function $f$ on $\mathbb{R}$. 
 Hence the local times $L_{\xi}$ can be obtained as the limit of $\{ L_{\xi^{(n)}}\}_{n\geq 1}$.
 This suggests us to approximate $L^{\xi}_\zeta$ by using $\big\{L^{\xi^{(n)}}_\zeta\big\}_{n\geq 1}$. Actually, it works in the sense of finite-dimensional distributions; see Theorem~2.4 in \cite{LambertSimatos2015}.

 \begin{lemma}\label{ConvergenceLocalTime}
 Under Condition~\ref{Main.Condition}, we have $L^{\xi^{(n)}}_\zeta \overset{\rm f.d.d.}\longrightarrow L^\xi_\zeta$ as $n\to\infty$.
 \end{lemma}
  
 At the end of this section, we construct another approximation for $L^\xi_\zeta $ by using solutions of the two equivalent stochastic Volterra equations \eqref{HawkesRep} and \eqref{HawkesRepresentation}. 
 For $n\geq 1$, let $R_{\varPi_n}$  be the resolvent of the second kind associated to the function $\gamma_n\cdot \overline\varPi_n$ and $R_n$  the corresponding two-parameter function defined as in \eqref{Resolvent.02}.  
 For $k\geq 1$, let $Z_{n,k}:=\{Z_{n,k}(t):t\geq 0 \}$ be the unique solution to 
 \beqlb\label{SVERepresentation}
 Z_{n,k}(t)\ar=\ar  
 \sum_{j=1}^{k} R_n(t,\ell_{n,j})  + \int_0^t \int_0^{Z_{n,k} (s-)} \int_0^\infty   R_n (t-s,y) \widetilde{N}_{\varPi_n} (ds,dz,dy) 
 \eeqlb
 with $\{ \ell_{n,j} \}_{j\geq 1}$ being a sequence of i.i.d $\mathscr{F}_0$-measurable random variables with common distribution $\varPi^*_n(dy)$ and $\widetilde{N}_{\varPi_n}(ds,dz,dy)$ being a $(\mathscr{F}_t)$-compensated Poisson random measure on $(0,\infty)^3$ with intensity $ \gamma_n \cdot ds\, dz\, \varPi_n(dy)$.  
 A direct consequence of Lemma~\ref{Lemma.SVR} tells that $L^{Y_n}_k \overset{\rm d}= Z_{n,k}$.
 For each $\zeta>0$, by \eqref{RescaledCP} and \eqref{TimeSpatialChange01} we have  \beqlb\label{eqn.28}
 L^{\xi^{(n)}}_\zeta \overset{\rm a.s.}= \bigg\{ \frac{n}{\kappa_n}\cdot L^{Y^{(n)}}_{[\zeta\cdot\kappa_n/n]}(nx): x\geq 0  \bigg\} \overset{\rm d}= \bigg\{ \frac{n}{\kappa_n} \cdot Z^{(n)}_{n,[\zeta\cdot\kappa_n/n]}(nt): t\geq 0 \bigg\} =: X^{(n)}_\zeta . 
 \eeqlb   
 Applying the change of variables to \eqref{SVERepresentation} induces that $X^{(n)}_\zeta$ is the unique strong solution to 
 \beqlb\label{Eqn.HR}
 X_{\zeta}^{(n)}(t)
 \ar=\ar \frac{n}{\kappa_n} \sum_{j=1}^{[\zeta\cdot\kappa_n/n]} R_n(nt,\ell_{n,j})  + \int_0^t\int_0^{X_{\zeta}^{(n)} (s-)} \int_0^\infty  \frac{n}{\kappa_n} \cdot R_n \big(n(t-s),y\big)\widetilde{N}^{(n)} (ds,dz,dy) \qquad 
 \eeqlb
 with $\widetilde{N}^{(n)} (ds,dz,dy) := \widetilde{N}_{\varPi_n} (n\cdot ds, n^{-1}\kappa_n \cdot dz ,dy)$. 
 The next corollary is a direct consequence of Lemma~\ref{ConvergenceLocalTime}. 
 
 \begin{corollary}\label{ConvergenceLocalTime02}
  Under Condition~\ref{Main.Condition}, we have $X^{(n)}_\zeta\overset{\rm f.d.d.}\longrightarrow L^\xi_\zeta$ as $n\to\infty$.
 \end{corollary}

   \section{Stochastic Volterra representation of $L^\xi_\zeta$}
 \label{Sec.SVETightness}
 \setcounter{equation}{0}

 In this section, we establish the stochastic Volterra representation for the process $L^\xi_\zeta$ by proving that the stochastic equation \eqref{Eqn.HR} converges weakly to \eqref{MainThm.SVE}. 
 As we have mentioned in the Introduction, the proof of tightness and the characterization of limit   processes are usually quite challenging. 
 To overcome difficulties encountered in proving the weak convergence of \eqref{Eqn.HR},  we need to specify the compound Poisson processes $\{Y_n \}_{n\geq 1}$ with  easy-to-handle jump-size distributions and also establish limit theorems for stochastic Volterra integrals driven by martingale measures. 
  
 \subsection{Asymptotic assumptions}

 Denote by  $\boldsymbol{e}_{c}$, $\boldsymbol{E}_c$ and $\overline{\boldsymbol{E}}_c$ the exponential probability density function, distribution function and tail-distribution function with mean $c$ respectively, i.e.,
 \beqnn
 \boldsymbol{e}_c(y):= \frac{1}{c}\cdot e^{-y/c},\quad 
 \boldsymbol{E}_c(y):= 1- e^{-y/c}
 \quad \mbox{and}\quad 
 \overline{\boldsymbol{E}}_c(y):=e^{-y/c}
 ,\quad y\geq 0.
 \eeqnn
 Consider two non-negative sequences $\{ \eta_n \}_{n\geq 1}\subset \mathbb{R}_+$ and $\{ \theta_n \}_{n\geq 1}\subset [0,1)$ satisfying  that 
 \begin{enumerate}
 	\item[$\bullet$]  $\eta_n=0$ and $\theta_n =  c \cdot \bar{\nu}(0)\cdot n^{-2}$ when $\bar{\nu}(0)<\infty$; 
 	
 	\item[$\bullet$]  $ n   \theta_n \rightarrow 0$, $n^2 \theta_n \to  \infty $ and $\bar{\nu}(\eta_n) = c^{-1}\cdot n^2 \theta_n $ when $\bar{\nu}(0)=\infty$. 
 \end{enumerate}
 For each $n\geq 1$, associated to $\eta_n$ and $\nu(dy)$ we define a probability law $  \Lambda_n(dy)  $ on $(0,\infty)$ by
 \beqlb\label{Lambda.n}
 \Lambda_n(dy):=\frac{1}{\bar{\nu}(\eta_n)}\cdot \nu\Big(\eta_n+\frac{dy}{n}\Big)
 \quad\mbox{with}\quad 
 \overline{\Lambda}_n(y):=\Lambda_n\big([y,\infty)\big) =\frac{\bar\nu(\eta_n+y/n)}{\bar{\nu}(\eta_n)},\quad y\geq 0 . 
 \eeqlb 
 By the finite integral in \eqref{LevyTriplet}, the law $  \Lambda_n(dy)  $ has a finite mean 
 \beqlb\label{eqn.540} 
 \big\|  \overline{\Lambda}_n \big\|_{L^1} 
 \ar=\ar \int_0^\infty y \, \Lambda_n(dy) 
 = \frac{n}{\bar{\nu}(\eta_n)}\int_0^\infty \bar\nu(\eta_n+y) \,dy 
 = \frac{c}{n\theta_n} \int_0^\infty  \bar\nu(\eta_n+y) \,dy <\infty,
 \eeqlb
 and its size-biased distribution $\Lambda_n^*(dy)$ has density function 
 \beqlb\label{eqn.541}
 \boldsymbol{\lambda}_n^*(y):= \frac{\overline{\Lambda}_n(y)}{\|  \overline{\Lambda}_n \|_{L^1} } = \frac{1}{n}\cdot \frac{\bar\nu(\eta_n+y/n)}{\int_0^\infty \bar\nu(\eta_n+z) \,dz},\quad y>0 . 
 \eeqlb
 \begin{proposition} \label{Prop.502}
 As $n\to \infty$, we have 
  \beqnn 
   \int_0^\infty \frac{\bar\nu(\eta_n+y)}{\bar{\nu}(\eta_n)}\,dy
  \to
  \begin{cases} \displaystyle
  \doublebar\nu(0)/\bar{\nu}(0),& \mbox{if $\bar{\nu}(0)<\infty$};\vspace{3pt}\\
  	0, &  \mbox{if $\bar{\nu}(0)=\infty$},
  	\end{cases}
  	\qquad\mbox{and}\qquad 
  	 \frac{1}{n}\int_0^\infty  \bar\nu(\eta_n+y) \,dy \to 0.
  \eeqnn
 \end{proposition}
 \proof The second limit follows directly from the first one and the fact that $n^{-1}
 \cdot\bar{\nu}(\eta_n)\to0$. 
 For the first limit, it is obvious when $\bar{\nu}(0)<\infty$. 
 When $\bar{\nu}(0)=\infty$, the monotonicity of $\bar\nu$ induces that for $\epsilon>0$, 
 \beqnn
 \int_0^\infty \frac{\bar\nu(\eta_n+y)}{\bar{\nu}(\eta_n)}\,dy
 \ar=\ar \int_0^\epsilon \frac{\bar\nu(\eta_n+y)}{\bar{\nu}(\eta_n)}\,dy +\int_\epsilon^\infty \frac{\bar\nu(\eta_n+y)}{\bar{\nu}(\eta_n)}\,dy
 \leq  \epsilon  +  \int_\epsilon^\infty \frac{\bar\nu(y)}{\bar{\nu}(\eta_n)}\,dy ,
 \eeqnn
 which goes to $0$ by letting $n\to\infty$ and then $\epsilon\to 0+$. 
 \qed
 
 For each $n\geq 1$, assume that the $n$-th compound Poisson process $Y_n$ has arrival rate
 \beqlb\label{Gamma.n}
 \gamma_n=\frac{1}{c} \cdot \Big( 1- \frac{b}{n} - \frac{1}{n}  \int_0^\infty y \,\nu(\eta_n+dy) \Big)^+ 
 = \frac{1}{c} \cdot \Big( 1- \frac{b}{n} - \frac{1}{n}  \int_0^\infty \bar\nu(\eta_n+y)\, dy \Big)^+
 \leq \frac{1}{c},
 \eeqlb
 which goes to $1/c$ as $n\to\infty$,
 and its jump-size distribution $\varPi_n(dy)$ on $(0,\infty)$ has density function
 \beqlb \label{Pi.n}
 \pi_n(y) \ar:=\ar (1-\theta_n ) \cdot \boldsymbol{e}_{c}(y)+\theta_n\cdot \boldsymbol{e}_{c}*d\Lambda_n(y) ,\quad y\geq 0 .  
 \eeqlb  
 Integrating both sides over $[x,\infty)$ and then using Fubini's theorem show that $\varPi_n(dy)$ has tail-distribution 
 \beqlb \label{eqn.504}
  \overline{\varPi}_n(x):= \varPi_n\big([x,\infty)\big)= \overline{\boldsymbol{E}}_c(x) +\theta_n \cdot  \boldsymbol{e}_{c}*\overline{\Lambda}_n(x),\quad x\geq 0 . 
 \eeqlb
 In addition, by the change of variables and \eqref{eqn.540}, 
 \beqlb\label{eqn.505}
  \big\|\overline{\varPi}_n\big\|_{L^1}
  = c + \theta_n \cdot \big\| \overline{\Lambda}_n \big\|_{L^1}
  \ar=\ar c+  \frac{n\theta_n}{\bar{\nu}(\eta_n)} \int_0^\infty \bar\nu(\eta_n+y) \,dy \cr
  \ar=\ar c\cdot \Big(1+ \frac{1}{n} \int_0^\infty \bar\nu(\eta_n+y)\, dy  \Big) , 
 \eeqlb
 which along with  \eqref{Gamma.n} induces that 
 \beqlb\label{Asym.Mean}
 \gamma_n\cdot \big\|\overline{\varPi}_n\big\|_{L^1} =  1-\frac{b}{n}  -\frac{1}{n^2} 
 \cdot \Big(b+ \int_0^\infty \bar\nu(\eta_n+y)\, dy \Big)\cdot \int_0^\infty \bar\nu(\eta_n+z)\, dz \in (0,1),
 \eeqlb
 and further as $n\to\infty$,
 \beqlb\label{Asym.Critical}
 \gamma_n\cdot \big\|\overline{\varPi}_n\big\|_{L^1} \to 1
 \quad \mbox{and}\quad 
 n\cdot \big( 1- \gamma_n\cdot \big\|\overline{\varPi}_n\big\|_{L^1} \big)\to b.
 \eeqlb 
 Recall the rescaled process $ \xi^{(n)} $ defined by \eqref{RescaledCP} with  $\kappa_n=n^2$. 
 It has Laplace exponent $ \varPhi^{(n)} $ as follows
 \beqlb\label{LapXi.n01}
 \varPhi^{(n)}(\lambda) 
 \ar=\ar n\cdot \big( 1-\gamma_n \cdot \big\|\overline\varPi_n\big\|_{L^1} \big)\cdot \lambda + n^2\gamma_n    \int_0^\infty \Big( e^{-\frac{\lambda}{n} y}-1+ \frac{\lambda}{n} y\Big)   \pi_n (y)\, dy .
 \eeqlb 
 
 \begin{proposition}\label{Prop.501}
 The sequence $\big\{\varPhi^{(n)} \big\}_{n\geq 1}$ satisfies	Condition~\ref{Main.Condition}.  
 \end{proposition}
 \proof
 For each $\lambda \geq 0$, plugging \eqref{Pi.n} into \eqref{LapXi.n01}   and then using \eqref{Asym.Critical} induce that as $n\to\infty$,
 \beqlb\label{eqn.501}
 \varPhi^{(n)}(\lambda)  
 \ar\sim\ar b\cdot \lambda + \frac{n^2}{c} \int_0^\infty \Big( e^{-\frac{\lambda}{n} y}-1+ \frac{\lambda}{n} y\Big)   \boldsymbol{e}_{c}(y)\, dy\cr
 \ar\ar +\bar{\nu}(\eta_n) \cdot \int_0^\infty \Big( e^{-\frac{\lambda}{n} y}-1+ \frac{\lambda}{n}  y\Big) \boldsymbol{e}_{c}*d\Lambda_n(y) \,dy  .
 \eeqlb 
 Applying the dominated convergence theorem and the fact that $e^{-z}-1+z\sim z^2/2$ as $z\to0$ to the second term on the right side of \eqref{eqn.501},  
 \beqnn
 \lim_{n\to\infty} \frac{n^2}{c} \int_0^\infty \Big( e^{-\frac{\lambda}{n} y}-1+ \frac{\lambda}{n} y\Big)   \boldsymbol{e}_{c}(y)\, dy 
 = \frac{\lambda^2}{2c} \int_0^\infty   y^2  \boldsymbol{e}_{c}(y)\, dy  
 = c\cdot \lambda^2. 
 \eeqnn
 An application of the basic property of convolution to the third term shows that it equals to 
 \beqnn
  \bar{\nu}(\eta_n) \cdot \bigg(\int_0^\infty   e^{-\frac{\lambda}{n} y}\boldsymbol{e}_{c}(y)\, dy \cdot \int_0^\infty   e^{-\frac{\lambda}{n} z}\Lambda_n(dz) -1+ c\cdot \frac{\lambda}{n} +\int_0^\infty \frac{\lambda}{n} y\, \Lambda_n(dy)  \bigg),
 \eeqnn 
 which can be written as the sum of  the following three terms 
 \beqnn
 A_1^{(n)}(\lambda)\ar:=\ar  \frac{\bar{\nu}(\eta_n)}{n} \cdot  c \cdot\lambda,\cr
 \ar\ar\cr
 A_2^{(n)}(\lambda)\ar:=\ar \bar{\nu}(\eta_n) \cdot   \int_0^\infty \Big(  e^{-\frac{\lambda}{n}  y}-1+\frac{\lambda}{n} y\Big)\, \Lambda_n(dy) , \cr
 A_3^{(n)}(\lambda)\ar:=\ar   \bar{\nu}(\eta_n) \cdot  \int_0^\infty  \Big( e^{-\frac{\lambda}{n}  y}-1\Big)\boldsymbol{e}_{c}(y)\, dy \cdot \int_0^\infty   e^{-\frac{\lambda}{n} z}\Lambda_n(dz) .     
 \eeqnn 
 The limit $ n^{-1}\cdot \bar{\nu}(\eta_n) \to 0$ immediately yields that $|A_1^{(n)}(\lambda)| \to 0$ as $n\to\infty$. 
 Also, the inequality $|e^{-z}-1|\leq z$ for any $z\geq 0$ induces that 
 \beqnn
 \big|A_3^{(n)} (\lambda)\big|\ar\leq \ar   \lambda\cdot  
 \int_0^\infty y\, \boldsymbol{e}_{c}(y)\, dy \cdot \frac{\bar{\nu}(\eta_n)}{n} \to 0. 
 \eeqnn
 By using (\ref{Lambda.n}), the change of variables and then the vague convergence of $\nu(\eta_n+dy)$ to $\nu(dy)$, 
 \beqnn
  A_2^{(n)}(\lambda)\ar=\ar    \int_0^\infty \Big(  e^{-\lambda\frac{y}{n}}-1+\lambda\frac{y}{n} \Big)\, \nu\Big(\eta_n+\frac{dy}{n}\Big)\cr
  \ar=\ar \int_0^\infty \big(  e^{-\lambda  y}-1+\lambda y\big)\, \nu(\eta_n+dy)\cr
  \ar\to\ar \int_0^\infty \big(  e^{-\lambda  y}-1+\lambda y\big)\, \nu(dy),
 \eeqnn
 as $n\to\infty$. 
 Taking all preceding limits back into \eqref{eqn.501}
 induces the limit $\varPhi^{(n)}(\lambda)  \to \varPhi (\lambda) $ for all $\lambda \geq 0$.
 \qed 
 
 By Proposition~\ref{Prop.501} and Corollary~\ref{ConvergenceLocalTime02}, the process $L^\xi_\zeta$ can be well-approximated in the sense of finite-dimensional distributions by the solution $X^{(n)}_\zeta $ of \eqref{Eqn.HR} with $\kappa_n=n^2$. 
 Before establishing the weak convergence of $\{X^{(n)}_\zeta \}_{n\geq 1}$, we first provide a stochastic Volterra equation that is equivalent to \eqref{Eqn.HR} and much easy to be dealt with. 
 By (\ref{eqn.504}) and (\ref{eqn.505}), the size-biased distribution $\varPi^*_n(dy)$ on $\mathbb{R}_+$ of $\varPi_n(dy)$ has probability density function
 \beqlb\label{eqn.506}
  \pi^*_n(y)  := \frac{\overline\varPi_n(y)}{\big\|\overline{\varPi}_n\big\|_{L^1}} 
  \ar=\ar \frac{\overline{\boldsymbol{E}}_c(y)}{\big\|\overline{\varPi}_n\big\|_{L^1}}  + \frac{\theta_n  \cdot \boldsymbol{e}_{c}*\overline{\Lambda}_n(y)}{\big\|\overline{\varPi}_n\big\|_{L^1}}
  = p_n \cdot  \boldsymbol{e}_{c}(y) +(1-p_n)\cdot    \boldsymbol{e}_{c}*\boldsymbol{\lambda}^*_n(y),
 \eeqlb
 which is a weighted sum of the two probability density functions $\boldsymbol{e}_{c}$ and $\boldsymbol{e}_{c}*\boldsymbol{\lambda}^*_n$ with weight  
 \beqlb\label{eqn.pn}
 p_n:=\Big(1+ \frac{1}{n} \int_0^\infty \bar\nu(\eta_n+y)\, dy  \Big)^{-1} 
 = 1-\frac{1}{n} \int_0^\infty \bar\nu(\eta_n+y)\, dy + o\big(1/n\big),
 \eeqlb
 as $n\to\infty$. 
 Let $B_n^c(\zeta)$ be a binomial random variable with number of trials $[n\zeta]$ and success probability $p_n$.
 Moreover, let $B_n^{\Lambda}(\zeta):=[n\zeta]-B^c_n(\zeta)$, which is also binomially distributed with parameters $([n\zeta],1-p_n)$. 
 Consider two independent sequences of i.i.d. non-negative random variables $\{ \ell_{n,j}^c \}_{j\geq 1}$ and $\{ \ell_{n,j}^\Lambda \}_{j\geq 1}$ satisfying that 
 \begin{enumerate}
 	\item[$\bullet$] $\ell_{n,j}^c$ has probability density function $\boldsymbol{e}_{c}$;
 	
 	\item[$\bullet$] $\ell_{n,j}^\Lambda$ has probability density function $ \boldsymbol{e}_{c}*\boldsymbol{\lambda}^*_n$. 
 \end{enumerate} 
 For each $j\geq 1$, the random variable $\ell_{n,j}$ can be realized by setting  $\ell_{n,j}=\ell_{n,j}^c$ with probability $p_n$ and $\ell_{n,j}=\ell_{n,j}^\Lambda$ with probability $1-p_n$.
 Consequently, one can identify that 
 \beqnn
 \sum_{j=1}^{[n\zeta]} R_n(\cdot,\ell_{n,j})  \overset{\rm d}=   \sum_{i=1}^{B_n^c(\zeta)} R_n(\cdot,\ell_{n,i}^c) +  \sum_{j=1}^{B_n^\Lambda(\zeta)} R_n(\cdot,\ell_{n,j}^\Lambda) . 
 \eeqnn
  
 Note that the compensated Poisson random measure $\widetilde{N}^{(n)} (ds,dz,dy)=\widetilde{N}_{\varPi_n}(n\cdot ds, n\cdot dz ,dy)$  on $(0,\infty)^3$ has intensity $ n^2\gamma_n \cdot ds\,dz\, \varPi_n(dy)$, which, by \eqref{Pi.n}, can be written as 
 \beqlb\label{Intensities}
 n^2\gamma_n(1-\theta_n ) \cdot \boldsymbol{e}_{c}(y)\,\cdot ds\,dz\,dy
 + 
 n^2\gamma_n\theta_n \cdot  \boldsymbol{e}_{c}*d\Lambda_n(y)\,\cdot ds\,dz\,dy.
 \eeqlb
 This allows us to define on an extension of the original probability space two orthogonal  compensated Poisson random measures $\widetilde{N}^{(n)}_c (ds,dz,dy)$ and $\widetilde{N}^{(n)}_\Lambda (ds,dz,dy)$ on $(0,\infty)^3$ with intensities being the two terms in \eqref{Intensities} respectively such that 
 \beqnn
 \widetilde{N}^{(n)} (ds,dz,dy)
 =\widetilde{N}^{(n)}_c (ds,dz,dy)+\widetilde{N}^{(n)}_\Lambda (ds,dz,dy).
 \eeqnn 
 
 In conclusion, the preceding arguments and notation allow us to write the stochastic Volterra equation \eqref{Eqn.HR} into the following equivalent form
 \beqlb\label{Eqn.HR01}
 X_{\zeta}^{(n)}(t) \ar=\ar \sum_{i=1}^4 \mathbf{I}^{(n)}_i(t) ,\quad t\geq 0 ,
 \eeqlb
 in which the four summands on the right side are given by 
 \beqlb
  \mathbf{I}^{(n)}_1(t)  \ar:=\ar \frac{1}{n} \sum_{j=1}^{B_n^c(\zeta)} R_n(nt,\ell_{n,j}^c) ,
  \quad \ \ \,  
  \mathbf{I}^{(n)}_2(t)  :=  \frac{1}{n} \sum_{j=1}^{B_n^\Lambda(\zeta)} R_n(nt,\ell_{n,j}^\Lambda),\label{eqn.I} \\
  \mathbf{I}^{(n)}_3(t) \ar:=\ar   \int_0^t\int_0^{X_{\zeta}^{(n)} (s-)} \int_0^\infty  \frac{1}{n}\cdot R_n\big(n(t-s),y\big)\widetilde{N}^{(n)}_c (ds,dz,dy), \label{eqn.Jc}\\
  \mathbf{I}^{(n)}_4(t) \ar:=\ar \int_0^t\int_0^{X_{\zeta}^{(n)} (s-)} \int_0^\infty  \frac{1}{n}\cdot R_n\big(n(t-s),y\big)\widetilde{N}^{(n)}_{\Lambda} (ds,dz,dy).  \label{eqn.JLambda}
 \eeqlb
 In the following subsections, we will prove the jointly weak convergence of the five terms in \eqref{Eqn.HR01} to the corresponding terms in \eqref{MainThm.SVE}. 
   
 \subsection{Auxiliary lemmas}\label{Sec.Auxiliarylemmas}
 
 
 \subsubsection{Asymptotics of rescaled resolvent}
 
 In view of \eqref{eqn.I}-\eqref{eqn.JLambda}, the two time-scaled functions $R_{\varPi_n}(n\cdot) $ and $R_n(n\cdot,y)$ play an important role in the convergence of $\{ X^{(n)}_\zeta \}_{n\geq 1}$. 
 By \eqref{Resolvent.01} and \eqref{Resolvent.02}, we have for $t\geq 0$ and $y>0$, 
 \beqlb
 R_{\varPi_n}(t)\ar=\ar\gamma_n \cdot \overline\varPi_n(t)+  \gamma_n \overline\varPi_n * R_{\varPi_n} (t),\label{Eqn.RPi.n} \\
 R_n(t, y)\ar=\ar 1_{\{y> t\}}+\int_{t-y}^t R_{\varPi_n} (s) \, ds.  \label{Eqn.R.n}
 \eeqlb
 In this section, we mainly study their asymptotic properties including uniform upper bounds and $L^p$-convergence with the help of the following upper-bound estimates for probability density functions of geometric random sums. 
 
 \begin{proposition}\label{Prop.GeometricSum}
 Consider a geometrically distributed random variable $N_q$ with parameter $q\in(0,1)$ and two independent sequences of  i.i.d. non-negative random variables $\{e_{i}\}_{i\geq 1}$ and $\{U_{i}\}_{i\geq 1}$. 
 If $e_{i}$ is exponentially distributed with mean $\lambda>0$, then the geometric sum $\sum_{i=1}^{N_q} (e_i+U_i)$ has a density function uniformly bounded by $q/\lambda$.  
 \end{proposition}
 \proof 
 Let $F_U(du)$ be the probability law of $U_1$ and $e_{k,\lambda}$ the density function of $\sum_{i=1}^k e_i$ for $k\geq 1$.  
 The density function of geometric sum $\sum_{i=1}^{N_q} (e_i+U_i)$ admits the following representation 
 \beqnn
 \lefteqn{\sum_{k=1}^\infty q(1-q)^{k-1} \int_{\mathbb{R}_+^k} e_{k,\lambda}(x-u_1-\cdots-u_k)  \prod_{i=1}^k F_U(du_i) }\ar\ar\cr
 \ar=\ar \int_{\mathbb{R}_+^\infty} \sum_{k=1}^\infty q(1-q)^{k-1}  e_{k,\lambda}(x-u_1-\cdots-u_k)  \prod_{i=1}^\infty F_U(du_i) .
 \eeqnn
 By (4.14) and (4.15) in \cite[p.156]{Kalashnikov1997}, 
 \beqnn
 \sup_{x\in\mathbb{R}}\, \sum_{k=1}^\infty q(1-q)^{k-1}  e_{k,\lambda}(x-u_1-\cdots-u_k) 
 \ar\leq\ar \sup_{x\in\mathbb{R}}\, \sum_{k=1}^\infty q(1-q)^{k-1}  e_{k,\lambda}(x)   
  = \frac{q}{\lambda},
 \eeqnn
 and hence the  density function of the geometric sum  is uniformly bounded by $q/\lambda$.
 \qed

 \begin{lemma}\label{Lemma.UpperBoundR}
 	The following two upper bounds hold uniformly in $t\geq 0$ and $y> 0$,
 	\beqlb\label{eqn.509}
 	  \sup_{n\geq 1} R_{\varPi_n}(t) \leq \frac{1}{c}
 	 \quad \mbox{and}\quad 
 	 \sup_{n\geq 1} R_n(t, y) \leq 1+ \frac{t\wedge y}{c}.
 	\eeqlb 
 \end{lemma}
 \proof Here we just prove the first upper bound. 
 The second one is a direct consequence of the first one. 
 By \eqref{Resolvent.0101} and the first equality in \eqref{eqn.506}, the resolvent $ R_{\varPi_n}$ has the following expansion 
 \beqlb\label{eqn.508}
  R_{\varPi_n}(t)  
  \ar=\ar \frac{\gamma_n \cdot \big\|\overline{\varPi}_n\big\|_{L^1}}{1-\gamma_n \cdot \big\|\overline{\varPi}_n\big\|_{L^1}} \cdot \sum_{k=1}^\infty \big(1-\gamma_n \cdot \big\|\overline{\varPi}_n\big\|_{L^1} \big) \big( \gamma_n \cdot \big\|\overline{\varPi}_n\big\|_{L^1} \big)^{k-1}\cdot \big(  \pi^*_n \big)^{(*k)}(t). 
  \eeqlb
 By \eqref{Asym.Mean} and \eqref{eqn.506},  we have $\gamma_n \cdot \big\|\overline{\varPi}_n\big\|_{L^1}<1$ and $\pi^*_n(y)= \boldsymbol{e}_{c} * d\mu_n(y)$ with 
 $\mu_n(dy):=p_n \cdot  \delta_0(dy) +(1-p_n)\cdot  \boldsymbol{\lambda}^*_n(y)\,dy$ for $y\geq0$. 
 Hence the sum on the right side of \eqref{eqn.508} is the probability density function of the geometric sum in Proposition~\ref{Prop.GeometricSum} with $q=1-\gamma_n \cdot \big\|\overline{\varPi}_n\big\|_{L^1}$, $\lambda=c$ and $ U_i$ distributed as $\mu_n(dy)$. 
 Hence 
 \beqnn
 \sup_{t\geq 0}R_{\varPi_n}(t) 
 \leq  \frac{\gamma_n \cdot \big\|\overline{\varPi}_n\big\|_{L^1}}{1-\gamma_n \cdot \big\|\overline{\varPi}_n\big\|_{L^1}} \cdot \frac{1-\gamma_n \cdot \big\|\overline{\varPi}_n\big\|_{L^1}}{c}  = \frac{\gamma_n \cdot \big\|\overline{\varPi}_n\big\|_{L^1}}{c}\leq \frac{1}{c}, 
 \eeqnn 
  uniformly in $n\geq 1$. 
 \qed 
  
 In the next lemma, we prove the weak convergence of measures with densities $R_{\varPi_n}(n\cdot)$ and $R_n(n\cdot,y)$ by considering the aysmptotics of their Laplace transforms.   
 
 \begin{lemma}\label{Lemma.ConR}
 For each $\beta\geq 0$ and $y\geq 0$, we have  as $n\to\infty$, 
  \beqnn
  \int_0^t e^{-\beta s}R_{\varPi_n}(ns) \,ds \overset{\rm u.c.}\to \int_0^t e^{-\beta s} W'(s)\, ds
  \quad \mbox{and}\quad 
  \int_0^t e^{-\beta s} R_{n}(ns,y)\, ds \overset{\rm u.c.}\to y\cdot \int_0^t e^{-\beta s} W'(s)\, ds.
  \eeqnn
 \end{lemma}
 \proof Denote by $\mathcal{L}_{R_{\varPi_n}}$ and $\mathcal{L}_{\overline{\varPi}_n}$ the Laplace transforms of $R_{\varPi_n}$ and $\overline{\varPi}_n$ respectively.  
 Taking Laplace transforms on both sides of \eqref{Eqn.RPi.n} shows that 
 \beqnn
  \mathcal{L}_{R_{\varPi_n}}(\lambda)= \gamma_n \cdot \mathcal{L}_{\overline{\varPi}_n}(\lambda) + \gamma_n \cdot \mathcal{L}_{\overline{\varPi}_n}(\lambda)\cdot \mathcal{L}_{R_{\varPi_n}}(\lambda)
  \quad \mbox{and}\quad 
  \mathcal{L}_{R_{\varPi_n}}(\lambda)=\frac{\gamma_n \cdot \mathcal{L}_{\overline{\varPi}_n}(\lambda) }{1-\gamma_n\cdot \mathcal{L}_{\overline{\varPi}_n}(\lambda)}, \quad \lambda >0,
 \eeqnn
 which along with the change of variables induces that 
  \beqlb\label{eqn.410}
  \int_0^{\infty} e^{-\lambda t}  R_{\varPi_n}(nt)\,  dt  
  = \frac{1}{n}\cdot\mathcal{L}_{R_{\varPi_n}}(\lambda/n) 
  = \frac{\gamma_n \cdot \mathcal{L}_{\overline{\varPi}_n}(\lambda/n)}{n\big(1-\gamma_n\cdot \mathcal{L}_{\overline{\varPi}_n}(\lambda/n)\big)}.
  \eeqlb
 By using the dominated convergence theorem and then \eqref{Asym.Mean}, we first have as $n\to\infty$, 
 \beqnn
  \gamma_n \cdot \mathcal{L}_{\overline{\varPi}_n}(\lambda/n)\sim \gamma_n\cdot \big\|\overline{\varPi}_n \big\|_{L^1} \to 1.
 \eeqnn 
 Moreover, by using integration by parts to the denominator of the last fraction in \eqref{eqn.410}, 
 \beqlb\label{eqn.553}
  n\big(1-\gamma_n\cdot \mathcal{L}_{\overline{\varPi}_n}(\lambda/n)\big) 
  \ar=\ar \frac{n^2}{\lambda} \Big(\frac{\lambda}{n}-\gamma_n\cdot  \int_0^{\infty} \frac{\lambda}{n} e^{- \frac{\lambda}{n} t } \overline{\varPi}_n (t)\,  dt\Big)\cr
  \ar=\ar \frac{n^2}{\lambda} \Big(\frac{\lambda}{n}+ \gamma_n\cdot\int_0^\infty \big( 1- e^{- \frac{\lambda}{n} t }\big)\,  \overline{\varPi}_n(dt) \Big)= \frac{\varPhi^{(n)}(\lambda)}{\lambda}, 
 \eeqlb
 which goes to $\varPhi(\lambda)/\lambda$ as $n\to\infty $ by Proposition~\ref{Prop.501}. 
 Plugging these back into \eqref{eqn.410} gives that 
 \beqlb\label{eqn.515}
  \lim_{n\to\infty} \int_0^{\infty} e^{-\lambda t}  R_{\varPi_n}(nt)\,  dt = \frac{\lambda}{\varPhi(\lambda)} = \int_0^\infty e^{-\lambda t} W'(t)\, dt, 
 \eeqlb
 which induces the weak convergence of measures with density functions $\{R_{\varPi_n}(n\cdot)\}_{n\geq 1}$ to the measure with density $W'$ and hence the first desired locally uniform convergence holds. 
 Similarly, we also have
 \beqnn
  \int_0^{\infty} e^{-\lambda t}  R_n(nt,y)\,  dt  
  = \frac{\int_0^\infty e^{-\frac{\lambda}{n} t} \cdot \mathbf{1}_{\{ y>t\}}\, dt }{n\big(1-\gamma_n\cdot \mathcal{L}_{\overline{\varPi}_n}(\lambda/n)\big)} 
  \to   y\cdot \frac{\lambda}{\varPhi(\lambda)} = \int_0^\infty e^{-\lambda t}\cdot y\cdot W'(t)\, dt,
 \eeqnn
  as $n\to\infty$ and the second desired locally uniform convergence also holds.  
 \qed
 
 \begin{corollary}\label{Coro.UpperBound}
 	We have $W'(t)\leq 1/c$ for any $t\in\mathbb{R}$.
 \end{corollary}
 \proof The definition of $W'$ first tells that $W'(t)=0$ for any $t<0$. 
 By Lemma~\ref{Lemma.ConR}, the $\sigma$-finite measure on $\mathbb{R}_+$ with density $R_{\varPi_n}(nt)$ converges vaguely as $n\to\infty$ to the $\sigma$-finite measure with density $W'(t)$, which along with Lemma~\ref{Lemma.UpperBoundR} induces that the $\sigma$-finite measure on $\mathbb{R}_+$ with density $1/c-R_{\varPi_n}(nt)$ converges vaguely as $n\to\infty$ to the $\sigma$-finite measure with density $1/c -W'(t)$ that is non-negative almost everywhere. 
 Finally, the continuity of $W'$ on $\mathbb{R}_+$ yields that $W'(t)\leq 1/c$ for all $t\geq 0$.
 \qed 
 
 The lack of monotonicity makes the proof of point-wise convergence of $R_{\varPi_n}$ and $R_n$ to the corresponding limits quite challenging or even impossible. 
 The next best thing we can do is to establish their $L^p$-convergence with the help of the following two uniform bound estimates for the Fourier transforms of $\{\overline{\varPi}_n\}_{n\geq 1}$. 
 
 \begin{proposition}\label{Thm.UpperBound}
  There exists a constant $C>0$ such that for any $\beta,y \geq 0$ and $z\in \mathbb{R}$, 
  \beqlb\label{eqn.514}
  \sup_{n\geq 1}\bigg|\int_0^\infty e^{\mathtt{i}zt -\frac{\beta}{n}  t} \cdot \gamma_n   \overline{\varPi}_n(t)\, dt \bigg|
	\leq  \frac{C}{|z|}\wedge 1 
	\quad \mbox{and}\quad 
	\sup_{n\geq 1}\bigg|\int_0^\infty e^{\mathtt{i}zt -\frac{\beta}{n}  t} \cdot \mathbf{1}_{\{y>t\}}\, dt \bigg|
	\leq  \frac{C}{|z|}\wedge y .
  \eeqlb
 \end{proposition}
 \proof 
 By \eqref{Asym.Mean}, we first have uniformly in $n\geq 1$, $z\in \mathbb{R}$ and $\beta \geq 0$,
 \beqlb\label{eqn.511}
  \bigg|\int_0^\infty e^{\mathtt{i}zt -\frac{\beta}{n} t} \cdot \gamma_n   \overline{\varPi}_n(t)\, dt \bigg| 
 \leq \gamma_n \cdot \big\|\overline{\varPi}_n\big\|_{L^1}
 \leq  1 .
 \eeqlb 
 Moreover, the differentiation property of Fourier transform induces that 
 \beqnn
 \int_0^\infty e^{\mathtt{i}zt -\frac{\beta}{n} t} \cdot \gamma_n  \overline{\varPi}_n(t)\, dt 
 \ar=\ar  \frac{\gamma_n}{z} \int_0^\infty e^{\mathtt{i}zt } \cdot \frac{\partial}{\partial t} \big( e^{-\frac{\beta}{n} t}   \overline{\varPi}_n(t)\big)\, dt \cr
 \ar=\ar  -\frac{\gamma_n}{z} \int_0^\infty e^{\mathtt{i}zt  -\frac{\beta}{n} t} \Big(   \frac{\beta}{n} +\pi_n(t)\Big)\, dt ,
 \eeqnn
 which along with $\gamma_n\leq 1/c$ induces that
 \beqnn
  \bigg|  \int_0^\infty e^{\mathtt{i}zt -\frac{\beta}{n} t} \cdot \gamma_n   \overline{\varPi}_n(t)\, dt \bigg| 
  \leq \frac{\gamma_n}{|z|} \int_0^\infty e^{ -\frac{\beta}{n} t} \cdot \Big(   \frac{\beta}{n} +\pi_n(t)\Big)\, dt 
   \leq \frac{2}{c\cdot |z|} . 
 \eeqnn 
 This together with \eqref{eqn.511} yields the first upper bound in \eqref{eqn.514}. 
 For the second one, we also have 
 \beqnn
 \sup_{n\geq 1}\bigg|\int_0^\infty e^{\mathtt{i}zt -\frac{\beta}{n}  t} \cdot \mathbf{1}_{\{y>t\}}\, dt \bigg| 
 \leq \int_0^\infty   \mathbf{1}_{\{y>t\}}\, dt =y .
 \eeqnn
  Additionally, integration by parts gives that 
  \beqnn
  \int_0^\infty e^{\mathtt{i}zt -\frac{\beta}{n}  t} \cdot \mathbf{1}_{\{y>t\}}\, dt = \frac{e^{(\mathtt{i}z -\frac{\beta}{n})y}-1}{\mathtt{i}z  - \beta/n} ,
  \eeqnn
 and then 
 \beqnn
  \sup_{n\geq 1}\bigg|\int_0^\infty e^{\mathtt{i}zt -\frac{\beta}{n}  t} \cdot \mathbf{1}_{\{y>t\}}\, dt \bigg| 
  \leq \sup_{n\geq 1} \frac{2}{|\mathtt{i}z  - \beta/n|} \leq \frac{2}{|z|} . 
  \eeqnn
  Consequently, the second upper bound in \eqref{eqn.514} holds. 
 \qed
  
 \begin{proposition}
  There exists a constant $\beta_0\geq 0$ such that for each $\beta \geq \beta_0$, we can always find two constants $n_0\geq 1$ and $T_0>0$ satisfying that 
 \beqlb\label{eqn.412}
 \sup_{n\geq n_0} \int_{T_0}^\infty  t e^{-\frac{\beta}{n} t} \cdot \gamma_n  \overline\varPi_n(t)\, dt \leq   \frac{c}{8}  
  \quad \mbox{and} \quad
  \inf_{n\geq n_0}\int_0^{T_0} t e^{-\frac{\beta}{n} t} \cdot \gamma_n \overline\varPi_n(t)\,dt \geq \frac{3\cdot c}{4} .
 \eeqlb 
 \end{proposition}
 \proof 
 Firstly, a simple calculation along with \eqref{eqn.504} induces that 
 \beqnn
 \int_0^{T} t e^{-\frac{\beta}{n} t} \cdot \gamma_n  \overline\varPi_n(t)\, dt 
 \geq \gamma_n  c \int_0^{T} t e^{-\frac{\beta}{n} t} \boldsymbol{e}_{c}(t)\, dt
 \geq \frac{\gamma_n}{(1/c+\beta/n)^2} \Big( 1-e^{-\frac{T}{c}}- \frac{T}{c} e^{-\frac{T}{c}} \Big) \to c,
 \eeqnn
 as $n,T\to\infty$, and also by the fact that $\gamma_n\leq 1/c$,
 \beqlb\label{eqn.610}
 \int_T^\infty t e^{-\frac{\beta}{n} t} \cdot \gamma_n \overline\varPi_n(t)\, dt 
 \ar\leq\ar   \int_T^\infty t \cdot  \boldsymbol{e}_{c}(t)\, dt
 + \frac{ \theta_n }{c}   \int_0^\infty t e^{-\frac{\beta}{n} t} \cdot  \boldsymbol{e}_{c}*\overline{\Lambda}_n(t)\, dt\cr
 \ar=\ar (c+T)e^{-T/c} + \frac{ \theta_n \big\|\overline{\Lambda}_n\big\|_{L^1}}{c}   \int_0^\infty t e^{-\frac{\beta}{n} t} \cdot  \boldsymbol{e}_{c}*\frac{\overline{\Lambda}_n}{\big\|\overline{\Lambda}_n\big\|_{L^1}}(t)\, dt . 
 \eeqlb
 Let $X$ and $Y$ be two independent random variables with density $\boldsymbol{e}_{c}$ and $\big\|\overline{\Lambda}_n\big\|_{L^1}^{-1}\cdot \overline{\Lambda}_n$ respectively. The second term on the right side of equality in \eqref{eqn.610} equals to
 \beqnn
 \frac{ \theta_n \big\|\overline{\Lambda}_n\big\|_{L^1}}{c}\mathbf{E}\Big[(X+Y)e^{-\frac{\beta}{n} (X+Y)} \Big] 
 \ar\leq\ar \frac{ \theta_n \big\|\overline{\Lambda}_n\big\|_{L^1}}{c} \mathbf{E}\big[X\big] +\frac{ \theta_n \big\|\overline{\Lambda}_n\big\|_{L^1}}{c}\mathbf{E}\big[Ye^{-\frac{\beta}{n} Y} \big]\cr
 \ar=\ar \theta_n \big\|\overline{\Lambda}_n\big\|_{L^1}  
 + \frac{ \theta_n }{c}\int_0^\infty te^{-\frac{\beta}{n}t} \overline{\Lambda}_n(t)\, dt .   
 \eeqnn  
 By \eqref{Lambda.n}, the change of variables, the identity $\bar\nu(\eta_n)= c^{-1}\cdot n^2\theta_n$ and then integration by parts,  
 \beqnn 
 \frac{ \theta_n }{c}\int_0^\infty te^{-\frac{\beta}{n}t} \overline{\Lambda}_n(t)\, dt
 \leq \int_0^\infty t e^{- \beta  t}  \overline{\nu}(t)\, dt 
 = \int_0^\infty \int_0^ts e^{- \beta s} \, ds \, \nu(dt)
 \leq \int_0^\infty  \Big( \frac{t^2}{2} \wedge \frac{1}{\beta^2} \Big)\, \nu(dt).
 \eeqnn 
 Taking these estimates back into the right side of the equality in \eqref{eqn.610} yields that 
 \beqlb\label{eqn.413}
 \int_T^\infty t e^{-\frac{\beta}{n} t} \gamma_n \cdot \overline\varPi_n(t)\, dt 
 \ar\leq\ar 
 (c+T)e^{-T/c}
 + \theta_n \big\|\overline{\Lambda}_n\big\|_{L^1}
 + \int_0^\infty  \Big( \frac{t^2}{2} \wedge \frac{1}{\beta^2} \Big)\, \nu(dt) \to 0,
 \eeqlb
 as $n\to \infty$ and then   $\beta,T\to \infty$ by using the dominated convergence theorem along with \eqref{LevyTriplet}, \eqref{eqn.540} and Proposition~\ref{Prop.502}.   
 Hence the desired uniform upper and lower bounds hold.
 \qed 
 
 \begin{proposition}\label{Prop.LowerBound}
 For each $\beta>\beta_0$, there exist constants $n_0\geq 1$ and $C\geq 0$ such that for any $z\in \mathbb{R}$, 
 \beqlb\label{LowerBound}
  \inf_{n\geq n_0} \bigg|1-\int_0^\infty e^{\mathtt{i}zt -\frac{\beta}{n}  t} \cdot \gamma_n \overline\varPi_n(t)\, dt\bigg| \geq C  \cdot \big( |z|\wedge 1 \big).
 \eeqlb 
 \end{proposition}
 \proof  
 Note that $\cos(x)\geq 1/2$ for any $|x|\leq 1$. 
 The two bound estimates in \eqref{eqn.412} yield that
 \beqnn
 \frac{\partial}{\partial z}\int_0^\infty \sin(z t) e^{-\frac{\beta}{n} t}  \cdot \gamma_n   \overline\varPi_n(t)\,dt  \ar=\ar \int_0^\infty \cos(z t)\cdot t e^{-\frac{\beta}{n} t} \cdot \gamma_n  \overline\varPi_n(t)\,dt  \cr
 \ar\geq\ar \frac{1}{2} \int_0^{T_0}  t e^{-\frac{\beta}{n} t}  \cdot\gamma_n   \overline\varPi_n(t)\,dt -\int_{T_0}^\infty t e^{-\frac{\beta}{n} t}  \cdot \gamma_n   \overline\varPi_n(t)\,dt,
 \eeqnn 
 which is larger or equal to $c/4$ for $|z|\leq 1/T_0$. 
 This along with the mean value theorem induces that
 \beqnn
  \bigg|1-\int_0^\infty e^{\mathtt{i}zt -\frac{\beta}{n} t} \cdot \gamma_n   \overline\varPi_n(t)\, dt \bigg|  
  \geq  \bigg|\int_0^\infty \sin(z t) e^{-\frac{\beta}{n}  t} \cdot \gamma_n   \overline\varPi_n(t)\,dt\bigg|
  \geq \frac{c}{4} \cdot |z|.
 \eeqnn
 Hence the desired inequality \eqref{LowerBound} holds for any $|z|\leq 1/T_0$. 
 \medskip
 
 For $|z|>1/T_0$, from Proposition~\ref{Thm.UpperBound}, there exists a constant $z_0>0$ such that for any $|z|\geq z_0$,
 \beqnn
 \sup_{n\geq 1} \bigg|\int_0^\infty e^{\mathtt{i}zt -\frac{\beta}{n} t} \cdot \gamma_n   \overline\varPi_n(t)\, dt \bigg|\leq \frac{1}{2}
 \quad \mbox{and hence}\quad 
 \inf_{n\geq 1} \bigg|1-\int_0^\infty e^{\mathtt{i}zt -\frac{\beta}{n}  t} \cdot  \gamma_n  \overline\varPi_n(t)\, dt \bigg|\geq \frac{1}{2}. 
 \eeqnn 
 Obviously, the desired inequality \eqref{LowerBound} holds for any $z\in\mathbb{R}$ if  $1/T_0\geq z_0$ and then the proof ends. 
 It remains to consider the case $1/T_0<z_0$. 
 Actually, it suffices to prove that
 \beqlb\label{eqn.415}
 \inf_{|z|>1/T_0}\inf_{n\geq n_0} \Big|1-\int_0^\infty e^{\mathtt{i}zt -\frac{\beta}{n} t}  \cdot \gamma_n \overline\varPi_n(t)\, dt\Big| \geq C_0,
 \eeqlb
  for some constant $C_0>0$. 
 Notice that
 \beqnn
 \bigg|1-\int_0^\infty e^{\mathtt{i}zt -\frac{\beta}{n}  t} \cdot \gamma_n  \overline\varPi_n(t)\, dt\bigg|  \geq 1-  F_n(z)
 \quad \mbox{with}\quad 
  F_n(z) = \int_0^\infty \cos(z t)e^{-\frac{\beta}{n} t} \cdot \gamma_n \overline\varPi_n(t)\, dt. 
 \eeqnn
  The continuity of $F_n$ guarantees that its local maximum on the interval $[1/T_0,z_0]$ can be attained at some point $z_n \in [1/T_0,z_0]$.
 For any $T>0$, the fact that $\cos(x)\leq 1$ induces that
 \beqlb\label{UpperboundF} 
 F_n(z_n)
 \ar\leq\ar \int_0^T \cos(z_n t)e^{-\frac{\beta}{n} t} \cdot \gamma_n \overline\varPi_n(t)\, dt + \int_T^\infty e^{-\frac{\beta}{n} t} \cdot \gamma_n \overline\varPi_n(t)\, dt.
 \eeqlb
 Using \eqref{eqn.413} again, we can choose $T>0$ large enough such that
 \beqnn
 \sup_{n\geq n_0} \int_T^\infty  e^{-\frac{\beta}{n} t} \cdot \gamma_n \overline\varPi_n(t)\, dt \leq \frac{1}{T} \cdot  \sup_{n\geq n_0} \int_T^\infty  t e^{-\frac{\beta}{n} t} \cdot \gamma_n \overline\varPi_n(t)\, dt \leq \frac{1}{2}
 \eeqnn
 and hence 
 \beqnn
 \inf_{n\geq n_0}   \int_0^T  e^{-\frac{\beta}{n} t} \cdot \gamma_n \overline\varPi_n(t)\, dt \geq \frac{1}{2}.
 \eeqnn
 By the periodicity of $\cos(z_nt)$, we have
 \beqlb\label{Sum}
 \int_0^T \cos(z_n t) e^{-\frac{\beta}{n} t} \cdot \gamma_n \overline\varPi_n(t)\, dt
 \ar\leq\ar \sum_{k=0}^{ [Tz_n /(2\pi)] } \int_{(2k\pi-\pi/2)/z_n}^{{(2k\pi+\pi/2)/z_n}} \cos(z_n t) e^{-\frac{\beta}{n} t} \cdot \gamma_n \overline\varPi_n(t)\, dt.
 \eeqlb
 We now start to analyze the maximum of the sum above.
 Notice that $\cos(z_nt)$ is unimodal on each interval $[(2k\pi-\pi/2)/z_n, (2k\pi+\pi/2)/z_n]$ for any $k\geq 0$ with the maximum arrived at the point $2k\pi /z_n$.
 Thus the more weight of $e^{-\frac{\beta}{n}  t} \gamma_n \cdot \overline\varPi_n(t)$ is distributed around the local maximum points, the larger the sum above will be.
 To obtain the maximum of the summation in (\ref{Sum}) we should split the weight of $\int_0^T e^{-\frac{\beta}{n} t} \gamma_n \cdot \overline\varPi_n(t)dt$ uniformly around these maximum points.
 More precisely, we can choose $T>0$ large enough such that
 \beqnn
 R_{z_n}:= \frac{z_n\int_0^T e^{-\frac{\beta}{n} t} \cdot \gamma_n \overline\varPi_n(t)\, dt}{2\gamma_n\cdot ([Tz_n /(2\pi)]+1)}<1.
 \eeqnn
 From the previous observation and the fact that $\cos(z_n t) \leq 1$ and $  \overline\varPi_n(t)\leq 1$, we have for any $k\geq 0$,
 \beqnn
 \int_{(2k\pi-\pi/2)/z_n}^{{(2k\pi+\pi/2)/z_n}} \cos(z_n t)e^{-\frac{\beta}{n} t}\cdot \gamma_n  \overline\varPi_n(t)\, dt
 \ar\leq \ar\gamma_n \int_{(2k\pi-R_{z_n})/z_n}^{{(2k\pi+R_{z_n})/z_n}} \cos(z_n t)\, dt \cr
 \ar=\ar  \frac{ \int_0^T e^{-\frac{\beta}{n} t} \cdot \gamma_n \overline\varPi_n(t)dt}{ [Tz_n /(2\pi)]+1}\cdot\frac{\sin(R_{z_n})}{R_{z_n}}.
 \eeqnn
 Taking this back into (\ref{Sum}) and then (\ref{UpperboundF}), we have
 \beqnn
 F_n(z_n) \leq \frac{\sin(R_{z_n})}{R_{z_n}} \int_0^T e^{-\frac{\beta}{n} t} \cdot \gamma_n \overline\varPi_n(t)\, dt  + \int_T^\infty e^{-\frac{\beta}{n} t} \cdot \gamma_n \overline\varPi_n(t)\, dt
 \eeqnn 
 and hence
 \beqnn
 \inf_{|z|>1/T_0}  \Big|1-\int_0^\infty e^{\mathtt{i}zt}e^{-\frac{\beta}{n} t} \gamma_n \cdot \overline\varPi_n(t)\, dt\Big| 
 \ar\geq\ar 1- \int_0^\infty e^{-\frac{\beta}{n} t} \cdot \gamma_n \overline\varPi_n(t)\, dt \cr
 \ar\ar +  \Big(1-\frac{\sin(R_{z_n})}{R_{z_n}} \Big)\cdot \int_0^T e^{-\frac{\beta}{n} t} \cdot \gamma_n \overline\varPi_n(t)\, dt \cr
 \ar\geq\ar  \frac{1}{2}\Big(1-\frac{\sin(R_{z_n})}{R_{z_n}} \Big) .
 \eeqnn
 By the fact that $z_n\in[1/T_0,z_0]$ and $\gamma_n<1$ for any $n\geq n_0$, we  have 
 \beqnn
 \inf_{n\geq n_0}R_{z_n}>0
 \quad\mbox{and hence}\quad
 \sup_{n\geq n_0}\frac{\sin(R_{z_n})}{R_{z_n}}<1 .
 \eeqnn 
 Consequently, the uniform lower bound \eqref{eqn.415} holds. 
 \qed
 
 \begin{lemma}\label{Lemma.509} 
  For any $\beta \geq\beta_0$ and $y> 0$, we have as $n\to\infty$, 
 \beqnn
   \int_0^\infty e^{-2\beta t}\big|R_{\varPi_n}(nt) - W'(t)\big|^2dt 
  + \int_0^\infty e^{-2\beta t}\big|R_n(nt,y) - y\cdot W'(t)\big|^2dt \to 0
 \eeqnn
 \end{lemma}
 \proof Here we just prove the convergence to $0$ of the first integral. 
 The second one can be proved in the same way. 
 We first give some upper bounds for the Fourier transform of $R^{(n)}_{\varPi}(n\cdot)$.  
 By \eqref{eqn.515}, 
 \beqnn
 \sup_{z\in\mathbb{R}}\bigg|\int_0^\infty e^{(\mathrm{i}z  -\beta) t} R_{\varPi_n}(nt)\,dt \bigg| 
 \leq \int_0^{\infty} e^{-\beta t} R_{\varPi_n}(nt)\, dt 
 \to \frac{\beta}{\varPhi(\beta)} >0,
 \eeqnn 
 as $n\to\infty$, which immediately yields that
 \beqnn
 \sup_{n\geq 1}\sup_{z\in\mathbb{R}}\Big|\int_0^\infty e^{(\mathrm{i}z  -\beta)  t} R_{\varPi_n}(nt)\,dt\Big| <\infty. 
 \eeqnn
 Similarly as in the proof of Lemma~\ref{Lemma.ConR}, we also have 
 \beqnn  
 \int_0^\infty e^{(\mathrm{i}z  -\beta)  t} R_{\varPi_n}(t)\, dt 
 =\frac{ \int_0^{\infty} e^{(\mathrm{i}z  -\beta) t} \overline{\varPi}_n (t) \, dt }{ 1- \int_0^{\infty} e^{(\mathrm{i}z  -\beta) t} \cdot \gamma_n\overline{\varPi}_n (t) \, dt},\quad z\in\mathbb{R}.
 \eeqnn
 By the change of variables and then using Proposition~\ref{Thm.UpperBound} and \ref{Prop.LowerBound},
 there exists a constant $C>0$ such that for any  $z\in\mathbb{R}$ and large $n\geq 1$, 
 \beqnn
 \bigg|\int_0^\infty e^{(\mathrm{i}z  -\beta) t} R_{\varPi_n}(nt)\,dt\bigg|
 = \frac{ \big|\int_0^{\infty} e^{(\mathrm{i}z  -\beta)\cdot \frac{t}{n}} \overline{\varPi}_n (t) \,dt\big| }{n\cdot\big| 1- \int_0^{\infty} e^{(\mathrm{i}z  -\beta)\cdot \frac{t}{n}}\cdot \gamma_n \overline{\varPi}_n (t)\, dt\big|} 
 \leq C \cdot \frac{\frac{n}{|z|}\wedge 1}{|z|\wedge n}= \frac{C}{|z|}.
 \eeqnn
 Putting these estimates together, there exist two constants constants $C >0$ and $n_0\geq 1$ such that
 \beqlb\label{UpperBoundRii}
 \sup_{n\geq n_0} \Big|\int_0^\infty e^{(\mathrm{i}z  -\beta) t} R_{\varPi_n}(nt)dt\Big|\ar \leq\ar C\Big(\frac{1}{|z|}\wedge 1\Big),
 \eeqlb
 uniformly in $z\in\mathbb{R}$. 
 By using the Fourier isometry along with the square integrability of  $e^{-\beta t}\cdot R_{\varPi_n}(nt)$ and $e^{-\beta t}\cdot  W'(t)$ as well as \eqref{eqn.2002}, 
 \beqnn
 \int_0^\infty e^{-2\beta t}\big|R_{\varPi_n}(nt) - W'(t)\big|^2dt   
 = \int_\mathbb{R} \bigg|\int_0^\infty e^{(\mathrm{i}z  -\beta) t} R_{\varPi_n}(nt)dt- \int_0^\infty e^{(\mathrm{i}z  -\beta) t} W'(t)dt \bigg|^2 dz,
 \eeqnn 
 which goes to $0$ as $n\to\infty$ by the dominated convergence theorem and  \eqref{eqn.515} with $\lambda=\beta - \mathrm{i}z  $.  
 \qed
 
 \begin{corollary}\label{Coro.410}
 For any $T\geq 0$ and $y>0$, we have as $n\to\infty$, 
 \beqlb\label{eqn.546} 
 	\int_0^T \big|R_{\varPi_n}(nt) - W'(t)\big|^2\,dt + \int_0^T\big|R_n(nt,y) - y\cdot W'(t)\big|^2\, dt \to 0. 
 \eeqlb 
 
 \end{corollary}
   
 \begin{lemma}\label{Lemma.510}
 For any $\beta\geq \beta_0$, there exist two constants $C>0$ and $n_0\geq 1$ such that for any  $\delta \geq 0$,
 \beqnn
  \sup_{n\geq n_0}\int_\mathbb{R}  \big|  e^{-\beta(t+\delta)} R_{\varPi_n}\big(n(t+\delta)\big) - e^{-\beta t} R_{\varPi_n}(nt)\big|^2\, dt \leq C\cdot \delta . 
 \eeqnn 
 \end{lemma}
 \proof 
 By the Fourier isometry, the square integrability of  $e^{-\beta t}\cdot R_{\varPi_n}(nt)$ and then the change of variables, 
 \beqnn
  \int_\mathbb{R}  \big|  e^{-\beta(t+\delta)} R_{\varPi_n}\big(n(t+\delta)\big) - e^{-\beta t} R_{\varPi_n}(nt)\big|^2\, dt
  =\int_\mathbb{R} \Big|\big(e^{\mathrm{i}z\delta}-1\big)\int_\mathbb{R}  e^{(\mathrm{i}z-\beta)t} R_{\varPi_n}(nt)\, dt   \Big|^2 dz . 
 \eeqnn 
 An application of \eqref{UpperBoundRii} and the inequality $\big|e^{\mathrm{i}z}-1\big| \leq |z|\wedge 2$ for all $z\in\mathbb{R}$ to the right-hand side yields that 
 \beqnn
 \sup_{n\geq n_0}\int_\mathbb{R}  \big|  e^{-\beta(t+\delta)} R_{\varPi_n}\big(n(t+\delta)\big) - e^{-\beta t} R_{\varPi_n}(nt)\big|^2\, dt 
 \ar\leq\ar  C\cdot \int_\mathbb{R}\big((z\delta)^2\wedge 1\big)   \Big( \frac{1}{z^2}\wedge 1 \Big) dz 
 \leq  C\cdot \delta , 
 \eeqnn
 uniformly in $\delta>0$.
 \qed

  The next corollary can be proved similarly as Corollary~\ref{Coro.51101} by using the preceding lemma.
  
 \begin{corollary}\label{Coro.511}
 For each $T\geq 0$,  there exist two constants $C>0$ and $n_0\geq 1$ such that for any  $\delta \in (0,1)$,
 \beqnn
  \sup_{n\geq n_0}\int_{-\infty}^T \big|  R_{\varPi_n}\big(n(t+\delta)\big) -  R_{\varPi_n}(nt)\big|^2\, dt  
  \leq C\cdot \delta.
 \eeqnn
 \end{corollary}
%

 \subsubsection{Uniform moment estimates}
  
 In this section, we mainly establish a uniform upper bound for all moments of the sequence $\{ X^{(n)}_\zeta \}_{n\geq 1}$ with the help of the following moment results for the measures $\Lambda_n(dy)$ and $\boldsymbol{e}_c* d\Lambda_n(dy)$ as well as the random variable $\ell^\Lambda_{n,1}$. 
  
 \begin{proposition}\label{Proposition.510}
 For each $p\geq 2$, there exists a constant $C>0$ independent of $n $ and $t $ such that 
 \beqlb
  n^2\theta_n  \int_0^\infty    \Big(t\wedge \frac{y}{n}  \Big)^p    \Lambda_n(dy) 
  \ar\leq\ar C\cdot (1+t)^{p-1}, \label{eqn.543}\\
  n^2\theta_n   \int_0^\infty  \Big(t\wedge \frac{y}{n}  \Big)^p  \boldsymbol{e}_{c}*d\Lambda_n(y)\, dy 
  \ar\leq\ar C\cdot \Big(\frac{\theta_n}{n^{p-2}}  +  (1+t)^{p-1} \Big). \label{eqn.544}
 \eeqlb 
 \end{proposition}
 \proof 
 By \eqref{Lambda.n} and then the change of variables and Fubini's theorem, we first have 
 \beqnn
  n^2\theta_n  \int_0^\infty    \Big(t\wedge \frac{y}{n}  \Big)^p  \,   \Lambda_n(dy)
  \ar=\ar c  \int_0^\infty    \Big(t\wedge \frac{y}{n}  \Big)^p \,  \nu\Big(\eta_n+\frac{dy}{n}\Big) \cr
  \ar=\ar c \int_0^\infty   (t\wedge y)^p   \,   \nu (\eta_n+ dy  )\cr
  \ar=\ar cp \int_0^t  y^{p-1} \cdot   \bar\nu (\eta_n+ y  )\,dy\cr
  \ar\leq\ar cp \int_0^t  y^{p-1} \cdot   \bar\nu (y)\,dy,
 \eeqnn
 uniformly in $n\geq1$ and $t\geq 0$, which, by \eqref{LevyTriplet}, can be uniformly bounded by
 \beqnn
  C \Big(\int_0^1 y^{p-1}   \cdot   \bar\nu (y)dy  +\int_1^{1\vee t} y^{p-1}   \cdot   \bar\nu (y)dy \Big)
  \leq C  \Big(1  +t^{p-1}\cdot \int_1^\infty   \bar\nu (y)dy \Big)
  \leq  C \cdot (1+t)^{p-1}
 \eeqnn
 and then \eqref{eqn.543} holds. 
 By the inequality $(x+y)\wedge t \leq x+ y\wedge t$ for  $x,y\geq 0$ and the power inequality, 
 \beqnn
  n^2\theta_n   \int_0^\infty  \Big(t\wedge \frac{y}{n}  \Big)^p  \cdot  \boldsymbol{e}_{c}*d\Lambda_n(y)\, dy
  \ar=\ar n^2\theta_n   \int_0^\infty   \int_0^\infty  \Big(t\wedge \frac{x+y}{n}  \Big)^p   \boldsymbol{e}_{c}(x)\, dx \Lambda_n(dy) \cr
  \ar\leq\ar C\cdot n^2\theta_n   \int_0^\infty  \Big(  \frac{x}{n}  \Big)^p   \boldsymbol{e}_{c}(x)\, dx
   +  C\cdot n^2\theta_n     \int_0^\infty  \Big(t\wedge \frac{y}{n}  \Big)^p  \,  \Lambda_n(dy),
 \eeqnn
 which along with the first result induces \eqref{eqn.544} immediately. 
 \qed 
  
 \begin{corollary}\label{Corollary.511}
 For any $p\geq 1$,  there exists a constant $C>0$ independent of $n $ and $t $ such that 
 \beqnn
  \int_0^\infty  \big|R_n(nt,y)\big|^{p} \cdot \boldsymbol{e}_{c}(y)\, dy   \leq C 
  \quad \mbox{and}\quad 
  \frac{\theta_n}{n^{p-2}}  \int_0^\infty  \big|R_n(nt,y)\big|^{p} \cdot   \boldsymbol{e}_{c}*d\Lambda_n(y)\, dy 
  \leq C\cdot (1+t)^p. 
 \eeqnn 
 \end{corollary} 
 \proof 
 By using the second inequality in \eqref{eqn.509} and then the power inequality,  
 \beqnn
  \sup_{n\geq 1}\sup_{t\geq 0} \int_0^\infty  \big|R_n(nt,y)\big|^{p} \cdot \boldsymbol{e}_{c}(y)\, dy  
  \ar\leq\ar C \cdot  \int_0^\infty   \boldsymbol{e}_{c}(y)\, dy + C \cdot  \int_0^\infty y^p  \cdot \boldsymbol{e}_{c}(y)\, dy <\infty 
 \eeqnn
 and the first desired inequality holds. 
 Similarly, the second one can be proved by using \eqref{eqn.544}, i.e.,
 \beqlb\label{eqn.545}
  \frac{\theta_n}{n^{p-2}}  \int_0^\infty  \big|R_n(nt,y)\big|^{p} \cdot   \boldsymbol{e}_{c}*d\Lambda_n(y)\, dy 
  \ar\leq\ar  \frac{\theta_n}{n^{p-2}}  \int_0^\infty  \Big(1+ \frac{(nt)\wedge y}{c}   \Big)^p  \cdot  \boldsymbol{e}_{c}*d\Lambda_n(y)\, dy\cr
  \ar\leq\ar C \cdot \frac{\theta_n}{n^{p-2}}    + C\cdot n^2\theta_n  \int_0^\infty    \Big(t\wedge  \frac{y}{n}\Big)^p    \cdot  \boldsymbol{e}_{c}*d\Lambda_n(y)\, dy , 
 \eeqlb
 which can be bounded by $C\cdot(1+t)^p$ uniformly in $n\geq 1$ and $t\geq 0$.  
 \qed
  
 \begin{proposition} \label{Prop.503}
 For each $p\geq 1$, there exists a constant $C>0$ independent of $n $ and $t $ such that 
 \beqnn
  \mathbf{E}\bigg[\Big| \frac{\ell^\Lambda_{n,1}}{n} \wedge t\Big|^p \bigg]  \leq C\cdot \bigg(\frac{1}{n^p}+ \frac{t^{p-1}}{n(1-p_n)} \cdot \int_0^\infty   (y\wedge t) \cdot  \bar\nu(y) \, dy    \bigg)   .
 \eeqnn  
 \end{proposition}
 \proof  
 Similarly as in \eqref{eqn.545}, we have uniformly in $n\geq 1$ and $t\geq 0$,
 \beqnn
  \mathbf{E}\bigg[\Big| \frac{\ell^\Lambda_{n,1}}{n} \wedge t\Big|^p \bigg] 
  \ar\leq\ar C\cdot \int_0^\infty \Big( \frac{y}{n} \Big)^p \boldsymbol{e}_{c}(y)\, dy
  + C\cdot \int_0^\infty \Big( \frac{y}{n} \wedge t \Big)^p \boldsymbol{\lambda}^*_n(y)\, dy.
 \eeqnn
 The first term on right side of this equality can be uniformly bounded by $C\cdot n^{-p}$.
 For the second term, by using \eqref{eqn.541} and then the change of variables and \eqref{eqn.pn} we have 
 \beqnn
  \int_0^\infty \Big( \frac{y}{n} \wedge t \Big)^p \boldsymbol{\lambda}^*_n(y)\, dy
  \ar=\ar \frac{1}{\int_0^\infty \bar\nu(\eta_n+z) \,dz}\int_0^\infty    ( y \wedge t )^p \cdot \bar\nu(\eta_n+y) dy \cr
  \ar=\ar \frac{p_n}{n(1-p_n)}\int_0^\infty   ( y \wedge t )^p \cdot \bar\nu(\eta_n+y) \, dy\cr
  \ar\leq\ar  
  \frac{1}{n(1-p_n)}  \int_0^\infty   (y\wedge t)^p \cdot \bar\nu(y) \, dy \cr
  \ar\leq\ar \frac{t^{p-1}}{n(1-p_n)}  \int_0^\infty   (y\wedge t) \cdot \bar\nu(y) \, dy . 
 \eeqnn
 Here the last integral is finite because of \eqref{LevyTriplet}. 
 The desired inequality holds. 
 \qed

 \begin{corollary}\label{Prop.BoundR}
 For each $p\geq 1$, there exists a constant $C>0$ independent of $n$ and $t$ such that
 \beqnn
  \mathbf{E}\Big[ \big| R_n(nt,\ell^c_{n,1})\big|^p \Big]\leq C
  \quad \mbox{and}\quad 
  \mathbf{E}\Big[ \big| R_n(nt,\ell^\Lambda_{n,1})\big|^p \Big]\leq C\cdot \frac{n^{p-1}}{1-p_n}\cdot (1+t)^p . 
 \eeqnn
 \end{corollary}
 \proof 
 By \eqref{eqn.509} and the power inequality, we have uniformly in $n\geq 1$ and $t\geq 0$, 
 \beqnn
  \mathbf{E}\Big[ \big| R_n(nt,\ell^c_{n,1})\big|^p \Big] 
  \leq C\cdot \Big( 1+\mathbf{E}\Big[\big| \ell^c_{n,1} \big|^p \Big] \Big) 
 \quad \mbox{and}\quad 
  \mathbf{E}\Big[ \big| R_n(nt,\ell^\Lambda_{n,1})\big|^p \Big] 
  \leq C \Big( 1+\mathbf{E}\Big[\big| \ell^\Lambda_{n,1} \wedge (nt)\big|^p \Big] \Big). 
 \eeqnn
 The two desired inequalities follow respectively from 
 \beqnn 
  \mathbf{E}\Big[ \big| \ell^c_{n,1} \big|^p \Big]<\infty
  \quad \mbox{and}\quad
  \mathbf{E}\Big[\big| \ell^\Lambda_{n,1} \wedge (nt)\big|^p \Big] \leq C\cdot \frac{n^{p-1}}{1-p_n}\cdot (1+t)^p,
 \eeqnn
 uniformly in $n\geq 1$ and $t\geq 0$;  
 see Proposition~\ref{Prop.503}. 
 \qed 
 
 \begin{proposition} \label{Prop.I}
 For each $p\geq 1$, there exists a constant $C>0$ such that for any $\zeta, t\geq 0$,
 \beqlb\label{eqn.416}
  \sup_{n\geq 1}  \mathbf{E}\Big[\big|\mathbf{I}^{(n)}_1(t)\big|^p  \Big]  \leq C \cdot \zeta^p
  \quad \mbox{and}\quad 
  \sup_{n\geq 1}	\mathbf{E}\Big[\big|\mathbf{I}^{(n)}_2(t)\big|^p  \Big] \leq C\cdot (\zeta \vee \zeta^p) \cdot \big(1+t\big)^{p}. 
 \eeqlb
 \end{proposition}
 \proof 
 By H\"older's inequality, it suffices to prove this  proposition for $p=k \in \mathbb{N}$. By using the fact $B_n^c(\zeta)\leq [n\zeta]$ a.s., Jensen's inequality and then the mutual-independence among $ \{\ell_{n,j}^c\}_{j\geq 1}$ 
 \beqnn
  \sup_{n\geq 1} \mathbf{E}\Big[\big|\mathbf{I}^{(n)}_1(t)\big|^k  \Big]
  \ar\leq\ar \zeta^k\cdot \sup_{n\geq 1}\cdot \mathbf{E}\Bigg[  \bigg| \frac{1}{[n\zeta]} \sum_{j=1}^{[n\zeta]} R_n(nt,\ell_{n,j}^c)\bigg|^k \Bigg]\cr
  \ar\leq \ar \zeta^k\cdot \sup_{n\geq 1}\cdot  \frac{1}{[n\zeta]} \sum_{j=1}^{[n\zeta]} \mathbf{E}\Big[  \big|R_n(nt,\ell_{n,j}^c)\big|^k \Big]\cr
  \ar=\ar \zeta^k\cdot \sup_{n\geq 1}\cdot    \mathbf{E}\Big[  \big|R_n(nt,\ell_{n,1}^c)\big|^k \Big], 
 \eeqnn
 which is bounded by $C\cdot \zeta^k$ uniformly in $t\geq 0$; see Corollary~\ref{Prop.BoundR}, and hence the first desired uniform upper bound holds.  
 For the second one, without loss of generality we assume that $[n\zeta]>k$.  
 By repeating the binomial expansion,
 \beqnn
  \bigg|  \sum_{j=1}^{B_n^\Lambda(\zeta)} R_n(nt,\ell_{n,j}^\Lambda)\bigg|^k 
  \ar=\ar \sum_{k_1=0}^k \sum_{k_2=0}^{k-k_1} \cdots \sum_{k_{B_n^\Lambda(\zeta)}=0}^{k-\sum_{i=1}^{B_n^\Lambda(\zeta)-1}k_i }  
  \prod_{j=1}^{B_n^\Lambda(\zeta)}{{k-\sum_{i=1}^{j-1}k_i }\choose{k_j}}  \Big| R_n(nt,\ell_{n,j}^\Lambda)\Big|^{k_j} .
 \eeqnn
 Taking expectations on  both sides of this equality conditionally on $B_n^\Lambda(\zeta)$ and then using Corollary~\ref{Prop.BoundR}, there exists a constant $C>0$ independent of $n$ and $t$ such that 
 \beqnn
  \lefteqn{ \mathbf{E}\Bigg[\bigg|  \sum_{j=1}^{B_n^\Lambda(\zeta)} R_n(nt,\ell_{n,j}^\Lambda)\bigg|^k\,\Bigg|\, B_n^\Lambda(\zeta)\Bigg] }\ar\ar\cr
  \ar=\ar \sum_{k_1=0}^k \sum_{k_2=0}^{k-k_1} \cdots \sum_{k_{B_n^\Lambda(\zeta)}=0}^{k-\sum_{i=1}^{B_n^\Lambda(\zeta)-1}k_i } 
  \prod_{j=1}^{B_n^\Lambda(\zeta)} {{k-\sum_{i=1}^{j-1}k_i }\choose{k_j}} \mathbf{E}\Big[  \big| R_n(nt,\ell_{n,j}^\Lambda)\big|^{k_j} \Big] \cr
  \ar\leq\ar C \cdot \sum_{k_1=0}^k \sum_{k_2=0}^{k-k_1} \cdots \sum_{k_{B_n^\Lambda(\zeta)}=0}^{k-\sum_{i=1}^{B_n^\Lambda(\zeta)-1}k_i }  
  \prod_{j=1}^{B_n^\Lambda(\zeta)} {{k-\sum_{i=1}^{j-1}k_i }\choose{k_j}} \frac{\big(n(1-p_n)\big)^{(k_j-1)^+}}{|1-p_n|^{k_j}}(1+t)^{k_j}\cr
  \ar=\ar  C \cdot  \sum_{k_1=0}^k \sum_{k_2=0}^{k-k_1} \cdots \sum_{k_{B_n^\Lambda(\zeta)}=0}^{k-\sum_{i=1}^{B_n^\Lambda(\zeta)-1}k_i }  
  \cdot \frac{\big(n(1-p_n)\big)^{k-J}}{|1-p_n|^k} \cdot (1+t)^k \cdot \prod_{j=1}^{B_n^\Lambda(\zeta)} {{k-\sum_{i=1}^{j-1}k_i }\choose{k_j}}
 \eeqnn
 with $J:= \# \{1\leq i\leq B_n^\Lambda(\zeta) : k_i\geq 1 \} \leq k $. 
 Note that all products of binomial coefficients above are uniformly bounded by a constant depending only on $k$, we have 
 \beqnn
  \mathbf{E}\Bigg[\bigg|  \sum_{j=1}^{B_n^\Lambda(\zeta)} R_n(nt,\ell_{n,j}^\Lambda)\bigg|^k\,\Bigg|\, B_n^\Lambda(\zeta)\Bigg] 
  \leq C \cdot (1+t)^{k} \cdot  \sum_{k_1=0}^k \sum_{k_2=0}^{k-k_1} \cdots \sum_{k_{B_n^\Lambda(\zeta)}=0}^{k-\sum_{i=1}^{B_n^\Lambda(\zeta)-1}k_i }   \frac{\big(n(1-p_n)\big)^{k-J}}{|1-p_n|^k} .
 \eeqnn
 Using the multinomial distribution and then the combination formula to the last multiple sum,
 \beqnn
  \mathbf{E}\Bigg[\bigg|  \sum_{j=1}^{B_n^\Lambda(\zeta)} R_n(nt,\ell_{n,j}^\Lambda)\bigg|^k\,\Bigg|\, B_n^\Lambda(\zeta)\Bigg]
  \ar\leq\ar C \cdot (1+t)^{k} \cdot  \sum_{j=1}^{k} {{B_n^\Lambda(\zeta)}\choose{j}} \cdot \frac{\big(n(1-p_n)\big)^{k-j}}{|1-p_n|^k} \cr
  \ar\leq\ar C \cdot (1+t)^{k} \cdot  \sum_{j=1}^{k} \big|B_n^\Lambda(\zeta) \big|^j  \cdot \frac{\big(n(1-p_n)\big)^{k-j}}{|1-p_n|^k}.
 \eeqnn
 Taking expectations on both sides of this inequality,  
 \beqnn
  \mathbf{E}\Big[\big|\mathbf{I}^{(n)}_1 (t)\big|^k  \Big]
  \ar\leq\ar  C\cdot (1+t)^{k} \cdot  \sum_{j=1}^{k} \mathbf{E}\Big[\big|B_n^\Lambda(\zeta) \big|^j \Big] \cdot \big(n(1-p_n)\big)^{-j}.
 \eeqnn
 Notice that $n(1-p_n) \to\doublebar{\nu}(0)\in(0,\infty]$ as $n\to\infty$, we have uniformly in $j\leq k$,
 \beqnn
  \mathbf{E}\Big[\big|B_n^\Lambda(\zeta)\big|^j\Big] = \sum_{i=1}^j {{j}\choose{i}} \big( [n\zeta]\cdot (1-p_n) \big)^i \leq C\cdot(\zeta\vee \zeta^j) \cdot \Big(\big(n(1-p_n)\big)\vee\big(n(1-p_n)\big)^j\Big).
 \eeqnn
 The second inequality in (\ref{eqn.416}) follows by combining the preceding two inequality together. 
 \qed 
 
 We are now ready to establish the uniform upper bounds for all moments of the processes $\{ X^{(n)}_\zeta\}_{n\geq1}$ with the help of the preceding results and the moment estimates for stochastic Volterra integrals given in Proposition~\ref{Prop.BDG}. 
  
 \begin{lemma}\label{Lemma.Moment}
 For each $p\geq 1$, there exists a constant $C>0$  such that for any $\zeta,t\geq 0$,
  \beqlb\label{eqn.4188}
   \sup_{n\geq 1}\mathbf{E}\Big[\big| X_\zeta^{(n)}(t) \big|^p  \Big] 
   \leq C\cdot (\zeta \vee \zeta^p) \cdot (1+t)^{2p-1}. 
  \eeqlb 
 \end{lemma}
 \proof 
 For $p=1$, taking expectations on both sides of \eqref{Eqn.HR01}, the non-negativity of $X_\zeta^{(n)}$ and Proposition~\ref{Prop.I} induce that uniformly in $t\geq 0$,
 \beqnn
  \sup_{n\geq 1}\mathbf{E}\Big[\big| X_\zeta^{(n)}(t) \big|  \Big]
  \ar=\ar  \sup_{n\geq 1}\mathbf{E}\Big[\mathbf{I}^{(n)}_1(t)  \Big] 
  +\sup_{n\geq 1}\mathbf{E}\Big[\mathbf{I}^{(n)}_2(t) \Big] 
  \leq C\cdot \zeta \cdot (1+t).
 \eeqnn
 By induction, we now proceed under the assumption that \eqref{eqn.4188} holds for some $p\geq 1$ and prove that it also holds for $2p$. 
 By the power inequality, 
 \beqlb\label{eqn.547}
  \mathbf{E}\Big[\big|X_{\zeta}^{(n)}(t)\big|^{2p}\big]
  \ar\leq\ar C\cdot \sum_{i=1}^4 \mathbf{E}\Big[\big|\mathbf{I}^{(n)}_i(t)\big|^{2p}  \Big] , 
 \eeqlb
 for some constant $C>0$ independent of $n$ and $t$.
 By Proposition~\ref{Prop.I}, 
 \beqlb\label{eqn.549}
  \mathbf{E}\Big[\big|\mathbf{I}^{(n)}_1(t)\big|^{2p}  \Big]
  + \mathbf{E}\Big[\big|\mathbf{I}^{(n)}_2 (t)\big|^{2p}  \Big] 
  \leq C\cdot(\zeta\vee \zeta^{2p})\cdot (1+t)^{2p},
 \eeqlb
 uniformly in $n\geq 1$ and $\zeta,t\geq 0$. 
 Applying \eqref{BDG} to $ \mathbf{E}\big[ |\mathbf{I}^{(n)}_3(t) |^{2p}  \big]$, there exists a constant $C>0$ depending only on $p$ such that  
 \beqnn
  \mathbf{E}\Big[\big|\mathbf{I}^{(n)}_3(t)\big|^{2p}  \Big]
  \ar\leq\ar C \sup_{r\in[0,t]} \mathbf{E}\Big[ \big|X_{\zeta}^{(n)}(r)\big|^p \Big] \cdot  \bigg|\int_0^t  \int_0^\infty  \frac{|R_n\big(ns,y\big)|^2}{n^2}\cdot  n^2\gamma_n(1-\theta_n ) \cdot \boldsymbol{e}_{c}(y)\,  ds\, dy\bigg|^{p} \cr
  \ar\ar + C \sup_{r\in[0,t]} \mathbf{E}\Big[ \big|X_{\zeta}^{(n)}(r)\big|\Big] \cdot   \int_0^t  \int_0^\infty  \frac{|R_n\big(ns,y\big)|^{2p}}{n^{2p}}\cdot  n^2\gamma_n(1-\theta_n ) \cdot \boldsymbol{e}_{c}(y)\,  ds\, dy.
 \eeqnn
 By \eqref{eqn.4188}, Corollary~\ref{Corollary.511} as well as the two facts that $\gamma_n\leq 1/c $ and $\theta_n\in (0,1)$,  
 \beqlb\label{eqn.548}
  \sup_{n\geq 1} \mathbf{E}\Big[\big|\mathbf{I}^{(n)}_3(t)\big|^{2p}  \Big] 
  \ar\leq\ar C\cdot (\zeta \vee \zeta^p)\cdot   (1+t)^{p} + C \cdot \zeta \cdot (1+t) 
  \leq
  C \cdot (\zeta \vee \zeta^p)\cdot   (1+t)^{p} ,
 \eeqlb 
 uniformly in $\zeta,t\geq 0$. 
 Similarly, we also have  
 \beqnn
  \mathbf{E}\Big[\big|\mathbf{I}^{(n)}_4(t)\big|^{2p}  \Big]
  \ar\leq\ar C \sup_{r\in[0,t]} \mathbf{E}\Big[ \big|X_{\zeta}^{(n)}(r)\big|^p \Big] \cdot     \bigg|\int_0^t  \int_0^\infty  \frac{|R_n\big(ns,y\big)|^2}{n^2}\cdot  n^2\gamma_n\theta_n \cdot  \boldsymbol{e}_{c}*d\Lambda_n(y)\,\cdot ds\cdot  dy\bigg|^{p} \cr
  \ar\ar + C \sup_{r\in[0,t]} \mathbf{E}\Big[ \big|X_{\zeta}^{(n)}(r)\big| \Big] \cdot    \int_0^t  \int_0^\infty  \frac{|R_n\big(ns,y\big)|^{2p}}{n^{2p}}\cdot  n^2\gamma_n\theta_n \cdot  \boldsymbol{e}_{c}*d\Lambda_n(y)\,\cdot ds\cdot dy,
 \eeqnn
 which also can be bounded by 
 $
 C\cdot(\zeta \vee \zeta^p)\cdot (1+t)^{4p-1} + C\cdot \zeta \cdot (1+t)^{2p+1}
 \leq C\cdot(\zeta \vee \zeta^p)\cdot (1+t)^{4p-1},
 $
 uniformly in $n\geq 1$ and $\zeta,t\geq 0$ by using \eqref{eqn.4188} and Corollary~\ref{Corollary.511} again. 
 Taking this and \eqref{eqn.549}-\eqref{eqn.548} back into \eqref{eqn.547} induces that
 the inequality \eqref{eqn.4188} holds for $2p$. 
 \qed  
    
 \subsection{$C$-tightness}\label{Sec.CTightness}
 
 In the next four subsections, we will prove the $C$-tightness of four sequences $ \{\mathbf{I}^{(n)}_i \}_{n\geq 1}$, $i=1,2,3,4,$,
 separately, which together with Corollary~3.33 in \cite[p.353]{JacodShiryaev2003} induces the $C$-tightness of the sequence 
 $$\big\{ (X^{(n)}_\zeta,\mathbf{I}^{(n)}_1,\mathbf{I}^{(n)}_2,\mathbf{I}^{(n)}_3,\mathbf{I}^{(n)}_4 ) \big\}_{n\geq 1}. $$ 
 
 \subsubsection{Convergence in probability of $\{\mathbf{I}^{(n)}_1\}_{n\geq 1}$}
 
 Note that $\mathbf{E}\big[\mathbf{1}_{\{\ell_{n,i}^c>t\}} \big]=  \overline{\boldsymbol{E}}_c(t)$ for $t\geq 0$.  
 In view of \eqref{eqn.I},
 centralizing each summand in $\mathbf{I}^{(n)}_1$ allows us to write it into
 \beqlb\label{eqn.550}
  \mathbf{I}^{(n)}_1(t) = \mathbf{I}^{(n)}_{1,1}(nt) +  R_{\varPi_n}* \mathbf{I}^{(n)}_{1,1}(nt) + \frac{B_n^c(\zeta)}{n}\cdot \mathbf{I}^{(n)}_{1,2}(nt) ,\quad t\geq 0,
 \eeqlb
 where 
 \beqnn
  \mathbf{I}^{(n)}_{1,1}(t) := \frac{1}{n} \sum_{i=1}^{B_n^c(\zeta)} \big( \mathbf{1}_{\{\ell_{n,i}^c>t \}} -  \overline{\boldsymbol{E}}_c(t)\big) 
  \quad \mbox{and}\quad  
  \mathbf{I}^{(n)}_{1,2}(t) := \overline{\boldsymbol{E}}_c(t) +R_{\varPi_n}*\overline{\boldsymbol{E}}_c(t).    
 \eeqnn 
 In the next two propositions, we prove the convergence in probability of the three terms on the right side of \eqref{eqn.550} separately. 
 
 \begin{proposition}
 We have $\sup_{t\geq 0}\big| \mathbf{I}^{(n)}_{1,1}(t)\big|+ \big\| \mathbf{I}^{(n)}_{1,1}\big\|_{L^1} + \sup_{t\geq 0}\big| R_{\varPi_n}* \mathbf{I}^{(n)}_{1,1}(t)\big|  \overset{\rm p}\to 0$ as $n\to\infty$.  
 \end{proposition}
 \proof 
 We first prove $\sup_{t\geq 0}\big| \mathbf{I}^{(n)}_{1,1}(t)\big|\to 0$ in probability by using the Glivenko-Cantelli theorem. 
 Notice that random variables $ \ell_{n,i}^c $, $i,n\geq 1$, are i.i.d., without loss of generality we may omit the subcript $n$ and write $\{ \ell_{n,i}^c\}_{i\geq 1}$ into $\{ \ell_{i}^c\}_{i\geq 1}$. 
 For each $\epsilon>0$, we have 
 \beqlb\label{eqn.519}
  \mathbf{P}\bigg( \sup_{t\geq 0}\big| \mathbf{I}^{(n)}_{1,1}(t)\big| \geq \epsilon \bigg)  
  \ar=\ar \mathbf{P}\bigg(  \sup_{t\geq 0}\big| \mathbf{I}^{(n)}_{1,1}(t)\big| \geq \epsilon ; B_n^c(\zeta)< [n\zeta]\cdot p_n -[n\zeta]^{2/3}\bigg)\cr
  \ar\ar +\mathbf{P}\bigg( \sup_{t\geq 0}\big| \mathbf{I}^{(n)}_{1,1}(t)\big| \geq \epsilon ; B_n^c(\zeta)\geq [n\zeta]\cdot p_n -[n\zeta]^{2/3}\bigg).
 \eeqlb
 By the Chernoff inequality for binomial random variables\footnote{For a binomial random variable $X$ with parameter $(k,p)$, we have $\mathbf{P}(X\leq kp-x) \leq \exp\{-\frac{x^2}{2kp}\}$ for any $x\geq 0$.}, the first probability on the right side of \eqref{eqn.519} can be bounded by
 \beqnn
  \mathbf{P}\Big(B_n^c(\zeta)< [n\zeta]\cdot p_n -[n\zeta]^{2/3}\Big)
  \leq \exp\Big\{ - \frac{[n\zeta]^{4/3}}{2[n\zeta]\cdot p_n}\Big\} ,
 \eeqnn
 which goes to $0$ as $n\to\infty$. 
 Further, conditioned on $B_n^c(\zeta)\geq  [n\zeta]\cdot p_n -[n\zeta]^{2/3}$ we have 
 \beqnn
  \sup_{t\geq 0}\big| \mathbf{I}^{(n)}_{1,1}(t)\big|\ar\leq\ar \sup_{t\geq 0} \bigg| \frac{1}{n}  \sum_{i=1}^{[n\zeta]\cdot p_n -[n\zeta]^{2/3}} \big( \mathbf{1}_{\{\ell_{i}^c>t \}} - \overline{\boldsymbol{E}}_c(t)  \big) \bigg| 
  +  \sup_{t\geq 0}\bigg| \frac{1}{n} \sum_{i=[n\zeta]\cdot p_n -[n\zeta]^{2/3}}^{B_n^c(\zeta)} \big( \mathbf{1}_{\{\ell_{i}^c>t \}} - \overline{\boldsymbol{E}}_c(t)  \big) \bigg|\cr
  \ar\leq\ar  \sup_{t\geq 0} \bigg| \frac{1}{n} \sum_{i=1}^{[n\zeta]\cdot p_n -[n\zeta]^{2/3}} \big( \mathbf{1}_{\{\ell_{i}^c>t \}} - \overline{\boldsymbol{E}}_c(t)  \big) \bigg| + \frac{2}{n}\cdot\Big|B_n^c(\zeta)- [n\zeta]\cdot p_n  \Big| +\frac{2\cdot [n\zeta]^{2/3}}{n} . 
 \eeqnn
 Hence for $n\geq 64\cdot \zeta^2/\epsilon^3$, the second probability on the right side of \eqref{eqn.519} can be bounded by
 \beqnn
  \mathbf{P}\Bigg( \bigg| \frac{1}{n} \sum_{i=1}^{[n\zeta]\cdot p_n -[n\zeta]^{2/3}} \big( \mathbf{1}_{\{\ell_{i}^c>t \}} - \overline{\boldsymbol{E}}_c(t)  \big) \bigg| \geq \frac{\epsilon}{2} \Bigg)
  + \mathbf{P}\Bigg( \frac{2}{n} \cdot \Big|B_n^c(\zeta)- [n\zeta]\cdot p_n \Big| \geq \frac{\epsilon}{4} \Bigg)
 \eeqnn
 By the Glivenko-Cantelli theorem and the fact that $[n\zeta]\cdot p_n -[n\zeta]^{2/3} \sim \zeta \cdot n$, the first probability vanishes as $n\to\infty$. 
 By Chebyshev's inequality, 
 \beqlb\label{eqn.660}
  \mathbf{P}\Bigg( \frac{2}{n} \cdot \Big|B_n^c(\zeta)- [n\zeta]\cdot p_n \Big| \geq \frac{\epsilon}{4} \Bigg)
  \leq \frac{32}{\epsilon^2} \cdot \frac{{\rm Var}\big(B_n^c(\zeta)\big)}{n^2} 
  = \frac{32}{\epsilon^2} \cdot \frac{[n\zeta]\cdot p_n (1-p_n)}{n^2}
  \to 0 ,
 \eeqlb 
 as $n\to\infty$.
 Putting these estimates together and then taking them back into \eqref{eqn.519} induce that 
 \beqnn
 \lim_{n\to\infty}\mathbf{P}\bigg(  \sup_{t\geq 0}\big| \mathbf{I}^{(n)}_{1,1}(t)\big| \geq \epsilon \bigg) =0. 
 \eeqnn
 
 We now prove  $\big\| \mathbf{I}^{(n)}_{1,1}\big\|_{L^1} \to 0$ in probability, which along with  Lemma~\ref{Lemma.UpperBoundR} induces that as $n\to\infty$, 
 \beqnn
  \sup_{t\geq 0}\big| R_{\varPi_n}* \mathbf{I}^{(n)}_{1,1}(t)\big| \leq \frac{1}{c} \cdot \big\| \mathbf{I}^{(n)}_{1,1} \big\|_{L^1} \to 0. 
 \eeqnn 
 Similarly as in the previous argument, for each $\epsilon>0 $ we also have
 \beqnn
  \mathbf{P}\Big( \big\| \mathbf{I}^{(n)}_{1,1}\big\|_{L^1}\geq \epsilon \Big)  
  \ar\sim\ar \mathbf{P}\Big(  \big\| \mathbf{I}^{(n)}_{1,1}\big\|_{L^1} \geq \epsilon ; B_n^c(\zeta)\geq [n\zeta]\cdot p_n -[n\zeta]^{2/3} \Big),
 \eeqnn
 as $n\to\infty$.
 It is easy to see that conditioned on $B_n^c(\zeta)\geq  [n\zeta]\cdot p_n -[n\zeta]^{2/3}$,
 \beqnn
  \big\| \mathbf{I}^{(n)}_{1,1}\big\|_{L^1}
  \ar\leq\ar \int_0^\infty \bigg| \frac{1}{n} \sum_{i=1}^{[n\zeta]\cdot p_n -[n\zeta]^{2/3}} \big( \mathbf{1}_{\{\ell_{i}^c>s \}} - \overline{\boldsymbol{E}}_c(s) \big) \bigg|\, ds \cr
  \ar\ar + \int_0^\infty \bigg| \frac{1}{n} \sum_{i=[n\zeta]\cdot p_n -[n\zeta]^{2/3}+1}^{B_n^c(\zeta)} \big( \mathbf{1}_{\{\ell_{i}^c>s \}} - \overline{\boldsymbol{E}}_c(s) \big) \bigg|\, ds .
 \eeqnn
 Here the second integral on the right side can be uniformly bounded by 
 \beqnn
  \frac{1}{n} \sum_{i=[n\zeta]\cdot p_n -[n\zeta]^{2/3}+1}^{[n\zeta] }    \ell_{i}^c  + c\cdot \frac{|B_n^c(\zeta)-[n\zeta]\cdot p_n|}{n} + c\cdot \frac{[n\zeta]^{2/3}}{n}.
 \eeqnn
 The last two terms vanish in probability as $n\to\infty$; see \eqref{eqn.660}. 
 Moreover, an application of the law of large number to the first term induces that almost surely 
 \beqnn
  \frac{1}{n} \sum_{i=[n\zeta]\cdot p_n -[n\zeta]^{2/3}+1}^{[n\zeta] }    \ell_{i}^c \sim \frac{[n\zeta](1-p_n)+[n\zeta]^{2/3} }{n} \cdot \mathbf{E}\big[  \ell_1^c \big] \to 0,
 \eeqnn
 as $n\to\infty$. 
 In conclusion, we can obtain that $ \| \mathbf{I}^{(n)}_{1,1} \|_{L^1} \overset{\rm p}\to 0   $ as $n\to\infty$ if 
 \beqlb\label{eqn.618}
  \int_0^\infty \bigg| \frac{1}{n} \sum_{i=1}^{[n\zeta]\cdot p_n -[n\zeta]^{2/3}} \big( \mathbf{1}_{\{\ell_{i}^c>s \}} - \overline{\boldsymbol{E}}_c(s)  \big) \bigg|\, ds \overset{\rm p}\to 0 .
 \eeqlb
 To prove this convergence in probability, we define a sequence of i.i.d. $L^1(\mathbb{R}_+;\mathbb{R})$-valued random variables $\{X_i\}_{i\geq 1}$ with 
 \beqnn
  X_i(t):=\mathbf{1}_{\{\ell_{i}^c>t \}} -\overline{\boldsymbol{E}}_c(t),\quad t\geq 0.
 \eeqnn 
 Note that $\mathbf{E}\big[  \|X_i \|_{L^1} \big] \leq 2/c$ and $ \mathbf{E}\big[   X_i(t) \big]=0$ for any $t\geq 0$. 
 The left-hand side of \eqref{eqn.618} can be bounded by
 \beqnn
  \zeta\cdot \bigg\| \frac{1}{[n\zeta]\cdot p_n -[n\zeta]^{2/3}} \sum_{i=1}^{[n\zeta]\cdot p_n -[n\zeta]^{2/3}} X_i  \bigg\|_{L^1} , 
 \eeqnn 
  which goes to $0$ in probability as $n\to\infty$ by Corollary~7.10 \cite[p.189]{LedouxTalagrand1991}.
 The whole proof ends.  
 \qed

 \begin{proposition}
 For each $T\geq 0$, we have $\sup_{t\in[0,T]}\big|\mathbf{I}^{(n)}_{c,2}(nt) - c\cdot W'(t)\big| \to 0$ as $n\to\infty$. 
 \end{proposition}
 \proof 
 For some constant $\beta> 0$, multiplying both sides of \eqref{Eqn.RPi.n} by the function $\frac{1}{\gamma_n}e^{-\frac{\beta}{n}\cdot t}$ and then integrating them over $(t,\infty)$, 
 \beqlb \label{eqn.661}
  \frac{1}{\gamma_n}\int_t^\infty e^{-\frac{\beta}{n} x}R_{\varPi_n}(x)dx 
  \ar=\ar \int_t^\infty e^{-\frac{\beta}{n} x}\overline\varPi_n(x)dx +  \int_t^\infty \int_0^x  e^{-\frac{\beta}{n} (x-y)}\overline\varPi_n(x-y)    e^{-\frac{\beta}{n} y} R_{\varPi_n} (y)dydx. \quad 
 \eeqlb
 Note that the last double integral equals to 
 \beqnn
  \int_0^\infty \int_0^x e^{-\frac{\beta}{n} (x-y)}\overline\varPi_n(x-y)    e^{-\frac{\beta}{n} y} R_{\varPi_n} (y)dydx
  -\int_0^t \int_0^x e^{-\frac{\beta}{n} (x-y)}\overline\varPi_n(x-y)    e^{-\frac{\beta}{n} y} R_{\varPi_n} (y)dydx.
 \eeqnn
 Applying Fubini's theorem to these two double integrals, they equal to  
 \beqnn
  \lefteqn{ \int_0^\infty   e^{-\frac{\beta}{n} y} R_{\varPi_n} (y)dy \cdot \int_0^\infty   e^{-\frac{\beta}{n} x}\overline\varPi_n(x) dx - \int_0^t  e^{-\frac{\beta}{n} (t-x)} R_{\varPi_n}(t-x)    \int_0^x e^{-\frac{\beta}{n} y}  \overline\varPi_n(y)dydx}\ar\ar\cr 
  \ar=\ar  \int_t^\infty e^{-\frac{\beta}{n} y} R_{\varPi_n}(y) dy \cdot \int_0^\infty   e^{-\frac{\beta}{n} x} \overline\varPi_n(x)dx  +  \int_0^t e^{-\frac{\beta}{n} (t-x)} R_{\varPi_n}(t-x) \int_x^\infty   e^{-\frac{\beta}{n} y} \overline\varPi_n(y)dy  dx.
 \eeqnn
 Plugging this back into the right side of \eqref{eqn.661} and then making some simple transform, 
 \beqnn
  \lefteqn{\int_t^\infty e^{-\frac{\beta}{n} x}\overline\varPi_n(x)\,dx +\int_0^t   e^{-\frac{\beta}{n} (t-x)}\overline\varPi_n(t-x) \int_x^\infty   e^{-\frac{\beta}{n} y} R_{\varPi_n} (y)\, dy  dx}
  \qquad \ar\ar\cr
  \ar=\ar \frac{1}{\gamma_n}\bigg( 1- \gamma_n \int_0^\infty   e^{-\frac{\beta}{n} z}\overline\varPi_n (z)\, dz \bigg) \cdot \int_t^\infty e^{-\frac{\beta}{n} x}R_{\varPi_n}(x)\,dx , 
 \eeqnn
 which along with the second equality in \eqref{eqn.504} induces that 
 \beqnn
  \lefteqn{ \int_t^\infty  e^{-\frac{\beta}{n} x} \overline{\boldsymbol{E}}_c(x) dx  + \int_0^t  e^{-\frac{\beta}{n} (t-y)}  R_{\varPi_n} (t-y) \int_y^\infty   e^{-\frac{\beta}{n}  x}\overline{\boldsymbol{E}}_c(x) dx dy}\ar\ar\cr
  \ar=\ar \frac{1}{\gamma_n} \Big( 1- \int_0^\infty \gamma_n \cdot e^{-\frac{\beta}{n}  x}\overline\varPi_n(x) dx   \Big)\cdot  \int_t^\infty e^{-\frac{\beta}{n}  x}R_{\varPi_n}(x)dx\cr
  \ar\ar - \theta_n  \int_t^\infty  e^{-\frac{\beta}{n}  x}  \boldsymbol{e}_{c}*\overline{\Lambda}_n(x) dx 
  - \theta_n\int_0^t  e^{-\frac{\beta}{n} (t-y)}  R_{\varPi_n} (t-y) \int_y^\infty   e^{-\frac{\beta}{n} x}  \boldsymbol{e}_{c}*\overline{\Lambda}_n(x) dx dy . 
 \eeqnn
 Since $\int_t^\infty  e^{-\frac{\beta}{n}\cdot x} \overline{\boldsymbol{E}}_c(x) dx =  (\beta/n+1/c)^{-1} e^{-(\beta/n+1/c)  t} $, the left-hand side of this equality equals to 
 \beqnn 
  \frac{1}{\beta/n+1/c}\Big( e^{-(\frac{\beta}{n} +\frac{1}{c}) t} +  \int_0^t  e^{-\frac{\beta}{n} (t-y)}  R_{\varPi_n} (t-y)   e^{-(\frac{\beta}{n} +\frac{1}{c}) y} \, dy\Big) 
  \ar=\ar \frac{e^{-\frac{\beta}{n} t}}{\beta/n+1/c}\big( \overline{\boldsymbol{E}}_c(t)+    R_{\varPi_n} * \overline{\boldsymbol{E}}_c(t)  \big) .  
 \eeqnn
 These two equations together with the fact that $c\gamma_n\sim 1$ induce  that  
 \beqlb\label{eqn.525}
  \mathbf{I}^{(n)}_{1,2}(nt)= \overline{\boldsymbol{E}}_c(nt)+    R_{\varPi_n} * \overline{\boldsymbol{E}}_c(nt)    \sim e^{\beta t} \cdot\Big( I_1^{(n)}(\beta,t) -I_2^{(n)}(\beta,t) -I_3^{(n)}(\beta,t)\Big), 
 \eeqlb
 locally uniformly in $t,\beta\geq 0$ as $n\to\infty$ with 
 \beqnn
 I_1^{(n)}(\beta,t)\ar:=\ar \Big( 1-  \gamma_n \int_0^\infty e^{-\frac{\beta}{n}  z}\overline\varPi_n(z) dz  \Big)  \cdot \int_{nt}^\infty e^{-\frac{\beta}{n}  x}R_{\varPi_n}(x)\,dx,\cr
 I_2^{(n)}(\beta,t)\ar:=\ar \frac{\theta_n}{c}   \int_{nt}^\infty  e^{-\frac{\beta}{n} y}  \boldsymbol{e}_{c}*\overline{\Lambda}_n(y)\, dy , \cr
 I_3^{(n)}(\beta,t)\ar:=\ar \frac{\theta_n}{c}   \int_0^{nt}  e^{-\frac{\beta}{n} (nt-x)}  R_{\varPi_n} (nt-x) \int_x^\infty   e^{-\frac{\beta}{n} y}  \boldsymbol{e}_{c}*\overline{\Lambda}_n(y)\, dy dx . 
 \eeqnn
 By using the change of variables and then \eqref{eqn.553} with $\lambda =\beta$ and Lemma~\ref{Lemma.ConR},  
 \beqlb\label{eqn.523}
 I_1^{(n)}(\beta,t)\ar= \ar n\bigg( 1-  \gamma_n  \int_0^\infty e^{-\frac{\beta}{n}  z}\overline\varPi_n(z) dz   \bigg)\cdot \bigg( \int_0^\infty e^{- \beta x}R_{\varPi_n}(nx)dx - \int_0^t e^{- \beta x}R_{\varPi_n}(nx)dx\bigg)\cr
 \ar\to\ar \frac{\varPhi(\beta)}{\beta} \Big( \int_0^\infty e^{- \beta  x}W'(x)dx-  \int_0^t e^{- \beta  x}W'(x)dx \Big) \cr
 \ar=\ar  1-\frac{\varPhi(\beta)}{\beta}  \cdot   \int_0^t e^{- \beta  x}W'(x)dx , 
 \eeqlb
 uniformly in $t$ and $\beta $ on compacts as $n\to\infty$.  
 Additionally, by Proposition~\ref{Prop.502} we also have 
 \beqlb\label{eqn.524}
  \lim_{n\to\infty}\sup_{\beta \geq 0}\sup_{t\geq 0}I_2^{(n)}(\beta,t)
  \ar\leq\ar  \lim_{n\to\infty}\frac{\theta_n}{c}  \int_0^\infty   \boldsymbol{e}_{c}*\overline{\Lambda}_n(x)\, dx  
  = \lim_{n\to\infty}\frac{1}{n} \cdot\int_0^\infty  \bar\nu(\eta_n +y)\, dy = 0.
 \eeqlb 
 For $I_3^{(n)}(\beta,t)$, we can split it into the following two terms
 \beqlb 
  I_{3,1}^{(n)}(\beta,t) \ar:=\ar \frac{\theta_n}{c}   \int_0^{nt}  e^{-\frac{\beta}{n} (nt-x)}  R_{\varPi_n} (nt-x) dx \cdot \int_0^\infty   e^{-\frac{\beta}{n} y}  \boldsymbol{e}_{c}*\overline{\Lambda}_n(y) dy, \cr
  I_{3,2}^{(n)}(\beta,t) \ar:=\ar -\frac{\theta_n}{c}   \int_0^{nt}  e^{-\frac{\beta}{n} (nt-x)}  R_{\varPi_n} (nt-x) \int_0^x   e^{-\frac{\beta}{n} y}  \boldsymbol{e}_{c}*\overline{\Lambda}_n(y) dy dx.\label{eqn.5222}
 \eeqlb
 By Fubini's theorem and the change of variables,   
 \beqlb\label{eqn.5223}
  I_{3,1}^{(n)}(\beta,t) \ar=\ar\frac{\theta_n}{c}   \int_0^{nt}  e^{-\frac{\beta}{n} x}  R_{\varPi_n} (x) \,dx \cdot \int_0^\infty   e^{-\frac{\beta}{n} y}  \overline{\Lambda}_n(y) \,dy \cdot  \int_0^\infty   e^{-\frac{\beta}{n} z}  \boldsymbol{e}_{c}(z)\, dz \cr
  \ar= \ar \frac{n\theta_n}{c}   \int_0^{t}  e^{-\beta x}  R_{\varPi_n} (nx)\, dx \cdot   \frac{n}{\bar\nu(\eta_n)} \int_0^\infty   e^{-\beta y}\bar\nu(\eta_n+y)\,  dy  \cdot \frac{n/c}{n/c+\beta}\cr
  \ar=\ar  \frac{n/c}{n/c+\beta}  \int_0^{t}  e^{-\beta (t-x)}  R_{\varPi_n} \big(n(t-x)\big)      \int_0^\infty   e^{-\beta y}\bar\nu(\eta_n+y) \, dy dx  .
 \eeqlb
 For $I_{3,2}^{(n)}(\beta,t) $, by using Fubini's theorem again we have
 \beqnn
  \int_0^x   e^{-\frac{\beta}{n}\cdot y}  \boldsymbol{e}_{c}*\overline{\Lambda}_n(y) dy
  \ar=\ar \int_0^x  e^{-\frac{\beta}{n}\cdot y} \overline{\Lambda}_n(y) \int_0^{x-y}e^{-\frac{\beta}{n}\cdot z}  \boldsymbol{e}_{c}(z)dz  dy \cr
  \ar=\ar \int_0^x \frac{n/c}{n/c+\beta} \big( 1- e^{-(\beta/n+1/c)\cdot (x-y)}   \big)  e^{-\frac{\beta}{n}\cdot y} \overline{\Lambda}_n(y)   dy  .
 \eeqnn 
 Plugging this back into the right side of \eqref{eqn.5222} gives that 
 \beqnn
  I_{3,2}^{(n)}(\beta,t) \ar=\ar - \frac{\theta_n}{c}   \int_0^{nt}  e^{-\frac{\beta}{n}\cdot (nt-x)}  R_{\varPi_n} (nt-x) \int_0^x \frac{n/c}{n/c+\beta} \big( 1- e^{-(\beta/n+1/c)\cdot (x-y)}   \big)  e^{-\frac{\beta}{n}\cdot y} \overline{\Lambda}_n(y) \,  dy dx \cr
  \ar=\ar -\frac{n^2\theta_n}{c}   \int_0^{t}  e^{- \beta  (t-x)}  R_{\varPi_n} \big(n(t-x)\big) \int_0^{x} \frac{n/c}{n/c+\beta} \big( 1- e^{-(\beta +n/c)\cdot (x-y)}   \big)  e^{- \beta y} \overline{\Lambda}_n(ny)   dy dx \cr
  \ar=\ar  -\frac{n/c}{n/c+\beta}  \int_0^{t}  e^{- \beta  (t-x)}  R_{\varPi_n} \big(n(t-x)\big) \int_0^{x}    e^{- \beta y} \bar\nu(\eta_n+y) \, dy dx\cr
  \ar\ar +\frac{n/c}{n/c+\beta} e^{- \beta t}\int_0^{t} R_{\varPi_n} \big(n(t-x)\big) \int_0^{x}    e^{-\frac{n}{c}\cdot (x-y)}     \bar\nu(\eta_n+y) \, dy dx . 
 \eeqnn
 Combing this together with \eqref{eqn.5223}, we have  as $n\to\infty$,
 \beqlb \label{eqn.557}
  I_3^{(n)}(\beta,t) 
  \ar\sim\ar \int_0^{t} e^{-\beta (t-x)}  R_{\varPi_n} \big(n(t-x)\big) \int_x^\infty e^{-\beta y}\bar\nu(\eta_n+y) \, dy dx \cr
  \ar\ar + e^{- \beta t}\int_0^{t} R_{\varPi_n} \big(n(t-x)\big) \int_0^{x} e^{-\frac{n}{c} (x-y)} \bar\nu(\eta_n+y) \, dy dx , 
 \eeqlb
 locally uniformly in $t,\beta\geq 0$.  
 Applying Fubini's theorem to the second double integral on the right side and then using Lemma~\ref{Lemma.UpperBoundR}, it equals to
 \beqnn
  e^{-\beta t}\int_0^{t}\bar\nu(\eta_n+x) \, dx \int_0^{t-x}    R_{\varPi_n} \big(n(t-x-y)\big) e^{-\frac{n}{c}  y}\, dy 
  \ar\leq\ar \frac{1}{c}  \int_0^{t}    \bar\nu(\eta_n+x) \, dx \cdot \int_0^{t}  e^{-\frac{n}{c}\cdot y} \, dy \cr
  \ar\leq\ar  \frac{1}{n}  \int_0^\infty    \bar\nu(\eta_n+x) \, dx ,
 \eeqnn 
 which vanishes uniformly in $t,\beta\geq 0$ as $n\to\infty$; see Proposition~\ref{Prop.502}. 
 For the first term on the right side of \eqref{eqn.557}, an application of Fubini's theorem again induces that 
 \beqnn 
  \int_0^{t}  e^{-\beta (t-x)}  R_{\varPi_n} \big(n(t-x)\big)   \int_x^\infty   e^{-\beta y}\bar\nu(\eta_n+y) \, dy dx= \int_0^\infty e^{-\beta y}\bar\nu(\eta_n+y) \int^{t}_{t-y }   e^{-\beta x}  R_{\varPi_n} (nx)\, dx  dy,
 \eeqnn
 which can be written as the sum of the following two terms
 \beqnn
  A^{(n)}_1(\beta,t) \ar:=\ar \int_0^\infty \Big(\int^{t}_{t-y }  e^{-\beta z}  R_{\varPi_n} (nz)\,dz -\int^{t}_{t-y } e^{-\beta z}  W'(z)\, dz\Big) e^{-\beta y} \bar\nu(\eta_n+y) \, dy, \cr
  A^{(n)}_2(\beta,t)\ar:=\ar \int_0^\infty  \Big(  \int^{t}_{t-y }  e^{-\beta  z}  W'(z)\, dz \Big)\cdot e^{-\beta y}\bar\nu(\eta_n+y) \, dy.
 \eeqnn
 By using the dominated convergence theorem along with \eqref{LevyTriplet} and the fact that 
 \beqnn 
  \sup_{n\geq 1}\sup_{t\in[0,T]}\Big|\int^{t}_{t-y }  e^{-\beta z}  R_{\varPi_n} (nz)\,dz \Big| + \sup_{t\in[0,T]} \Big|\int^{t}_{t-y } e^{-\beta z}  W'(z)\, dz\Big| \leq C\cdot (T\wedge y) ,
 \eeqnn
 uniformly in $T,y,\beta\geq 0$; see Lemma~\ref{Lemma.UpperBoundR} and \eqref{eqn.201},  
 and then using Lemma~\ref{Lemma.ConR}, 
 \beqnn
  \lim_{n\to\infty}\big|  A^{(n)}_1(\beta,t) \big|
  \ar\leq\ar \int_0^\infty \lim_{n\to\infty}\sup_{t\in[0,T]}\Big|\int^{t}_{t-y }  e^{-\beta z}  R_{\varPi_n} (nz)\,dz -\int^{t}_{t-y } e^{-\beta z}  W'(z)dz\Big| e^{-\beta y} \bar\nu(y) \, dy
  =0.
 \eeqnn
 Moreover, applying the monotone convergence theorem to $ A^{(n)}_2(\beta,t)$ we also have 
 \beqnn
  A^{(n)}_2(\beta,t) \to \int_0^\infty  \Big(  \int^{t}_{t-y }  e^{-\beta z}  W'(z)\, dz \Big)\cdot e^{-\beta y}\bar\nu(y) \, dy,
 \eeqnn
 uniformly in $\beta$ and $t$ on compacts as $n\to\infty$. 
 Taking all these limits back into \eqref{eqn.557}, 
 \beqnn
  I_3^{(n)}(\beta,t)  
  \ar\to\ar \int_0^\infty  \Big(  \int^{t}_{t-y }  e^{-\beta z}  W'(z)dz \Big)\cdot e^{-\beta y}\bar\nu(y)  dy,
 \eeqnn
 locally uniformly in $t$ and $\beta $ as $n\to\infty$.  
 Taking this and \eqref{eqn.523}-\eqref{eqn.524} back into \eqref{eqn.525} induces that 
 \beqnn
  \mathbf{I}^{(n)}_{1,2}(nt)
  \ar\to \ar e^{\beta t}\cdot \Big(1-\frac{\varPhi(\beta)}{\beta}  \cdot \int_0^t e^{- \beta  x}W'(x)\, dx  -  \int_0^\infty  \Big(  \int^{t}_{t-y }  e^{-\beta z}  W'(z) \, dz \Big)\cdot e^{-\beta y}\bar\nu(y) \, dy \Big),
 \eeqnn
 uniformly in $t$ and $\beta $ on compacts. 
 Note that $\varPhi(\beta)/\beta \to b$ as $\beta \to 0+$. 
 Finally, by \eqref{eqn.103} and \eqref{eqn.10301},  
 \beqnn
  \mathbf{I}^{(n)}_{1,2}(nt) \overset{\rm u.c.}
  \to 1-b W(t)  -  \int_0^\infty  \int^{t}_{t-y }  W'(z) \, dz \bar\nu(y) \, dy =c \cdot W'(t),
 \eeqnn 
 as $n\to\infty$ and the whole proof ends. 
 \qed 
 
 The next lemma follows immediately by applying the locally uniform convergence results in the preceding two propositions and Lemma~\ref{Lemma.UpperBoundR} to the right side of \eqref{eqn.550}. 
 
 \begin{lemma}   \label{Lemma.I1}
 For each $T\geq 0$, we have $\sup_{t\in[0,T]}\big|\mathbf{I}^{(n)}_1(t) - \zeta \cdot c\cdot W'(t) \big| \overset{\rm p}\to 0$ as $n\to\infty$.
 \end{lemma}
  
 \subsubsection{Weak convergence of $\{\mathbf{I}^{(n)}_2\}_{n\geq 1}$}
  
 By \eqref{Eqn.R.n} and the change of variables, we can split $\mathbf{I}^{(n)}_2(t)  $ into the following two terms
 \beqnn
  \mathbf{I}^{(n)}_{2,1}(t): = \frac{1}{n} \sum_{j=1}^{B_n^\Lambda(\zeta)} \mathbf{1}_{\{\ell_{n,j}^\Lambda\geq nt\}} 
  \quad \mbox{and}\quad 
  \mathbf{I}^{(n)}_{2,2}(t):= \sum_{j=1}^{B_n^\Lambda(\zeta)} \int_{t-\frac{\ell_{n,j}^\Lambda}{n}}^{t} R_{\varPi_n}(ns) \, ds  . 
 \eeqnn
 The next proposition is a direct consequence of the fact that $B_n^\Lambda(\zeta)/n\sim \zeta\cdot p_n \overset{\rm p}\to 0$ as $n\to\infty$. 
 
 \begin{proposition}\label{Prop.519}
  We have $\sup_{t\geq 0} \big| \mathbf{I}^{(n)}_{2,1}(t) \big| \overset{\rm p}\to 0$ as $n\to\infty$.
 \end{proposition}
 
 \begin{proposition}\label{Prop.520}
  The sequence $ \big\{ \mathbf{I}^{(n)}_{2,2} \big\}_{n\geq 1}$ is $C$-tight. 
 \end{proposition}
 \proof  
 By Aldous's tightness criterion; see Theorem~1 in \cite{Aldous1978}, it suffices to prove that for each $T\geq 0$ and $\varepsilon>0$, there exist constants $n_0\geq 1$ and $\delta_0>0$ such that 
 \beqlb\label{eqn.527}
  \sup_{n\geq n_0}\sup_{0<\delta \leq \delta_0} \sup_{\tau \in \mathcal{S}_T^{(n)}} \mathbf{E}\Big[  \big| \mathbf{I}^{(n)}_{2,2}(\tau+\delta)-\mathbf{I}^{(n)}_{2,2}(\tau) \big| \Big] 
  \leq \varepsilon,
 \eeqlb
 where $\mathcal{S}_T^{(n)}$ is the collection of all $(\mathscr{F}_{nt})$-stopping times bounded by $T$.  
 We first have  
 \beqnn
  \Big| \mathbf{I}^{(n)}_{2,2}(\tau+\delta)-\mathbf{I}^{(n)}_{2,2}(\tau) \Big|
  \ar=\ar \Bigg| \sum_{j=1}^{B_n^\Lambda(\zeta)} \bigg(\int_{\tau+\delta-\frac{\ell_{n,j}^\Lambda}{n}}^{\tau+\delta} R_{\varPi_n}(ns) ds - \int_{\tau-\frac{\ell_{n,j}^\Lambda}{n}}^{\tau} R_{\varPi_n}(ns) ds\bigg)\Bigg| \cr
  \ar\leq\ar    \sum_{j=1}^{B_n^\Lambda(\zeta)} \bigg|\int_{\tau+\delta-\frac{\ell_{n,j}^\Lambda}{n}}^{\tau+\delta} R_{\varPi_n}(ns) ds - \int_{\tau-\frac{\ell_{n,j}^\Lambda}{n}}^{\tau} R_{\varPi_n}(ns) ds\bigg|.
 \eeqnn
 When $\delta\geq \ell_{n,j}^\Lambda/n$, the disjunction of the two intervals $(\tau+\delta-\frac{\ell_{n,j}^\Lambda}{n},\tau+\delta ]$ and $(\tau-\frac{\ell_{n,j}^\Lambda}{n},\tau]$ along with Lemma~\ref{Lemma.UpperBoundR} induces  that
 \beqlb\label{eqn.591}
  \bigg|\int_{\tau+\delta-\frac{\ell_{n,j}^\Lambda}{n}}^{\tau+\delta} R_{\varPi_n}(ns) ds - \int_{\tau-\frac{\ell_{n,j}^\Lambda}{n}}^{\tau} R_{\varPi_n}(ns) ds\bigg| 
  \leq \frac{2}{c} \cdot \frac{\ell_{n,j}^\Lambda}{n}.
 \eeqlb
 When $ \delta<\ell_{n,j}^\Lambda/n$, we also have $(\tau+\delta-\frac{\ell_{n,j}^\Lambda}{n},\tau+\delta ]\cap (\tau-\frac{\ell_{n,j}^\Lambda}{n},\tau]= (\tau+\delta-\frac{\ell_{n,j}^\Lambda}{n},\tau]$ and then
 \beqlb\label{eqn.592}
  \lefteqn{\bigg|\int_{\tau+\delta-\frac{\ell_{n,j}^\Lambda}{n}}^{\tau+\delta} R_{\varPi_n}(ns) ds - \int_{\tau-\frac{\ell_{n,j}^\Lambda}{n}}^{\tau} R_{\varPi_n}(ns) ds\bigg|}\ar\ar\cr  
  \ar=\ar \bigg|\int_{\tau}^{\tau+\delta} R_{\varPi_n}(ns) ds+ \int_{\tau-\frac{\ell_{n,j}^\Lambda}{n}}^{\tau-\frac{\ell_{n,j}^\Lambda}{n}+\delta} R_{\varPi_n}(ns) ds\bigg|\leq \frac{2}{c} \cdot \delta.
 \eeqlb 
 In conclusion, we have uniformly in $n\geq1 $, $\tau \in \mathcal{S}_T^{(n)}$ and $\delta >0$, 
 \beqlb\label{eqn.593}
  \Big| \mathbf{I}^{(n)}_{2,2}(\tau+\delta)-\mathbf{I}^{(n)}_{2,2}(\tau) \Big|
  \leq \frac{2}{c}\sum_{j=1}^{B_n^\Lambda(\zeta)} \Big(\frac{\ell_{n,j}^\Lambda}{n}\wedge \delta \Big).
 \eeqlb
 By using the independence between $B_n^\Lambda(\zeta)$ and $\{  \ell_{n,j}^\Lambda\}_{j\geq 1}$ and then Proposition~\ref{Prop.503} with $p=1$,  
 \beqnn
  \sup_{\tau \in \mathcal{S}_T^{(n)}} \mathbf{E}\Big[  \big| \mathbf{I}^{(n)}_{2,2}(\tau+\delta)-\mathbf{I}^{(n)}_{2,2}(\tau) \big| \Big]
  \ar\leq\ar \frac{2}{c} \cdot \mathbf{E}\Bigg[   \sum_{j=1}^{B_n^\Lambda(\zeta)} \Big(\frac{\ell_{n,j}^\Lambda}{n}\wedge \delta \Big) \Bigg]\cr
  \ar=\ar \frac{2}{c}\cdot \mathbf{E}\big[    B_n^\Lambda(\zeta) \big] \cdot  \mathbf{E}\Big[\frac{\ell_{n,1}^\Lambda}{n}\wedge \delta   \Big]\cr
  \ar=\ar \frac{2}{c}\cdot   [n\zeta]\cdot (1-p_n) \cdot \mathbf{E}\Big[\frac{\ell_{n,1}^\Lambda}{n}\wedge \delta   \Big]\cr
  \ar\leq\ar C \cdot \Big(1-p_n+   \int_0^\infty   (y\wedge \delta) \bar\nu(y) \, dy \Big) ,
 \eeqnn
 which goes to $0$ as $n\to\infty$ and $\delta \to 0+$.  Hence \eqref{eqn.527} holds.
 \qed

 \begin{lemma}\label{Lemma.I2}
  As $n\to\infty$, we have that $\mathbf{I}^{(n)}_2 $ converges weakly in $D(\mathbb{R}_+;\mathbb{R}_+)$ to  
 \beqnn
  \mathbf{I}_2(t):= \int_0^\zeta \int_0^\infty \big( W(t)-W(t-y) \big) N_0(dz,dy),\quad t\geq 0.  
 \eeqnn  
 \end{lemma}
 \proof 
 By Proposition~\ref{Prop.519} and \ref{Prop.520}, the sequence $\{\mathbf{I}^{(n)}_2\}_{n\geq 1}$ is $C$-tight and it suffices to  prove that 
 \beqlb \label{eqn.529}
  \lim_{n\to\infty} \mathbf{E}\Bigg[ \exp\bigg\{ -\sum_{i=1}^d\lambda_i \cdot \mathbf{I}_{2,2}^{(n)}(t_i) \bigg\} \Bigg]
  = \mathbf{E}\Bigg[ \exp\bigg\{ -\sum_{i=1}^d\lambda_i \cdot \mathbf{I}_2(t_i) \bigg\} \Bigg]
 \eeqlb
 for any $d \in \mathbb{Z}_+$, $0\leq t_1<t_2<\cdots <t_d$ and $ \lambda_1,\cdots,\lambda_d \in \mathbb{R}_+$.
 By the exponential formula of stochastic integral with respect to Poisson random measure; see \cite[p.8]{Bertoin1996}, we first have
 \beqlb\label{eqn.528}
 \lefteqn{ \mathbf{E}\Bigg[ \exp\bigg\{ -\sum_{i=1}^d\lambda_i\cdot \mathbf{I}_{2}(t_i) \bigg\} \Bigg] }\ar\ar\cr
  \ar=\ar \mathbf{E}\Bigg[ \exp\bigg\{ - \int_0^\zeta \int_0^\infty \sum_{i=1}^d\lambda_i\big( W(t_i)-W(t_i-y) \big) N_0(dz,dy) \bigg\} \Bigg]\cr
  \ar=\ar \exp\Bigg\{- \zeta \cdot  \int_0^\infty\bigg(1- \exp\bigg\{ - \sum_{i=1}^d \lambda_i  \big( W(t_i)-W(t_i-y) \big) \bigg\} \bigg) \bar\nu(y) \,dy   \Bigg\}.
 \eeqlb
 On the other hand, by the mutual-independence among $B_n^\Lambda(\zeta)$ and $\{\ell_{n,j}^\Lambda\}_{j\geq 1}$ and then the Laplace transform of binomial distribution, 
 \beqnn
  \mathbf{E}\Bigg[ \exp\bigg\{ -\sum_{i=1}^d\lambda_i \cdot \mathbf{I}^{(n)}_{2,2}(t_i) \bigg\} \Bigg]
  \ar=\ar \mathbf{E}\Bigg[ \exp\bigg\{ -  \sum_{j=1}^{B_n^\Lambda(\zeta)}\sum_{i=1}^d  \lambda_i  \int_{t_i-\ell_{n,j}^\Lambda/n}^{t_i} R_{\varPi_n}(ns)\, ds \bigg\} \Bigg]\cr
  \ar=\ar  \mathbf{E}\Bigg[  \mathbf{E}\bigg[\exp\bigg\{ - \sum_{i=1}^d \lambda_i  \int_{t_i-\ell_{n,1}^\Lambda/n}^{t_i} R_{\varPi_n}(ns)\, ds \bigg\}\bigg]^{B_n^\Lambda(\zeta)} \Bigg] \cr
  \ar=\ar   \Bigg[ p_n +\big(1-p_n\big) \cdot \mathbf{E}\bigg[\exp\bigg\{ - \sum_{i=1}^d \lambda_i  \int_{t_i-\ell_{n,1}^\Lambda/n}^{t_i} R_{\varPi_n}(ns)\, ds \bigg\}\bigg] \Bigg]^{[n\zeta]}\cr
  \ar=\ar \Bigg[ 1 +\big(1-p_n\big) \cdot \mathbf{E}\bigg[\exp\bigg\{ - \sum_{i=1}^d \lambda_i  \int_{t_i-\ell_{n,1}^\Lambda/n}^{t_i} R_{\varPi_n}(ns)\, ds \bigg\}-1\bigg] \Bigg]^{[n\zeta]},
 \eeqnn
 which is asymptotically equivalent to 
 \beqnn  
  \exp\Bigg\{- [n\zeta]\cdot  \big(1-p_n\big) \cdot \mathbf{E}\bigg[1-\exp\bigg\{ - \sum_{i=1}^d \lambda_i  \int_{t_i-\ell_{n,j}^\Lambda/n}^{t_i} R_{\varPi_n}(ns)\, ds \bigg\}\bigg] \Bigg\} ,
 \eeqnn 
 since both $1-p_n$ and the last expectation go to $0$ as $n\to\infty$.
 By the definition of $\ell_{n,j}^\Lambda$ and then using Fubini's  theorem as well as the change of variables, 
 \beqnn
  \lefteqn{\mathbf{E}\bigg[1- \exp\bigg\{ - \sum_{i=1}^d \lambda_i  \int_{t_i-\ell_{n,j}^\Lambda/n}^{t_i} R_{\varPi_n}(ns)\, ds \bigg\} \bigg] }\quad\ar\ar\cr
  \ar=\ar 
  \int_0^\infty \bigg(1- \exp\bigg\{ - \sum_{i=1}^d \lambda_i  \int_{t_i-x/n}^{t_i} R_{\varPi_n}(ns)\, ds \bigg\} \bigg)\boldsymbol{e}_{c}*\boldsymbol{\lambda}^*_n(x)\,dx\cr 
  \ar=\ar 
  \int_0^\infty \boldsymbol{e}_{c}(x)\, dx \int_0^\infty\bigg(1- \exp\bigg\{ - \sum_{i=1}^d \lambda_i  \int_{t_i-(y+x)/n}^{t_i} R_{\varPi_n}(ns)\, ds \bigg\} \bigg) \boldsymbol{\lambda}^*_n(y) \, dy \cr
  \ar=\ar 
  \int_0^\infty \boldsymbol{e}_{c}(x)\,dx\int_0^\infty\bigg(1- \exp\bigg\{ - \sum_{i=1}^d \lambda_i  \int_{t_i-(y+x)/n}^{t_i} R_{\varPi_n}(ns)\, ds \bigg\} \bigg) \frac{\bar\nu(\eta_n+y/n)}{n\int_{\eta_n}^\infty \bar\nu(z)dz}\, dy \cr
  \ar=\ar 
  \int_0^\infty \boldsymbol{e}_{c}(x)\,dx \int_0^\infty\bigg(1- \exp\bigg\{ - \sum_{i=1}^d \lambda_i  \int_{t_i-y-x/n}^{t_i} R_{\varPi_n}(ns) \,ds \bigg\} \bigg) \frac{\bar\nu(\eta_n+y)}{\int_{\eta_n}^\infty \bar\nu(z)dz}\, dy . 
 \eeqnn
 Combing this together with the fact that $ [n\zeta]\cdot  \big(1-p_n\big) \sim \zeta\cdot \int_{\eta_n}\bar\nu(z)dz$; see \eqref{eqn.pn}, we have as $n\to\infty$,  
 \beqnn
  \lefteqn{\mathbf{E}\Bigg[ \exp\bigg\{ -\sum_{i=1}^d\lambda_i \cdot \mathbf{I}^{(n)}_{2,2}(t_i) \bigg\} \Bigg]}\ar\ar\cr
  \ar\sim\ar 
  \exp\Bigg\{- \zeta \cdot \int_0^\infty \boldsymbol{e}_{c}(x)\, dx \int_0^\infty\bigg(1- \exp\bigg\{ - \sum_{i=1}^d \lambda_i  \int_{t_i-y-x/n}^{t_i} R_{\varPi_n}(ns) \, ds \bigg\} \bigg) \bar\nu(\eta_n+y)\, dy  \Bigg\}\cr
  \ar\to\ar 
  \exp\Bigg\{- \zeta \cdot  \int_0^\infty\bigg(1- \exp\bigg\{ - \sum_{i=1}^d \lambda_i   \big( W(t_i)-W(t_i-y) \big) \bigg\} \bigg) \bar\nu(y)\, dy   \Bigg\}, 
 \eeqnn
 which together with \eqref{eqn.528}  immediately induces \eqref{eqn.529}. 
 \qed
  
 \subsubsection{$C$-tightness of $\{\mathbf{I}^{(n)}_3\}_{n\geq 1}$}

 In this section we prove the $C$-tightness of the sequence $\{\mathbf{I}^{(n)}_3\}_{n\geq 1}$.
 Plugging \eqref{Eqn.R.n} into \eqref{eqn.Jc}, we can split $\mathbf{I}^{(n)}_3(t)$ into the following two terms.  
 \beqlb
  \mathbf{I}^{(n)}_{3,1}(t) \ar:=\ar  \int_0^t\int_0^{X_{\zeta}^{(n)} (s-)} \int_0^\infty  \frac{1}{n}\cdot \mathbf{1}_{\{y\geq n(t-s) \}}  \widetilde{N}^{(n)}_c (ds,dz,dy) , \label{eqn.5331}\\
  \mathbf{I}^{(n)}_{3,2}(t) \ar:=\ar  \int_0^t\int_0^{X_{\zeta}^{(n)} (s-)} \int_0^\infty  \frac{1}{n}\cdot \Big(\int_{n(t-s)-y}^{n(t-s)} R_{\varPi_n}(r) dr\Big)\widetilde{N}^{(n)}_c (ds,dz,dy).\label{eqn.5332}
 \eeqlb
 By Corollary~3.33 in \cite[p.353]{JacodShiryaev2003}, it suffices to prove the $C$-tightness of the two sequences $\{ \mathbf{I}^{(n)}_{3,1}\}_{n\geq 1}$ and $\{\mathbf{I}^{(n)}_{3,2}\}_{n\geq 1}$ separately.
 We first prove the weak convergence of the first sequence to $0$ by using the following $C$-tightness criterion for c\`adl\`ag stochastic processes established in \cite{HorstXuZhang2023a}.
 
 \begin{proposition} \label{tightness condition}
 For a sequence of c\`adl\`ag stochastic processes $\{X^{(n)}\}_{n\geq 1}$ defined on a common probability space with $\sup_{n \geq 1}  \mathbf{E} \big[|X^{(n)}(0)|^q\big] < \infty$ for some $q>0$, it is $C$-tight if the following two conditions hold for any $T\geq 0$ and some constant $\beta>2$.
 \begin{enumerate}
  \item[(1)] 
  There exist some constants $C>0 $, $p\geq 1$ and $\rho,\alpha>0$ such that for any $\delta\in(0,1)$ and $n\geq 1$,
  \beqnn
   \sup_{t\in[0,T]} \mathbf{E}\Big[\big| X^{(n)}(t+\delta)-X^{(n)}(t)\big|^{p}\Big] 
   \leq C\cdot \Big(  \delta^{1+\rho}+ \frac{\delta}{n^\alpha}\Big).
  \eeqnn
 		
  \item[(2)] 
  $\displaystyle\sup_{k=0,1,\cdots,[Tn^\beta]} \sup_{\delta\in[0,1/n^{\beta}]} \big|  X^{(n)}(k/n^\beta+\delta) -X^{(n)}(k/n^\beta)\big|  \overset{\rm p}\to 0$  as $n\to\infty$. 
 \end{enumerate}  
 \end{proposition}
   
 \begin{lemma}\label{Lemma.526}
  We have $\mathbf{I}^{(n)}_{3,1} \to 0$ weakly in $ D (\mathbb{R}_+;\mathbb{R})$ as $n\to\infty$.  
 \end{lemma}
 \proof 
 For each $t\geq 0$, by using \eqref{BDG} along with Lemma~\ref{Lemma.Moment} as well as the two facts that $c\gamma_n\leq 1 $ and $\theta_n\in (0,1)$, there exists a constant $C>0$ such that for any $n\geq 1$,
 \beqnn
  \mathbf{E}\Big[ \big|\mathbf{I}^{(n)}_{3,1}(t) \big|^2 \Big] 
  \ar\leq \ar C\cdot   \int_0^t \int_0^\infty  \mathbf{1}_{\{y\geq n(t-s) \}}   \cdot \boldsymbol{e}_{c}(y)\, ds dy
  \leq  C\cdot   \int_0^t e^{-ns/c} ds 
  \leq \frac{C}{n} , 
 \eeqnn
 which goes to $0$ as $n\to\infty$. This induces that 
 $
  \mathbf{I}^{(n)}_{3,1} \overset{\rm f.f.d.}\longrightarrow 0.
 $
 For the tightness, it suffices to prove that $\{\mathbf{I}^{(n)}_{3,1}\}_{n\geq 1}$ satisfies the two conditions in Proposition~\ref{tightness condition}. 
 For each $t\geq 0$ and $\delta\in(0,1)$, we have
 \beqlb\label{eqn.565}
  \mathbf{I}^{(n)}_{3,1}(t+\delta)-\mathbf{I}^{(n)}_{3,1}(t) = J^{(n)}_{1,1}(t,\delta)-J^{(n)}_{1,2}(t,\delta)
 \eeqlb
 with 
 \beqnn 
  J^{(n)}_{1,1}(t,\delta)\ar:=\ar \int_t^{t+\delta}\int_0^{X_{\zeta}^{(n)} (s-)} \int_{n(t+\delta-s)}^\infty  \frac{1}{n} \, \widetilde{N}^{(n)}_c (ds,dz,dy) , \cr
  J^{(n)}_{1,2}(t,\delta)\ar:=\ar   \int_0^t\int_0^{X_{\zeta}^{(n)} (s-)} \int_{n(t-s)}^{n(t+\delta-s)}  \frac{1}{n} \, \widetilde{N}^{(n)}_c (ds,dz,dy) . 
 \eeqnn
 \medskip
 
 \textit{\textbf{Condition (1).}}
 For any $p\geq 1$, by using \eqref{BDG} along with Lemma~\ref{Lemma.Moment} and then the change of variables, there exists a constant $C>0$ independent of $n$ and  $\delta$ such that
 \beqlb\label{eqn.530}
  \sup_{t\in[0,T]} \mathbf{E}\Big[ \big| J^{(n)}_{1,1}(t,\delta) \big|^{2p}  \Big] 
  \ar\leq\ar C \cdot \sup_{t\in[0,T]}\mathbf{E} \Big[\big|X^{(n)}_\zeta(t)\big|^p \Big] \cdot  \Big| \int_0^\delta ds  \int_{ ns}^\infty   \boldsymbol{e}_{c}(y)\,  dy \Big|^{p} \cr
  \ar\ar  +  \frac{C}{n^{2p-2}} \cdot \sup_{t\in[0,T]}\mathbf{E} \Big[\big|X^{(n)}_\zeta(t)\big| \Big] \cdot  \int_0^{\delta}ds \int_{ ns}^\infty    \boldsymbol{e}_{c}(y)\, dy \cr  
  \ar\leq\ar C \cdot \Big(\delta^p + \frac{\delta}{n^{2p-2}}\Big) . 
 \eeqlb 
 Similarly, we also have uniformly in $n\geq 1$ and $\delta\in (0,1)$,
 \beqlb\label{eqn.531}
  \sup_{t\in[0,T]}\mathbf{E}\Big[ \big| J^{(n)}_{1,2}(t,\delta) \big|^{2p}  \Big] 
  \ar\leq\ar C\cdot \sup_{t\in[0,T]}\mathbf{E} \Big[\big|X^{(n)}_\zeta(t)\big|^p \Big] \cdot \Big| \int_0^T ds \int_{ ns}^{n(s+\delta)}       \boldsymbol{e}_{c}(y)\, dy \Big|^{p} \cr
  \ar\ar +   \frac{C}{n^{2p-2}} \sup_{t\in[0,T]}\mathbf{E} \Big[\big|X^{(n)}_\zeta(t)\big|^p \Big]  \int_0^T ds\int_{ ns}^{n(s+\delta)}      \boldsymbol{e}_{c}(y)\,   dy \cr
  \ar\leq\ar  C  \Big| \int_0^T ds \int_{ns}^{n(s+\delta)}    \boldsymbol{e}_{c}(y)\,  dy \Big|^{p} 
  + \frac{C}{n^{2p-2}}  \int_0^T ds \int_{ns}^{n(s+\delta)}    \boldsymbol{e}_{c}(y)\, dy. 
 \eeqlb
 A simple calculation along with the inequality that $  x (1-e^{-1/x})\leq 1$ for any $x\geq 0$ shows  
 \beqnn
  \int_0^T  ds \int_{ns}^{n(s+\delta)}    \boldsymbol{e}_{c}(y)\,  dy
  \ar=\ar \int_0^T e^{-\frac{ns}{c}} \, ds\cdot \big( 1-e^{- n\delta/c} \big)
  \leq \frac{c}{n} \cdot \big( 1-e^{- n\delta/c} \big) 
  \leq \delta . 
 \eeqnn  
 Plugging this back into \eqref{eqn.531} gives that uniformly in $n\geq 1$ and $\delta\in (0,1)$, 
 \beqlb\label{eqn.532}
  \sup_{t\in[0,T]} \mathbf{E}\Big[ \big| J^{(n)}_{1,2}(t,\delta) \big|^{2p}  \Big] \leq C\cdot \Big(\delta^p +\frac{\delta}{n^{2p-2}}\Big). 
 \eeqlb
 Armed with the two upper-bound estimates \eqref{eqn.530} and \eqref{eqn.532}, we use the power inequality to get 
 \beqnn
  \sup_{t\in[0,T]} \mathbf{E}\Big[ \big| \mathbf{I}^{(n)}_{3,1}(t+\delta)-\mathbf{I}^{(n)}_{3,1}(t)\big|^{2p}  \Big] 
  \ar\leq\ar  C\cdot \sup_{t\in[0,T]} \mathbf{E}\Big[ \big|J^{(n)}_{1,1}(t,\delta)\big|^{2p}  \Big] \cr
  \ar\ar +C\cdot \sup_{t\in[0,T]} \mathbf{E}\Big[ \big|J^{(n)}_{1,2}(t,\delta)\big|^{2p}  \Big] \cr
  \ar\leq\ar  C\cdot \Big(\delta^p +\frac{\delta}{n^{2p-2}}\Big),
 \eeqnn
 uniformly in $n\geq 1$ and $\delta \in(0,1)$ and hence Condition (1) in Proposition~\ref{tightness condition} is satisfied.
 \medskip
 
 \textit{\textbf{Condition (2).}} For $\beta>2$ and each $T> 0$, by \eqref{eqn.565} it holds that
 \beqlb\label{eqn.559}
  \lefteqn{\sup_{k=0,1,\cdots, [n^\beta T]} \sup_{\delta\in[0, 1/n^{\beta}]} \big| \mathbf{I}^{(n)}_{3,1}(k/n^{\beta}+\delta)-\mathbf{I}^{(n)}_{3,1}(k/n^{\beta})\big| }\ar\ar\cr
  \ar\leq\ar \sup_{k=0,1,\cdots, [n^\beta T]}\sup_{\delta\in[0, 1/n^{\beta}]} \big| J^{(n)}_{1,1}(k/n^{\beta},\delta)\big| 
  + \sup_{k=0,1,\cdots, [n^\beta T]} \sup_{\delta\in[0, 1/n^{\beta}]} \big|J^{(n)}_{1,2}(k/n^{\beta},\delta) \big|.
 \eeqlb
 It suffices to prove that both of the two terms on the right side converge in probability to $0$.
 \medskip
 
 \textbf{(i)}  By the facts that $c\gamma_n\leq 1$ and $\theta_n\in (0,1)$, we have uniformly  in $n\geq 1$, $k\geq 1$ and $\delta>0$, 
 \beqlb\label{eqn.558}
 \lefteqn{ \sup_{\delta\in[0, 1/n^{\beta}]} \big| J^{(n)}_{1,1}(k/n^{\beta},\delta)\big|  }\ar\ar\cr
  \ar\leq\ar  \sup_{\delta\in[0, 1/n^{\beta}]}\int_{k/n^{\beta}}^{k/n^{\beta}+\delta} \frac{n}{c} \cdot X_{\zeta}^{(n)} (s) \,  ds \int_{ n(t+\delta-s)}^\infty  \boldsymbol{e}_{c}(y)\, dy \cr
  \ar\ar +  \sup_{\delta\in[0, 1/n^{\beta}]}\int_{k/n^{\beta}}^{k/n^{\beta}+\delta}\int_0^{X_{\zeta}^{(n)} (s-)} \int_{ n(t+\delta-s)}^\infty  \frac{1}{n} \, N^{(n)}_c (ds,dz,dy) \cr
  \ar\leq\ar  \int_{k/n^{\beta}}^{(k+1)/n^{\beta}} \frac{n}{c} \cdot X_{\zeta}^{(n)} (s) \,  ds
  + \int_{k/n^{\beta}}^{(k+1)/n^{\beta}}\int_0^{X_{\zeta}^{(n)} (s-)} \int_0^\infty  \frac{1}{n} \, N^{(n)}_c (ds,dz,dy) . 
 \eeqlb  
 By the definition of compensated Poisson random measure, the last stochastic integral on the right side of the second inequality can be bounded by
 \beqnn
  \int_{k/n^{\beta}}^{(k+1)/n^{\beta}} \frac{n}{c} \cdot X_{\zeta}^{(n)} (s) \,  ds 
  + \Bigg|\int_{k/n^{\beta}}^{(k+1)/n^{\beta}} \int_0^{X_{\zeta}^{(n)} (s-)} \int_0^\infty  \frac{1}{n} \, \widetilde{N}^{(n)}_c (ds,dz,dy)\Bigg|  .
 \eeqnn
 Therefore, we have uniformly in $n\geq 1$,
 \beqlb\label{eqn.560}
  \lefteqn{\sup_{k=0,1,\cdots, [n^\beta T]}\sup_{\delta\in[0, 1/n^{\beta}]} \big| J^{(n)}_{c,1}(k/n^{\beta},\delta)\big| }\ar\ar\cr
  \ar\leq\ar 
  \sup_{k=0,1,\cdots, [n^\beta T]} \int_{k/n^{\beta}}^{(k+1)/n^{\beta}} \frac{2n}{c} \cdot X_{\zeta}^{(n)} (s) \,  ds \cr
  \ar\ar + \sup_{k=0,1,\cdots, [n^\beta T]}\Bigg|\int_{k/n^{\beta}}^{(k+1)/n^{\beta}}\int_0^{X_{\zeta}^{(n)} (s-)} \int_0^\infty  \frac{1}{n} \, \widetilde{N}^{(n)}_c (ds,dz,dy)\Bigg|  . 
 \eeqlb
 For any $\epsilon>0$, by using Chebyshev's inequality, H\"older's inequality and then Lemma~\ref{Lemma.Moment} we have 
 \beqnn
 \lefteqn{ \mathbf{P}\bigg(\sup_{k=0,1,\cdots, [n^\beta T]}
  \int_{k/n^{\beta}}^{(k+1)/n^{\beta}} \frac{2n}{c} \cdot X_{\zeta}^{(n)} (s) \,  ds\geq \epsilon  \bigg)}\ar\ar\cr 
  \ar\leq\ar 
  \sum_{k=0}^{[n^\beta T]}\mathbf{P}\bigg( \int_{k/n^{\beta}}^{(k+1)/n^{\beta}} \frac{2n}{c} \cdot X_{\zeta}^{(n)} (s) \,  ds\geq \epsilon  \bigg)\cr
  \ar\leq\ar 
  \frac{1}{\epsilon^{2p}}\sum_{k=0}^{[n^\beta T]} \mathbf{E}\Bigg[\bigg|\int_{k/n^{\beta}}^{(k+1)/n^{\beta}} \frac{2n}{c} \cdot X_{\zeta}^{(n)} (s) \,  ds \bigg|^{2p} \Bigg]\cr
  \ar\leq\ar 
  \frac{C}{\epsilon^{2p}} \cdot n^{2p- (2p-1)\beta} \sum_{k=0}^{[n^\beta T]}\int_{k/n^{\beta}}^{(k+1)/n^{\beta}}   \mathbf{E}\Big[\big|X_{\zeta}^{(n)} (s)\big|^{2p} \Big] \,  ds , 
 \eeqnn
 which can be bounded by $C \cdot n^{2p- (2p-1)\beta} \to 0$ as $n\to\infty$ for all $p\geq 2$. 
 Similarly, we also have 
 \beqlb \label{eqn.563}
  \lefteqn{ \mathbf{P}\bigg(\sup_{k=0,1,\cdots, [n^\beta T]}\int_{k/n^{\beta}}^{(k+1)/n^{\beta}}\int_0^{X_{\zeta}^{(n)} (s-)} \int_0^\infty  \frac{1}{n} \, \widetilde{N}^{(n)}_c (ds,dz,dy)\geq \epsilon  \bigg)}\ar\ar\cr
  \ar\leq\ar 
  \sum_{k=0}^{[n^\beta T]}\mathbf{P}\bigg(\int_{k/n^{\beta}}^{(k+1)/n^{\beta}}\int_0^{X_{\zeta}^{(n)} (s-)} \int_0^\infty  \frac{1}{n} \, \widetilde{N}^{(n)}_c (ds,dz,dy)\geq \epsilon  \bigg)\cr
  \ar\leq\ar 
  \frac{1}{\epsilon^{2p}}\sum_{k=0}^{[n^\beta T]} \mathbf{E}\Bigg[\bigg| \int_{k/n^{\beta}}^{(k+1)/n^{\beta}}\int_0^{X_{\zeta}^{(n)} (s-)} \int_0^\infty  \frac{1}{n} \, \widetilde{N}^{(n)}_c (ds,dz,dy) \bigg|^{2p} \Bigg]. 
 \eeqlb
 Similarly as in \eqref{eqn.530}, the last expectation can be bounded by $C\cdot  \big( n^{-p\beta} + n^{2-2p-\beta} \big)$ uniformly in $n\geq 1$ and $k\geq 0$. 
 Consequently, there exists a constant $C>0$ such that for any $n\geq 1$, 
 \beqlb \label{eqn.564}
  \mathbf{P}\bigg(\sup_{k=0,1,\cdots, [n^\beta T]}\int_{k/n^{\beta}}^{(k+1)/n^{\beta}}\int_0^{X_{\zeta}^{(n)} (s-)} \int_0^\infty  \frac{1}{n} \,  \widetilde{N}^{(n)}_c (ds,dz,dy)\geq \epsilon  \bigg)
  \leq  C\cdot  \big( n^{\beta(1-p)} + n^{2(1-p)} \big),
 \eeqlb
 which also goes to $0$ as $n\to\infty $ for any $p>1$. 
 Taking the preceding two limits back into \eqref{eqn.560}, we have as $n\to\infty$, 
 \beqnn
  \sup_{k=0,1,\cdots, [n^\beta T]}\sup_{\delta\in[0, 1/n^{\beta}]} \big| J^{(n)}_{1,1}(k/n^{\beta},\delta)\big| \overset{\rm p}\to 0.
 \eeqnn
 
 \medskip
 
 \textbf{(ii)} We now prove that the second term on the right-hand side of \eqref{eqn.559} converges in probability to $0$.  
 Similarly as in \eqref{eqn.558}, we also have 
 \beqnn
  \sup_{\delta\in[0, 1/n^{\beta}]} \big|J^{(n)}_{1,2}(k/n^{\beta},\delta) \big| 
  \ar\leq\ar   \int_0^{k/n^{\beta}} X_{\zeta}^{(n)} (s )ds  \int_{ n(k/n^{\beta}-s)}^{n((k+1)/n^{\beta}-s)}  \frac{2n}{c}\cdot    \boldsymbol{e}_{c}(y)\, dy \cr
  \ar\ar + \Bigg| \int_0^{k/n^{\beta}}\int_0^{X_{\zeta}^{(n)} (s-)}  \int_{ n(k/n^{\beta}-s)}^{n((k+1)/n^{\beta}-s)}  \frac{1}{n} \, \widetilde{N}^{(n)}_c (ds,dz,dy)  \Bigg| ,
 \eeqnn
 and then 
 \beqlb\label{eqn.562}
  \lefteqn{ \sup_{k=0,1,\cdots, [n^\beta T]}\sup_{\delta\in[0, 1/n^{\beta}]} \big| J^{(n)}_{1,2}(k/n^{\beta},\delta)\big| }\ar\ar\cr
  \ar\leq\ar \sup_{k=0,1,\cdots, [n^\beta T]} \int_0^{k/n^{\beta}} X_{\zeta}^{(n)} (s )ds  \int_{ n(k/n^{\beta}-s)}^{n((k+1)/n^{\beta}-s)}  \frac{2n}{c}\cdot    \boldsymbol{e}_{c}(y)\, dy\cr
  \ar\ar + \sup_{k=0,1,\cdots, [n^\beta T]}  \bigg| \int_0^{k/n^{\beta}}\int_0^{X_{\zeta}^{(n)} (s-)}  \int_{ n(k/n^{\beta}-s)}^{n((k+1)/n^{\beta}-s)}  \frac{1}{n}  \widetilde{N}^{(n)}_c (ds,dz,dy)  \bigg|. \qquad
 \eeqlb
 As before, it suffices to prove the convergence in probability of the two terms on the right-hand side to $0$.  
 Firstly, taking expectation of the first term and then using the next inequality 
 \beqnn
  \sup_{k\geq 0}\sup_{s\geq 0} \int_{ n(k/n^{\beta}-s)}^{n((k+1)/n^{\beta}-s)}  \frac{2n}{c}\cdot    \boldsymbol{e}_{c}(y)  \,dy
  \leq  \frac{2n^{2-\beta}}{c^2},
 \eeqnn
 uniformly in $n\geq 1$ and Lemma~\ref{Lemma.Moment}, we can get that
 \beqnn
  \mathbf{E}\Bigg[ \sup_{k=0,1,\cdots, [n^\beta T]} \int_0^{k/n^{\beta}} X_{\zeta}^{(n)} (s )ds  \int_{ n(k/n^{\beta}-s)}^{n((k+1)/n^{\beta}-s)}  \frac{2n}{c}\cdot    \boldsymbol{e}_{c}(y)\, dy \Bigg]
  \ar\leq\ar  \frac{C}{n^{\beta-2}}  \int_0^{T}\mathbf{E}\big[X_{\zeta}^{(n)} (s )\big] ds
  \leq \frac{C}{n^{\beta-2}}, 
 \eeqnn
 which vanishes as $n\to\infty$ since $\beta>2$.  
 For the second term, similarly as in \eqref{eqn.563} and \eqref{eqn.564} we also have 
 \beqnn
  \lefteqn{ \mathbf{P}\bigg(\sup_{k=0,1,\cdots, [n^\beta T]} \Bigg| \int_0^{k/n^{\beta}}\int_0^{X_{\zeta}^{(n)} (s-)}  \int_{ n(k/n^{\beta}-s)}^{n((k+1)/n^{\beta}-s)}  \frac{1}{n}  \widetilde{N}^{(n)}_c (ds,dz,dy)  \Bigg|\geq \epsilon  \bigg)}\ar\ar\cr
  \ar\leq\ar \frac{1}{\epsilon^{2p}} \sum_{k=0}^{[n^\beta T]}\mathbf{E}\Bigg[ \bigg| \int_0^{k/n^{\beta}}\int_0^{X_{\zeta}^{(n)} (s-)}  \int_{ n(k/n^{\beta}-s)}^{n((k+1)/n^{\beta}-s)}  \frac{1}{n}  \widetilde{N}^{(n)}_c (ds,dz,dy)  \bigg|^{2p} \Bigg]\cr
  \ar\leq\ar 
  \frac{C}{\epsilon^{2p}}\sum_{k=0}^{[n^\beta T]} \Bigg( \bigg|  \int_0^{k/n^{\beta}} ds \int_{ n(k/n^{\beta}-s)}^{n((k+1)/n^{\beta}-s)}     \boldsymbol{e}_{c}(y)\, dy \bigg|^{p} 
  +  \int_0^{k/n^{\beta}} ds \int_{ n(k/n^{\beta}-s)}^{n((k+1)/n^{\beta}-s)}    \frac{ \boldsymbol{e}_{c}(y)}{n^{2p-2}} \, dy \Bigg),  
 \eeqnn
 which can be bounded by $ C\cdot  \big( n^{p(1- \beta)+\beta} + n^{3-2p} \big) \to 0$ as $n\to\infty$ for all $p>2$. 
 Taking the preceding two limits back into \eqref{eqn.562}, we have as $n\to\infty$,
 \beqnn
  \sup_{k=0,1,\cdots, [n^\theta T]}\sup_{\delta\in[0, 1/n^{\theta}]} \big| J^{(n)}_{1,2}(k/n^{\theta},\delta)\big| \overset{\rm p}\to 0. 
 \eeqnn  
 \qed

 We now prove the $C$-tightness of the continuous processes $\{\mathbf{I}^{(n)}_{3,2}\}_{n\geq 1}$ by using the well-known Kolmogorov-Chentsov tightness criterion for continuous processes; see Problem~4.11 in \cite[p.64]{KaratzasShreve1988}.  
 The Kolmogorov-Chentsov tightness criterion states that a sequence of continuous processes $\{X^{(n)}\}_{n\geq 1}$ is tight if for each $T\geq 0$, there exist constants $C,\beta,p,\rho>0 $ such that for any $\delta\in(0,1)$, 
 \beqlb\label{eqn.KolTight}
  \sup_{n\geq 1}\mathbf{E}\Big[\big|X^{(n)}(0)\big|^{\beta}\Big]<\infty
  \quad\mbox{and}\quad  
  \sup_{n\geq 1} \sup_{t\in[0,T]} \mathbf{E}\Big[\big|X^{(n)}(t+\delta)-X^{(n)}(t)\big|^{p}\Big]\leq C\cdot  \delta^{1+\rho} . 
 \eeqlb
   
 \begin{lemma}
  The sequence $\{\mathbf{I}^{(n)}_{3,2}\}_{n\geq 1}$ is $C$-tight.
 \end{lemma}
 \proof 
 It suffices to prove that the two inequalities in \eqref{eqn.KolTight} hold for $\{\mathbf{I}^{(n)}_{3,2}\}_{n\geq 1}$.  
 The first one is obvious because $\mathbf{I}^{(n)}_{3,2}(0)\overset{\rm a.s.}=0$ for all $n\geq 1$. We now prove the second one. 
 For each $t\geq 0$ and $\delta\in(0,1)$, by  the change of variables we have
 \beqnn
  \mathbf{I}^{(n)}_{3,2}(t+\delta)-\mathbf{I}^{(n)}_{3,2}(t) = J^{(n)}_{2,1}(t,\delta)+J^{(n)}_{2,2}(t,\delta), 
 \eeqnn
 where the two terms on the right side are given by
 \beqnn 
  J^{(n)}_{2,1}(t,\delta)\ar:=\ar \int_t^{t+\delta}\int_0^{X_{\zeta}^{(n)} (s-)} \int_0^\infty  \frac{1}{n}\cdot \Big(\int_{n(t+\delta-s) -y}^{n(t+\delta-s)} R_{\varPi_n}( r ) dr\Big)  \widetilde{N}^{(n)}_c (ds,dz,dy) , \cr
  J^{(n)}_{2,2}(t,\delta)\ar:=\ar   \int_0^t\int_0^{X_{\zeta}^{(n)} (s-)} \int_0^\infty  \frac{1}{n}\cdot \Big(\int_{n(t-s) -y}^{n(t-s)}\big( R_{\varPi_n}( r+ n\delta ) - R_{\varPi_n} ( r )\big) dr\Big)   \widetilde{N}^{(n)}_c (ds,dz,dy) . 
 \eeqnn
 For any $p\geq 2$, by the power inequality we have 
 \beqlb \label{eqn.572}
  \mathbf{E}\Big[ \big|\mathbf{I}^{(n)}_{3,2}(t+\delta) - \mathbf{I}^{(n)}_{3,2}(t) \big|^{2p} \Big]  
  \leq 
  C\cdot \mathbf{E}\Big[ \big| J^{(n)}_{2,1}(t,\delta)\big|^{2p} \Big] +C\cdot\mathbf{E}\Big[ \big|J^{(n)}_{2,2}(t,\delta) \big|^{2p} \Big] , 
 \eeqlb
 for some constant $C>0$ relying only on $p$. 
 By using \eqref{BDG} along with Lemma~\ref{Lemma.Moment} as well as the two facts that $c\gamma_n\leq 1$ and $\theta_n\in (0,1)$, and then using the change of variables, 
 \beqlb\label{eqn.670}
  \sup_{t\in[0,T]} \mathbf{E}\Big[ \big|J^{(n)}_{2,1}(t,\delta)\big|^{2p}  \Big] 
  \ar\leq \ar C\cdot \Big( \big|A^{(n)}_2(\delta)\big|^p + A^{(n)}_{2p}(\delta) \Big),
 \eeqlb
 for some constant $C>0$ independent of $n$ and  $\delta$, where
 \beqnn
  A^{(n)}_{2k}(\delta):=   n^2\int_0^{\delta} ds  \int_0^\infty \Big(\int_{s -y/n}^{s} R_{\varPi_n}(nr) dr\Big)^{2k}   \boldsymbol{e}_{c}(y)\,  dy,\quad k\geq 1.
 \eeqnn
 By the non-negativity and uniform boundedness of $R_{\varPi_n}$; see Lemma~\ref{Lemma.UpperBoundR},
 \beqnn
  A^{(n)}_{2k}(\delta) 
  \ar\leq\ar  \sup_{t\in[0,\delta]} \sup_{y\geq 0}\Big(\int_{t-y/n}^t R_{\varPi_n}(nr) dr\Big)^{2k-2}   \int_0^{\delta} ds  \int_0^\infty  n^2\Big(\int_{s -y/n}^{s} R_{\varPi_n}(nr) dr\Big)^2   \boldsymbol{e}_{c}(y)\,  dy   \cr
  \ar\leq\ar \Big(\frac{\delta}{c}\Big)^{2k-2}\cdot \int_0^\delta ds \int_0^\infty  \Big(\frac{y}{c}\Big)^2  \boldsymbol{e}_{c}(y)\,  dy, 
 \eeqnn
 which can be bounded by $C\cdot \delta^{2k-1}$ uniformly in $n\geq 1$ and $\delta>0$.  
 Plugging this result with $k=1$ and $p$ back into \eqref{eqn.670}, there exists a constant $C>0$ such that for any $\delta>0$,
 \beqlb\label{eqn.673}
  \sup_{n\geq 1}\sup_{t\in[0,T]} \mathbf{E}\Big[ \big|J^{(n)}_{2,1}(t,\delta)\big|^{2p}  \Big] 
  \leq C\cdot \big( \delta^p+ \delta^{2p-1} \big).
 \eeqlb
 For $J^{(n)}_{2,2}(t,\delta)$, similarly as in \eqref{eqn.670} we also have 
 \beqlb\label{eqn.573}
  \sup_{t\in[0,T]} \mathbf{E}\Big[ \big|J^{(n)}_{2,2}(t,\delta)\big|^{2p}  \Big]  
  \ar\leq\ar C\cdot \Big( \big| B^{(n)}_2(\delta)\big|^p + B^{(n)}_{2p}(\delta) \Big),
 \eeqlb
 uniformly in $n\geq 1$ and $\delta\in (0,1)$, where 
 \beqlb\label{eqn.672}
  B^{(n)}_{2k}(\delta):= n^2\int_0^\infty  \boldsymbol{e}_{c}(y)\, dy \int_0^T \Big(\int_{s-y/n}^{s}\big( R_{\varPi_n}\big( n(r+\delta)\big) - R_{\varPi_n} ( nr )\big) dr\Big)^{2k} \,ds,
  \quad k\geq1.
 \eeqlb
 We first consider the inner integrals in $B^{(n)}_{2k}(\delta)$.  
 By using H\"older's inequality and Corollary~\ref{Coro.511}, 
 \beqlb\label{eqn.671}
  \lefteqn{\int_0^T \Big(\int_{s-y/n}^{s}\big( R_{\varPi_n}\big( n(r+\delta)\big) - R_{\varPi_n} ( nr )\big) dr\Big)^{2k} \,ds}\ar\ar\cr
  \ar\leq\ar  \int_0^T \Big(\frac{y}{n}\int_{s-y/n}^{s}\big( R_{\varPi_n}\big( n(r+\delta)\big) - R_{\varPi_n} ( nr )\big)^2 dr\Big)^{k} \,ds\cr
  \ar\leq\ar  \sup_{t\in[0,T]} \sup_{y\geq 0} \Big( \int_{t-y/n}^{t}\big( R_{\varPi_n}\big( n(r+\delta)\big) - R_{\varPi_n} ( nr )\big)^2 dr\Big)^{k-1} \cr
  \ar\ar \times \int_0^T \Big(\frac{y}{n}\Big)^{k}\int_{s-y/n}^{s}\big( R_{\varPi_n}\big( n(r+\delta)\big) - R_{\varPi_n} ( nr )\big)^2 dr \,ds \cr
  \ar\leq\ar  \Big( \int_{-\infty}^{T}\big( R_{\varPi_n}\big( n(r+\delta)\big) - R_{\varPi_n} ( nr )\big)^2 dr\Big)^{k-1} \cr
  \ar\ar \times \int_0^T \Big(\frac{y}{n}\Big)^{k}\int_{s-y/n}^{s}\big( R_{\varPi_n}\big( n(r+\delta)\big) - R_{\varPi_n} ( nr )\big)^2 dr \,ds \cr
  \ar\leq\ar C\cdot \delta^{k-1}\cdot \Big(\frac{y}{n}\Big)^{k}\cdot \int_{-\infty}^T \int_{s-y/n}^{s}\big( R_{\varPi_n}\big( n(r+\delta)\big) - R_{\varPi_n} ( nr )\big)^2 dr \,ds,
 \eeqlb
 uniformly in $n\geq 1$ and $\delta\in(0,1)$. 
 By the fact that $R_{\varPi_n}(t)=0$ for any $t<0$ and then using Fubini's theorem to the last double integral, 
 \beqnn
  \int_{-\infty}^T \int_{s-y/n}^{s}\big( R_{\varPi_n}\big( n(r+\delta)\big) - R_{\varPi_n} ( nr )\big)^2 dr \,ds 
  \ar=\ar \int_{-\delta}^T \int_{(s-y/n)\vee (-\delta)}^{s}\big( R_{\varPi_n}\big( n(r+\delta)\big) - R_{\varPi_n} ( nr )\big)^2 dr \,ds\cr
  \ar\leq\ar \int_{-\delta}^T \big| R_{\varPi_n}\big(n(s+ \delta)\big) - R_{\varPi_n} (ns)\big|^2 \,ds  \int_{s}^{s+y/n}  \, dr\cr
  \ar\leq\ar \frac{y}{n} \cdot \int_{-\infty}^T \big| R_{\varPi_n}\big(n(s+ \delta)\big) - R_{\varPi_n} (ns)\big|^2 \,ds , 
 \eeqnn
 which can be bounded by $C\cdot \frac{y}{n} \cdot \delta$ uniformly in $n\geq 1$, $y>0$ and $\delta\in(0,1)$. 
 Plugging this into the last term in \eqref{eqn.671} and then taking it back into \eqref{eqn.672},  there exists  a constant $C>0$ such that for any $n\geq 1$ and $\delta\in(0,1)$,
 \beqnn
  B^{(n)}_{2k}(\delta)\leq C\cdot n^2\int_0^\infty  \Big(\frac{y}{n}\Big)^{k+1} \boldsymbol{e}_{c}(y)\, dy \cdot \delta^k
  \leq \frac{C}{n^{k-1}}  \cdot \delta^k,
 \eeqnn 
 Taking this result with $k=1$ and $p$  back into \eqref{eqn.573} yields that uniformly in $\delta\in(0,1)$, 
 \beqnn
  \sup_{n\geq 1}\sup_{t\in[0,T]} \mathbf{E}\Big[ \big|J^{(n)}_{2,2}(t,\delta)\big|^{2p}  \Big] \leq C\cdot  \delta^p  .
 \eeqnn
 Finally, plugging this and \eqref{eqn.673} back into \eqref{eqn.572} shows that the sequence $\{\mathbf{I}^{(n)}_{3,2}\}_{n\geq 1}$ satisfies the second inequality in \eqref{eqn.KolTight} and hence is $C$-tight.  
 \qed
 
 \begin{corollary}\label{Corollary.528}
  The sequence $\{\mathbf{I}^{(n)}_{3}\}_{n\geq 1}$ is $C$-tight. 
 \end{corollary}
%
%
  \subsubsection{$C$-tightness of $\{\mathbf{I}^{(n)}_{4}\}_{n\geq 1}$} 
  
 In this section, we prove the $C$-tightness of the sequence $\{\mathbf{I}^{(n)}_{4}\}_{n\geq 1}$. 
 By \eqref{Eqn.R.n} and the change of variables, we can split the process $\mathbf{I}^{(n)}_{4}$ into the following two terms 
 \beqlb
  \mathbf{I}^{(n)}_{4,1}(t)\ar:=\ar \int_0^t\int_0^{X_{\zeta}^{(n)} (s-)} \int_0^\infty  \frac{1}{n}\cdot \mathbf{1}_{\{y\geq n(t-s) \}}  \widetilde{N}^{(n)}_{\Lambda} (ds,dz,dy), \label{eqn.667}\\
  \mathbf{I}^{(n)}_{4,2}(t)\ar:=\ar \int_0^t\int_0^{X_{\zeta}^{(n)} (s-)} \int_0^\infty  \Big(\int_{t-s-y/n }^{t-s} R_{\varPi_n}(nr) dr\Big)\widetilde{N}^{(n)}_{\Lambda} (ds,dz, dy). \label{eqn.668} 
 \eeqlb
 
 \begin{lemma}\label{Lemma.5291}
 For each $T\geq 0$, we have $\sup_{t\in[0,T]}\big|\mathbf{I}^{(n)}_{4,1}(t) \big| \overset{\rm p}\to 0$ as $n\to\infty$.   
 \end{lemma} 
 \proof  
 By the fact that $\widetilde{N}^{(n)}_{\Lambda} (ds,dz,dy)= N^{(n)}_{\Lambda} (ds,dz,dy)- \gamma_n\cdot n^2\theta_n \cdot  \boldsymbol{e}_{c}*d\Lambda_n(y)\,  ds \, dz \, dy $,  
 \beqnn
  \sup_{t\in[0,T]}\big| \mathbf{I}^{(n)}_{4,1}(t) \big| 
  \ar\leq\ar   \int_0^T\int_0^{X_{\zeta}^{(n)} (s-)} \int_0^\infty  \frac{1}{n}\cdot  N^{(n)}_{\Lambda} (ds,dz,dy)  +   \gamma_n\cdot n\theta_n  \int_0^T X_{\zeta}^{(n)} (s) \,  ds  .
 \eeqnn
 Taking expectations on both sides of this inequality and then using Lemma~\ref{Lemma.Moment} as well  as the fact that $c\gamma_n\leq 1$, there exists a constant $C>0$ such that for any $n\geq 1$, 
 \beqnn
  \mathbf{E}\bigg[   \sup_{t\in[0,T]}\big| \mathbf{I}^{(n)}_{4,1}  (t) \big| \bigg]  
  \ar\leq\ar C\cdot  n\theta_n ,
 \eeqnn
 which goes to $0$ as $n\to\infty$. 
 The desired locally uniform convergence in probability holds. 
 \qed

 \begin{lemma}\label{Lemma.529}
  The sequence $\{\mathbf{I}^{(n)}_{4,2}\}_{n\geq 1}$ is $C$-tight.
 \end{lemma}
 \proof 
 It suffices to prove the two inequalities in \eqref{eqn.KolTight} hold for $\{\mathbf{I}^{(n)}_{4,2}\}_{n\geq 1}$. 
 The first inequality  holds obviously since $ \mathbf{I}^{(n)}_{4,2}(0)\overset{\rm a.s.}=0$ for each $n\geq 1$. 
 For the second inequality, by the power inequality we have for $p>2$ and any $\delta\in(0,1)$,
 \beqlb\label{eqn.666}
  \mathbf{E}\Big[ \big|\mathbf{I}^{(n)}_{4,2}(t+\delta)-\mathbf{I}^{(n)}_{4,2}(t) \big|^{2p} \Big] 
  \leq C\cdot \mathbf{E}\Big[ \big| J_{4,1}^{(n)}(t,\delta) \big|^{2p} \Big] + C\cdot \mathbf{E}\Big[ \big|J_{4,2}^{(n)}(t,\delta) \big|^{2p} \Big]
 \eeqlb
 for some constant $C>0$ depending only on $p$, where 
 \beqnn
  J_{4,1}^{(n)}(t,\delta)\ar:=\ar \int_t^{t+h}\int_0^{X_{\zeta}^{(n)} (s-)} \int_0^\infty  \Big(\int_{t+\delta-s-y/n }^{t+\delta-s} R_{\varPi_n}(nr) dr\Big)\widetilde{N}^{(n)}_{\Lambda} (ds,dz, dy),\cr
  J_{4,2}^{(n)}(t,\delta)\ar:=\ar \int_0^t\int_0^{X_{\zeta}^{(n)} (s-)} \int_0^\infty  \Big(\int_{t+\delta-s-y/n }^{t+\delta-s} R_{\varPi_n}(nr) dr-\int_{t-s-y/n }^{t-s} R_{\varPi_n}(nr) dr\Big)\widetilde{N}^{(n)}_{\Lambda} (ds,dz, dy)\cr
  \ar=\ar \int_0^t\int_0^{X_{\zeta}^{(n)} (s-)} \int_0^\infty  \Big(\int_{t-s-y/n }^{t-s} \big(R_{\varPi_n}\big(n(r+\delta)\big)- R_{\varPi_n}(nr)\big) dr\Big)\widetilde{N}^{(n)}_{\Lambda} (ds,dz, dy) . 
 \eeqnn 
 By using \eqref{BDG}, Lemma~\ref{Lemma.Moment} and then the change of variables, we have 
 \beqnn
  \sup_{t\in[0,T]}\mathbf{E}\Big[ \big|J_{4,1}^{(n)}(t,\delta)  \big|^{2p} \Big] 
  \ar\leq\ar C \cdot \bigg|n^2\theta_n \int_0^\delta \int_0^\infty  \Big(\int_{s-y/n }^{s} R_{\varPi_n}(nr) dr\Big)^2   \boldsymbol{e}_{c}*d\Lambda_n(y)\, ds\, dy \bigg|^p\cr
  \ar\ar + C\cdot n^2 \theta_n  \int_0^\delta   \int_0^\infty  \Big(\int_{s-y/n }^{s} R_{\varPi_n}(nr) dr\Big)^{2p}   \boldsymbol{e}_{c}*d\Lambda_n(y)\,  ds\, dy .
 \eeqnn
 By Lemma~\ref{Lemma.UpperBoundR}, we have $ \int_{s-y/n }^{s} R_{\varPi_n}(nr) dr \leq \frac{1}{c}(s\wedge \frac{y}{n})$ for any $s,y\geq 0$ and then
 \beqlb\label{eqn.665}
  \sup_{t\in[0,T]} \mathbf{E}\Big[ \big|J_{4,1}^{(n)}(t,\delta)  \big|^{2p} \Big] 
  \ar\leq\ar  C \cdot \bigg|n^2\theta_n \int_0^\delta \int_0^\infty   \Big(s\wedge \frac{y}{n}\Big)^2   \boldsymbol{e}_{c}*d\Lambda_n(y)\, ds\,dy \bigg|^p\cr
  \ar\ar + C\cdot n^2\theta_n  \int_0^\delta   \int_0^\infty  \Big(s\wedge \frac{y}{n}\Big)^{2p}   \boldsymbol{e}_{c}*d\Lambda_n(y)\,  ds\, dy\cr
  \ar\leq\ar   C \cdot \bigg|n^2\theta_n   \int_0^\infty   \Big(\delta\wedge \frac{y}{n}\Big)^2   \boldsymbol{e}_{c}*d\Lambda_n(y)\, ds\,dy \bigg|^p\cdot \delta^p\cr
  \ar\ar + C\cdot n^2\theta_n   \int_0^\infty  \Big(\delta\wedge \frac{y}{n}\Big)^{2}   \boldsymbol{e}_{c}*d\Lambda_n(y)\,dy \cdot \delta^{2p-1} \cr
  \ar\ar\cr
  \ar\leq\ar  C\cdot \big( \delta^p + \delta^{2p-1} \big) ,
 \eeqlb
 uniformly in $n\geq 1$ and $\delta\in (0,1)$. 
 Here the last inequality follows from \eqref{eqn.544}. 
 Similarly,  
 \beqlb\label{eqn.664}
 \sup_{t\in[0,T]}\mathbf{E}\Big[ \big| J_{4,2}^{(n)}(t,\delta) \big|^{2p} \Big] 
 \leq C\cdot \Big( \big|I^{(n)}_2(\delta)\big|^p + I^{(n)}_{2p}(\delta) \Big),
 \eeqlb
 for some constant $C>0$ independent of $n$, where 
 \beqnn
 I^{(n)}_{2k}(\delta)
 :=  n^2 \theta_n  \int_0^T ds \int_0^\infty  \Big(\int_{s-y/n }^{s}\big( R_{\varPi_n}\big(n(r+\delta)\big)  - R_{\varPi_n}(nr)\big) dr\Big)^{2k}  \boldsymbol{e}_{c}*d\Lambda_n(  y)\, dy,
 \quad k\geq 1.   
 \eeqnn
 By H\"older's inequality and Corollary~\ref{Coro.511}, we have uniformly in $s\in[0,T]$, $\delta\in (0,1)$ and $y\geq 0$,
 \beqnn
  \lefteqn{\Big(\int_{s-y/n }^{s}\big( R_{\varPi_n}\big(n(r+\delta)\big)  - R_{\varPi_n}(nr)\big) dr\Big)^{2k}}\ar\ar\cr
  \ar\leq\ar \Big( \int_{s-y/n }^{s}\big| R_{\varPi_n}\big(n(r+\delta)\big)  - R_{\varPi_n}(nr)\big|  dr\Big)^{2}  \Big(\int_{-1}^T\big| R_{\varPi_n}\big(n(r+\delta)\big)  - R_{\varPi_n}(nr)\big| dr\Big)^{2k-2}\cr
  \ar\leq\ar \Big(\int_{s-y/n }^{s}\big| R_{\varPi_n}\big(n(r+\delta)\big)  - R_{\varPi_n}(nr)\big| dr\Big)^{2} \Big((T+1)\int_{-1}^T\big| R_{\varPi_n}\big(n(r+\delta)\big)  - R_{\varPi_n}(nr)\big|^2 dr\Big)^{k-1}\cr
  \ar\leq\ar C \cdot \Big(\int_{s-y/n }^{s}\big| R_{\varPi_n}\big(n(r+\delta)\big)  - R_{\varPi_n}(nr)\big| dr\Big)^{2}\cdot \delta^{k-1}. 
 \eeqnn 
 Plugging this back into $ I^{(n)}_{2k}$, we have 
 \beqnn
  I^{(n)}_{2k}(\delta)
  \ar\leq\ar C\cdot \delta^{k-1}  \int_0^\infty   n^2 \theta_n\cdot \boldsymbol{e}_{c}*d\Lambda_n(  y)\, dy \cr
  \ar\ar\quad \times \int_0^T \Big|\int_{s-y/n }^{s}\big( R_{\varPi_n}(n(r+\delta))  - R_{\varPi_n}(nr)\big) dr\Big|^{2}ds  \cr
  \ar\ar\cr
  \ar=\ar   C\cdot \delta^{k-1} \cdot \big( A^{(n)}_1(\delta) + A^{(n)}_2(\delta)\big),
 \eeqnn
 where
 \beqnn
  A^{(n)}_1(\delta) \ar:=\ar n^2 \theta_n  \int_{nT}^\infty \,   \boldsymbol{e}_{c}*d\Lambda_n(y)\,dy  \int_0^T  \Big|\int_{s-y/n }^{s}\big( R_{\varPi_n}(n(r+\delta)) -R_{\varPi_n}(nr)\big) dr\Big|^2\, ds ,\cr
  A^{(n)}_2(\delta) \ar:=\ar n^2 \theta_n \int_0^{nT} \boldsymbol{e}_{c}*d\Lambda_n(y)\, dy \int_0^T  \Big|\int_{s-y/n }^{s}\big( R_{\varPi_n}(n(r+\delta)) -R_{\varPi_n}(nr)\big) dr\Big|^2\, ds .  
 \eeqnn
 By the fact that $R_{\varPi_n}(x)=0$ for $x<0$ and H\"older's inequality,
 \beqnn
  A^{(n)}_1(\delta) \ar\leq\ar n^2 \theta_n \int_{nT}^\infty  T  \cdot \boldsymbol{e}_{c}*d\Lambda_n(y)\, dy \cdot   \Big|\int_{-1}^T\big( R_{\varPi_n}(n(r+\delta)) -R_{\varPi_n}(nr)\big) dr\Big|^2 \cr
  \ar\leq\ar   n^2 \theta_n  \int_{nT}^\infty  (T+1)^2  \boldsymbol{e}_{c}*d\Lambda_n(y)\, dy \cdot \int_{-1}^T\big| R_{\varPi_n}(n(r+\delta)) -R_{\varPi_n}(nr)\big|^2 dr  \cr
  \ar\leq\ar n^2 \theta_n  \int_0^\infty  \Big((T+1)\wedge \frac{y}{n}\Big)^2  \boldsymbol{e}_{c}*d\Lambda_n(y)\, dy \cdot \int_{-1}^T\big| R_{\varPi_n}(n(r+\delta)) -R_{\varPi_n}(nr)\big|^2 dr,
 \eeqnn
 which can be bounded by $C\cdot \delta$ uniformly in $n\geq 1$ and $\delta\in (0,1)$ because of  Corollary~\ref{Coro.511} and \eqref{eqn.544}. 
 For $ A^{(n)}_2(\delta)$, using H\"older's inequality again gives that
 \beqnn
  A^{(n)}_2(\delta)\ar\leq\ar n^2 \theta_n \cdot\int_0^{nT} \Big(\int_0^T ds \int_{s-y/n }^{s}\big( R_{\varPi_n}(n(r+\delta)) -R_{\varPi_n}(nr)\big)^2 \,dr   \Big)  \cdot \frac{y}{n}\cdot \boldsymbol{e}_{c}*d\Lambda_n(y)\, dy ,
 \eeqnn
 which can be further divided into the following two parts
 \beqnn
  A^{(n)}_{2,1}(\delta)\ar:= \ar n^2\theta_n \cdot \int_0^{nT}  \Big(\int_0^Tds \int_{(s-y/n)\wedge 0 }^{0}\big| R_{\varPi_n}(n(r+\delta)) -R_{\varPi_n}(nr)\big|^2\, dr   \Big)\cdot \frac{y}{n} \cdot \boldsymbol{e}_{c}*d\Lambda_n(y) \, dy , \cr
  A^{(n)}_{2,2}(\delta)\ar:= \ar n^2\theta_n \cdot \int_0^{nT}  \Big(\int_0^T ds \int_{(s-y/n)\vee 0 }^{s}\big| R_{\varPi_n}(n(r+\delta)) -R_{\varPi_n}(nr)\big|^2 \,dr  \Big)\cdot \frac{y}{n} \cdot \boldsymbol{e}_{c}*d\Lambda_n(y) \, dy. 
 \eeqnn
 The fact that $R_{\varPi_n}(x)=0$ for $x<0$ induces that 
 \beqnn
  \int_0^Tds \int_{(s-y/n)\wedge 0 }^{0}\big| R_{\varPi_n}(n(r+\delta)) -R_{\varPi_n}(nr)\big|^2\, dr
  \ar\leq\ar \int_0^{y/n} \int_{-(s\wedge \delta) }^{0}\big( R_{\varPi_n}(n(r+\delta))  \big)^2 dr  ds \cr 
  \ar\leq\ar  \int_0^{y/n} \frac{s\wedge \delta}{c^2}ds
  \leq \frac{\delta  }{c^2 } \cdot \frac{y}{n} ,
 \eeqnn
 which along with Proposition~\ref{Proposition.510} yields that
 \beqlb\label{eqn.663}
  A^{(n)}_{2,1}(\delta) 
  \ar\leq\ar  \frac{\delta}{c^2} \cdot n^2\theta_n \cdot \int_0^{nT}    \Big(\frac{y}{n}\Big)^2 \cdot \boldsymbol{e}_{c}*d\Lambda_n(y)  dy  \cr
  \ar\leq\ar \frac{\delta}{c^2} \cdot n^2\theta_n \cdot \int_0^\infty   \Big(T\wedge \frac{y}{n}\Big)^2 \cdot \boldsymbol{e}_{c}*d\Lambda_n(y)  dy  
  \leq C\cdot \delta.
 \eeqlb
 Moreover, an application of Fubini's theorem to the inner integral of $A^{(n)}_{2,2}(\delta)$ induces that
 \beqnn
  \int_0^T ds \int_{(s-y/n)\vee 0}^{s}\big( R_{\varPi_n}(n(r+\delta)) -R_{\varPi_n}(nr)\big)^2 \, dr  
  \ar\leq\ar \frac{y}{n}  \int_0^T  \big| R_{\varPi_n}(n(s+\delta)) -R_{\varPi_n}(ns)\big|^2\, ds .
 \eeqnn
 Similarly as in \eqref{eqn.663}, we also have $A^{(n)}_{2,2}(\delta) \leq C\cdot \delta$ uniformly in $n\geq 1$ and $\delta\in(0,1)$.  
 Consequently,
 \beqnn
  I^{(n)}_{2k}(\delta)
  \ar\leq\ar C\cdot \delta^{k-1} \cdot \Big(A^{(n)}_{1}(\delta)+A^{(n)}_{2,1}(\delta)+A^{(n)}_{2,2}(\delta) \Big) \leq C\cdot \delta^k.
 \eeqnn
 Plugging this result with $k=1$ and $p$ into \eqref{eqn.664} yields that uniformly in $\delta\in(0,1)$, 
 \beqnn
  \sup_{n\geq 1}\sup_{t\in[0,T]}\mathbf{E}\Big[ \big| J_{4,2}^{(n)}(t,\delta) \big|^{2p} \Big] 
  \leq C\cdot  \delta^p  ,
 \eeqnn
 Taking this and \eqref{eqn.665} back into \eqref{eqn.666}, we have 
 \beqnn
  \sup_{n\geq 1}\sup_{t\in[0,T]} \mathbf{E}\Big[ \big| \mathbf{I}^{(n)}_{4,2}(t+\delta)-\mathbf{I}^{(n)}_{4,2} \big|^{2p} \Big] 
  \leq  C\cdot \delta^p 
 \eeqnn
 and then the sequence $\{\mathbf{I}^{(n)}_{4,2}\}_{n\geq 1}$ is $C$-tight. 
 \qed 
  
 \begin{corollary}\label{Coro.431}
  The sequence $\{\mathbf{I}^{(n)}_{4}\}_{n\geq 1}$ is $C$-tight. 
 \end{corollary} 
 
 \subsection{Limit characterization}
 
 Armed with the $C$-tightness results proved in the last section, we now turn to characterize their limits, which will be used to derive the stochastic Volterra representation for the process $L^\xi_\zeta$.  
 By Lemma~\ref{Lemma.I1} and \ref{Lemma.I2}, cluster points of $\{\mathbf{I}^{(n)}_{1}\}_{n\geq 1}$ and $\{\mathbf{I}^{(n)}_{2}\}_{n\geq 1}$ are clear and it remains to characterize the limits of $\{\mathbf{I}^{(n)}_{3}\}_{n\geq 1}$ and $\{\mathbf{I}^{(n)}_{4}\}_{n\geq 1}$. 
 Note that integrands in \eqref{eqn.Jc} and \eqref{eqn.JLambda} vary in $n$, which makes their limit characterizations much challenging.  
 To overcome this difficulty, we provide some good approximations for them in the next section.  
 
  \subsubsection{Good approximations for $\{ \mathbf{I}_3^{(n)}\}_{n\geq 1}$ and $\{\mathbf{I}_4^{(n)}\}_{n\geq 1}$}
 
 The second limit in Lemma~\ref{Lemma.ConR} tells that $R_{\varPi_n}(nr,y) \sim y\cdot W'(r)$ as $n\to\infty$. Hence it is sensible to conjecture that the sequence $\{\mathbf{I}_{3,2}^{(n)}\}_{n\geq 1}$ can be asymptotically approximated by  $\{\boldsymbol{\mathcal{I}}_3^{(n)}\}_{n\geq 1}$ with 
 \beqnn
 \boldsymbol{\mathcal{I}}^{(n)}_{3}(t):= \int_0^t\int_0^{X_{\zeta}^{(n)} (s-)} \int_0^\infty  \frac{y}{n}\cdot W'(t-s) \widetilde{N}^{(n)}_c (ds,dz,dy) ,\quad t\geq 0.
 \eeqnn
 For each $n\geq 1$, we consider the following two integrated processes
 \beqnn
 I_{ \mathbf{I}_{3}}^{(n)}(t) := \int_0^t  \mathbf{I}^{(n)}_{3}(s)\, ds 
 \quad \mbox{and}\quad 
 I_{\boldsymbol{\mathcal{I}}_{3}}^{(n)}(t):=\int_0^t \boldsymbol{\mathcal{I}}^{(n)}_{3}(s)\, ds ,\quad t\geq 0,
 \eeqnn

 \begin{lemma}\label{Lemma.432}
 	The following hold.
 	\begin{enumerate}
 		\item[(1)] The sequence $\{(I_{ \mathbf{I}_{3}}^{(n)},I_{\boldsymbol{\mathcal{I}}_{3}}^{(n)}) \}_{n\geq 1}$ 
 		is $C$-tight. 
 		
 		\item[(2)]  For each $T\geq 0$, we have $	\sup_{t\in[0,T]}\big|I_{ \mathbf{I}_{3}}^{(n)}(t)-I_{\boldsymbol{\mathcal{I}}_{3}}^{(n)}(t)\big|   \overset{\rm p}\to  0$ as $n\to\infty$. 
 		
 	\end{enumerate}

 \end{lemma}
 \proof 
 The first claim follows directly from the second one and Corollary~\ref{Corollary.528}. We now prove the second claim. 
 For convention, we denote by $\{\boldsymbol{\varepsilon}_3^{(n)}\}_{n\geq 1}$ the error process, i.e.,
 \beqnn
 \boldsymbol{\varepsilon}_3^{(n)}(t)
 := \mathbf{I}_3^{(n)}(t)- \boldsymbol{\mathcal{I}}^{(n)}_{3}(t),
 \quad t\geq 0.
 \eeqnn
 By the triangle inequality and the decomposition of $\mathbf{I}_{3}^{(n)}$; see \eqref{eqn.5331}-\eqref{eqn.5332}, 
 \beqnn
 \sup_{t\in[0,T]}\big|I_{ \mathbf{I}_{3}}^{(n)}(t)-I_{\boldsymbol{\mathcal{I}}_{3}}^{(n)}(t)\big|  
 \leq \big\| \boldsymbol{\varepsilon}_3^{(n)} \big\|_{L^1_T} 
 \leq \big\| \mathbf{I}_{3,1}^{(n)}\big\|_{L^1_T} +\big\| \mathbf{I}_{3,2}^{(n)}- \boldsymbol{\mathcal{I}}^{(n)}_{3}\big\|_{L^1_T}.
 \eeqnn
 Firstly, Lemma~\ref{Lemma.526} along with the continuous mapping theorem induces that 
 \beqnn
  \int_0^t \mathbf{I}_{3,1}^{(n)}(s)\, ds \to 0,
 \eeqnn
 weakly in $C(\mathbb{R}_+;\mathbb{R})$ as $n\to\infty$, which directly yields that 
 $\big\| \mathbf{I}_{3,1}^{(n)}\big\|_{L^1_T}\to 0 $ in distribution and hence in probability.
 Secondly, by using \eqref{BDG} along with Lemma~\ref{Lemma.UpperBoundR} and then the change of variables, there exists a constant $C>0$ independent of $n$ such that
 \beqnn
 \sup_{t\in[0,T]}\mathbf{E}\Big[ \big| \mathbf{I}_{3,2}^{(n)} (t)- \boldsymbol{\mathcal{I}}^{(n)}_{3}(t) \big|^2 \Big]
 \ar\leq \ar C  \sup_{t\in[0,T]} \int_0^t  \int_0^\infty  \Big(   R_n\big(n(t-s),y\big)-y\cdot W'(t-s) \Big)^2   \boldsymbol{e}_{c}(y)\, ds \,dy\cr
 \ar\leq\ar C  \int_0^\infty \boldsymbol{e}_{c}(y) \, dy  \int_0^T\big(   R_n\big(ns,y\big)-y\cdot W'(s) \big)^2 \,ds ,
 \eeqnn
 which goes to $0$ by the dominated convergence theorem and \eqref{eqn.546}.  By H\"older's inequality,  Fubini's theorem and the dominated convergence theorem, we have 
 \beqnn
 \mathbf{E}\Big[ \big\| \mathbf{I}_{3,2}^{(n)}- \boldsymbol{\mathcal{I}}^{(n)}_{3}\big\|_{L^1_T}^2 \Big]
 \ar\leq\ar T\cdot \mathbf{E}\Big[ \big\| \mathbf{I}_{3,2}^{(n)}- \boldsymbol{\mathcal{I}}^{(n)}_{3}\big\|_{L^2_T}^2 \Big]
 = T\cdot  \int_0^T \mathbf{E}\Big[ \big| \mathbf{I}_{3,2}^{(n)} (t)- \boldsymbol{\mathcal{I}}^{(n)}_{3}(t) \big|^2 \Big]dt \to 0,
 \eeqnn
 as $n\to\infty$ and hence $\big\| \mathbf{I}_{3,2}^{(n)}- \boldsymbol{\mathcal{I}}^{(n)}_{3}\big\|_{L^1_T} \to 0$ in distribution and hence in probability.
 \qed 
 
 By the first limit in Lemma~\ref{Lemma.ConR}, we also see that $R_{\varPi_n}(nr) $ can be well-approximated by $ W'(r)$ as $n\to\infty$, which allows us to approximate  the sequence $\{\mathbf{I}_{4,2}^{(n)}\}_{n\geq 1}$  by  $\{\boldsymbol{\mathcal{I}}_4^{(n)}\}_{n\geq 1}$ with 
 \beqlb\label{eqn.669}
 \boldsymbol{\mathcal{I}}_4^{(n)}(t) \ar:=\ar \int_0^t\int_0^{X_{\zeta}^{(n)} (s-)} \int_0^\infty   \Big(\int_{t-s-y/n}^{t-s}W'(r)dr\Big) \widetilde{N}^{(n)}_{\Lambda} (ds,dz,dy),
 \quad t\geq 0 .  
 \eeqlb
 
 \begin{lemma}\label{Lemma.433}
 	The following hold.
 	\begin{enumerate}
 		\item[(1)] The sequence $\{\boldsymbol{\mathcal{I}}_4^{(n)} \}_{n\geq 1}$ is $C$-tight.
 		
 		\item[(2)] For each $T\geq 0$, we have $\sup_{t\in[0,T]}\big|\mathbf{I}_{4}^{(n)} (t)- \boldsymbol{\mathcal{I}}^{(n)}_{4}(t)\big| \overset{\rm p}\to 0$  as $n\to\infty$. 
 	\end{enumerate} 
 \end{lemma}
 \proof  
 The first claim can be proved by repeating the proof of Lemma~\ref{Lemma.529} with $R_{\varPi_n}(nr) $ replaced by $W'(r)$. 
 For the second claim, by the decomposition of $\mathbf{I}_4^{(n)}$ (see \eqref{eqn.667}-\eqref{eqn.668}) and Lemma~\ref{Lemma.5291}, it remains to prove that as $n\to\infty$,
 \beqnn
 \sup_{t\in[0,T]}\big|\mathbf{I}_{4,2}^{(n)}(t) -\boldsymbol{\mathcal{I}}_4^{(n)}(t)\big| \overset{\rm p}\to 0 . 
 \eeqnn 
 By the first claim and Lemma~\ref{Lemma.529}, it suffices to prove that for any $t\geq 0$,
 \beqnn
 \lim_{n\to\infty} \mathbf{E}\Big[ \big| \mathbf{I}_{4,2}^{(n)}(t) -\boldsymbol{\mathcal{I}}_4^{(n)}(t) \big|^2  \Big] = 0. 
 \eeqnn 
 Indeed, by \eqref{eqn.667} and \eqref{eqn.668}, 
 \beqnn
 \mathbf{I}_{4,2}^{(n)}(t) -\boldsymbol{\mathcal{I}}_4^{(n)}(t)
 \ar=\ar \int_0^t\int_0^{X_{\zeta}^{(n)} (s-)} \int_0^\infty   \Big(\int_{t-s-y/n}^{t-s}\big(R_{\varPi_n}(nr)-W'(r)\big)dr\Big) \widetilde{N}^{(n)}_{\Lambda} (ds,dz,dy). 
 \eeqnn
 By using \eqref{BDG}, Lemma~\ref{Lemma.Moment} and the change of variables, 
 \beqnn
 \mathbf{E}\Big[ \big|\mathbf{I}_{4,2}^{(n)}(t) -\boldsymbol{\mathcal{I}}_4^{(n)}(t) \big|^2  \Big]
 \ar\leq\ar  C\cdot  n^2 \theta_n  \int_0^\infty  \boldsymbol{e}_{c}*d\Lambda_n(y)\, dy  \int_0^t \Big(\int_{s-y/n}^{s}\big(R_{\varPi_n}(nr)-W'(r)\big)dr\Big)^2 \, ds .
 \eeqnn
 Similar to  establishing the estimates of $I^{(n)}_{2k}$ in the proof of Lemma~\ref{Lemma.529}, we have
 \beqnn
 \mathbf{E}\Big[ \big| \mathbf{I}_{4,2}^{(n)}(t) -\boldsymbol{\mathcal{I}}_4^{(n)}(t) \big|^2  \Big] 
 \ar\leq\ar n^2 \theta_n  \int_0^\infty  \Big( t\wedge \frac{y}{n} \Big)^2  \boldsymbol{e}_{c}*d\Lambda_n(y)\, dy \cdot \int_0^T \big(R_{\varPi_n}(nr)-W'(r)\big)^2 dr , 
 \eeqnn
 which vanishes as $n\to\infty$; see \eqref{eqn.546} and \eqref{eqn.544}. 
 \qed

 \subsubsection{Convergence of stochastic Volterra integrals}
 
 By Lemma~\ref{Lemma.432} and \ref{Lemma.433}, we can characterize the limits of  $\{I_{\mathbf{I}_{3}}^{(n)}\}_{n\geq 1}$ and  $\{\mathbf{I}_{4}^{(n)}\}_{n\geq 1}$ by  
 identifying the weak limit of the two sequences $\{I_{\boldsymbol{\mathcal{I}}_{3}}^{(n)} \}_{n\geq 1}$ and $\{\boldsymbol{\mathcal{I}}_4^{(n)} \}_{n\geq 1}$.
 Applying the stochastic Fubini's theorem along with \eqref{eqn.201} to $I_{\boldsymbol{\mathcal{I}}_3}^{(n)}$, we have
 \beqlb\label{eqn.442}
 I_{\boldsymbol{\mathcal{I}}_3}^{(n)}(t)
 \ar=\ar \int_0^t dr\int_0^r\int_0^{X_{\zeta}^{(n)} (s-)} \int_0^\infty  \frac{y}{n}\cdot W'(r-s) \widetilde{N}^{(n)}_c (ds,dz,dy)\cr
 \ar=\ar \int_0^t\int_0^{X_{\zeta}^{(n)} (s-)} \int_0^\infty  \frac{y}{n}\cdot W(t-s) \widetilde{N}^{(n)}_c (ds,dz,dy) ,\quad t\geq 0, 
 \eeqlb
 which together with \eqref{eqn.669} motivates us to establish the weak convergence theorem for the following stochastic Volterra integral 
 \beqnn
  Z^{(n)} (t):= \int_0^t \int_0^{X^{(n)} (s-)}\int_0^\infty  G \big(t-s,y/n \big)\widetilde{N}^{(n)} (ds,dz,dy), \quad t \geq 0.
 \eeqnn 
 for some process $X^{(n)} \in D(\mathbb{R}_+;\mathbb{R}_+)$, function $G$ on $(0,\infty)^2$ and compensated Poisson random measure $\widetilde N^{(n)}(ds,dz,dy)$ on $(0,\infty)^3$ with intensity $ds\,dz \, \mu_n (dy)$ and $\mu_n (dy)$ being a $\sigma$-finite measure on $(0,\infty)$. 
    
 \begin{lemma}\label{Lemma.434}
 Assume that $G(t,y)=F(t)\cdot y$ for some locally bounded function $F$ on $\mathbb{R}_+$. 
 If  
 \beqlb\label{eqn.627}
   \int_0^\infty \frac{y^2}{n^2}\, \mu_n(dy) \to \sigma^2  \in (0,\infty)
 \quad \mbox{and}\quad 
 X^{(n)} \to X^* \in C(\mathbb{R}_+;\mathbb{R}_+) \mbox{ weakly in $D(\mathbb{R}_+;\mathbb{R}_+)$},
 \eeqlb 
 as $n\to\infty$,   we have $Z^{(n)}  \overset{\rm f.d.d.}\longrightarrow Z^*$ with the limit process $Z^* $ given by
 \beqnn
 Z^*(t):= \int_0^t \int_0^{X^* (s)} F(t-s) B^* (ds,dz), \quad t\geq 0,
 \eeqnn  
 where $B^* (ds,dz)$ is Gaussian white noise with intensity $ \sigma^2\cdot ds\, dz$. 
 \end{lemma} 
 \proof It suffices to prove that  for any $d\in \mathbb{Z}_+$ and $0\leq t_1<\cdots<t_d$,
 \beqlb\label{eqn.624}
 \big( Z^{(n)}(t_1),\cdots, Z^{(n)}(t_d) \big) \overset{\rm d}\to \big( Z^*(t_1),\cdots, Z^*(t_d) \big),
 \eeqlb
 as $n\to\infty$. 
 However, the dependence of the integrand $G$ on the time variable does not allow us prove it by using the convergence theorem for It\^o's integrals directly. 
 To ``drop'' the dependence of  integrand on the time variable, we consider instead the integral processes 
 \beqlb
 Z_i^{(n)}(t)
 \ar :=\ar \int_0^t\int_0^{X^{(n)}(s-)} \int_0^\infty F(t_i-s  )\cdot \frac{y}{n}\, \widetilde N^{(n)}(ds,dz,dy), \label{eqn.625} \\
 Z^*_i(t)
 \ar: =\ar \int_0^t \int_0^{X(s)} F(t_i-s) B^*(ds,dz), \quad t\geq 0,\, i=1,\cdots, d.\label{eqn.626}
 \eeqlb
 It is obvious that $Z^{(n)}(t_i)\overset{\rm a.s.}=Z_i^{(n)}(t_i)$ and $Z^*(t_i)\overset{\rm a.s.}=Z^*_i(t_i)$ for all $i=1,\cdots,d$. 
 By Proposition~3.14 in \cite[p.349]{JacodShiryaev2003} and the continuity of $\big( Z^*_{1}, ..., Z^*_d \big)$, we can obtain \eqref{eqn.624} by proving that 
 \beqlb\label{eqn.401}
 \big( Z^{(n)}_1, ..., Z^{(n)}_d \big) 
 \to
 \big( Z^*_{1}, ..., Z^*_d \big),
 \eeqlb
 weakly in $D(\mathbb{R}_+;\mathbb{R}^d)$ as $n\to\infty$. 
 
 We now prove \eqref{eqn.401} by using Lemma~\ref{kurztheorem} with  $\mathbb{H}=L^2(\mathbb{R}_+;\mathbb{R})$. 
 We first write \eqref{eqn.625} and \eqref{eqn.626} into the form of stochastic integrals driven by $(L^2(\mathbb{R}_+;\mathbb{R}))^\#$-martingales. 
 For $i=1,\cdots ,d$, we define a mapping $\mathbf{F}_i$ from $ D(\mathbb{R}_+;\mathbb{R}_+) \times [0,\infty) $ into $ L^2(\mathbb{R}_+;\mathbb{R}) $ by
 \beqlb\label{eqn.629}
  (x,s) \mapsto \mathbf{F}_i(x,s)(z) := F(t_i-s)\cdot \mathbf{1}_{\{0<z\leq x(s)\}}.
 \eeqlb
 and introduce two standard $(L^2(\mathbb{R}_+;\mathbb{R}))^\#$-martingales  
 \beqlb\label{eqn.420}
 {\bf \widetilde N}^{(n)}(t) :=   \int_0^\infty \frac{y}{n}\,\widetilde N^{(n)}\big((0,t],dz,dy\big) 
 \quad \mbox{and} \quad 
 \mathbf{B}^*(t) := B^*\big((0,t], dz\big),
 \quad t \geq 0.
 \eeqlb
 On terms of these notation, we can represent the two integral processes $Z^{(n)}_i$ and $Z^*_i$ as 
 \beqnn
 Z^{(n)}_i(t)= \mathbf{F}_i(X^{(n)},-)\cdot d{\bf \widetilde N}^{(n)}(t)
 \quad\mbox{and}\quad 
 Z^*_i(t)= \mathbf{F}_i(X^*,-)\cdot d\mathbf{B}^*(t),\quad t\geq 0. 
 \eeqnn
 By Lemma~\ref{kurztheorem}, the weak convergence \eqref{eqn.401} holds if we can prove the following two claims.
 \begin{enumerate}
 	\item[$\bullet$] The sequence $\{{\bf \widetilde N}^{(n)}\}_{n\geq 1}$ is uniformly tight;
 	
 	\item[$\bullet$] $( \mathbf{F}_1(X^{(n)},\cdot),\cdots, \mathbf{F}_d(X^{(n)},\cdot),{\bf \widetilde N}^{(n)})\Rightarrow (\mathbf{F}_1(X^*,\cdot),\cdots, \mathbf{F}_d(X^*,\cdot),\mathbf{B}^*)$.
 \end{enumerate}
 
 \textit{Uniform tightness.} Recall Definition~\ref{Definition.A1}(1). It suffices to prove that for any $T\geq 0$, there exists a constant $C>0$ such that for any c\`adl\`ag $L^2(\mathbb{R}_+;\mathbb{R})$-valued process $U$ with $\sup_{t\in[0,T]}\|U(t)\|_{L^2}\leq 1$ a.s.,
 \beqnn
 \sup_{n\geq 1} \mathbf{E}\Big[\sup_{t\in[0,T]} \big|U(-)\cdot d {\bf \widetilde N}^{(n)}(t) \big|^2 \Big] \leq C. 
 \eeqnn
 By \eqref{eqn.420} and then the Burkholder-Davis-Gundy inequality, 
 \beqnn
 U(-)\cdot d {\bf \widetilde N}^{(n)}(t) \ar=\ar \int_0^t \int_0^\infty U(s-,z) \int_0^\infty \frac{y}{n}\,\widetilde N^{(n)}\big(ds,dz,dy\big)
 \eeqnn
 and there exists a constant $C>0$ independent of $n,T,U$ such that
 \beqnn
  \mathbf{E}\bigg[\sup_{t\in[0,T]} \big|U(-)\cdot d {\bf \widetilde N}^{(n)}(t) \big|^2 \bigg]
  \ar\leq\ar C \cdot \mathbf{E}\bigg[\int_0^T \int_0^\infty  \big|U(s-,z)\big|^2 \int_0^\infty  \frac{y^2}{n^2}\, N^{(n)}\big(ds,dz,dy\big)\bigg]\cr
  \ar=\ar C \cdot \mathbf{E}\bigg[\int_0^T ds \int_0^\infty  \big|U(s-,z)\big|^2 dz \cdot  \int_0^\infty \frac{y^2}{n^2}  \, \mu_n(dy)\bigg],
 \eeqnn
 which is bounded uniformly in $n$ and $U$ because of the assumption of $U$ and the first limit in \eqref{eqn.627}. 
 Hence the first claim holds. 
 
 \textit{Weak convergence.} By Definition~\ref{Definition.A1}(2), it suffices to prove that  for any $k\geq 1$ and $f_1,\cdots,f_k\in L^2(\mathbb{R}_+;\mathbb{R})$, 
 \beqnn
 \big( \mathbf{F}_1(X^{(n)},\cdot),\cdots, \mathbf{F}_d(X^{(n)},\cdot),{\bf \widetilde N}^{(n)}(f_1),\cdots , {\bf \widetilde N}^{(n)}(f_k)\big) 
 \ar\to\ar \big(\mathbf{F}_1(X^*,\cdot),\cdots, \mathbf{F}_d(X^*,\cdot),\mathbf{B}^*(f_1),\cdots, \mathbf{B}^*(f_k) \big),
 \eeqnn
 weakly in $D(\mathbb{R}_+;(L^2(\mathbb{R}_+;\mathbb{R}))^d\times \mathbb{R}^k)$ as $n\to\infty$.
 For simplicity, we just prove it with $d=k=1$.  
 The general case can be proved in the same way.
 Note that $\mathbf{B}^*(f_1) \in C(\mathbb{R}_+;\mathbb{R}) $ and the continuity of $X^*$ is inherited by $\mathbf{F}_1(\xi^*,\cdot)$.
 Moreover, the function $f_1$ can be write as $f_1^+-f_1^-$  with $f_1^+,f_1^- \in L^2(\mathbb{R}_+;\mathbb{R}_+)$.  
 By Corollary~3.33 in \cite[p.353]{JacodShiryaev2003}, it suffices to prove that
 \beqlb \label{eqn.630}
 \mathbf{F}_1 \big(X^{(n)},\cdot \big)\to \mathbf{F}_1 \big(X^*,\cdot \big)
 \quad \mbox{and}\quad
 {\bf \widetilde N}^{(n)}(f_1) \to \mathbf{B}^*(f_1),
 \eeqlb
 weakly in $D(\mathbb{R}_+;L^2(\mathbb{R}_+;\mathbb{R}))$ and  $D(\mathbb{R}_+;\mathbb{R})$ separately for $f_1\in L^2(\mathbb{R}_+;\mathbb{R}_+)$. 
 By  the Skorokhod representation theorem, we may assume $X^{(n)} \overset{\rm a.s.}\to X^*$ in $D(\mathbb{R}_+;\mathbb{R}_+)$. 
 The continuity of $X^*$ induces that
 \beqnn
 \sup_{t\in[0,T]} \big|X^{(n)}(t)-X^*(t) \big| \overset{\rm a.s.}\to 0
 \eeqnn
 as $n \to 0$  for any $T\geq 0$; see Proposition~1.17(b) in \cite[p.328]{JacodShiryaev2003}. This along with \eqref{eqn.629} induces that
 \beqlb\label{eqn.423}
  \sup_{t\in[0,T]} \big\|\mathbf{F}_1(X^{(n)},t) - \mathbf{F}_1(X^*,t) \big\|_{L^2}^2 
 \ar=\ar \sup_{t\in[0,T]} |F(t_i-t)|^2\cdot \int_0^\infty \big| \mathbf{1}_{\{0<z\leq X^{(n)}(t)\}}-   \mathbf{1}_{\{0<z\leq X^*(t)\}}\big|^2\,dz \cr
 \ar=\ar \sup_{t\in[0,T]} |F(t_i-t)|^2\cdot \big|X^{(n)}(t)-X^*(t) \big|
 \eeqlb
 which goes to $0$ almost surely as $n\to\infty$ and hence the first limit in \eqref{eqn.630} holds. 
 For the second one, note that 
 \beqnn
 {\bf \widetilde N}^{(n)}(f_1,t)= \int_0^t \int_0^\infty \int_0^\infty f_1(z) \cdot \frac{y}{n} \, \widetilde N^{(n)}\big(ds,dz,dy\big), \quad t\geq 0,
 \eeqnn
 is a spectrally positive L\'evy process with Laplace exponent 
 \beqnn
 \int_0^\infty \int_0^\infty \bigg(\exp\Big\{-\lambda \cdot f_1(z) \cdot \frac{y}{n} \Big\}-1+ \lambda \cdot f_1(z) \cdot \frac{y}{n}\bigg)\, dz\, \mu_n(dy),\quad \lambda \geq 0,
 \eeqnn
 which, by the first limit in \eqref{eqn.627}, converges as $n\to\infty$ to 
 \beqnn
 \frac{\sigma^2}{2}\int_0^\infty \big|f_1(z)  \big|^2 \,dz \cdot \lambda^2.
 \eeqnn
 This together with Corollary~4.3 in \cite[p.440]{JacodShiryaev2003} and Theorem III-7 in \cite{KarouiMeleard1990} induces that 
 \beqnn
 {\bf \widetilde N}^{(n)}(f_1)\to \sigma\cdot \Big(\int_0^\infty \big|f_1(z)  \big|^2 \,dz \Big)^{1/2}\cdot B  \overset{\rm d}= \mathbf{B}^*(f_1),
 \eeqnn
 weakly in $D(\mathbb{R}_+;\mathbb{R})$ as $n\to\infty$, where $B$ is a standard Brownian motion. 
 The second claim holds. 
 \qed

 \begin{lemma}\label{Lemma.435}
  Assume that there exists a $\sigma$-finite measure $m(dy)$ on $(0,\infty)$ such that for any $T>0$ and non-negative measurable function $f$ on $\mathbb{R}_+$, 
  \beqlb\label{eqn.6271}
   \sup_{t\in[0,T]}\int_0^\infty \big| G(t,y) \big|^2 \, m(dy) < \infty 
    \quad \mbox{and} \quad 
   \sup_{n\geq 1} \int_0^\infty f(y) \,\mu_n(n\cdot dy) \leq  \int_0^\infty f(y) \, m(dy)  . 
  \eeqlb 
  If $X^{(n)}  \to X^* \in C(\mathbb{R}_+;\mathbb{R}_+) $ weakly in $D(\mathbb{R}_+;\mathbb{R}_+)$ and $\mu_n(n\cdot dy)\to \mu^*(dy) $ vaguely as $n\to\infty$, we have $Z^{(n)} \overset{\rm f.d.d.}\longrightarrow Z^*$ with the limit process $Z^*$ given by
  \beqnn
  Z^*(t):= \int_0^t \int_0^{X^*(s)} \int_0^\infty G(t-s,y)\, N^*(ds,dz,dy), \quad t\geq 0,
  \eeqnn  
  where $N^*(ds,dz,dy)$ is a Poisson random measure on $(0,\infty)^3$ with intensity $ ds\,dz\,\mu^*(dy)$.
 \end{lemma}
 \proof 
  Similarly as in the proof of Lemma~\ref{Lemma.434}, for any $d\geq 1$ and $0\leq t_1<\cdots<t_d$ it suffices to prove that 
 $
 \big( Z^{(n)}_1, ..., Z^{(n)}_d \big) 
 \to
 \big( Z^*_{1}, ..., Z^*_d \big) 
 $
 weakly in $D(\mathbb{R}_+;\mathbb{R}^k)$ as $n\to\infty$, where
 \beqlb 
 Z_i^{(n)}(t)
 \ar :=\ar \int_0^t\int_0^{X^{(n)}(s-)} \int_0^\infty G(t_i-s,y ) \, \widetilde N^{(n)}(ds,dz,n\cdot dy), \label{eqn.6251} \\ 
 Z^*_i(t)
 \ar: =\ar \int_0^t \int_0^{X(s)} \int_0^\infty G(t_i-s,y) N^*(ds,dz,dy), \quad t\geq 0,\, i=1,\cdots, d.\label{eqn.6261}
 \eeqlb
 Let $\mathbf{m}(dz,dy) :=dz\, m(dy) $ and $L^2_{\mathbf{m}}(\mathbb{R}_+^2;\mathbb{R})$ be the Hilbert space of all square integrable functions on $(0,\infty)^2$ with respect to $\mathbf{m}(dz,dy)$ endowed with the norm $\|\cdot\|_{L^2_{\mathbf{m}}}$. 
 For $i=1,\cdots ,d$, we define a mapping $\mathbf{G}_i$ from $ D(\mathbb{R}_+;\mathbb{R}_+) \times [0,\infty) $ into $L^2_{\mathbf{m}}(\mathbb{R}_+^2;\mathbb{R}) $ by
 \beqnn
 (x,s) \mapsto \mathbf{G}_i(x,s)(z,y) := G(t_i-s,y)\cdot \mathbf{1}_{\{0<z\leq x(s)\}}.
 \eeqnn
 and introduce two standard $(L^2_{\mathbf{m}}(\mathbb{R}_+^2;\mathbb{R}))^\#$-martingales  
 \beqlb\label{eqn.421}
 {\bf \widetilde N}^{(n)}(t) :=  \widetilde N^{(n)}\big((0,t],dz,n\cdot dy\big) 
 \quad \mbox{and} \quad 
 {\bf \widetilde N}^* (t) :=  \widetilde N^* \big((0,t],dz,dy\big),
 \quad t \geq 0.
 \eeqlb 
 Similarly, we can write the two equations \eqref{eqn.6251} and \eqref{eqn.6261} as 
 \beqnn
 Z^{(n)}_i(t)= \mathbf{G}_i(X^{(n)},-)\cdot d{\bf \widetilde N}^{(n)}(t)
 \quad\mbox{and}\quad 
 Z^*_i(t)= \mathbf{G}_i(X^*,-)\cdot d {\bf \widetilde N}^*(t),\quad t\geq 0. 
 \eeqnn
 By Lemma~\ref{kurztheorem} with  $\mathbb{H}=L^2_{\mathbf{m}}(\mathbb{R}_+^2;\mathbb{R})$, it suffices to prove the uniform tightness of $\{{\bf \widetilde N}^{(n)}\}_{n\geq 1}$ and 
 \beqnn
 ( \mathbf{G}_1(X^{(n)},\cdot),\cdots, \mathbf{G}_d(X^{(n)},\cdot),{\bf \widetilde N}^{(n)})\Rightarrow (\mathbf{G}_1(X^*,\cdot),\cdots, \mathbf{G}_d(X^*,\cdot),{\bf \widetilde N}^*).
 \eeqnn 
 
 \textit{Uniform tightness.} Recall Definition~\ref{Definition.A1}(1). It suffices to prove that for any $T\geq 0$, there exists a constant $C>0$ such that for any c\`adl\`ag $L^2_{\mathbf{m}}(\mathbb{R}_+^2;\mathbb{R})$-valued process $U$ with $\sup_{t\in[0,T]}\|U(t)\|_{L^2_{\mathbf{m}}}\leq 1$ a.s.,
 \beqnn
 \sup_{n\geq 1}\mathbf{E}\Big[\sup_{t\in[0,T]} \big|U(-)\cdot d {\bf \widetilde N}^{(n)}(t) \big|^2 \Big] \leq C. 
 \eeqnn
 By \eqref{eqn.421} and then the Burkholder-Davis-Gundy inequality,  
 \beqnn
  U(-)\cdot d {\bf \widetilde N}^{(n)}(t) \ar=\ar \int_0^t \int_0^\infty  \int_0^\infty U(s-,z,y)\,\widetilde N^{(n)}\big(ds,dz,n\cdot dy\big)
 \eeqnn
 and there exists a constant $C>0$ independent of $n,T,U$ such that
 \beqnn
 \mathbf{E}\Big[\sup_{t\in[0,T]} \big|U(-)\cdot d {\bf \widetilde N}^{(n)}(t) \big|^2 \Big]
 \ar\leq\ar C \cdot \mathbf{E}\Big[\int_0^T \int_0^\infty   \int_0^\infty  \big|U(s-,z,y)\big|^2\, N^{(n)}\big(ds,dz,n\cdot dy\big)\Big]\cr
 \ar=\ar C \cdot \mathbf{E}\Big[\int_0^T ds \int_0^\infty dz  \int_0^\infty \big|U(s-,z,y)\big|^2 \,  \mu_n(n\cdot dy)\Big].
 \eeqnn
 which, by the first inequality in \eqref{eqn.6271} and the fact that $\sup_{t\in[0,T]}\|U(t)\|_{L^2_{\mathbf{m}}}\leq 1$ a.s, can be bounded uniformly in $n$ and $U$ by
 \beqnn
 C \cdot \mathbf{E}\Big[\int_0^T ds \int_0^\infty  dz \int_0^\infty \big|U(s-,z,y)\big|^2 \, m(dy)\Big] \leq C\cdot T.
 \eeqnn

 \textit{Weak convergence.} 
 Similarly as in the end of the proof of Lemma~\ref{Lemma.434}, it suffices to prove that
 \beqnn
 \mathbf{G}_1(\xi^{(n)},\cdot)\to \mathbf{G}_1(\xi^*,\cdot)
 \quad \mbox{and}\quad
 {\bf \widetilde N}^{(n)}(f_1) \to {\bf \widetilde N}^*(f_1)
 \eeqnn
 weakly in $D(\mathbb{R}_+;L^2(\mathbb{R}_+^2;\mathbb{R}))$ and  $D(\mathbb{R}_+;\mathbb{R})$ separately for $f_1\in L^2_{\mathbf{m}}(\mathbb{R}_+^2;\mathbb{R}_+)$. 
 The first limit can be proved similarly as in \eqref{eqn.423}. 
 For the second one, note that 
 \beqnn
 {\bf \widetilde N}^{(n)}(f_1,t)= \int_0^t \int_0^\infty \int_0^\infty f_1(z,y)   \, \widetilde N^{(n)}\big(ds,dz,n\cdot dy\big), \quad t\geq 0,
 \eeqnn
 is a spectrally positive L\'evy process with Laplace exponent 
 \beqnn
 \int_0^\infty \int_0^\infty \big(e^{-\lambda \cdot f_1(z,y) }-1+ \lambda \cdot f_1(z,y) \big)\, dz\, \mu_n(n\cdot dy),\quad \lambda \geq 0,
 \eeqnn
 which, by the assumption that $\mu_n(n\cdot dy)\to \mu^*(dy) $ vaguely as $n\to\infty$, converges as $n\to\infty$ to 
 \beqnn
 \int_0^\infty \int_0^\infty \big(e^{-\lambda \cdot f_1(z,y) }-1+ \lambda \cdot f_1(z,y) \big)\, dz\, \mu^*( dy).
 \eeqnn
 This together with Theorem III-7 in \cite{KarouiMeleard1990} induces that 
 ${\bf \widetilde N}^{(n)}(f_1)\to  \mathbf{N}^*(f_1)$.
 \qed

 \subsubsection{Proof of Theorem~\ref{Main.Thm01}}
 The proof is carried out in the next three steps. 
 The uniqueness of solution to \eqref{MainThm.SVE} will be proved in the end of Section~\ref{Sec.LaplaceFunctionals}. 
 
 \textit{Step 1: $C$-tightness.} An application of Corollary~3.33 in \cite[p.353]{JacodShiryaev2003} along with Lemma~\ref{Lemma.I1},
 \ref{Lemma.I2}, \ref{Lemma.432}, \ref{Lemma.433} and Corollary~\ref{Corollary.528},  
 \ref{Coro.431} immediately induces the $C$-tightness of the sequence 
 \beqlb\label{eqn.690}
 \big\{ \big(X_\zeta^{(n)}, \mathbf{I}_1^{(n)}, \mathbf{I}_2^{(n)}, \mathbf{I}_3^{(n)}, I_{\mathbf{I}_3}^{(n)}, I_{\boldsymbol{\mathcal{I}}_3}^{(n)},    \mathbf{I}_4^{(n)}, \boldsymbol{\mathcal{I}}^{(n)}_{4} \big)\big\}_{n\geq 1}. 
 \eeqlb
 
 \textit{Step 2: Limit characterization.} Assume that $\big(X_\zeta , \mathbf{I}_1 , \mathbf{I}_2, \mathbf{I}_3, I_{\mathbf{I}_3}, I_{\boldsymbol{\mathcal{I}}_3} ,  \mathbf{I}_4, \boldsymbol{\mathcal{I}}_{4} \big)$ is a cluster point of the sequence \eqref{eqn.690}. 
 We now characterize  the limit processes one-by-one.
 
 \begin{enumerate}
 	\item[$\bullet$]  By Lemma~\ref{Lemma.I1} and \ref{Lemma.I2}, we first have 
 	\beqlb\label{eqn.707}
 	\mathbf{I}_1(t)\ar=\ar \zeta\cdot c\cdot W(t) 
 	\quad \mbox{and}\quad 
 	\mathbf{I}_2(t) =\int_0^\zeta\int_0^\infty \big(W(t)-W(t-y) \big) N_0(dz,dy),\quad t\geq 0. \qquad 
 	\eeqlb
 	
 	\item[$\bullet$] An application of Lemma~\ref{Lemma.434} along with the two facts that $X_{\zeta}^{(n)}\to X_\zeta$ weakly in $D(\mathbb{R}_+;\mathbb{R}_+)$ and 
 	\beqnn
 	\int_0^\infty \frac{y^2}{n^2} \cdot n^2\gamma_n(1-\theta_n ) \cdot \boldsymbol{e}_{c}(y)\,\cdot ds\,dz\,dy \sim \frac{1}{c}\int_0^\infty  y^2   \boldsymbol{e}_{c}(y)\,  dy \to 2c
 	\eeqnn
 	to the stochastic integral \eqref{eqn.442} induces that the limit process $I_{\boldsymbol{\mathcal{I}}_3}$ admits the  representation
 	\beqnn
 	I_{\boldsymbol{\mathcal{I}}_3}(t)= \int_0^t\int_0^{X_{\zeta}(s)}   W(t-s) B_c (ds,dz),
 	\quad t\geq 0. 
 	\eeqnn  
 	Using the stochastic Fubini's theorem \eqref{eqn.SFT01} along with \eqref{eqn.201} to the preceding stochastic integral, we have 
 	\beqnn
 	I_{\boldsymbol{\mathcal{I}}_3}(t)
 	\ar=\ar \int_0^t\int_0^{X_{\zeta}(s)} \int_0^{t-s}  W(r)dr B_c (ds,dz) 
 	= \int_0^t dr \int_0^r\int_0^{X_{\zeta}(s)}   W'(r-s) B_c (ds,dz) . 
 	\eeqnn 
 	Moreover, applying the continuous mapping theorem to the weak convergence $(I_{\mathbf{I}_3}^{(n)}, I_{\boldsymbol{\mathcal{I}}_3}^{(n)})\to (I_{\mathbf{I}_3} , I_{\boldsymbol{\mathcal{I}}_3} )$ and then using Lemma~\ref{Lemma.432}(2), 
 	\beqnn
 	I_{\mathbf{I}_3}(t)= \int_0^t  \mathbf{I}_3(s)ds=  I_{\boldsymbol{\mathcal{I}}_3}(t),\quad t\geq 0,
 	\eeqnn
 	which along with the continuity of $\mathbf{I}_3$ induces that 
 	\beqlb\label{eqn.708}
 	\mathbf{I}_3(t) = \int_0^t\int_0^{X_{\zeta}(s)}   W'(t-s) B_c (ds,dz),\quad t\geq 0. 
 	\eeqlb
 	
 	\item[$\bullet$] By Lemma~\ref{Lemma.433}, we have $\mathbf{I}_4=\boldsymbol{\mathcal{I}}_{4} $. 
 	By the proof of Proposition~\ref{Prop.501}, 
 	\beqnn
 	n^3\gamma_n\theta_n \cdot  \boldsymbol{e}_{c}*d\Lambda_n(ny)\,\cdot ds\,dz\,dy \to ds\,dz\,\nu(dy),
 	\eeqnn
 	vaguely $n\to\infty$. 
 	By Lemma~\ref{Lemma.435} and the weak convergence $X_\zeta^{(n)}\to X_\zeta$, 
 	\beqlb\label{eqn.709}
 	\mathbf{I}_4(t)=\boldsymbol{\mathcal{I}}_{4}(t) = \int_0^t \int_0^{X_\zeta(s)}\int_0^\infty \big(W(t-s)-W(t-s-y) \big) \widetilde{N}_\nu(ds,dz,dy),\quad t\geq 0.
 	\eeqlb 
 \end{enumerate}

 \medskip
  \textit{Step 3: Stochastic Volterra equation.}
 Notice that the  limit process $(X_\zeta,\mathbf{I}_1,\mathbf{I}_2,\mathbf{I}_3,\mathbf{I}_4)$ is determined uniquely in distribution by $X_\zeta$, which, by Corollary~\ref{ConvergenceLocalTime02}, always equals in distribution to $L^\xi_\zeta$.
 In conclusion, we have as $n\to\infty$,
 \beqnn
 \big(X_\zeta^{(n)}, \mathbf{I}_1^{(n)}, \mathbf{I}_2^{(n)}, \mathbf{I}_3^{(n)}, \mathbf{I}_4^{(n)} \big) \to \big(X_\zeta, \mathbf{I}_1, \mathbf{I}_2, \mathbf{I}_3, \mathbf{I}_4\big) \in C(\mathbb{R}_+;\mathbb{R}_+^3\times \mathbb{R}^2),
 \eeqnn
 weakly in $D(\mathbb{R}_+;\mathbb{R}_+^3\times \mathbb{R}^2)$. 
 By Proposition 2.4 in \cite[p.339]{JacodShiryaev2003} and the continuous mapping theorem, 
 \beqnn
 \sup_{t\in[0,T]} \bigg| X_\zeta(t)- \sum_{i=1}^4\mathbf{I}_i(t)  \bigg| 
 \ar \overset{\rm d}=\ar \lim_{n\to\infty} \sup_{t\in[0,T]} \bigg| X_\zeta^{(n)}(t)-\sum_{i=1}^4 \mathbf{I}_i^{(n)}(t)  \bigg| \overset{\rm a.s.}=0,
 \eeqnn
 for any $T\geq 0$, which yields that 
 \beqlb\label{eqn.710}
  X_\zeta(t)= \sum_{i=1}^4\mathbf{I}_i(t),\quad t\geq 0.
 \eeqlb
 Plugging \eqref{eqn.707}-\eqref{eqn.709} into the right-hand side, we can get \eqref{MainThm.SVE} immediately.

 \subsubsection{Proof of Theorem~\ref{Thm.Comparison}} 
 
 Without loss of generality, we may assume that $\zeta_1\leq \zeta_2$. 
 Consider two sequences of compound Poission processes $\{ Y_{1,n} \}_{n\geq 1}$ and $\{ Y_{2,n} \}_{n\geq 1}$, in which $Y_{i,n}$, $i=1,2$, has a drift $-1$, jump-size distribution $\varPi_n(dy)$ with density $\pi_n$ given by $\eqref{Pi.n}$ and arrival rate $\gamma_{i,n}$ defined by \eqref{Gamma.n} with $b = b_i$. 
 For $t\geq 0$, let $ \xi^{(n)}_{i}(t):= Y_{i,n}(n^2t)/n $ that is a spectrally positive L\'evy process with Laplace exponent denoted by $ \varPhi^{(n)}_i $ given by \eqref{LapXi.n01} with $\gamma_{n}=\gamma_{i,n}$.  
 Repeating the proof of Proposition~\ref{Prop.501} shows that $\varPhi^{(n)}_{i}\to \varPhi_{i}$ pointwisely with $\varPhi_{i}$ being the Laplace exponent of $\xi_i$ and then $ \xi^{(n)}_{i}\to \xi_i$ weakly in $D(\mathbb{R}_+;\mathbb{R})$ as $n\to\infty$. 
 
 For each $n\geq 1$, let $N^{n}_{c}(ds,dz,dy,du)$ and $N^{n}_{\Lambda} (ds,dz,dy,du)$ be two orthogonal $(\mathscr{F}_{nt})$-compensated Poisson random measures on $(0,\infty)^4$ with intensities 
  \beqnn
  (1-\theta_n) \cdot \boldsymbol{e}_{c}(y)\,\cdot ds\,dz\,dy\,du
 \quad \mbox{and}\quad 
  \theta_n \cdot  \boldsymbol{e}_{c}*d\Lambda_n(y)\,\cdot ds\,dz\,dy\,du.
 \eeqnn
 Let $B_n^c(\zeta_1)$ and $B_n^c(\zeta_2-\zeta_1)$ be two independent binomial random variables with common success probability $p_n$ and number of trials $[n\zeta_1]$ and  $[n\zeta_2]-[n\zeta_1]$ respectively. 
 Moreover,  let $B^c_n(\zeta_2):=B^c_n(\zeta_1)+B_n^c(\zeta_2-\zeta_1)$,  $B_n^{\Lambda}(\zeta_1):=[n\zeta_1]-B^c_n(\zeta_1)$ and $B_n^{\Lambda}(\zeta_2):=[n\zeta_2]-B^c_n(\zeta_2)$. 
 For $i=1,2$, let 
 $$N_{\varPi_i}^n (ds,dz,dy):=N_{c}^{n}(ds,dz,dy,[0,\gamma_{i,n}])+N_{\Lambda}^{n} (ds,dz,dy,[0,\gamma_{i,n}]).$$
 By Lemma~\ref{Lemma.HR}, the process $L^{Y_{n,i}}_{[n\zeta_i]}$ equals in distribution to the unique strong solution $Z^{i,n}_{[n\zeta_i]}$ to 
 \beqnn
 Z^{i,n}_{[n\zeta_i]}(t)
 \ar=\ar  \sum_{j=1}^{B^c_n(\zeta_i)} \mathbf{1}_{\{\ell^c_{n,j}> t\}} 
 + \int_0^t \int_0^{Z^{i,n}_{[n\zeta_i]} (s-)} \int_0^\infty   \mathbf{1}_{\{y>t-s\}} \, N_{c}^{n} (ds,dz,dy,[0,\gamma_{i,n}])\cr
 \ar\ar +\sum_{j=1}^{B^\Lambda_n(\zeta_i)} \mathbf{1}_{\{\ell^\Lambda_{n,j}> t\}}+  \int_0^t \int_0^{Z^{i,n}_{[n\zeta_i]} (s-)} \int_0^\infty   \mathbf{1}_{\{y>t-s\}} \, N_{\Lambda}^{n} (ds,dz,dy,[0,\gamma_{i,n}]), \quad t\geq 0.
 \eeqnn
 Similarly as in \eqref{eqn.28}, we also have 
 $ L^{\xi^{(n)}_i}_{\zeta_i} \overset{\rm d}= X^{b_i,(n)}_{\zeta_i} = \big\{ n^{-1}\cdot Z^{i,n}_{[n\zeta_i]}(nt):t\geq 0 \big\}.$
  The proof of the next lemma is elementary and omitted.
 
 \begin{lemma}\label{Lemma.535}
  If $b_1\geq b_2$, we have for any $n\geq 1$,
  \beqnn
 	\mathbf{P}\Big( X^{b_1,(n)}_{\zeta_1}(t)\leq X^{b_2,(n)}_{\zeta_2}(t) , \, t\geq 0 \Big)=\mathbf{P}\Big( Z^{1,n}_{[n\zeta_1]}(t)\leq Z^{2,n}_{[n\zeta_2]}(t) , \, t\geq 0 \Big)=1.
  \eeqnn   
 \end{lemma}
 
 Similarly as in \eqref{Eqn.HR01}-\eqref{eqn.JLambda}, the process $ (X^{b_1,(n)}_{\zeta_1}, X^{b_2,(n)}_{\zeta_2})$ is the unique solution to the following two dimensional stochastic Volterra equation
 \beqnn
 X_{\zeta}^{b_i(n)}(t) \ar=\ar \sum_{j=1}^4 \mathbf{I}^{(n)}_{i,j}(t) ,\quad t\geq 0 ,\, i=1,2,
 \eeqnn
 with the summands are given as in \eqref{eqn.I}-\eqref{eqn.JLambda}, where the resolvent $R_{\varPi_n,i}$ and two parameter function $R_{i,n}$ are defined as in \eqref{Eqn.RPi.n} and \eqref{Eqn.R.n} respectively with $\gamma_n=\gamma_{i,n}$, and
 \beqnn
 \widetilde{N}^{(n)}_{i,c}(ds,dz,dy,[0,\gamma_{i,n}])
 \ar:=\ar N_{c}^{n} (n\cdot ds,n\cdot dz,dy,[0,\gamma_{i,n}])- n^2\gamma_{i,n}(1-\theta_n) \cdot \boldsymbol{e}_{c}(y)\,\cdot ds\,dz\,dy,\cr 
 \widetilde{N}^{(n)}_{i,\Lambda} (ds,dz,dy,[0,\gamma_{i,n}])
 \ar:=\ar N^{n}_{\Lambda} (n\cdot ds,n\cdot dz,dy,[0,\gamma_{i,n}])- n^2\gamma_{i,n}\theta_n \cdot  \boldsymbol{e}_{c}*d\Lambda_n(y)\,\cdot ds\,dz\,dy.
 \eeqnn  
 Repeating all arguments in Section~\ref{Sec.Auxiliarylemmas}-\ref{Sec.CTightness} and the proof of Theorem~\ref{Main.Thm01} proves that 
 $$ \Big( 	X^{b_1,(n)}_{\zeta_1},	X^{b_2,(n)}_{\zeta_2} \Big) \to \Big( 	X^{b_1}_{\zeta_1},	X^{b_2}_{\zeta_2} \Big),$$
 weakly in $D(\mathbb{R}_+;\mathbb{R}_+^2)$ as $n\to\infty$ with $X_{\zeta_i}^{b_i}\overset{\rm d}=L^{\xi_i}_{\zeta_i}$ for  $i=1,2$, and there are three common driving noises $(N_0,B_c,\widetilde{N}_\nu)$ defined as in Theorem~\ref{Main.Thm01} such that the quadruple $(X^{b_i}_{\zeta_i},N_0,B_c,\widetilde{N}_\nu)$ satisfies \eqref{MainThm.SVE} with characteristics $(\zeta_i\,;\, b_i,c,\nu)$. 
 By Skorokhod's representation theorem and then Proposition 1.17(b) in \cite[p.328]{JacodShiryaev2003}, we may assume that 
 $$ \big(X^{b_1,(n)}_{\zeta_1},	X^{b_2,(n)}_{\zeta_2} \big) \overset{\rm a.s.}\to \big(X^{b_1}_{\zeta_1},	X^{b_2}_{\zeta_2} \big),$$ 
 uniformly on compacts as $n\to\infty$, which
 along with Lemma~\ref{Lemma.535} yields that 
 \beqnn
 	\mathbf{P}\Big( 	X^{b_1}_{\zeta_1}(t)\leq 	X^{b_2}_{\zeta_2}(t) , \, t\geq 0 \Big)= 	\lim_{n\to\infty}\mathbf{P}\Big( 	X^{b_1,(n)}_{\zeta_1}(t)\leq 	 X^{b_2,(n)}_{\zeta_2}(t) , \, t\geq 0 \Big) =1 
 \eeqnn
 and the whole proof of Theorem~\ref{Thm.Comparison} ends.

  \section{Strong existence and uniqueness}
 \label{Sec.PathUniqueness}
 \setcounter{equation}{0}
 
 In this section, we prove the  existence and uniqueness of strong solution to \eqref{MainThm.SVE} under the assumption that $\doublebar{\nu}(0)<\infty $. 
 Our proof replies on the next lemma, which is a direct consequence of Theorem~1.5 in \cite{Kurtz2014} that generalizes the Yamada-Watanabe theorems established in \cite{YamadaWatanabe1971}. The detailed proof is omitted. 
 We say \textsl{pathwise uniqueness} holds for \eqref{MainThm.SVE}  if any two solutions defined on the same filtrated probability space endowed with the same driving noises are distinguishable.

 \begin{lemma} \label{Lemma.YamadaWatanabe}
  If the existence of solutions and pathwise uniqueness hold for \eqref{MainThm.SVE}, then the strong solution exists uniquely. 
 	
 \end{lemma}

 \textit{\textbf{Proof of Theorem~\ref{Thm.StrongUniqueness}.}} 
 The existence of solutions to \eqref{MainThm.SVE} follows from Theorem~\ref{Main.Thm01}. 
 Meanwhile, the next two lemmas respectively show that  \eqref{MainThm.SVE} is equivalent to \eqref{eqn.1004} and  pathwise uniqueness holds for \eqref{eqn.1004}. 
 In conclusion, there exists a unique strong solution to \eqref{MainThm.SVE}.
 \qed

 \begin{lemma}
 The stochastic equation \eqref{MainThm.SVE} is equivalent to \eqref{eqn.1004}. 
 \end{lemma}
 \proof By the change of variables and the fact that $W'(x)=0$ for $x<0$, 
 \beqnn
  W(t-s)-W(t-s-y) = \int_{t-s-y}^{t-s} W'(r)dr = \int_0^t W'(t-r)\cdot \mathbf{1}_{\{ s<r\leq s+y \}} \, dr.
 \eeqnn
 Plugging this into the last stochastic integral in \eqref{MainThm.SVE}, we have 
 \beqnn
 \lefteqn{\int_0^t \int_0^{X_\zeta(s)} \int_0^\infty  \big(W(t-s)-W(t-s-y) \big) \widetilde{N}_\nu(ds,dz,dy)}\ar\ar\cr 
 \ar=\ar \int_0^t \int_0^{X_\zeta(s)} \int_0^\infty \Big( \int_0^t W'(t-r)\cdot \mathbf{1}_{\{ s<r\leq s+y \}} \,dr \Big) \widetilde{N}_\nu(ds,dz,dy).
 \eeqnn
 The stochastic Fubini theorem \eqref{eqn.SFT02} along with
 \beqnn
 \int_0^t  W'(t-r) \int_0^r ds\int_0^\infty \mathbf{1}_{\{ s<r\leq s+y \}}  \,\nu(dy)\,dr <\infty 
 \eeqnn
 induces that 
 \beqnn
 \lefteqn{\int_0^t \int_0^{X_\zeta(s)} \int_0^\infty  \big(W(t-s)-W(t-s-y) \big) \widetilde{N}_\nu(ds,dz,dy)}\ar\ar\cr 
 \ar=\ar \int_0^t W'(t-r) \int_0^r \int_0^{X_\zeta(s)} \int_0^\infty   \mathbf{1}_{\{ s<r\leq s+y \}} \widetilde{N}_\nu(ds,dz,dy) \, dr\cr
 \ar=\ar \int_0^x W'(x-r) \int_0^r \int_0^{X_\zeta(s)} \int_{r-s}^\infty    \widetilde{N}_\nu(ds,dz,dy) \, dr.
 \eeqnn
 Taking this back into \eqref{MainThm.SVE} and then using the second equality in \eqref{eqn.1004}, it can be written as
 \beqlb\label{eqn.1002}
 X_\zeta(x)\ar=\ar \zeta \cdot c\cdot W'(t)  
 +  \int_0^\zeta\int_0^\infty  \big(W(t-s)-W(t-s-y) \big) N_0(dz, dy) + \int_0^t  W'(t-s)  dM(s).\quad 
 \eeqlb
 Since $W' $ is differentiable on $(0,\infty)$ with $ W'(0)=1/c$, an application of the Fubini's theorem to the last stochastic integral shows that 
 \beqlb\label{eqn.1003}
 \int_0^t  W'(t-s)dM(s) 
 \ar=\ar \frac{1}{c}\cdot M(t)+  \int_0^t \int_0^{t-s} W''(r)\,dr\, dM(s)\cr
 \ar=\ar \frac{1}{c} \cdot M(t) + \int_0^t   W^{''}(t-s)M(s)\, ds. 
 \eeqlb
 Further, the differentiability of $W'' $ on $(0,\infty)$ and $ W''(0+)<\infty$ induces that
 \beqnn
  \int_0^t   W^{''}(t-s)M(s)\, ds  \ar=\ar  W''(0+)\cdot M(t) + \int_0^t M(s) \int_0^{t-s} W'''(r)\,dr\, ds \cr
 \ar=\ar W''(0+) \cdot  M(t) + \int_0^t  W'''* M(s)\, ds .
 \eeqnn
 The representation \eqref{eqn.1004} follows immediately by taking this back into \eqref{eqn.1003} and then \eqref{eqn.1002}.
 \qed

 \begin{lemma} \label{Lemma.PathUnique}
 	The pathwise uniqueness holds for \eqref{eqn.1004}.
 \end{lemma}
 \proof Assume that $X_1$ and $X_2$ are two solutions to \eqref{MainThm.SVE} with common driving noises $(N_0,B_c,\widetilde{N}_\nu)$. 
 Let $Z(t):= X_1(t)-X_2(t)$ for $t\geq 0$.
 By \eqref{eqn.1004},  we have 
 \beqlb\label{eqn.1005}
 Z(t) =   \int_0^t \Big( W''(0+) \cdot \overline{M}(s) +  W'''*\overline{M}(s) \Big) \, ds + \frac{1}{c}\cdot \overline{M}(t),
 \eeqlb
 is a $(\mathscr{F}_t)$-semimartiangle, where
 \beqlb\label{eqn.1006}
 \overline{M}(t)\ar=\ar \int_0^t   \int_0^r \int_0^\infty \int_{r-s}^\infty  \big( \mathbf{1}_{\{ z\leq X_1(s) \}}-\mathbf{1}_{\{ z\leq X_2(s) \}} \big)  \widetilde{N}_\nu(ds,dz,dy)\, dr \cr
 \ar\ar + \int_0^t \int_0^\infty \big( \mathbf{1}_{\{ z\leq X_1(s) \}}-\mathbf{1}_{\{ z\leq X_2(s) \}} \big)   B_c(ds,dz)  . 
 \eeqlb
 Consider a strictly decreasing positive sequence $\{a_n \}_{n\geq 1}$ such that $\int_{a_n}^{a_{n-1}}z^{-1}dz=n$ and $a_n\to0$ as $n\to\infty$, e.g.,
 $a_n=\exp\big\{-n(n+1)/2\big\}$. 
 Let $x\mapsto g_n(x)$ be a non-negative continuous function on $\mathbb{R}$ which has support in $(a_n,a_{n-1})$ and satisfies 
 \beqnn
 \int_{a_n}^{a_{n-1}}g_n(x)\,dx=1
 \quad \mbox{and}\quad
  x\cdot g_n(x)\leq\frac{2}{n}.
 \eeqnn
  Moreover, we define the non-negative and twice continuously differentiable function
 \beqnn
 f_n(z)= \int_0^{|z|}dy\int_0^yg_n(x)dx , \qquad z\in \mathbb{R}.
 \eeqnn
 It is easy to see that $f_n(z)\rightarrow |z|$ non-decreasingly as $n\to\infty$ and 
 \beqlb\label{eqn.1007}
 \sup_{z\in\mathbb{R}}|f_n'(z)|\leq 1
 \quad \mbox{and}\quad
 \sup_{z\in\mathbb{R}}|zf''_n(z)|\leq \frac{2}{n},\quad n\geq 1. 
 \eeqlb
 By \eqref{eqn.1005}-\eqref{eqn.1006} and using It\^o's formula to $f_n\big(Z(t)\big)$; see Theorem 5.1 in \cite[p.66]{IkedaWatanabe1989}, 
 \beqnn
  f_n\big(Z(t)\big)\ar=\ar  \int_0^t f'_n\big(Z(s)\big)  \Big( W''(0+) \cdot \overline{M}(s) +  W'''*\overline{M}(s) \Big) \, ds   + \frac{1}{c}  \int_0^t \big|Z(s)\big|\cdot f''_n\big(Z(s)\big)  \, ds\cr
 \ar \ar + \frac{1}{c}\int_0^t   f'_n\big(Z(s)\big)  \int_0^s \int_0^\infty \int_{s-r}^\infty  \big( \mathbf{1}_{\{ z\leq X_1(r) \}}-\mathbf{1}_{\{ z\leq X_2(r) \}} \big)  \widetilde{N}_\nu(dr,dz,dy)\, ds\cr 
 \ar\ar + \frac{1}{c}\int_0^t \int_0^\infty f'_n\big(Z(s)\big)\big( \mathbf{1}_{\{ z\leq X_1(s) \}}-\mathbf{1}_{\{ z\leq X_2(s) \}} \big)  B_c(ds,dz).
 \eeqnn
 By the non-negativity of $f_n$, \eqref{eqn.1007} and the definition of compensated Poisson random measure,  
 \beqnn
  f_n\big(Z(t)\big)\ar\leq\ar  \int_0^t   \Big| W''(0+) \cdot  \overline{M}(s) +  W'''*\overline{M}(s) \Big| \, ds   + \frac{t}{c}\cdot \frac{2}{n}  \cr
 \ar \ar + \frac{1}{c}\int_0^t  \bigg| \int_0^s \int_0^\infty \int_{s-r}^\infty  \big( \mathbf{1}_{\{ z\leq X_1(r) \}}-\mathbf{1}_{\{ z\leq X_2(r) \}} \big)  \widetilde{N}_\nu(dr,dz,dy)\bigg| \,ds\cr 
 \ar\ar + \frac{1}{c}\int_0^t \int_0^\infty f'_n\big(Z(s)\big)\big( \mathbf{1}_{\{ z\leq X_1(s) \}}-\mathbf{1}_{\{ z\leq X_2(s) \}} \big)  B_c(ds,dz)\cr
 \ar\leq\ar \big|W''(0+) \big| \cdot \int_0^t  \big|\overline{M}(s)\big|\, ds +  \int_0^t  \big|W'''\big|*\big|\overline{M}\big|(s) \, ds  + \frac{t}{c}\cdot \frac{2}{n}   + \frac{1}{c}\int_0^t    \big|Z\big|* \bar\nu(s)\,ds\cr
 \ar \ar   + \frac{1}{c}\int_0^t   \int_0^s \int_0^\infty \int_{s-r}^\infty  \big| \mathbf{1}_{\{ z\leq X_1(r) \}}-\mathbf{1}_{\{ z\leq X_2(r) \}} \big|  N_\nu(dr,dz,dy)  \,ds\cr 
 \ar\ar + \frac{1}{c}\int_0^t \int_0^\infty f'_n\big(Z(s)\big)\big( \mathbf{1}_{\{ z\leq X_1(s) \}}-\mathbf{1}_{\{ z\leq X_2(s) \}} \big)  B_c(ds,dz) . 
 \eeqnn
 Taking expectations on both sides and then using Fubini's theorem and Young's inequality for convolution, 
 \beqnn
 \mathbf{E}\Big[ f_n\big(Z(t)\big) \Big] 
 \ar\leq\ar \bigg( \big|W''(0+) \big| + \int_0^t  \big|W'''\big|(r)dr\bigg) \cdot  \int_0^t \mathbf{E}\Big[ \big|\overline{M}(s)\big|\Big] \, ds   + \frac{2 \doublebar\nu(0)}{c} \cdot \int_0^t  \mathbf{E}\Big[\big|Z(s)\big|\Big]\,ds + \frac{t}{c}\cdot \frac{2}{n} . 
 \eeqnn
 By using the monotone convergence theorem along with the fact that $f_n(z)\rightarrow |z|$ non-decreasingly as $n\to\infty$, we have uniformly in $x\in[0,T]$, 
 \beqlb\label{eqn.1008}
 \mathbf{E}\Big[ \big|Z(t)\big| \Big] 
 \ar\leq\ar C\cdot  \int_0^t \mathbf{E}\Big[ \big|\overline{M}(s)\big| +\big|Z(s)\big|\Big] \, ds . 
 \eeqlb
 On the other hand, by \eqref{eqn.1005} we also have 
 \beqnn
 \overline{M}(t)= c\cdot Z(t) -c\cdot \int_0^t \Big(  W''(0+)   \cdot \overline{H}(s) +  W'''*\overline{M}(s) \Big) \, ds.
 \eeqnn
 By using the triangle inequality and then Young's inequality for convolution,
 \beqnn
 \big|\overline{M}(t)\big|\leq c\cdot \big|Z(t)\big| + c\cdot\Big( \big|W''(0+) \big| + \int_0^t \big| W'''(r)\big|\,dr\Big) \int_0^t    \big|\overline{M}(s)\big|\,ds .
 \eeqnn
 Taking expectations on both sides induces that for any $T\geq 0$,
 \beqnn
 \mathbf{E}\Big[\big|\overline{M}(t)\big|\Big] 
 \leq c\cdot \mathbf{E}\Big[\big|Z(t)\big|\Big] +C \int_0^t  \mathbf{E}\Big[  \big|\overline{M}(s)\big| \Big]\,ds 
 \leq C \int_0^t \mathbf{E}\Big[ \big|\overline{M}(s)\big| +\big|Z(s)\big|\Big] \, ds ,
 \eeqnn
  uniformly in $t\in [0,T]$. 
  Here the second inequality follows from  \eqref{eqn.1008}. 
  Finally, combining this together with \eqref{eqn.1008} yields that
 \beqnn
 \mathbf{E}\Big[  \big|Z(t)\big| + \big|\overline{M}(t)\big| \Big] \leq C \int_0^t \mathbf{E}\Big[ \big|\overline{M}(s)\big| +\big|Z(s)\big|\Big] \, ds,\quad t\in[0,T]. 
 \eeqnn
  By Gr\"onwall's inequality, we have  $\mathbf{E}\big[ |Z(t) | +  |\overline{M}(t) | \big]=0$
 and hence $X_1(t)\overset{\rm a.s.}=X_2(t)$ for any $t\in[0,T]$.
 Finally, the continuity of  $X_1$ and $X_2$ induces that they are indistinguishable. 
 \qed

 \section{Moment estimates and equivalent representations}
 \label{Sec.MomentEstimate}
 \setcounter{equation}{0}
 
 In this section we first prove the moment estimates given in Theorem~\ref{Thm.Moment} and then establish a general equivalent representation of \eqref{MainThm.SVE} that yields Theorem~\ref{MainThm.EquavilentRep} as a corollary.
 
 
 \subsection{Proof of Theorem~\ref{Thm.Moment}} 
 
  The proof is carried out based on the equivalent equation \eqref{MainThm.SVENew} that, for convention, is written as
 \beqlb\label{eqn.5000}
 X_\zeta(t)= \sum_{i=1}^4\widetilde{\mathbf{I}}_i(t),\quad t\geq 0.
 \eeqlb
 Taking expectations on both sides and then by using the non-negativity of $X_\zeta$, we first have
 \beqlb \label{eqn.5001}
 \mathbf{E}\big[  X_\zeta(t)   \big] = \widetilde{\mathbf{I}}_i(t) = \zeta \cdot \big(1-bW(t)  \big) \leq \zeta , 
 \eeqlb
 uniformly in $\zeta,t\geq 0$ and hence \eqref{eqn.Moment} hold for $p=1$. 
 By induction, we now proceed under the assumption that \eqref{eqn.Moment} holds for some $p\geq 1$ and prove that it also holds for $2p$. 
 By the power inequality and \eqref{eqn.5000},  
 \beqlb\label{eqn.5003}
 \mathbf{E}\Big[\big| X_\zeta(t) \big|^{2p}  \Big] 
 \ar\leq\ar C \cdot \sum_{i=1}^4 \mathbf{E}\Big[\big| \widetilde{\mathbf{I}}_i(t) \big|^{2p}  \Big] ,
 \eeqlb 
 for some constant $C>0$ depending only on $p$. 
 The inequality in \eqref{eqn.5001} tells that $$\sup_{t\geq 0} \mathbf{E}\Big[\big| \widetilde{\mathbf{I}}_1(t) \big|^{2p}  \Big]= \sup_{t\geq 0} \big| \widetilde{\mathbf{I}}_1(t) \big|^{2p} \leq \zeta^{2p}.$$ 
 Applying \eqref{BDG} to $\mathbf{E}\big[ | \mathbf{I}_2(t) |^{2p}  \big]$ and then using \eqref{eqn.103} and \eqref{eqn.10301}, there exists a constant $C>0$ depending only on $p$ such that 
 \beqnn
 \mathbf{E}\Big[\big| \mathbf{I}_2(t) \big|^{2p}  \Big]
 \ar\leq\ar C\cdot \bigg|  \zeta \cdot \int_0^\infty \big| W(t)-W(t-y) \big|^2 \, \bar\nu(y)\,dy \bigg|^{p}\cr
 \ar\ar +  C\cdot   \zeta \cdot \int_0^\infty \big| W(t)-W(t-y) \big|^{2p} \, \bar\nu(y)\,dy \cr
 \ar\leq\ar  C\cdot \big|\zeta\cdot W(t)\big|^{p}\cdot \bigg|  \int_0^\infty \big( W(t)-W(t-y) \big)\, \bar\nu(y)\,dy \bigg|^{p}\cr
 \ar\ar +  C\cdot \zeta \cdot \big|W(t)\big|^{2p-1} \int_0^\infty \big( W(t)-W(t-y) \big) \, \bar\nu(y)\,dy \cr
 \ar\ar\cr
 \ar\leq\ar  C\cdot |\zeta|^p \cdot \big| W(t)\big|^{p} +  C\cdot \zeta \cdot \big|W(t)\big|^{2p-1}. 
 \eeqnn
 Applying \eqref{BDG0} to $\mathbf{E}\big[ | \mathbf{I}_3(t) |^{2p}  \big]$ and then using \eqref{eqn.Moment} give that uniformly in $\zeta,t \geq 0$, 
 \beqnn
 \mathbf{E}\Big[\big| \mathbf{I}_3(t) \big|^{2p}  \Big]
 \ar\leq\ar C \cdot \sup_{s\in [0,t]} \mathbf{E}\Big[ \big| X_\zeta(s) \big|^{p} \Big] \cdot  \bigg| c\int_0^t   \big|W'(s)\big|^2 ds  \bigg|^{p}\cr
 \ar\leq\ar C \cdot  \sup_{s\in [0,t]} \mathbf{E}\Big[ \big| X_\zeta(s) \big|^{p} \Big] \cdot  \bigg| \int_0^t  W'(s)  ds  \bigg|^{p}\cr
 \ar\leq\ar C \cdot  \sup_{s\in [0,t]} \mathbf{E}\Big[ \big| X_\zeta(s) \big|^{p} \Big] \cdot  \big| W(t) \big|^{p}\cr
 \ar\ar\cr
 \ar\leq\ar C \cdot \big(\zeta\vee \zeta^p\big) \cdot \big( 1+W(t)\big)^{3p-2} . 
 \eeqnn
 Finally, an application of \eqref{BDG} to $\mathbf{E}\big[ | \mathbf{I}_4(t) |^{2p}  \big]$ and then using \eqref{eqn.Moment} gives that
 \beqlb\label{eqn.5002}
 \mathbf{E}\Big[\big| \mathbf{I}_4(t) \big|^{2p}  \Big] 
 \ar\leq\ar C\cdot \sup_{s\in [0,t]}\mathbf{E}\Big[ \big| X_\zeta(s) \big|^{p} \Big]  \cdot \bigg|  \int_0^t ds  \int_0^\infty  \big( W(s)-W(s-y) \big)^2 \, \nu(dy) \bigg|^p \cr
 \ar\ar + C\cdot \sup_{s\in [0,t]}\mathbf{E}\big[  X_\zeta(s)   \big]  \cdot   \int_0^t ds  \int_0^\infty  \big( W(s)-W(s-y) \big)^{2p} \, \nu(dy)  \cr
 \ar\leq\ar C\cdot \sup_{s\in [0,t]}\mathbf{E}\Big[ \big| X_\zeta(s) \big|^{p} \Big]  \cdot \bigg|  \int_0^t ds  \int_0^\infty  \big( W(s)-W(s-y) \big)^2 \, \nu(dy) \bigg|^p \cr
 \ar\ar + C\cdot  \sup_{s\in [0,t]}\mathbf{E}\big[  X_\zeta(s)   \big] \cdot  \big|W(t) \big|^{2p-2}  \cdot   \int_0^t ds  \int_0^\infty  \big( W(s)-W(s-y) \big)^{2} \, \nu(dy).
 \eeqlb
 By using \eqref{eqn.201} and the monotonicity of $W$, 
 \beqnn
 \int_0^t ds  \int_0^\infty  \big( W(s)-W(s-y) \big)^{2} \, \nu(dy)
 \ar\leq\ar \int_0^t ds  \int_0^\infty \Big(W(s)\wedge \frac{y}{c} \Big) \big( W(s)-W(s-y) \big)  \, \nu(dy)\cr
 \ar\leq\ar   \int_0^\infty \Big(W(t)\wedge \frac{y}{c} \Big) \int_0^t \big( W(s)-W(s-y)\big)\,ds   \, \nu(dy)\cr
 \ar=\ar  \int_0^\infty \Big(W(t)\wedge \frac{y}{c} \Big) \int_{t-y}^t  W(s) \,ds   \, \nu(dy)\cr 
 \ar\leq\ar \frac{W(t)}{c}\int_0^\infty \Big(\big(cW(t)y\big)\wedge y^2 \Big)  \, \nu(dy) ,
 \eeqnn
 which can be uniformly bounded by $c^{-2}\cdot\int_0^\infty (y \wedge y^2)   \, \nu(dy)$ if $cW(t)\leq 1$, or by
 \beqnn
 \frac{W(t)}{c}\int_0^1  y^2    \, \nu(dy) + \frac{W(t)}{c}\int_1^\infty  cW(t)y   \, \nu(dy)  
 \ar\leq\ar |W(t)|^2\int_0^\infty  (y\wedge y^2)  \, \nu(dy) . 
 \eeqnn
 if $cW(t)>1$.
 Plugging this into \eqref{eqn.5002} and then using \eqref{eqn.Moment} induces that 
 \beqnn
 \mathbf{E}\Big[\big| \mathbf{I}_4(t) \big|^{2p}  \Big] 
 \ar\leq\ar C\cdot \sup_{s\in [0,t]}\mathbf{E}\Big[ \big| X_\zeta(s) \big|^{p} \Big]  \cdot \Big(1+W(t)  \Big)^{2p}  \cr
 \ar\ar + C\cdot  \sup_{s\in [0,t]}\mathbf{E}\big[  X_\zeta(s)   \big] \cdot |W(t)|^{2p-2}  \cdot  \Big(1+W(t)  \Big)^{2}\cr
 \ar\leq\ar C\cdot \big(\zeta\vee \zeta^p\big)  \cdot \Big(1+W(t)  \Big)^{4p-2}  + C\cdot  \zeta \cdot   \Big(1+W(t)  \Big)^{2p} 
 \eeqnn
 Taking all preceding estimates back into \eqref{eqn.5003}, we see that \eqref{eqn.Moment} also holds for $2p$.

  \subsection{Equivalent representations}
 
 We now provide an equivalent representation of \eqref{MainThm.SVE} by using the stochastic Fubini's theorems given in Proposition~\ref{StoFubiniThm} and the next proposition that can be found in e.g.  \cite[p.36]{GripenbergLondenStaffans1990} or Lemma~3 in \cite{BacryDelattreHoffmannMuzy2013}. 
 
 \begin{proposition}\label{Prop.601}
 	Consider two functions $g,h\in L^1_{\rm loc}(\mathbb{R}_+;\mathbb{R}) $. 
 	Let $R_g \in L^1_{\rm loc}(\mathbb{R}_+;\mathbb{R})$ be the resolvent of the second kind of $g$ defined as in \eqref{Resolvent.01}. 
 	The linear Volterra equation 
 	\beqnn
 	f(t)= h(t) +g*f(t) ,\quad t\geq 0,
 	\eeqnn
 	has a unique locally integrable solution.
 	Moreover, the solution has the representation
 	\beqnn
 	f(t)= h(t) + R_g*h (t), \quad  t\geq 0.
 	\eeqnn 
 \end{proposition}

 For any $\beta\in\mathbb{R}$, recall the scaled function $W_\beta$ defined above Lemma~\ref{Lemma.Identity02}. 
 By \eqref{eqn.201},
 \beqlb \label{eqn.20101}
 W'_\beta(0)=\frac{1}{c} 
 \quad \mbox{and}\quad
 \sup_{x\in \mathbb{R}} W'_\beta(x)\leq \frac{1}{c}.
 \eeqlb

 \begin{lemma}
 	For any $\beta\in\mathbb{R}$, the stochastic Volterra equation \eqref{MainThm.SVE} is equivalent to 
 	\beqlb \label{Main.SVE021} 
 	X_\zeta(t)\ar=\ar \zeta \cdot c\cdot W'_\beta(t)  +  \int_0^\zeta\int_0^\infty  \big( W_\beta(t)-W_\beta(t-y) \big) \, N_0(dz, dy)   \cr
 	\ar\ar  + \int_0^t (\beta-b)\cdot W'_\beta(t-s) X_\zeta(s)ds + \int_0^t \int_0^{X_\zeta(s)}   W'_\beta(t-s) \, B_c(ds,dz) \cr
 	\ar\ar + \int_0^t \int_0^{X_\zeta(s)} \int_0^\infty \big( W_\beta(t-s)-W_\beta(t-s-y) \big) \, \widetilde{N}_\nu(ds,dz,dy) ,\quad t\geq 0 . 
 	\eeqlb

 \end{lemma}
 \proof 
 We just prove that any solution of \eqref{Main.SVE021} also solves \eqref{MainThm.SVE}.  The converse can be proved similarly. 
 For convention, we write \eqref{Main.SVE021} as 
 \beqlb\label{eqn.602}
 X_\zeta(t)\ar=\ar  \sum_{i=1}^4 \mathbf{I}_{\beta,i}(t) + (\beta-b)  W'_\beta* X_\zeta(t) ,
 \quad t\geq 0, 
 \eeqlb
 with  $\mathbf{I}_{\beta,i}(t)$, $i=1,\cdots,4$ representing the other four terms on the right-hand side of \eqref{Main.SVE02} sequentially.
 Multiplying both sides of the identity in Lemma~\ref{Lemma.Identity02} by $\beta-b $ gives that  
 \beqnn
 (\beta-b) W'(t) = (\beta-b) W'_\beta(t)+ \big((\beta-b) W'\big)* \big((\beta-b) W'_\beta\big)(t), \quad t\geq 0.
 \eeqnn
 Armed with this equality, we apply Proposition~\ref{Prop.601} to \eqref{eqn.602} and  get 
 \beqnn
 X_\zeta(t)\ar=\ar \sum_{i=1}^4 \Big( \mathbf{I}_{\beta,i}(t) +  (\beta-b) W'* \mathbf{I}_{\beta,i}(t) \Big),\quad t\geq 0. 
 \eeqnn
 Comparing this to \eqref{MainThm.SVE} or \eqref{eqn.710}, it suffices to prove that 
 \beqlb\label{eqn.601}
 \mathbf{I}_i(t)=\mathbf{I}_{\beta,i}(t) +  (\beta-b) W'* \mathbf{I}_{\beta,i} (t),\quad i=1,\cdots,4.
 \eeqlb
 
 \begin{enumerate}
 	\item[$\bullet$] For $i=1$, by Lemma~\ref{Lemma.Identity02} we have 
 	\beqnn
 	\mathbf{I}_{\beta,i} (t)+  (\beta-b) W'* \mathbf{I}_{\beta,i}(t) = \zeta\cdot c\cdot \big( W'_\beta(t) + (\beta-b) W'*  W'_\beta(t)\big)= \zeta\cdot c\cdot W'(t) = \mathbf{I}_1(t).
 	\eeqnn
 	
 	\item[$\bullet$] For $i=2$, the non-negativity of $W'$ and $W'_\beta$ allows us to use Fubini's theorem pathwisely to obtain that
 	\beqnn
 	(\beta-b) W'* \mathbf{I}_{\beta,2}(t) 
 	\ar=\ar \int_0^t (\beta-b) W^{\prime} (t-r) \int_0^\zeta\int_0^\infty \Big(\int_{r-y}^r W'_\beta(s)\,ds\Big)   \, N_0(dz, dy)\, dr\cr
 	\ar=\ar  \int_0^\zeta\int_0^\infty  \Big( \int_0^t(\beta-b) W^{\prime} (t-r) \int_{r-y}^r W'_\beta(s)\,ds  \, dr \Big) \, N_0(dz, dy).
 	\eeqnn
 	Additionally, using the change of variables and then Fubini's theorem again to the integrand in the last stochastic integral gives that  
 	\beqlb\label{eqn.613}
 	\int_0^t (\beta-b) W'(t-r) \int_{r-y}^r W'_\beta(s)\,ds\, dr
 	\ar=\ar \int_0^t (\beta-b) W'(t-r) \int_0^y W'_\beta(r-s)\,ds\, dr\cr
 	\ar=\ar \int_0^y \int_0^t (\beta-b) W'(t-r) W'_\beta(r-s)\,dr \,ds \cr
 	\ar=\ar \int_0^y   (\beta-b) W'* W'_\beta(t-s) \,ds \cr
 	\ar=\ar \int_{t-y}^t  (\beta-b) W'* W'_\beta(s) \,ds. 
 	\eeqlb
 	Combining these together and then using the identity in Lemma~\ref{Lemma.Identity02}, we have 
 	\beqnn
 	\mathbf{I}_{\beta,2}(t) + (\beta-b) W'* \mathbf{I}_{\beta,2}(t) 
 	\ar=\ar  \int_0^\zeta\int_0^\infty \int_0^t \Big( \int_{t-y}^t   W'_\beta(s) \,ds +  \int_{t-y}^t  (\beta-b) W'* W'_\beta(s) \,ds \Big) \, N_0(dz, dy) \cr
 	\ar=\ar  \int_0^\zeta\int_0^\infty  \Big( \int_{t-y}^t   W' (s) \,ds \Big) \, N_0(dz, dy),
 	\eeqnn 
 	and hence \eqref{eqn.601} holds for $i=2$.
 	
 	\item[$\bullet$] For $i=3$, by \eqref{eqn.201} and \eqref{eqn.20101} we first have 
 	\beqnn
 	\int_0^t |\beta-b| W'(t-r) \cdot \int_0^r\big| W'_\beta(s)\big|^2\, ds\, dr \leq \frac{|\beta-b|}{c^3}\cdot t^2,
 	\eeqnn
 	which along with Theorem~\ref{Thm.Moment} indicates that conditions in Proposition~\ref{StoFubiniThm} are satisfied. Using the identity \eqref{eqn.SFT01} gives that  
 	\beqnn
 	(\beta-b) W'* \mathbf{I}_{\beta,3}(t) 
 	\ar=\ar \int_0^t (\beta-b) W'(t-r) \int_0^r \int_0^{X_\zeta(s)}   W'_\beta(r-s)  \, B_c(ds,dz)\, dr\cr
 	\ar=\ar  \int_0^t \int_0^{X_\zeta(s)}    (\beta-b) W'* W'_\beta (t-s)\, B_c(ds,dz) .
 	\eeqnn
 	Similarly as in the preceding argument, we use the identity in Lemma~\ref{Lemma.Identity02} again to get
 	\beqnn
 	\mathbf{I}_{\beta,3}(t) + (\beta-b) W'* \mathbf{I}_{\beta,3}(t) 
 	\ar=\ar  \int_0^t \int_0^{X_\zeta(s)} \Big( W'_\beta (t-s) + (\beta-b) W'* W'_\beta (t-s)\Big)\, B_c(ds,dz) \cr
 	\ar=\ar \int_0^t \int_0^{X_\zeta(s)}  W' (t-s)  \, B_c(ds,dz),
 	\eeqnn
 	and hence \eqref{eqn.601} holds for $i=3$.
 	
 	\item[$\bullet$] For $i=4$, similarly as in the preceding case, by using \eqref{eqn.20101} and \eqref{LevyTriplet} we also have 
 	\beqnn
 	\int_0^t(\beta-b) W'(t-r) \int_0^\infty \Big| \int_{r-y}^{r} W'_\beta(u)\,du \Big|^2\nu(dy) \, dr
 	\leq \frac{t}{c^3}\int_0^\infty (t\wedge y)^2\nu(dy) <\infty,
 	\eeqnn
 	and hence conditions in Proposition~\ref{StoFubiniThm} are satisfied. 
 	By \eqref{eqn.SFT02}, we have $(\beta-b) W'* \mathbf{I}_{\beta,4}(t)$ equals to 
 	\beqnn 
 	\lefteqn{ \int_0^t  (\beta-b) W'(t-r) \int_0^r \int_0^{X_\zeta(s)} \int_0^\infty  \Big(\int_{r-s-y}^{r-s} W'_\beta(u)\,du\Big)  \, \widetilde{N}_\nu(ds,dz,dy)\, dr}\ar\ar\cr
 	\ar=\ar \int_0^t \int_0^{X_\zeta(s)} \int_0^\infty  \Big(  \int_0^{t-s}(\beta-b) W'(t-s-r)\int_{r-y}^{r} W'_\beta(u)\,du\, dr\Big) \, \widetilde{N}_\nu(ds,dz,dy)\cr
 	\ar=\ar \int_0^t \int_0^{X_\zeta(s)} \int_0^\infty  \Big( \int_{t-s-y}^{t-s} (\beta-b) W'* W'_\beta(r) dr\Big) \, \widetilde{N}_\nu(ds,dz,dy)
 	\eeqnn
 	Here the last equality is obtained similarly as in \eqref{eqn.613}. 
 	Again, combining this with $ \mathbf{I}_{\beta,4}(t) $ and then using the identity in Lemma~\ref{Lemma.Identity02} gives that
 	\beqnn
 	\mathbf{I}_{\beta,4}(t) + (\beta-b) W'* \mathbf{I}_{\beta,4}(t) 
 	\ar=\ar  \int_0^t \int_0^{X_\zeta(s)} \int_0^\infty  \Big( \int_{t-s-y}^{t-s}   W' (r) dr\Big)  \, \widetilde{N}_\nu(ds,dz,dy) ,
 	\eeqnn
 	and hence \eqref{eqn.601} holds for $i=4$.
 	\qed 
 \end{enumerate}
 
%
%

 \section{H\"older regularity and maximal inequality}
 \label{Sec.EquivalentSVE}
 \setcounter{equation}{0}
 

 This section devotes to the proof of Theorem~\ref{Thm.Regularity} by using the next equivalent representation of \eqref{MainThm.SVE} 
 	\beqlb \label{Main.SVE0212} 
 	X_\zeta(t) 
 	\ar=\ar \zeta \Big(1-\beta \cdot W_\beta(t)\Big)  +  \int_0^\zeta\int_0^\infty \big( W_\beta(t)-W_\beta(t-y) \big) \, \widetilde{N}_0(dz, dy)  \cr
 	\ar\ar  +\int_0^t 	(\beta-b) \cdot W'_\beta(t-s) X_\zeta(s)ds  + \int_0^t \int_0^{X_\zeta(s)}   W'_\beta(t-s) \, B_c(ds,dz)\cr
 	\ar\ar + \int_0^t \int_0^{X_\zeta(s)} \int_0^\infty  \big( W_\beta(t-s)-W_\beta(t-s-y) \big)\, \widetilde{N}_\nu(ds,dz,dy),\quad t\geq 0 . 
 	\eeqlb 
 It can be identified similarly as in Remark~\ref{Remark.EquivalentSVE}. 
 For convention, we denote by $\widetilde{\mathbf{I}}_{\beta,i}(t) $ with $i=1,\cdots,5$ the five terms on the right side of 
 \eqref{Main.SVE0212}.  
 By the Kolmogorov continuity theorem; see e.g. Theorem 2.1 in \cite[p.26]{RevuzYor2005}, 
 it suffices to identify that they satisfy certain constraints on the moments of their increments.
 
 \begin{proposition}\label{Proposition.603}
 	We have uniformly in $\zeta \geq 0$,  $\beta \in \mathbb{R}$ and $0\leq t_1<t_2$, 
 	\beqlb \label{eqn.622}
 	\Big| \widetilde{\mathbf{I}}_{\beta,1}(t_2) -\widetilde{\mathbf{I}}_{\beta,1}(t_1)\Big|
 	\leq \zeta\cdot \frac{|\beta|}{c} \cdot \big|t_2-t_1 \big| . 
 	\eeqlb
 \end{proposition}
 \proof  
 By \eqref{eqn.20101}, the scale function $W_\beta$ enjoys the global Lipschitz continuity,  i.e.,  
 \beqnn
 \big| W_\beta(t_2) -  W_\beta(t_1)\big| =\int_{t_1}^{t_2 } W'(r)dr 
 \leq \frac{|t_2-t_1|}{c},
 \eeqnn
 uniformly in  $0\leq t_1<t_2$ and hence the desired upper bound \eqref{eqn.622} follows.
 \qed 
 
 \begin{proposition}\label{Proposition.604}
 	For each $p \geq 1$ and $\beta >0$, there exists a constant $C>0$ such that for any $\zeta, t\geq 0$ and $0\leq t_1<t_2\leq t$,
 	\beqnn
 	\mathbf{E}\Big[ \big|\widetilde{\mathbf{I}}_{\beta,2}(t_2) - \widetilde{\mathbf{I}}_{\beta,2}(t_1)  \big|^{2p}  \Big] 
 	\leq C\cdot (\zeta \vee \zeta^p)  \cdot \big(|t_2-t_1|\wedge 1 \big)^p.
 	\eeqnn
 \end{proposition}
 \proof 
 By the definition of $\widetilde{\mathbf{I}}_{\beta,2}$, we first have 
 \beqnn
 \widetilde{\mathbf{I}}_{\beta,2}(t_2) - \widetilde{\mathbf{I}}_{\beta,2}(t_1)
 \ar=\ar \int_0^\zeta \int_0^\infty \Big( \int_{t_2-y}^{t_2} W'_\beta(r)dr  - \int_{t_1-y}^{t_1} W'_\beta(r)dr \Big) \widetilde{N}_0(dz,dy).
 \eeqnn 
 An application of \eqref{BDG} gives that for some constant $C>0$ depending only on $p$, 
 \beqlb \label{eqn.604}
 \mathbf{E}\Big[ \big|\widetilde{\mathbf{I}}_{\beta,2}(t_2) - \widetilde{\mathbf{I}}_{\beta,2}(t_1)  \big|^{2p}  \Big]  
 \ar\leq\ar C\cdot \Big(  \zeta^p \cdot \big|A_{2}(t_1,t_2)\big|^p + \zeta \cdot  A_{2p}(t_1,t_2) \Big),
 \eeqlb
 where 
 \beqnn
 A_{2k}(t_1,t_2):= \int_0^\infty \Big|\int_{t_2-y}^{t_2} W'_\beta(r)dr  - \int_{t_1-y}^{t_1} W'_\beta(r)dr \Big|^{2k} \bar\nu(y)\,dy,\quad k\geq 1. 
 \eeqnn
 The upper bound in \eqref{eqn.20101} tells that uniformly in $y,t\geq 0$ and $0\leq t_1<t_2\leq t$, 
 \beqnn
 \Big|\int_{t_2-y}^{t_2} W'_\beta(r)dr  - \int_{t_1-y}^{t_1} W'_\beta(r)dr \Big| 
 \leq \frac{ |t_2-t_1|\wedge y}{c}  \wedge W_\beta(t) .
 \eeqnn
 Then there exists a constant $C>0$ independent of $\zeta $ and $ t,t_1,t_2$ such that
 \beqnn
 A_{2k}(t_1,t_2) 
 \ar\leq\ar C \int_0^\infty  \big( |t_2-t_1|\wedge y \wedge W_\beta(t) \big)^{2k} \bar\nu(y)\,dy  \cr
 \ar\leq\ar C \cdot \int_0^\infty  \big( |t_2-t_1|\wedge y \wedge W_\beta(t) \big)  \bar\nu(y)\,dy \cdot \big||t_2-t_1|\wedge W_\beta(t) \big|^{2k-1} \cr
 \ar\ar\cr
 \ar\leq\ar C\cdot \big(1+|t_2-t_1|\wedge W_\beta(t) \big)\cdot \big||t_2-t_1|\wedge W_\beta(t) \big|^{2k-1} . 
 \eeqnn 
 Plugging this with $k=1 $ and $p$ back into \eqref{eqn.604} gives that 
 \beqnn
 \mathbf{E}\Big[ \big|\widetilde{\mathbf{I}}_{\beta,2}(t_2) - \widetilde{\mathbf{I}}_{\beta,2}(t_1)  \big|^{2p}  \Big]  
 \ar\leq\ar C\cdot \zeta^p \cdot \big(1+|t_2-t_1|\wedge W_\beta(t) \big)^p\cdot \big||t_2-t_1|\wedge W_\beta(t) \big|^p \cr
 \ar\ar\cr
 \ar\ar + C\cdot  \zeta \cdot  \big(1+|t_2-t_1|\wedge W_\beta(t) \big)\cdot \big||t_2-t_1|\wedge W_\beta(t) \big|^{2p-1} ,
 \eeqnn
 which along with the fact that $\sup_{t\in\mathbb{R}}W_\beta(t)\leq (\beta+b)^{-1}$  yields the desired upper bound follows. 
 \qed 
 
 \begin{proposition}\label{Proposition.607}
 	For each $p \geq 1$ and $\beta >0$,  there exists a constant $C>0$ such that for any $t\geq 0$ and $0\leq t_1<t_2\leq t$,
 	\beqnn
 	\mathbf{E}\Big[  \big| \widetilde{\mathbf{I}}_{\beta,3}(t_2) - \widetilde{\mathbf{I}}_{\beta,3}(t_1) \big|^{2p} \Big] 
 	\leq C\cdot   \big(\zeta\vee \zeta^{2p}\big)\cdot  \big(1+W(t)\big)^{4p-2}\cdot t^p   \cdot |t_2-t_1|^p. 
 	\eeqnn  
 \end{proposition}
 \proof By the definition of $\widetilde{\mathbf{I}}_{\beta,3}$ and H\"older's inequality, 
 \beqnn
 \Big| \widetilde{\mathbf{I}}_{\beta,3}(t_2) - \widetilde{\mathbf{I}}_{\beta,3}(t_1) \Big|^{2p}
 \ar=\ar |\beta|^{2p}\cdot \bigg|\int_0^{t_2}  \big( W'_\beta(t_2-s)-W'_\beta(t_1-s)\big) X_\zeta(s)ds  \bigg|^{2p}\cr
 \ar\leq\ar |\beta|^{2p}\cdot \bigg|\int_\mathbb{R} \big( W'_\beta(t_2-s)-W'_\beta(t_1-s)\big)^2\,ds\bigg|^{p} \cdot \bigg|\int_0^t  \big| X_\zeta(s)\big|^2 \, ds  \bigg|^{p}.
 \eeqnn
 Note that the claim in Lemma~\ref{Lemma.201} also holds for $W'_\beta$, then uniformly in $t_2>t_1\geq 0$,
 \beqlb\label{eqn.50004}
 \int_\mathbb{R} \big( W'_\beta(t_2-s)-W'_\beta(t_1-s)\big)^2\,ds 
 \leq C\cdot |t_2-t_1| .
 \eeqlb
 Consequently, there exists a constant $C>0$ depending only in $p$ and $\beta$ such that 
 \beqnn
 \mathbf{E}\Big[  \big| \widetilde{\mathbf{I}}_{\beta,3}(t_2) - \widetilde{\mathbf{I}}_{\beta,3}(t_1) \big|^{2p} \Big] 
 \ar\leq\ar C\cdot   |t_2-t_1|^p \cdot \mathbf{E}\bigg[ \Big|\int_0^{t_2}  \big| X_\zeta(s)\big|^2 \, ds  \Big|^{p} \bigg].
 \eeqnn 
 Using H\"older's inequality, Fubini's theorem and \eqref{eqn.Moment} to the last expectation, 
 \beqnn
 \mathbf{E}\bigg[ \Big|\int_0^{t_2}  \big| X_\zeta(s)\big|^2 \, ds  \Big|^{p} \bigg]
 \ar\leq\ar \mathbf{E}\bigg[ |t_2|^{p-1}  \int_0^{t_2}  \big| X_\zeta(s)\big|^{2p} \, ds   \bigg] \cr
 \ar\leq\ar t^{p}\cdot \sup_{s\in[0,t]} \mathbf{E}\Big[  \big| X_\zeta(s)\big|^{2p} \Big] 
 \leq C\cdot \big(\zeta\vee \zeta^{2p}\big)\cdot  \big(1+W(t)\big)^{4p-2}\cdot t^p. 
 \eeqnn
 The desired upper bounds follows by combining all estimates above together.
 \qed 
 
 \begin{proposition}\label{Proposition.605}
  For each $p \geq 1$ and $\beta >0$,  there exists a constant $C>0$ such that for any $\zeta, t\geq 0$ and $0\leq t_1<t_2\leq t$,
  \beqnn
   \mathbf{E}\Big[ \big|\widetilde{\mathbf{I}}_{\beta,4}(t_2) -  \widetilde{\mathbf{I}}_{\beta,4}(t_1)  \big|^{2p}  \Big] 
   \leq C\cdot (\zeta\vee \zeta^p)\cdot \big(1+W (t)\big)^{2p-2} \cdot   |t_2-t_1|^p.
  \eeqnn 
 \end{proposition}
 \proof 
 The fact that $W'_\beta(x)=0$ if $x<0$ allows us to write $\widetilde{\mathbf{I}}_{\beta,4}(t_2) - \widetilde{\mathbf{I}}_{\beta,4}(t_1)$ as 
 \beqnn
 \widetilde{\mathbf{I}}_{\beta,4}(t_2) - \widetilde{\mathbf{I}}_{\beta,4}(t_1)
 \ar=\ar \int_0^{t_2} \int_0^{X_\zeta(s)} \big( W'_\beta(t_2-s) -W'_\zeta(t_1-s)\big) \, B_c(ds,dz). 
 \eeqnn
 By using \eqref{BDG0}, there exists a constant $C>0$ depending only on $p$ such that
 \beqnn
 \mathbf{E}\Big[ \big|\widetilde{\mathbf{I}}_{\beta,4}(t_2) - \widetilde{\mathbf{I}}_{\beta,4}(t_1)  \big|^{2p}  \Big]  
 \ar\leq\ar C\cdot \sup_{s\in [0,t_2]}\mathbf{E}\Big[ \big|X_\zeta(s)\big|^p \Big] \cdot  \bigg|\int_{\mathbb{R}}\big( W'_\beta(t_2-s) - W'_\beta(t_1-s)\big)^2\,ds \bigg|^{p} , 
 \eeqnn
 which can be bounded by $C\cdot (\zeta\vee \zeta^p)\cdot \big(1+W(t)\big)^{2p-2} \cdot  |t_2-t_1|^p$ uniformly in $\zeta,t \geq 0$ and $0\leq t_1<t_2\leq t$; see \eqref{eqn.Moment} and \eqref{eqn.50004}. 
 \qed 
 
 \begin{proposition}\label{Proposition.606}
  For each $p \geq 1$ and $\beta > 0$, there exists a constant $C>0$ such that for any $\zeta, t\geq 0$ and $0\leq t_1<t_2\leq t$,
 \beqnn
  \mathbf{E}\Big[ \big|\widetilde{\mathbf{I}}_{\beta,5}(t_2) - \widetilde{\mathbf{I}}_{\beta,5}(t_1)  \big|^{2p}  \Big] 
  \leq C\cdot (\zeta \vee \zeta^p)\cdot \big(1+W (t)\big)^{2p-2} \cdot (1+t)^{2p}  \cdot  |t_2-t_1|^p. 
 \eeqnn  
 \end{proposition}
 \proof 
 Similarly as in the proof of Proposition~\ref{Proposition.605}, we have 
 \beqnn
 \widetilde{\mathbf{I}}_{\beta,5}(t_2) -  \widetilde{\mathbf{I}}_{\beta,5}(t_1)
 =\int_0^{t_2}\int_0^{X_\zeta(s)} \int_0^\infty \Big( \int_{t_2-s-y}^{t_2-s} W'_\beta(r)dr - \int_{t_1-s-y}^{t_1-s} W'_\beta(r)dr  \Big)\, \widetilde{N}_\nu(ds,dz,dy)  
 \eeqnn
 and there exists a constant $C>0$ depending only on $p$ such that
 \beqlb\label{eqn.620}
 \mathbf{E}\Big[ \big|\widetilde{\mathbf{I}}_{\beta,5}(t_2) - \widetilde{\mathbf{I}}_{\beta,5}(t_1)  \big|^{2p}  \Big]
 \ar\leq\ar C\cdot  \sup_{s\in[0,t_2]} \mathbf{E}\Big[ \big| X_\zeta(s)\big|^p \Big]\cdot \big| A_2(t_1,t_2) \big|^p \cr
 \ar\ar +  C\cdot  \sup_{s\in[0,t_2]} \mathbf{E}\Big[ \big| X_\zeta(s)\big|  \Big]\cdot  A_{2p}(t_1,t_2) ,
 \eeqlb
 where 
 \beqnn
 A_{2k}(t_1,t_2) 
 \ar:=\ar \int_{0}^{t_2} ds \int_0^\infty \Big(\int_{t_2-s-y}^{t_2-s} W'_\beta(r)dr - \int_{t_1-s-y}^{t_1-s} W'_\beta(r)dr \Big)^{2k} \,\nu(dy),\quad k\geq 1 . 
 \eeqnn
 For $k=1$, by the fact that $W'_\beta(x)=0$ if $x<0$ and the change of variables, 
 \beqnn
 A_2(t_1,t_2)
 \ar=\ar \int_{t_1}^{t_2} ds \int_0^\infty \Big(\int_{t_2-s-y}^{t_2-s} W'_\beta(r)dr  \Big)^{2} \,\nu(dy)\cr
 \ar\ar +\int_{0}^{t_1} ds \int_0^\infty \Big(\int_{t_2-s-y}^{t_2-s} W'_\beta(r)dr - \int_{t_1-s-y}^{t_1-s} W'_\beta(r)dr \Big)^{2} \,\nu(dy)\cr 
 \ar=\ar \int_{0}^{t_2-t_1} ds \int_0^\infty \Big(\int_{s-y}^{s} W'_\beta(r)dr  \Big)^{2} \,\nu(dy)\cr
 \ar\ar +\int_{0}^{t_1} ds \int_0^\infty \Big(\int_{t_2-t_1+s-y}^{t_2-t_1+s} W'_\beta(r)dr - \int_{s-y}^{s} W'_\beta(r)dr \Big)^{2} \,\nu(dy)
 \eeqnn
 By \eqref{eqn.20101} and  the inequality $\sup_{t\in\mathbb{R}}W_\beta(t)<\infty$, we have $\int_{s-y}^{s} W'_\beta(r)dr \leq C(1\wedge y)$ uniformly in $s,y\geq 0$, which along with \eqref{LevyTriplet} induces that
 \beqnn
 \int_0^{t_2-t_1} ds\int_0^\infty \Big( \int_{s-y}^{s} W'_\beta(r)dr  \Big)^2\,\nu(dy) 
 \ar\leq\ar C \int_0^{t_2-t_1}  ds\int_0^\infty (1\wedge y^2)\,\nu(dy)
 \leq   C\cdot   \big|t_2-t_1\big|,
 \eeqnn
 uniformly in $0\leq t_1<t_2<t$.  
 Additionally, by using the change of variables and then H\"older's inequality, 
 \beqnn
 \lefteqn{ \int_{0}^{t_1} ds \int_0^\infty \Big(\int_{t_2-t_1+s-y}^{t_2-t_1+s} W'_\beta(r)dr - \int_{s-y}^{s} W'_\beta(r)dr \Big)^{2} \,\nu(dy) }\ar\ar\cr 
 \ar=\ar \int_{0}^{t_1} \int_0^\infty \Big( \int_{s-y}^{s} \big( W'_\beta(t_2-t_1+r)- W'_\beta(r)\big) dr  \Big)^2\, ds\,\nu(dy)\cr
 \ar\leq\ar \int_0^\infty \nu(dy) \int_{0}^{t_1}  (s\wedge y) \int_{s-y}^{s} \big( W'_\beta(t_2-t_1+r)- W'_\beta(r)\big)^2 dr \, ds, 
 \eeqnn
 which can be divided into the following two terms
 \beqnn
 B_1(t_1,t_2)\ar:=\ar \int_0^\infty \nu(dy) \int_{0}^{t_1} (s\wedge y) \int_{(s-y)\wedge 0}^0 \big( W'_\beta(t_2-t_1+r) \big)^2 dr \, ds,\cr
 B_2(t_1,t_2)\ar:=\ar  \int_0^\infty \nu(dy) \int_{0}^{t_1} (s\wedge y) \int_{(s-y)\vee 0}^s \big( W'_\beta(t_2-t_1+r)- W'_\beta(r)\big)^2 dr \, ds .
 \eeqnn
 By \eqref{eqn.20101} as well as the two facts that $W'_\beta(x)=0$ if $x<0$ and $ W_\beta(\infty)<\infty$, we have
 \beqnn
 \int_{(s-y)\wedge 0}^0 \big( W'_\beta(t_2-t_1+r) \big)^2 dr \, ds
 \leq C\cdot \mathbf{1}_{\{y>s \}},
 \eeqnn
 uniformly in $x,y\geq 0$ and hence 
 \beqnn
 B_1(t_1,t_2)\ar\leq\ar 
C \int_0^\infty \nu(dy) \int_{0}^{t_1} (s\wedge y) \cdot \mathbf{1}_{\{y>s \}}  \, ds 
= C \int_{0}^{t_1} s\bar\nu(s)\, ds 
\leq C\cdot (1+t)^2 , 
 \eeqnn 
 uniformly in $0\leq t_1<t_2<t$. 
 For $B_2(t_1,t_2)$, by Fubini's theorem,
 \beqnn
 B_2(t_1,t_2)
 \ar=\ar \int_0^\infty \nu(dy) \int_{0}^{t_1}  \big( W'_\beta(t_2-t_1+s )- W'_\beta(s)\big)^2 \int_{s}^{(s+y)\wedge t_1} (r\wedge y) dr \, ds \cr 
 \ar\leq\ar \int_0^\infty (t_1\wedge y)^2 \nu(dy) \cdot \int_\mathbb{R} \big( W'_\beta(t_2-t_1+s)- W'_\beta(s)\big)^2   ds, 
 \eeqnn
 which can be bounded by $C\cdot (1+t)^2\cdot   |t_2-t_1|  $ uniformly in $0\leq t_1<t_2<t$; see \eqref{eqn.50004}. 
 Putting all preceding estimates together induces that 
 \beqlb \label{eqn.621}
 A_2(t_1,t_2) \leq C\cdot (1+t)^2\cdot |t_2-t_1|. 
 \eeqlb
 for some constant $C$ independent of $t,t_1,t_2$. 
 Similarly as in \eqref{eqn.591}-\eqref{eqn.593} we also have
 \beqnn
 \Big|\int_{t_2-s-y}^{t_2-s} W'_\beta(r)dr - \int_{t_1-s-y}^{t_1-s} W'_\beta(r)dr\Big|
 \ar\leq\ar \frac{2}{c}\cdot  |t_2-t_1| 
 \eeqnn
 and then $A_{2p}(t,\delta)\leq C \cdot A_2(t,\delta) \cdot |t_2-t_1|^{2p-2}
 \leq C\cdot(1+t)^2\cdot |t_2-t_1|^{2p-1}.$
 Taking this and \eqref{eqn.621} back into \eqref{eqn.620} and then using \eqref{eqn.Moment}, there exists constant $C>0$ such that for any $t\geq 0$ and $0\leq t_1<t_2\leq t$, 
 \beqnn
 \mathbf{E}\Big[ \big|\widetilde{\mathbf{I}}_{\beta,5}(t_2) - \widetilde{\mathbf{I}}_{\beta,5}(t_1)  \big|^{2p}  \Big] 
 \ar\leq\ar  C\cdot (\zeta \vee \zeta^p)\cdot (1+W(t))^{2p-2} (1+t)^{2p}  \cdot  |t_2-t_1|^p \cr
 \ar\ar +  C\cdot   \zeta  \cdot  (1+t)^2\cdot |t_2-t_1| \cdot (|t_2-t_1|\wedge W(t_2) )^{2p-2}, 
 \eeqnn
 and hence the desired upper bound holds. 
 \qed 
 
 Armed with the preceding propositions, we are now ready to prove Theorem~\ref{Thm.Regularity}. 
 To obtain the desired uniform upper bound for all moments of the H\"older coefficient, we need the Garsia-Rodemich-Rumsey inequality; see e.g. Lemma~1.1 in \cite{GarsiaRodemichRumseyRosenblatt1970}  with $\psi(u)=|u|^{2p}$ and $p(u)=|u|^{q+\frac{1}{2p}}$ for $p,q>0$ such that $2pq>1$ or Theorem~1.1 in \cite{Walsh1986}. 
 It states that for a continuous function $f$ on $\mathbb{R}_+$, there exists a constant $C_{p,q}>0$ depending only on $p$ and $q$ such that for any $x_2>x_1\geq 0$, 
 \beqlb\label{eqn.405}
 \big|f(x_2)-f(x_1) \big|^{2p} \leq C_{p,q} \cdot \big| x_2-x_1\big|^{2pq-1} \int_{x_1}^{x_2} ds \int_{x_1}^{x_2} \frac{|f(s)-f(r)|^{2p}}{|s-r|^{2pq+1}}dr.
 \eeqlb

 \textit{\textbf{Proof of Theorem~\ref{Thm.Regularity}.}} 
 The Kolmogorov continuity theorem along with Proposition~\ref{Proposition.603}-\ref{Proposition.606} yields that $ \widetilde{\mathbf{I}}_{\beta,1}$ is Lipschitz continuous and other four terms in \eqref{Main.SVE0212}, $ \widetilde{\mathbf{I}}_{\beta,i}$ with $i=2,\cdots,5$, are all H\"older continuous with exponent in $(0,1/2)$.  
 Consequently, the solution $X_\zeta$ is H\"older continuous with exponent strictly less than $1/2$. 
 
 We now prove the first upper bound in \eqref{Regularity.01}. 
 For $\kappa\in (0,1/2)$, by the power inequality, there exists constant $C>0$ depending only on $p$ such that for any $t\geq 0$, 
 \beqlb\label{eqn.623}
 \big\|X_\zeta\big\|_{C^{0,\kappa}_t}^{2p} 
 \leq C \sum_{i=1}^5 \big\| \widetilde{\mathbf{I}}_{\beta,i}\big\|_{C^{0,\kappa}_t}^{2p} 
 \quad \mbox{and}\quad 
 \mathbf{E}\Big[ \big\|X_\zeta\big\|_{C^{0,\kappa}_t}^{2p} \Big]
 \leq C \cdot \big\| \widetilde{\mathbf{I}}_{\beta,1}\big\|_{C^{0,\kappa}_t}^{2p} + C \sum_{i=2}^5  \mathbf{E}\Big[ \big\| \widetilde{\mathbf{I}}_{\beta,i}\big\|_{C^{0,\kappa}_t}^{2p}\Big]
 .
 \eeqlb
 By Lemma~\ref{Proposition.603} with $\beta>0$, there exists a constant $C>0$ such that for any $t\geq 0$, 
 \beqnn 
 \big\| \widetilde{\mathbf{I}}_{\beta,1}\big\|_{C^{0,\kappa}_t}^{2p}
 \ar=\ar  \zeta^{2p}\cdot \beta^{2p}\cdot \sup_{0\leq t_1< t_2\leq t} \frac{|W_\beta(t_2)-W_\beta(t_1)|^{2p}}{|t_2-t_1|^{2p\kappa}}
 \leq  
 \zeta^{2p}\cdot \frac{\beta^{2p}}{c^{2p\kappa}} \cdot \big|W_\beta(t)\big|^{2p(1-\kappa)} \leq C \cdot \zeta^{2p}. 
 \eeqnn
 We now provide the upper bound estimates for the last four expectations in \eqref{eqn.623} by using the inequality \eqref{eqn.405}. 
 More precisely, an applications of the inequality \eqref{eqn.405} to $\widetilde{\mathbf{I}}_{\beta,2}$ shows that 
 \beqnn
 \big| \widetilde{\mathbf{I}}_{\beta,2}(t_2)- \widetilde{\mathbf{I}}_{\beta,2}(t_1) \big|^{2p} 
 \leq C_{p,q}\cdot  \big| t_2-t_1\big|^{2pq-1} \int_{t_1}^{t_2} ds \int_{t_1}^{t_2} \frac{| \widetilde{\mathbf{I}}_{\beta,2}(s)- \widetilde{\mathbf{I}}_{\beta,2}(r)|^{2p}}{|s-r|^{2pq+1}}dr 
 \eeqnn
 and then 
 \beqnn
 \big\| \widetilde{\mathbf{I}}_{\beta,2}\big\|_{C^{0,\kappa}_t}^{2p}
 \ar\leq\ar C_{p,q}\cdot \sup_{0\leq t_1< t_2\leq t} \big| t_2-t_1\big|^{2p(q-\kappa)-1 } \int_{t_1}^{t_2} ds \int_{t_1}^{t_2} \frac{| \widetilde{\mathbf{I}}_{\beta,2}(s)- \widetilde{\mathbf{I}}_{\beta,2}(r)|^{2p}}{|s-r|^{2pq+1}}dr.
 \eeqnn
 In particular, by choosing $p>(1-2\kappa)^{-1}$ and $q=\frac{1}{2p}+\kappa$ we have
 \beqnn
 \big\| \widetilde{\mathbf{I}}_{\beta,2}\big\|_{C^{0,\kappa}_t}^{2p}
 \ar\leq\ar C_{p,q}  \int_{0}^{t} ds \int_{0}^{t} \frac{| \widetilde{\mathbf{I}}_{\beta,2}(s)- \widetilde{\mathbf{I}}_{\beta,2}(r)|^{2p}}{|s-r|^{2p\kappa+2}}dr.
 \eeqnn
 Taking expectations on both sides and then using Fubini's theorem as well as  Proposition~\ref{Proposition.604}, 
 \beqnn
 \mathbf{E}\Big[ \big\| \widetilde{\mathbf{I}}_{\beta,2}\big\|_{C^{0,\kappa}_t}^{2p} \Big]
 \ar\leq\ar C_{p,q}  \int_{0}^{t} ds \int_{0}^{t} \frac{\mathbf{E}\big[| \widetilde{\mathbf{I}}_{\beta,2}(s)- \widetilde{\mathbf{I}}_{\beta,2}(r)|^{2p}\big]}{|s-r|^{2p\kappa+2}}dr\cr
 \ar\leq\ar C \cdot (\zeta \vee \zeta^p)  \cdot \int_{0}^{t} ds \int_{0}^{t}  |s-r|^{p-2p\kappa-2}dr
 \leq  C \cdot (\zeta \vee \zeta^p)  \cdot t^{p-2p\kappa}, 
 \eeqnn
 uniformly in $\zeta, t \geq 0$. 
 Similarly, with the help of Lemma~\ref{Proposition.607}-\ref{Proposition.606}, we can also prove that 
 \beqnn
 \mathbf{E}\Big[ \big\|\widetilde{\mathbf{I}}_{\beta,3}\big\|_{C^{0,\kappa}_t}^{2p} \Big] 
 \ar\leq\ar C\cdot  (\zeta \vee \zeta^{2p})\cdot  \big(1+W(t)\big)^{4p-2} \cdot t^{2p(1-\kappa)},\cr
 \mathbf{E}\Big[ \big\|\widetilde{\mathbf{I}}_{\beta,4}\big\|_{C^{0,\kappa}_t}^{2p} \Big]
 \ar\leq\ar C \cdot (\zeta \vee \zeta^p) \cdot \big(1+W(t)\big)^{2p-2} \cdot t^{p-2p\kappa},\cr
 \mathbf{E}\Big[ \big\|\widetilde{\mathbf{I}}_{\beta,5}\big\|_{C^{0,\kappa}_t}^{2p} \Big] 
 \ar\leq\ar C \cdot (\zeta \vee \zeta^p) \cdot \big(1+W(t)\big)^{2p-2} \cdot (1+t)^{2p }\cdot t^{p-2p\kappa}.
 \eeqnn 
 Plugging all preceding upper bound estimates into the right side of the second inequality in \eqref{eqn.623} and then using the inequality $W(t)\leq t/c$ for all $t\geq 0$, we have uniformly in $\zeta ,t\geq 0$,
 \beqnn
 \mathbf{E}\Big[ \big\|X_\zeta\big\|_{C^{0,\kappa}_t}^{2p} \Big] \leq C \cdot (\zeta \vee \zeta^{2p}) \cdot (1+t)^{2p(3-\kappa)}.
 \eeqnn
 
 For the second upper bound in \eqref{Regularity.01},
 by the triangle inequality and H\"older continuity we have 
 \beqnn
 \sup_{s\in[0,t]} \big|X_\zeta(s)\big| \leq \zeta + \sup_{s\in[0,t]} \big|X_\zeta(s)-\zeta\big| 
 \leq \zeta + \big\|X_\zeta\big\|_{C^{0,\kappa}_t} \cdot t^\kappa ,  
 \eeqnn
 for any $\kappa\in (0,1/2)$.
 Then by the power inequality and the preceding results, 
 \beqnn
 \mathbf{E}\bigg[ \sup_{s\in[0,t]} \big|X_\zeta(s)\big|^{2p} \bigg] 
 \ar\leq\ar C\cdot \zeta^{2p} +  C\cdot \mathbf{E}\Big[\big\|X_\zeta\big\|_{C^{0,\kappa}_t}^{2p}\Big] \cdot t^{2p\kappa}
 \leq  C\cdot (\zeta \vee \zeta^{2p}) \cdot (1+t)^{6p} ,
 \eeqnn
 for some constant $C>0$ independent of $\zeta$ and $t$.
 \qed

   \section{Nonlinear Volterra equation and Laplace functionals}
 \label{Sec.NonlinearVE}
 \setcounter{equation}{0} 
 
  This section is devoted to provide detailed proof for Theorem~\ref{Thm.LaplaceF} as well as the proof of weak uniqueness of solution to \eqref{MainThm.SVE}. 
  It consists of two subsections in which the well-posedness of the nonlinear Volterra equation \eqref{NonlinearVolterra} and the affine representation for the Laplace functional of solutions to \eqref{MainThm.SVE} are established respectively.

 \subsection{Well-posedness of nonlinear Volterra equation}

 As a preparation, we first provide in the next two propositions an upper bounded estimate and Lipschitz continuity of the nonlinear operator $\boldsymbol{\mathcal{R}}$.
 For convention, we write $\mathcal{V} \circ f(t)$ for the second term on the right side of \eqref{OperatorR} and hence
 \beqlb\label{eqn.50001}
 \boldsymbol{\mathcal{R}} \circ f(t) =  c\cdot \big(f(t)\big)^2 + \mathcal{V} \circ f(t).
 \eeqlb
  
 \begin{proposition} \label{Prop.701}
 There exists a constant $C>0$ such that for any $f\in L^\infty_{\rm loc}(\mathbb{R}_+;\mathbb{R})$ and $t\geq 0$,
 \beqnn
 0\leq  \big(\boldsymbol{\mathcal{R}} \circ f\big)*W'(t)
  \leq C \cdot \exp\big\{\|f\|_{L^\infty_t}\cdot t \big\}    \cdot (1+t)^2 \cdot 
  \big\| f \big\|_{L^2_t}^2 .
 \eeqnn
 \end{proposition}
 \proof By \eqref{eqn.50001}, we have $\big(\boldsymbol{\mathcal{R}} \circ f\big)*W'(t) =c\cdot |f|^2*W'(t) + (\mathcal{V} \circ f)*W'(t) $. By \eqref{eqn.201},  
 \beqlb\label{eqn.701}
 0\leq c\cdot |f|^2*W'(t) = c\int_0^t |f(s)|^2 W'(t-s)ds \leq  \int_0^t |f(s)|^2 ds. 
 \eeqlb
 Moreover, applying the inequality $0\leq e^{-z}-1+z \leq |z|^2 \cdot e^{|z|} $ for all $z\in\mathbb{R}$ to $\mathcal{V} \circ f(t)$ induces that 
 \beqnn
 0\leq \mathcal{V} \circ f(t)  
 \ar\leq\ar  
 \int_0^\infty \exp\Big\{\int_{(t-y)^+}^{t} \big|f(r)\big|\,dr\Big\} \Big(\int_{(t-y)^+}^{t} f(r)\,dr\Big)^2 \nu ( dy ) \cr
 \ar\leq\ar \exp\big\{\|f\|_{L^\infty_t}\cdot t \big\}   \int_0^\infty   \Big(\int_{(t-y)^+}^{t} f(r)\,dr\Big)^2 \nu ( dy ).
 \eeqnn
 By using H\"older's inequality to the last inner integral, we further have  
 \beqnn
  0\leq \mathcal{V} \circ f(t)  
  \ar\leq\ar \exp\big\{\|f\|_{L^\infty_t}\cdot t \big\}    \int_0^\infty (t\wedge y) \nu ( dy )   \int_{(t-y)^+}^{t} |f(r)|^2\, dr , 
 \eeqnn  
 which along with \eqref{eqn.201} and Fubini's theorem yields that
 \beqnn
 0\leq \big(\mathcal{V} \circ f\big)*W'(t) 
 \ar\leq\ar \frac{1}{c} \int_0^t \exp\big\{\|f\|_{L^\infty_s}\cdot s \big\}     \int_0^\infty (s\wedge y) \nu ( dy )   \int_{(s-y)^+}^{s} |f(r)|^2\, dr \,ds \cr
 \ar\leq\ar \frac{1}{c} \cdot \exp\big\{\|f\|_{L^\infty_t}\cdot t \big\}   \int_0^\infty (t\wedge y) \nu ( dy )  \int_0^t \int_{(s-y)^+}^{s} |f(r)|^2\, dr \,ds \cr
 \ar=\ar \frac{1}{c} \cdot \exp\big\{\|f\|_{L^\infty_t}\cdot t \big\}    \int_0^t y\, \nu ( dy )  \int_0^t \int_{(s-y)^+}^{s} |f(r)|^2\, dr \,ds \cr
 \ar\ar + \frac{1}{c} \cdot \exp\big\{\|f\|_{L^\infty_t}\cdot t \big\}    \int_t^\infty t\, \nu ( dy )  \int_0^t \int_{(s-y)^+}^{s} |f(r)|^2\, dr \,ds.
 \eeqnn 
 The first term on the right side of the last equality can be bounded by 
 \beqnn
 \lefteqn{\frac{1}{c} \cdot \exp\big\{\|f\|_{L^\infty_t}\cdot t \big\}     \int_0^t y\, \nu ( dy )  \int_0^{t+y} \int_{(s-y)^+}^{s\wedge t} |f(r)|^2\, dr \,ds} \ar\ar\cr
 \ar=\ar \frac{1}{c} \cdot \exp\big\{\|f\|_{L^\infty_t}\cdot t \big\}   \int_0^t y\, \nu ( dy )  \int_0^t|f(s)|^2\int_s^{s+y}  \,dr\, ds \cr
 \ar=\ar \frac{1}{c} \cdot \exp\big\{\|f\|_{L^\infty_t}\cdot t \big\}   \int_0^t y^2\, \nu ( dy )  \int_0^t|f(s)|^2\, ds .
 \eeqnn 
 Here the first equality is obtained by using Fubini's theorem.  
 Moreover, the second term equals to
 \beqnn
 \frac{1}{c} \cdot \exp\big\{\|f\|_{L^\infty_t}\cdot t \big\}   \int_t^\infty t\, \nu ( dy )  \int_0^t \int_{0}^{s} |f(r)|^2\, dr \,ds 
 \leq   \frac{1}{c} \cdot \exp\big\{\|f\|_{L^\infty_t}\cdot t \big\}    \int_t^\infty t^2\, \nu ( dy )  \int_0^t |f(s)|^2\, ds .
 \eeqnn
 Putting these results together and then using \eqref{LevyTriplet}, we have 
 \beqnn
  0\leq \big(\mathcal{V} \circ f\big)*W'(t)
  \ar\leq\ar \frac{1}{c} \cdot \exp\big\{\|f\|_{L^\infty_t}\cdot t \big\}   \int_0^\infty (t\wedge y)^2\, \nu ( dy )  \cdot \int_0^t |f(s)|^2\, ds\cr
  \ar\leq\ar C \cdot \exp\big\{\|f\|_{L^\infty_t}\cdot t \big\}   \cdot (1+t)^2 \cdot \int_0^t |f(s)|^2\, ds,
 \eeqnn
 for some constant $C>0$ independent of $f$ and $t$. 
 The desired upper bound follows by combining this together with \eqref{eqn.701}.
 \qed 
  
%
  
 \begin{proposition}\label{Prop.703}
 There exists a constant $C>0$ such that for any $f_1,f_2\in L^\infty_{\rm loc}(\mathbb{R}_+;\mathbb{R})$ and $t\geq 0$,
 \beqlb \label{eqn.702}
 \Big|\big(\boldsymbol{\mathcal{R}} \circ f_1 -\boldsymbol{\mathcal{R}} \circ f_2 \big)*W'(t)\Big|
 \leq C\cdot \exp\big\{ ( \|f_1\|_{L^\infty_t} + \|f_2\|_{L^\infty_t} ) \cdot (1+t)  \big\} \cdot   (1+t)^2 \cdot \big\|f_1 -f_2 \big\|_{L^1_t} .
 \eeqlb
 \end{proposition}
 \proof By \eqref{eqn.50001}, we first have
 \beqnn
 \Big|\big(\boldsymbol{\mathcal{R}} \circ f_1 -\boldsymbol{\mathcal{R}} \circ f_2 \big)*W'(t)\Big| \leq c\cdot \Big| |f_1|^2*W'(t) -  |f_2|^2*W'(t)\Big| + \Big|\big(\mathcal{V} \circ f_1 - \mathcal{V} \circ f_1 \big)*W'(t)\Big|.
 \eeqnn
  By using \eqref{eqn.201} to the first term on the right side, 
 \beqnn
  c\cdot\Big|  |f_1|^2*W'(t) -  |f_2|^2*W'(t)\Big| 
  \ar\leq\ar \int_0^t \big||f_1(s)|^2-|f_2(s)|^2\big| \, ds\cr
  \ar=\ar \int_0^t |f_1(s)+f_2(s)|\cdot |f_1(s)-f_2(s)|\,  ds\cr
  \ar\leq\ar \big( \|f_1\|_{L^\infty_t} + \|f_2\|_{L^\infty_t} \big) \cdot  \int_0^t |f_1(s)-f_2(s)| \, ds.
 \eeqnn
 Additionally, an application of the next inequality that can be proved immediately by using the mean-value theorem
 \beqnn
  \big|(e^{-x}-1+x)-(e^{-z}-1+z) \big| \leq  e^{ |x|+|z| } \cdot \big(|x|+|z| \big) \cdot |x-z|,\quad x,z\in\mathbb{R},
 \eeqnn
 to $\big|\mathcal{V} \circ f_1(t)-\mathcal{V} \circ f_2(t)\big|$ induces that it can be bounded by
 \beqnn
  \lefteqn{\int_0^\infty e^{ \int_{(t-y)^+}^t \big(|f_1(r) | + |f_2(r) |\big) \,dr  }\cdot   \int_{(t-y)^+}^t \big(|f_1(r) | + |f_2(r) |\big) \,dr   \cdot \int_{(t-y)^+}^t \big|f_1(r)- f_2(r)\big|\,dr \,\nu(dy) } \quad \ar\ar\cr
 \ar\leq\ar \exp\big\{ ( \|f_1\|_{L^\infty_t} + \|f_2\|_{L^\infty_t} ) \cdot t \big\} \cdot \big(\|f_1\|_{L^\infty_t} + \|f_2\|_{L^\infty_t}\big)\cdot \int_0^\infty (t\wedge y) \int_{(t-y)^+}^t \big|f_1(r)- f_2(r)\big|\,dr \,\nu(dy).
 \eeqnn
 This together with \eqref{eqn.201} induces  that 
 \beqnn
 \Big|\big(\mathcal{V} \circ f_1 -\mathcal{V} \circ f_1 \big)*W'(t)\Big|  
 \ar\leq\ar \frac{1}{c} \int_0^t \big|\mathcal{V} \circ f_1(s)-\mathcal{V} \circ f_2(s)\big|\,ds \cr
 \ar\leq\ar \frac{1}{c}\cdot \exp\big\{ ( \|f_1\|_{L^\infty_t} + \|f_2\|_{L^\infty_t} ) \cdot (1+t)  \big\}\cr
 \ar\ar \times   \int_0^t\, ds \int_0^\infty (s\wedge y) \int_{(s-y)^+}^s\big|f_1(r)- f_2(r)\big|\,dr \,\nu(dy) . 
 \eeqnn
 Similarly as in the proof of Proposition~\ref{Prop.701}, there exists a constant $C>0$ independent of $t$ and $f_1,f_2$ such that
 \beqnn
 \int_0^t\, ds \int_0^\infty (s\wedge y) \int_{(s-y)^+}^s\big|f_1(r)- f_2(r)\big|\,dr \,\nu(dy) 
 \leq C\cdot (1+t)^2 \cdot \int_0^t \big|f_1(s)-f_2(s)\big|\, ds  . 
 \eeqnn
 The desired upper bound \eqref{eqn.702} follows directly by putting all preceding result together.
 \qed 
 
 A function $V_\mu$ on $\mathbb{R}_+$ is said to be a \textsl{$L^\infty_{\rm loc}(\mathbb{R}_+;\mathbb{R})$-noncontinuable solution} of \eqref{NonlinearVolterra} if there exists a constant $T_\infty \in (0,\infty]$ such that $V_\mu \in L^\infty ([0,T];\mathbb{R})$  for any $T\in(0,T_\infty)$ and $ \big\|V_\mu \big\|_{L^\infty_{T_\infty}}= \infty$ if $T_\infty<\infty$. 
 Moreover, it turns to be a $L^\infty_{\rm loc}(\mathbb{R}_+;\mathbb{R})$-global solution of \eqref{NonlinearVolterra} if $T_\infty=\infty$. 
 In the next lemma, we show that the existence of $L^\infty_{\rm loc}(\mathbb{R}_+;\mathbb{R})$-noncontinuable solution to \eqref{NonlinearVolterra} immediately induces that the global solution uniquely exists.


 
 \begin{lemma} \label{Lemma.704}
  If \eqref{NonlinearVolterra} has a $L^\infty_{\rm loc}(\mathbb{R}_+;\mathbb{R})$-noncontinuable solution, then it has a unique global solution  $V_\mu \in D(\mathbb{R}_+;\mathbb{R}_+)$ that satisfies  
  \beqnn
  V_\mu (t) \leq     \frac{\mu([0,t])}{c} ,\quad t\geq 0.  
  \eeqnn
 \end{lemma}
 \proof Assume that $(\widetilde{V}_\mu,T_\infty)$ is a $L^\infty_{\rm loc}(\mathbb{R}_+;\mathbb{R})$-noncontinuable solution. 
 By Proposition~\ref{Prop.701} and \eqref{eqn.201}, we have for almost every $t\geq 0$,
 \beqnn
 \widetilde{V}_\mu (t) \leq    W'*d\mu(t) \leq \frac{\mu([0,t])}{c}.  
 \eeqnn
 For any $T\in[0,T_\infty)$, by using Proposition~\ref{Prop.701} again there exists a constant $C_T>0$ such that 
 \beqnn
 -\widetilde{V}_\mu (t) = -  W'*d\mu(t) +  \big(\boldsymbol{\mathcal{R}} \circ\widetilde{V}_\mu\big)*W'(t) \leq C_T   \int_0^t |-\widetilde{V}_\mu|^2\, ds,
 \quad t\in[0,T].
 \eeqnn
 By the classic comparison theorem, we have $-\widetilde{V}_\mu(t)\leq  V(t)  $ for almost every $t\in [0,T]$, where $V\equiv 0$ is the unique solution to the Riccati equation $dV(t) = C_T \cdot |V(t)|^2$. 
 In conclusion, 
 \beqnn
 0\leq \widetilde{V}_\mu (t) \leq \frac{\mu([0,t])}{c},\quad t\in [0,T_\infty),
 \eeqnn
 which yields that $T_\infty =\infty$ and  $\widetilde{V}_\mu$ is a $L^\infty_{\rm loc}(\mathbb{R}_+;\mathbb{R})$-global solution. 
 
 Associated to $\widetilde{V}_\mu$ we define the function $ V_\mu$ as follows
 \beqnn
 V_\mu (t) :=  W'*d\mu (t)-\big( \boldsymbol{\mathcal{R}}\circ \widetilde{V}_\mu\big) *W'(t),\quad t\geq 0.
 \eeqnn
 It is obvious that $ V_\mu (t)= \widetilde{V}_\mu(t)$ for almost every $t\geq 0$. 
 The continuity of $\big( \boldsymbol{\mathcal{R}}\circ \widetilde{V}_\mu\big) *W'$ follows by the regularity property of convolution. 
 Moreover, by Corollary~6.2 in \cite[p.98]{GripenbergLondenStaffans1990},
 \beqnn
 W'*d\mu (t) = W'(0)\cdot \mu\big([0,t]\big) + \int_{[0,t]} \big(W'(t-s)-W'(0)\big) \mu(ds),
 \eeqnn
  is c\`ad\`ag in $t$.
  Hence $V_\mu \in D(\mathbb{R}_+;\mathbb{R}_+)$ satisfies \eqref{NonlinearVolterra} for all $t\geq 0$.  
  
 The uniqueness follows directly from Proposition~\ref{Prop.703} and Gr\"onwall's inequality. Indeed, assume that $V_\mu^{(1)}$ and $V_\mu^{(2)}$ are two global solutions in $D(\mathbb{R}_+;\mathbb{R}_+)$ to \eqref{NonlinearVolterra}. 
 For each $T\geq 0$, by  Proposition~\ref{Prop.703} there exists a constant $C>0$ such that for any $t\in[0,T]$,
 \beqnn
 \big|V_\mu^{(1)}(t)-V_\mu^{(2)} (t) \big| \ar=\ar \Big|\big(\boldsymbol{\mathcal{R}} \circ V_\mu^{(1)} - \boldsymbol{\mathcal{R}} \circ V_\mu^{(2)}\big)*W'(t)\Big|
 \leq C  \int_0^t \big| V_\mu^{(1)}(s)-V_\mu^{(2)} (s) \big|ds,
 \eeqnn
 which along with Gr\"onwall's inequality induces that $V_\mu^{(1)}(t)=V_\mu^{(2)} (t) $ for all $t\in[0,T]$. 
 \qed

 \textit{\textbf{Proof of Theorem~\ref{Thm.LaplaceF}: Part I.}} By Lemma~\ref{Lemma.704}, it suffices to prove the existence of $L^\infty_{\rm loc}(\mathbb{R}_+;\mathbb{R})$-noncontinuable solutions to \eqref{NonlinearVolterra}. 
 The proof is carried out in the following three steps. 
 
  \medskip

 {\it Step 1.} We first prove the existence of local solutions near $0$. 
 Consider a mapping $\mathcal{R}_0$ that acts on $f\in L^\infty_{\rm loc}(\mathbb{R}_+;\mathbb{R})$  according to
 \beqnn
 \boldsymbol{\mathcal{R}}_0 \circ f(t) :=  W'*d\mu(t)  - \big(\boldsymbol{\mathcal{R}} \circ  f\big)*W'(t),\quad t\geq 0.
 \eeqnn
 For each $T,K\geq 0$, let $\mathcal{B}_{T,K}$ be the collection of all functions $f\in L^\infty ([0,T];\mathbb{R})$ with $\|f\|_{L^\infty_T}\leq K$. 
 It can be easily identify that $\mathcal{B}_{T,K}$ is a closed, bounded and convex subset in $L^\infty([0,T];\mathbb{R})$.
 For any $f_1,f_2 \in \mathcal{B}_{T,K}$, by \eqref{eqn.201}, Proposition~\ref{Prop.701} and \ref{Prop.703} there exists a constant $C_0>0$ such that for any $K\geq0$ and $T\in(0,1]$,
 \beqnn
 \big\|\boldsymbol{\mathcal{R}}_0 \circ f_1\big\|_{L^\infty_T} 
 \leq C_0 + C_0 \cdot K^2 e^{K}   \cdot T
 \quad \mbox{and}\quad 
 \big\|\boldsymbol{\mathcal{R}}_0 \circ f_1- \boldsymbol{\mathcal{R}}_0 \circ f_2\big\|_{L^\infty_T} 
 \leq C_0\cdot T e^{4K} \cdot \big\|f_1-f_2\big\|_{L^\infty_T}.
 \eeqnn 
 Choosing $K_0\in (1,1+1/C_0)$ and $T_0\in(0,e^{-4K_0}]$, we have  
 \beqnn
 C_0 + C_0 \cdot K_0^2 e^{K_0}   \cdot T_0 \leq K_0
 \quad \mbox{and}\quad 
 C_0\cdot T_0 e^{4K_0}<1 . 
 \eeqnn 
 Consequently, the mapping $\boldsymbol{\mathcal{R}}_0$ is a contractive from $\mathcal{B}_{T_0,K_0}$ to itself. 
 By Banach's fixed point theorem, there exists a unique point $V_0 \in \mathcal{B}_{T_0,K_0}$ such that 
 \beqnn
 V_0 (t)= \boldsymbol{\mathcal{R}}_0 \circ V_0(t) = W'*d\mu(t) - \big(\boldsymbol{\mathcal{R}} \circ  V_0\big)*W'(t), 
 \eeqnn
 for almost every $t\in[0,T_0]$ and hence $V_0$ is a local solution of \eqref{NonlinearVolterra} on $[0,T_0]$.  
 
  \medskip 
 
 {\it Step 2.} We now extend the preceding local solution onto a larger interval.
 Denote by $\mathcal{T}$ the collection of all $T>0$ such that \eqref{NonlinearVolterra} has a $L^\infty([0,T];\mathbb{R})$-local solution.
 We assert that $\mathcal{T}$ is an open interval containing $[0,T_0]$.
 Indeed, for any $t_0\in \mathcal{T}$ and some $k_0>0$, assume that $v_0\in \mathcal{B}_{t_0,k_0}$ is a local solution of \eqref{NonlinearVolterra}.
 For $t\geq 0$, let
 \beqnn
 H_1(t):= W'*d\mu(t_0+t)-\int_0^{t_0} \boldsymbol{\mathcal{R}} \circ  v_0(s) W'(t_0+t-s)ds.
 \eeqnn
 Similarly as in the proof of Proposition~\ref{Prop.701}, there exists a constant $C>0$ such that for any $t\geq 0$,
 \beqnn
 \big|H_1(t) \big| \leq \frac{\mu([0,t_0+t])}{c} + C .
 \eeqnn
 We consider a mapping $\mathcal{R}_1$ acting on functions $f\in L^\infty_{\rm loc}(\mathbb{R}_+;\mathbb{R})$ by 
 \beqnn
 \boldsymbol{\mathcal{R}}_1\circ f(t):=  H_1(t) -  \big(\boldsymbol{\mathcal{R}} \circ  f\big)*W'(t),\quad t\geq 0 .
 \eeqnn 
 From Proposition~\ref{Prop.701} and \ref{Prop.703}, there exists a constant $C_1>0$ such that for any $T\in[0,1]$, $K>0$, $f_1,f_2\in \mathcal{B}_{T,K}$ and $t\in[0,T]$,
 \beqnn 
 \|\boldsymbol{\mathcal{R}}_1 \circ f_1\|_{L^\infty_T} 
 \leq C_1 + C_1\cdot K^2 e^{K}   \cdot T
 \quad \mbox{and}\quad 
 \|\boldsymbol{\mathcal{R}}_1 \circ f_1 - \boldsymbol{\mathcal{R}}_1 \circ f_2\|_{L^\infty_T} 
 \leq C_1\cdot T e^{4K} \cdot \|f_1-f_2\|_{L^\infty_T}.
 \eeqnn  
 Similarly as in Step 1, we choose $K_1\in (1,1+1/C_1)$ and $T_1\in(0,e^{-4K_1}]$, which induces that 
 \beqnn
 C_1 + C_1 \cdot K_1^2 e^{K_1}   \cdot T_1 \leq K_1
 \quad \mbox{and}\quad 
 C_1\cdot T_1 e^{4K_1}<1.
 \eeqnn 
 Hence $\boldsymbol{\mathcal{R}}_1$ is a contractive mapping from $\mathcal{B}_{T_1,K_1}$ to itself. 
 By Banach's fixed point theorem, there exists a unique point $v_1 \in \mathcal{B}_{T_1,K_1}$ such that 
 \beqnn
 v_1 (t)= \boldsymbol{\mathcal{R}}_1 \circ v_1(t) = H_1(t) - \big(\boldsymbol{\mathcal{R}} \circ  v_1\big)*W'(t), 
 \eeqnn
 for almost every $t\in[0, T_1]$. 
 It is easy to identify that the function
 \beqnn
 V_1(t) := v_0(t) \cdot \mathbf{1}_{\{t\in[0,t_0]\}} +v_1(t-t_0) \cdot \mathbf{1}_{\{t\in (t_0,t_0+T_1]\}},\quad t\in [0,t_0+T_1].
 \eeqnn 
 is a $L^\infty([0,t_0+T_1];\mathbb{R})$-local solution to \eqref{NonlinearVolterra} and hence the interval $\mathcal{T}$ is open.  
 
 \medskip
 
 {\it Step 3.} We now prove the existence of $L^\infty_{\rm loc}(\mathbb{R}_+;\mathbb{R})$-noncontinuable solutions. 
 Let $T_\infty := \sup \mathcal{T}$ and $V \in L^\infty([0,T_\infty);\mathbb{R})$ be a solution of \eqref{NonlinearVolterra} on $[0,T_\infty)$.
 To assert that $(V,T_\infty)$ is a $L^\infty_{\rm loc}(\mathbb{R}_+;\mathbb{R})$-noncontinuable solution to \eqref{NonlinearVolterra}, it remains to identify that $\|V\|_{L^\infty_{T_\infty}}=\infty$ if $T_\infty<\infty$. 
 If not, one can repeat Step 2 to find two constants $t_1\in(0,1)$ and $k_1>0$ such that the equation \eqref{NonlinearVolterra} has a local solution in $L^\infty([0,T_\infty +t_1];\mathbb{R})$. 
 Consequently, we have $T_\infty +t_1 \in \mathcal{T}$, which  contradicts to the definition of $T_\infty$.  
 \qed


 \subsection{Laplace functionals and weak uniqueness}
 \label{Sec.LaplaceFunctionals}
 
 Associated to the two solutions $X_\zeta$ of \eqref{MainThm.SVE} and $V_\mu$ of \eqref{NonlinearVolterra},  we first introduce in the next proposition an auxiliary process that plays an important role in the following argument. 
 For convention, we set $V_\mu(t)=0$ for $t<0$. 
  

 \begin{proposition} \label{MartingaleReProp}
 For any $T\geq 0$, the random variable 
 \beqlb\label{eqn.720}
  \boldsymbol{Y_T}
  :=  X_\zeta*d\mu(T) - \big(\boldsymbol{\mathcal{R}}\circ V_\mu\big)  *X_\zeta(T) 
 \eeqlb
 is integrable and has the following equivalent representation
 \beqlb \label{eqn.721}
 \boldsymbol{Y_T}
 \ar=\ar  \zeta\cdot c \cdot  V_\mu (T)+ \int_0^\zeta \int_0^\infty \Big( \int_{T-y}^{T} V_\mu (r)\,dr\Big) N_0(dz,dy)\cr
 \ar\ar + \int_0^T \int_0^{X_\zeta(s)}   V_\mu (T-s) \, B_c(ds,dz) \cr
 \ar\ar + \int_0^T\int_0^{X_\zeta(s)} \int_0^\infty \Big( \int_{T-s-y}^{T-s} V_\mu (r)\,dr\Big) \, \widetilde{N}_\nu(ds,dz,dy) .
 \eeqlb
 \end{proposition}
 \proof The local boundedness of $V_\mu$ induces that $\boldsymbol{\mathcal{R}}\circ V_\mu \in L^\infty_{\rm loc}(\mathbb{R}_+;\mathbb{R})$, which together with  Theorem~\ref{Thm.Moment} immediately yields that $\mathbf{E}\big[|\boldsymbol{Y_T}|\big]<\infty$. 
 Plugging \eqref{eqn.710} into the right side of \eqref{eqn.720},  
 \beqnn
 \boldsymbol{Y_T} = \sum_{i=1}^4  \mathbf{I}_i *\big( d\mu  - \boldsymbol{\mathcal{R}}\circ V_\mu  \big) (T).
 \eeqnn
 Thus it suffices to prove that the preceding four summands equal to the corresponding  terms on the right side of \eqref{eqn.721} respectively. 
 
 \begin{enumerate}
  \item[$\bullet$] For $i=1$, by \eqref{NonlinearVolterra} we have 
 	\beqnn
 	\mathbf{I}_1 *\big( d\mu - \boldsymbol{\mathcal{R}}\circ V_\mu  \big) (T)
 	\ar=\ar \zeta \cdot c \cdot \Big(  W'*d\mu(T) -  (\boldsymbol{\mathcal{R}}\circ V_\mu) *W' (T) \Big)
 	=\zeta \cdot c \cdot V_\mu(T). 
 	\eeqnn 
 
  \item[$\bullet$] For $i=2$, an application of Fubini's theorem along with the following two identities
  \beqnn
  \int_0^T  \int_{T-t-y}^{T-t} W'(r)dr \mu(dt) 
  \ar=\ar \int_0^T  \int_0^y W'(T-t-r)dr \mu(dt) \cr
  \ar=\ar \int_0^y \int_0^T W'(T-t-r)\mu(dt) dr \cr
  \ar=\ar \int_0^y \int_0^{T-r} W'(T-r-t)\mu(dt) dr\cr
  \ar=\ar \int_0^y   W'*d\mu(T-r) dr 
  \eeqnn
  and 
  \beqnn
  \int_0^T \boldsymbol{\mathcal{R}}\circ V_\mu (T-t) \int_{t-y}^t W'(r)dr  \, dt 
  \ar=\ar \int_0^T \boldsymbol{\mathcal{R}}\circ V_\mu (T-t) \int_0^y W'(t-r)dr \, dt \cr
  \ar=\ar  \int_0^y \int_0^T \boldsymbol{\mathcal{R}}\circ V_\mu (T-t) W'(t-r) \, dt\, dr\cr
  \ar=\ar \int_0^y   \big(\boldsymbol{\mathcal{R}}\circ V_\mu \big) *W'(T-r)  \, dr
  \eeqnn
  induces that 
  \beqnn
  \mathbf{I}_2 *\big( d\mu(T) - \boldsymbol{\mathcal{R}}\circ V_\mu  \big) (T)
  \ar=\ar \int_0^T  \int_0^\zeta\int_0^\infty \Big(\int_{T-t-y}^{T-t} W'(r)dr \Big)  \, N_0(dz, dy) \mu(dt)\cr
  \ar\ar - \int_0^T \boldsymbol{\mathcal{R}}\circ V_\mu (T-t) \int_0^\zeta\int_0^\infty \Big(\int_{t-y}^t W'(r)dr \Big) \, N_0(dz, dy)\, dt \cr
  \ar=\ar  \int_0^\zeta\int_0^\infty \Big( \int_0^T  \int_{t-y}^t W'(r)dr \mu(dt)  \Big) \, N_0(dz, dy) \cr
  \ar\ar - \int_0^\zeta\int_0^\infty  \Big(\int_0^T \boldsymbol{\mathcal{R}}\circ V_\mu (T-t) \int_{t-y}^t W'(r)dr  \, dt \Big) \, N_0(dz, dy)\cr
  \ar=\ar \int_0^\zeta\int_0^\infty  \Big( \int_0^y   W'*d\mu(T-r) dr \Big) \, N_0(dz, dy)  \cr
  \ar\ar - \int_0^\zeta\int_0^\infty  \Big( \int_0^y   \big(\boldsymbol{\mathcal{R}}\circ V_\mu \big) *W'(T-r)  \, dr \Big) \, N_0(dz, dy) \cr
  \ar=\ar   \int_0^\zeta\int_0^\infty \int_0^y\big(  d\mu  - \boldsymbol{\mathcal{R}}\circ V_\mu   \big) *W'(T-r) \, dr \, N_0(dz, dy) \cr
  \ar=\ar \int_0^\zeta\int_0^\infty \int_0^y  V_\mu (T-r) \, dr \, N_0(dz, dy)\cr
  \ar=\ar \int_0^\zeta\int_0^\infty \int_{T-y}^T  V_\mu (r) \, dr \, N_0(dz, dy).
  \eeqnn
  Here the last two equalities follow from \eqref{NonlinearVolterra} and the change of variables respectively.

  \item[$\bullet$] For $i=3$, by using \eqref{eqn.SFT01} along with \eqref{eqn.201} and the fact that $\boldsymbol{\mathcal{R}}\circ V_\mu \in L^\infty_{\rm loc}(\mathbb{R}_+;\mathbb{R})$,
  \beqnn
  \mathbf{I}_3 *\big( d\mu(T) - \boldsymbol{\mathcal{R}}\circ V_\mu  \big) (T)
  \ar=\ar \int_0^T \int_0^{T-t} \int_0^{X_\zeta(s)}  W'(T-t-s) \, B_c(ds,dz) \mu(dt) \cr
  \ar\ar - \int_0^T  \boldsymbol{\mathcal{R}}\circ V_\mu(T-t) \int_0^{t} \int_0^{X_\zeta(s)}  W'(t-s) \, B_c(ds,dz)dt \cr
  \ar=\ar \int_0^T  \int_0^{X_\zeta(s)} \int_0^{T-s} W'(T-s-r) \mu(dr)\, B_c(ds,dz)\cr
  \ar\ar -\int_0^T   \int_0^{X_\zeta(s)} \int_0^{T-s} \boldsymbol{\mathcal{R}}\circ V_\mu(T-s-r)  W'(r)\, dr \, B_c(ds,dz)\cr
  \ar=\ar \int_0^T \int_0^{X_\zeta(s)}   V_\mu (T-s) \, B_c(ds,dz) . 
  \eeqnn

 \item[$\bullet$] For $i=4$, it is easy to use \eqref{eqn.100} to identify that conditions in Proposition~\ref{StoFubiniThm} are satisfied.  
 By using \eqref{eqn.SFT02} along with 
 \beqnn
  \int_0^{T-s} \int_{T-s-t-y}^{T-s-t} W' (r)\,dr  \mu(dt)
  \ar=\ar  \int_0^{T-s} \int_{0}^{y} W' (T-s-t-r)\,dr  \mu(dt) \cr
  \ar=\ar \int_{0}^{y} \int_0^{T-s}  W' (T-s-r-t) \mu(dt) \,dr \cr
  \ar=\ar \int_{0}^{y} \int_0^{T-s-r}  W' (T-s-r-t) \mu(dt) \,dr\cr
  \ar=\ar \int_{0}^{y}  W' *d\mu(T-s-r) \,dr 
  \eeqnn
  and 
 \beqnn
  \int_0^{T-s}\boldsymbol{\mathcal{R}}\circ V_\mu(T-s-t)   \int_{t-y}^{t} W' (r)\,dr\,  dt 
  \ar=\ar \int_0^{T-s}\boldsymbol{\mathcal{R}}\circ V_\mu(T-s-t)  \int_0^y W' (t-r)\,dr\, dt\cr
  \ar=\ar \int_0^y dr \int_0^{T-s}\boldsymbol{\mathcal{R}}\circ V_\mu(T-s-t)   W' (t-r) \, dt \cr
  \ar=\ar \int_0^y dr \int_0^{T-s-r}\boldsymbol{\mathcal{R}}\circ V_\mu(T-s-r-t)   W' (t)  \,dt \cr
  \ar=\ar  \int_0^y (\boldsymbol{\mathcal{R}}\circ V_\mu)*W'(T-s-r) \, dr
 \eeqnn
 induces that 
 \beqnn
  \lefteqn{\mathbf{I}_4 *\big( d\mu(T) - \boldsymbol{\mathcal{R}}\circ V_\mu  \big) (T)}\ar\ar\cr
  \ar\ar\cr 
  \ar=\ar \int_0^T \int_0^t\int_0^{X_\zeta(s)} \int_0^\infty \Big( \int_{t-s-y}^{t-s} W' (r)\,dr\Big) \, \widetilde{N}_\nu(ds,dz,dy) \mu(dt) \cr
  \ar\ar - \int_0^T  \boldsymbol{\mathcal{R}}\circ V_\mu(T-t) \int_0^t\int_0^{X_\zeta(s)} \int_0^\infty \Big( \int_{t-s-y}^{t-s} W' (r)\,dr\Big) \, \widetilde{N}_\nu(ds,dz,dy)dt \cr
  \ar=\ar  \int_0^T \int_0^t\int_0^{X_\zeta(s)} \int_0^\infty \int_0^{T-s} \int_{T-s-t-y}^{T-s-t} W' (r)\,dr  \mu(dt) \, \widetilde{N}_\nu(ds,dz,dy) \cr
  \ar\ar -  \int_0^T  \int_0^{X_\zeta(s)} \int_0^\infty \int_0^{T-s}\boldsymbol{\mathcal{R}}\circ V_\mu(T-s-t)   \int_{t-y}^{t} W' (r)\,dr \, dt\, \widetilde{N}_\nu(ds,dz,dy)\cr
  \ar=\ar  \int_0^T \int_0^t\int_0^{X_\zeta(s)} \int_0^\infty  \int_{0}^{y}  W' *d\mu(T-s-r) \,dr  \, \widetilde{N}_\nu(ds,dz,dy) \cr
  \ar\ar -  \int_0^T  \int_0^{X_\zeta(s)} \int_0^\infty \int_0^y (\boldsymbol{\mathcal{R}}\circ V_\mu)*W'(T-s-r) \, dr\, \widetilde{N}_\nu(ds,dz,dy)\cr
  \ar=\ar \int_0^T \int_0^t\int_0^{X_\zeta(s)} \int_0^\infty  \int_{0}^{y} V_\mu(T-s-r) \,dr  \, \widetilde{N}_\nu(ds,dz,dy) \cr
  \ar=\ar \int_0^T \int_0^t\int_0^{X_\zeta(s)} \int_0^\infty  \int_{T-s-y}^{T-s} V_\mu(r) \,dr  \, \widetilde{N}_\nu(ds,dz,dy). 
 \eeqnn
 \qed
 \end{enumerate}

 Associated to $\boldsymbol{Y_T}$ we define a Doob's martingale $Y_T:=\big\{Y_T(t) : t \in[0,T] \big\}$ with 
 \beqlb\label{eqn.7201}
 Y_T(t):= \mathbf{E}\big[ \boldsymbol{Y_T} \,\big|\, \mathscr{F}_t\big] .
 \eeqlb
  Conditionally on $\mathscr{F}_t$, we take expectations on both sides of \eqref{eqn.721}  and obtain the following representation for the martingale $Y_T$:
 \beqlb  \label{eqn.9000} 
 Y_T(t) \ar=\ar  \zeta\cdot c \cdot  V_\mu (T)+ \int_0^\zeta \int_0^\infty \Big( \int_{T-y}^{T} V_\mu (r)\,dr\Big) N_0(dz,dy)  + \int_0^{t} \int_0^{X_\zeta(s)}   V_\mu(T-s) \, B_c(ds,dz)\cr
 \ar\ar  + \int_0^{t} \int_0^{X_\zeta(s)} \int_0^\infty \Big( \int_{T-s-y}^{T-s} V_\mu(r)\,dr\Big) \, \widetilde{N}_\nu(ds,dz,dy) ,\quad t\in[0,T]. 
 \eeqlb
 Additionally, we consider a  $(\mathscr{F}_t)$-process $Z_T:=\{Z_T(t) :t\in [0,T] \}$ defined by 
 \beqlb\label{eqn.Z}
 Z_T(t):= \mathbf{E}\big[  X_\zeta*d\mu(T) \,\big|\, \mathscr{F}_t \big] - \int_t^T \big(\boldsymbol{\mathcal{R}}\circ V_\mu\big) (T-s) \mathbf{E}\big[  X_\zeta(s) \,\big|\, \mathscr{F}_t \big] ds  . 
 \eeqlb
 
 \begin{lemma}
 	The process $Z_T $ is a c\`adl\`ag $(\mathscr{F}_t)$-semimartingale with the following decomposition
 	\beqlb\label{eqn.9001}
 	Z_T(t) \ar=\ar	Y_T(t) + \int_0^t \boldsymbol{\mathcal{R}} \circ V_\mu (T-s) X_\zeta(s)ds,
 	\quad t\in[0,T] . 
 	\eeqlb  
 \end{lemma} 
 \proof   
 Plugging \eqref{eqn.720} into \eqref{eqn.7201} induces that 
 \beqnn
 Y_T(t) 
 \ar=\ar  \mathbf{E}\big[    X_\zeta *d\mu(T) \big|  \mathscr{F}_t \big]  -\int_0^t \boldsymbol{\mathcal{R}}_0\circ V_\mu (T-s) X_\zeta(s)\, ds  - \int_t^T \boldsymbol{\mathcal{R}}_0\circ V_\mu (T-s) \mathbf{E}\big[ X_\zeta(s) \big| \mathscr{F}_t \big]\, ds  . 
 \eeqnn
 The desired representation \eqref{eqn.9001} follows immediately by moving the second term on the right-hand side to the left-hand side.  
 \qed  
 
 Note  that $Z_T(0)=Y_T(0)$. 
 Applying It\^o's formula (see Theorem 5.1 in \cite[p.66]{IkedaWatanabe1989}) along with \eqref{eqn.9000}-\eqref{eqn.9001} to $e^{-Z_T(t)}$ and then using \eqref{NonlinearVolterra}, we have  
 \beqlb \label{eqn.705}
 e^{-Z_T(t)}  
 \ar=\ar e^{-Y_T(0)} -\int_0^t e^{-Z_T(s)} \cdot \boldsymbol{\mathcal{R}} \circ V_\mu (T-s) X_\zeta(s)ds   + \int_0^t e^{-Z_T(s)}\cdot c\cdot\big|V_\mu (T-s)\big|^2 \cdot X_\zeta(s)ds \cr
 \ar\ar + \int_0^t e^{-Z_T(s)} \cdot \mathcal{V}\circ V_\mu (T-s) \cdot X_\zeta(s)ds  - \int_0^t \int_0^{X_\zeta(s)}  e^{-Z_T(s)} \cdot V_\mu(T-s) \, B_c(ds,dz) \cr
 \ar\ar  \cr 
 \ar\ar + \int_0^{t} \int_0^{X_\zeta(s)} \int_0^\infty e^{-Z_T(s-)}\cdot \Big( \exp\Big\{ -\int_{T-s-y}^{T-s} V_\mu(r)\,dr\Big\} -1 \Big) \, \widetilde{N}_\nu(ds,dz,dy) \cr 
 \ar=\ar e^{-Y_T(0)} - \int_0^t \int_0^{X_\zeta(s)}  e^{-Z_T(s)} \cdot V_\mu(T-s) \, B_c(ds,dz) \cr
 \ar\ar  \cr 
 \ar\ar + \int_0^{t} \int_0^{X_\zeta(s)} \int_0^\infty e^{-Z_T(s-)}\cdot \Big( \exp\Big\{ -\int_{T-s-y}^{T-s} V_\mu(r)\,dr\Big\} -1 \Big) \, \widetilde{N}_\nu(ds,dz,dy) . 
 \eeqlb
 Consider a process $U_T:=\big\{U_T(t):t\in[0,T] \big\}$ defined by 
 \beqnn
 U_T(t)\ar:=\ar  \int_0^{t} \int_0^{X_\zeta(s)} \int_0^\infty  \Big( \exp\Big\{-\int_{T-s-y}^{T-s} V_\mu (r)dr\Big\} -1 \Big) \, \widetilde{N}_\nu(ds,dz,dy) \cr
 \ar\ar -  \int_0^{t } \int_0^{ X_\zeta(s)}   V_\mu(T-s) \, B_c(ds,dz) ,
 \eeqnn
 which is a $(\mathscr{F}_t)$-martingale and allows us to rewrite \eqref{eqn.705} into 
 \beqlb\label{eqn.7051}
 e^{-Z_T(t)} \ar=\ar e^{-Y_T(0)} + \int_0^{t }   e^{-Z_T(s-)} dU_T(s), 
 \quad t\in[0,T],
 \eeqlb
 Multiplying both sides by $e^{Y_T(0)}$ shows that the process 
 $
  \mathcal{E}_{U_T}:= \big\{e^{Y_T(0)-Z_T(t)}:t \in [0,T] \big\}
 $
 is the Dol\'ean-Dade exponential associated to $U_T$.

 \begin{lemma}\label{Lemma.StochExp}
 	The Dol\'eans-Dade exponential $\mathcal{E}_{U_T}$ is a true $(\mathscr{F}_t)$-martingale.  
 \end{lemma}
 \proof  
 By It\^o's formula; see Theorem 5.1 in \cite[p.66]{IkedaWatanabe1989}, the stochastic exponential $\mathcal{E}_{U_T}$ admits the following representation
 \beqnn
 \mathcal{E}_{U_T}(t)
 \ar=\ar\exp\bigg\{ - \int_0^t \boldsymbol{\mathcal{R}}\circ V_\mu (T-s) X_\zeta(s) ds -\int_0^{t } \int_0^{ X_\zeta(s)}   V_\mu (T-s) \, B_c(ds,dz) \cr
 \ar\ar\ \  \quad\quad - \int_0^t \int_0^{X_\zeta(s)}\int_0^\infty  \Big(\int_{T-s-y}^{T-s} V_\mu (r)dr \Big)  \widetilde{N}_\nu (ds,dz,dy) \bigg\},\quad t\geq 0.
 \eeqnn 
 The non-negativity of $V_\mu$ induces that all jumps of $U_T$ are larger than $-1$, which induces that $ \mathcal{E}_{U_T}$ is a non-negative local martingale and hence a supermartingale.
 By Fatou's lemma, we have $\mathbf{E}[\mathcal{E}_{U_T}(t)] \leq 1$.
 Thus $\mathcal{E}_{U_T}$ is a true $(\mathscr{F}_t)$-martingale if we can identify that
 \beqlb\label{eqn.703}
 \mathbf{E}\big[\mathcal{E}_{U_T}(t)\big]=1,\quad t \in [0,T].
 \eeqlb

 For each $t_0 \in [0,T]$ and $n\geq 1$, let 
 \beqnn
 \tau_n:=\inf\bigg\{ r\geq 0: \int_0^r X_\zeta (s)ds \geq n \bigg\} \wedge t_0
 \quad\mbox{and}\quad 
 \mathcal{E}_{U_T}^n(t):=  \mathcal{E}_{U_T}(\tau_n \wedge t),\quad t\geq 0. 
 \eeqnn 
 Here we make the convention that $\inf \emptyset= \infty$.  
 It is obvious that as $n\to\infty$,
 \beqlb\label{eqn.50003}
 \tau_n \overset{\rm a.s.}\to t_0
 \quad \mbox{and}\quad 
 \mathcal{E}_{U_T}^n(\cdot)\overset{\rm a.s.}\to \mathcal{E}_{U_T}(t_0\wedge \cdot). 
 \eeqlb 
 The fact that $V_\mu\in L^\infty_{\rm loc}(\mathbb{R}_+;\mathbb{R}_+)$ induces that there exists a constant $C>0$ such that for any $n\geq 1$, 
 \beqnn
  \sup_{t\geq 0}\int_0^{t\wedge \tau_n} X_\zeta(s) \big|V_\mu(T-s)\big|^2\, ds \leq \int_0^{\tau_n} X_\zeta(s)\,ds \cdot   \big\|V_\mu\big\|^2_{L^\infty_T} < C\cdot n
 \eeqnn
 Also, by the inequality $|1-(1+z)e^{-z}| \leq |z|^2 $ for any $z\geq 0$, 
 \beqnn
 \lefteqn{  \sup_{t\geq 0}\int_0^{t\wedge \tau_n} X_\zeta (s)\, ds \int_0^\infty \bigg| 1-\Big( 1+\int_{T-s-y}^{T-s}V_\mu(r)\, dr \Big)     \exp\Big\{- \int_{T-s-y}^{T-s}V_\mu(r) \, dr\Big\}  \bigg| \,\nu(dy) }\ar\ar\cr
 \ar\leq\ar \int_0^{ \tau_n} X_\zeta (s)\, ds \int_0^\infty \bigg|  \int_{T-s-y}^{T-s}V_\mu(r) \, dr  \bigg|^2 \,\nu(dy)  
 \leq \int_0^{ \tau_n} X_\zeta (s)\, ds \cdot \big\|V_\mu\big\|^2_{L^\infty_T}\cdot\int_0^\infty \big( T\wedge y\big)^2 \,\nu(dy) , 
 \eeqnn 
 which can be uniformly bounded by $C\cdot n$.
 By using Theorem~IV.3 in \cite{LepingleMemin1978} with 
 \beqlb \label{eqn.704}
 M^c_t = -  \int_0^{t\wedge \tau_n } \int_0^{ X_\zeta(s)}   V_\mu(T-s) \, B_c(ds,dz)
 \eeqlb
 and 
 \beqnn
 y(s,z)=  \exp\Big\{- \int_{T-s-z}^{T-s}V_\mu(r) dr \Big\}-1 , \quad  \nu^M_t(dz)= \int_0^{t } \mathbf{1}_{\{ s\leq \tau_n \}} \cdot X_\zeta(s)\, ds\, \nu(dz),
 \eeqnn 
 the process $\mathcal{E}_{U_T}^n$ is a uniformly integrable martingale for each $n\geq 1$. Thus
 \beqnn
 1= \mathbf{E}\big[\mathcal{E}_{U_T}^n(t_0)\big]
 \ar=\ar \mathbf{E}\big[\mathcal{E}_{U_T}^n(t_0); \tau_n = t_0\big] + \mathbf{E}\big[\mathcal{E}_{U_T}^n(t_0);\tau_n<t_0\big] \cr
 \ar\ar\cr
 \ar=\ar \mathbf{E}\big[\mathcal{E}_{U_T}(t_0);\tau_n= t_0\big] + \mathbf{E}\big[\mathcal{E}_{U_T}^n(t_0);\tau_n<t_0\big] .
 \eeqnn
 By the monotone convergence theorem and the first limit in \eqref{eqn.50003}, 
 \beqnn
  \lim_{n\to\infty} \mathbf{E}\big[\mathcal{E}_{U_T}(t_0); \tau_n=t_0\big] = \mathbf{E}\big[\mathcal{E}_{U_T}(t_0)\big] .
 \eeqnn
 Therefore, to obtain \eqref{eqn.703} it suffices to prove that 
 \beqlb\label{eqn.706}
  \lim_{n\to\infty} \mathbf{E}\big[\mathcal{E}_{U_T}^n(t_0);\tau_n<t_0\big] =0.
 \eeqlb 

 Associate with the martingale $\mathcal{E}_{U_T}^n$, we define a probability law  $\mathbf{Q}^n$ on $(\Omega,\mathscr{F},\mathscr{F}_t)$ by 
 \beqnn
 \frac{d  \mathbf{Q}^n}{d \mathbf{P}}= \mathcal{E}_{U_T}^n(\tau_n). 
 \eeqnn 
  We now consider the random elements $(X_\zeta, N_0,B_c,\widetilde{N}_\nu)$ under this new probability law.
  \begin{enumerate}
  	\item[$\bullet$] Note that $\mathcal{E}_{U_T}^n(0)\overset{\rm a.s.}=1$, the Poisson random measure $N_0(dz,dy)$ is $\mathscr{F}_0$-measurable and has the same law under $\mathbf{P}$ and $\mathbf{Q}^n$.
  	
  	\item[$\bullet$] By the classical Girsanov's Theorem; see Theorem~3.11 in \cite[p.168]{JacodShiryaev2003}, under $\mathbf{Q}^n$ the continuous martingale $M^c$ defined by \eqref{eqn.704} has predictable quadratic variation 
  	\beqnn
  	\langle M^c\rangle_t= 2c\int_0^t \mathbf{1}_{\{ s\leq \tau_n \}} \cdot X_\zeta(s) \cdot  \big|V_\mu(T-s)\big|^2 ds,\quad t\geq 0,
  	\eeqnn
  	and the Gaussian white noise $B_c(ds,dz)$ has intensity $2c\cdot\mathbf{1}_{\{ s\leq \tau_n \}} \cdot ds\,dz$. 
  	
  	\item[$\bullet$] By Girsanov's theorem for random measure; see Theorem~3.17 in \cite[p.170]{JacodShiryaev2003}, the Poisson random measure $N_\nu(ds,dy,dz)$ is a random point measure under $\mathbf{Q}^n$ with intensity
  	\beqnn
  	\mathbf{1}_{\{ s\leq \tau_n \}}\cdot \exp\Big\{- \int_{T-s-y }^{T-s} V_\mu(r) dr\Big\} \, ds \, dz \, \nu (dy) ,
  	\eeqnn
  	
  	\item[$\bullet$] For each $t_1>0$, we consider the auxiliary process
  	\beqnn
  	X_{\zeta,t_1}(t)
  	\ar:=\ar \zeta\cdot c \cdot W'(t)  
  	+  \int_0^\zeta\int_0^\infty \big( W(t)-W(t-y) \big) \, N_0(dz, dy)   \cr 
  	\ar\ar + \int_0^t \int_0^{X_\zeta(s)}   W'(t_1-s) \, B_c(ds,dz) \cr
  	\ar\ar + \int_0^t \int_0^{X_\zeta(s)} \int_0^\infty  \big( W(t_1-s)-W(t_1-s-y) \big) \, \widetilde{N}_\nu(ds,dz,dy), \quad t\geq 0.
  	\eeqnn
 It is obvious that $X_{\zeta,t_1}$ is a $(\mathscr{F}_t)$-semimartingale with $X_{\zeta,t_1}(t_1)\overset{\rm a.s.}= X_{\zeta}(t_1)$ under $\mathbf{P}$ and $\mathbf{Q}^n$. 
 By Girsanov's theorem for semimartingales; see Theorem~3.24 in \cite[p.172]{JacodShiryaev2003}, the process $X_{\zeta,t_1}$ is also a semimartingale under $\mathbf{Q}^n$ with the following representation 
 \beqnn
  X_{\zeta,t_1}(t)
  \ar=\ar \zeta\cdot \big(1-b W(t)\big)
  +  \int_0^\zeta\int_0^\infty \big( W(t)-W(t-y) \big) \, \widetilde{N}_0(dz, dy)   \cr 
  \ar\ar +\int_0^{t\wedge \tau_n}    X_{\zeta}(s)\cdot \bigg[ \int_0^\infty \big( W(t_1-s)-W(t_1-s-y) \big) \cr
  \ar\ar \qquad \times \Big(  \exp\Big\{- \int_{T-s-y}^{T-s}V_\mu(r)\, dr\Big\} -1 \Big) \nu (dy) -2c \cdot W'(t_1-s)\bigg] ds \cr
  \ar\ar + \int_0^t \int_0^{X_\zeta(s)}   W'(t_1-s) \, B_c(ds,dz) \cr
  \ar\ar + \int_0^t \int_0^{X_\zeta(s)} \int_0^\infty  \big( W(t_1-s)-W(t_1-s-y) \big) \, \widetilde{N}_\nu(ds,dz,dy), \quad t\geq 0.
  \eeqnn
 By setting $t=t_1$ and then using the arbitrariness of $t_1$, the stochastic Volterra equation \eqref{MainThm.SVE} under $\mathbf{Q}^n$ turns to be 
 \beqnn
  X_{\zeta}(t)
 \ar=\ar \zeta\cdot \big(1-b W(t)\big)
 +  \int_0^\zeta\int_0^\infty \big( W(t)-W(t-y) \big) \, \widetilde{N}_0(dz, dy)    \cr 
 \ar\ar +\int_0^{t\wedge \tau_n}    X_{\zeta}(s) \cdot \bigg[ \int_0^\infty \big( W(t-s)-W(t-s-y) \big) \cr
 \ar\ar \qquad \times \Big(  \exp\Big\{- \int_{T-s-y}^{T-s}V_\mu(r)\, dr\Big\} -1 \Big) \nu (dy) -2c \cdot W'(t-s)\bigg] \, ds \cr
 \ar\ar + \int_0^t \int_0^{X_\zeta(s)}   W'(t-s) \, B_c(ds,dz) \cr
 \ar\ar + \int_0^t \int_0^{X_\zeta(s)} \int_0^\infty  \big( W(t-s)-W(t-s-y) \big) \, \widetilde{N}_\nu(ds,dz,dy), \quad t\geq 0.
 \eeqnn
  \end{enumerate}
  We write $\mathbf{E}^{\mathbf{Q}^n}$ for the expectation under the law $\mathbf{Q}^n$. 
 Taking expectations on both sides of this equation and then using Fubini's theorem,  
 \beqnn
 \mathbf{E}^{\mathbf{Q}^n}\big[ X_{\zeta}(t) \big]
 \ar=\ar \zeta\cdot \big(1-b W(t)\big) +\int_0^{t\wedge \tau_n}   \mathbf{E}^{\mathbf{Q}^n}\big[ X_{\zeta}(s)\big] \cdot \bigg[ \int_0^\infty \big( W(t-s)-W(t-s-y) \big) \cr
 \ar\ar \qquad \times \Big(  \exp\Big\{- \int_{T-s-y}^{T-s}V_\mu(r)\, dr\Big\} -1 \Big) \nu (dy) -2c \cdot W'(t-s)\bigg]\, ds.
 \eeqnn
 Since $W$ is non-decreasing and $V_\mu$, $W'$ are non-negative, we have 
 \beqnn
  \mathbf{E}^{\mathbf{Q}^n}\big[ X_{\zeta}(t) \big]
  \ar\leq\ar \zeta\cdot \big(1-b W(t)\big),\quad t\geq 0.
 \eeqnn 
 By the definition of $\tau_n$, Chebyshev's inequality and Fubini's theorem,
 \beqnn
 \mathbf{E}\big[\mathcal{E}_{U_T}^n(\tau_n); \tau_n<t_0 \big]  =\mathbf{Q}^n \big(\tau_n<t_0 \big) 
 \ar=\ar \mathbf{Q}^n\bigg(\int_0^{t_0}X_{\zeta} (s)ds \geq n\bigg)\cr 
 \ar\leq\ar \frac{1}{n}\mathbf{E}^{\mathbf{Q}^n}\bigg[\int_0^{t_0}X_{\zeta} (s)ds\bigg]\cr
 \ar=\ar \frac{1}{n}\int_0^{t_0}\mathbf{E}^{\mathbf{Q}^n}\big[X_{\zeta} (s)\big]\, ds \cr
 \ar\leq\ar \frac{\zeta}{n}\int_0^{t_0} \big(1-b W(s)\big)\, ds,
 \eeqnn
 which vanishes as $n\to\infty$.
 Hence both \eqref{eqn.706} and \eqref{eqn.703} hold.
 In conclusion, the local  martingale $\mathcal{E}_{U_T}$ is a true $(\mathscr{F}_t)$-martingale under $\mathbf{P}$.
 \qed


 \textit{\textbf{Proof of Theorem~\ref{Thm.LaplaceF}: Part II.}} Note that $e^{-Z_T(t)}= e^{-Y_T(0)}\cdot \mathcal{E}_{U_T}(t)$ for $t\in[0,T]$. 
 By Lemma~\ref{Lemma.StochExp}, the process $e^{-Z_T}$ is a true $(\mathscr{F}_t)$-martingale. 
 Moreover, by \eqref{eqn.Z} with $t=T$ we have
 \beqnn
 Z_T(T)=  X_\zeta *d\mu(T) .
 \eeqnn
 Taking expectations on both sides of \eqref{eqn.7051} yields that
 \beqnn
 \mathbf{E}\Big[ \exp\big\{-X_\zeta *d\mu(T)\big\} \Big]
 = \mathbf{E}\big[ e^{-Z_T(T)} \big]
 = \mathbf{E}\big[ e^{-Y_T(0)} \big].
 \eeqnn
 By \eqref{eqn.9000} with $t=0$ and then using the exponential formula of stochastic integral with respect to Poisson random measure; see \cite[p.8]{Bertoin1996},  
 \beqnn
 \mathbf{E}\big[ e^{-Y_T(0)} \big] 
 \ar=\ar \exp\Big\{  -\zeta \cdot c\cdot V_\mu(T)  \Big\}\cdot  \mathbf{E}\bigg[\exp\Big\{ - \int_0^\zeta\int_0^\infty \Big( \int_{T-y}^T V_\mu (r)dr \Big) \, N_0(dz, dy)  \Big\} \bigg]\cr
 \ar=\ar \exp\bigg\{  -\zeta \cdot c\cdot V_\mu(T) -\int_0^\infty \Big(1- \exp\Big\{-\int_{T-y}^{T} V_\mu(r)dr\Big\} \Big) \bar\nu(y)dy  \bigg\},
 \eeqnn
 and hence the representation \eqref{LaplaceFun01}-\eqref{LaplaceFun02} hold. 
 \qed

 \textit{\textbf{Proof of uniqueness of solution to \eqref{MainThm.SVE}.}} 
 Assume that $X_\zeta^{(1)}$ and $X_\zeta^{(2)}$ are two solutions of \eqref{MainThm.SVE}.  
 For any $T>0$, consider $X_\zeta^{(1)}$ and $X_\zeta^{(2)}$ as two $L^1([0,T];\mathbb{R}_+)$-valued random variables.
 By \eqref{LaplaceFun01},
 \beqnn 
 \mathbf{E}\bigg[ \exp\Big\{- \int_0^T X_\zeta^{(1)}(T-s)f(s)ds \Big\} \bigg]
 = \mathbf{E}\bigg[ \exp\Big\{- \int_0^T X_\zeta^{(2)}(T-s)f(s)ds \Big\} \bigg],
 \eeqnn
 for any $f \in L^\infty([0,T]:\mathbb{R}_+)$.
 Note that $L^\infty([0,T];\mathbb{R}_+)$ is the dual space of $L^1([0,T];\mathbb{R}_+)$, the preceding identity yields that $X_\zeta^{(1)}$ and $X_\zeta^{(2)}$ have the same law on $L^1([0,T];\mathbb{R}_+)$ and hence on $C([0,T];\mathbb{R}_+)$. In conclusion, the  uniqueness of solution holds for \eqref{MainThm.SVE}.
 \qed

 \appendix
 \renewcommand{\theequation}{A.\arabic{equation}}
  \section{Stochastic integrals driven by $\mathbb{H}^\#$-semimartingales} 
 \label{AppendixHMartinagle}
 \setcounter{equation}{0}

 We recall some basic theory of stochastic integrals with respect to infinite-dimensional semimartingales that were firstly studied by Kurtz and Protter \cite{KurtzProtter1996}.
 Let $\mathbb{H}$ be a separable Banach space endowed with a norm $\|\cdot\|_{\mathbb{H}}$.
 We first recall the definition of  $\mathbb{H}^\#$-semimartingales.

 \begin{definition}
 	We say $Y$ is an {\rm $\mathbb{H}^{\#}$-semimartingale} if it is $(\mathscr{F}_t)$-adapted and indexed by $\mathbb{H}\times [0,\infty)$ such that the following hold
 	\begin{enumerate} 
 \item[$\bullet$] For each $h\in \mathbb{H}$, the process $Y(h):=\{ Y(h,t):t\geq 0 \}$ is a c\`adl\`ag $\mathbb{R}$-valued $(\mathscr{F}_t)$-semimartingale starting from $0$;
 			
 \item[$\bullet$] For each $m\in\mathbb{Z}_+$, $h_1,\dots,h_m\in \mathbb{H}$, and $a_1,\dots,a_m\in\mathbb{R}$, the following finite additivity holds
 		\beqnn
 		Y\bigg(\sum_{k=1}^m a_k h_k,t\bigg)\overset{\rm a.s.}= \sum_{k=1}^m a_k Y\big(h_k,t \big),
 		\quad t\geq 0.
 		\eeqnn
 		
 	\end{enumerate}
 \end{definition}

 Let $\mathbb{H}_0$ be a dense subset of $\mathbb{H}$ and $\mathcal{S}_0$ the collection of $\mathbb{H}$-valued stochastic processes of the form
 \beqnn
 X(t):= \sum_{k=1}^m \xi_k(t)\varphi_k \quad \mbox{with}\quad \xi_k(t):=\sum_{i=0}^{\infty} \eta_i^k\cdot\mathbf{1}_{[\tau_i^k,\tau_{i+1}^k)}(t),
 \eeqnn
 where $m\geq 1$, $\varphi_1,\cdots,\varphi_m\in\mathbb{H}_0$, $\{\tau_i^k\}_{i\geq 0}$ is a sequence of non-decreasing $(\mathscr{F}_t)$-stopping times and $\eta_i^k \in\mathbb{R}^d$ is  $\mathscr{F}_{\tau_i^k}$-measurable.
 For any $X\in\mathcal{S}_0$,  define
 \beqnn
 X_-\cdot dY(t) =\sum_{k=1}^m \int_0^t \xi_k(s-)d Y(\varphi_k,t), \quad t\geq 0.
 \eeqnn
 \begin{definition}
 	The $\mathbb{H}^\#$-semimartingale $Y$ is {\rm standard} if
 	\beqlb\label{eqn.Ht}
 	\mathcal{H}_t:=  \Big\{ \sup_{s\leq t}|X_-\cdot dY(s)| : X \in\mathcal{S}_0,\, \sup_{s\leq t}\|X(s)\|_{\mathbb{H}}\leq 1  \Big\}
 	\eeqlb
 	is stochastically bounded for each $t\geq 0$.
 \end{definition}
 For any $\mathbb{H}$-valued c\`adl\`ag process $X$ and standard $\mathbb{H}^\#$-semimartingale $Y$,
 we can find a sequence $\{X^\epsilon\}_{\epsilon>0}\subset \mathcal{S}_0$ such that as $\epsilon \to0$,
 \beqnn
 \sup_{t\in[0,T]}\|X^\epsilon(t)-X(t)\|_\mathbb{H}\overset{\rm  a.s. }\to0
 \quad \mbox{and} \quad
 X_-\cdot dY := \lim_{\epsilon\to 0+} X^\epsilon_-\cdot dY
 \eeqnn
 exists a.s. in the sense that $\sup_{t\in[0,T]} |X_-\cdot dY(t)- X^\epsilon_-\cdot dY(t) | \overset{\rm P}\to 0$.
 Moreover, the limit process $X_-\cdot Y$ is c\`adl\`ag,  independent of $\{X^\epsilon\}_{\epsilon>0}$ and called the \textit{stochastic integral} of $X$ with respect to $Y$. 
 For any $(\mathscr{F}_t)$-stopping time $\sigma$, let  $X_-^\sigma(t):= X_-(t) \mathbf{1}_{[0,\sigma)}(t)$ for $t\geq 0$. We have the following identity
 \beqnn
 X_-\cdot dY(t\wedge \sigma)= X_-^\sigma\cdot dY(t).
 \eeqnn

 \begin{definition}\label{Definition.A1}
 Consider a sequence of $\mathbb{H}^\#$-semimartingales $\{Y_n\}_{n\geq 1}$.
 \begin{enumerate}
 	\item[(1)] It is {\rm uniformly tight} if for each  $t\geq 0$, the family $\{\mathcal{H}_{n,t}\}_{n\geq 1}$ is uniformly stochastically bounded, where $\mathcal{H}_{n,t}$ is defined as in (\ref{eqn.Ht}) with $Y$ replaced by $Y_n$.
 	
 	\item[(2)] We say $Y_n$ {\rm converges weakly} to $Y$ and write $Y_n \Rightarrow Y$ if for any $m\geq 1$ and $f_1,\cdots,f_m\in\mathbb{H}$,
 	\beqnn
 	(Y_n(f_1),\cdots, Y_n(f_m))\overset{\rm d}\to(Y(f_1),\cdots, Y(f_m))\quad \mbox{in }D([0,\infty),\mathbb{R}^m).
 	\eeqnn
 	In addition,
 	we also write $(X_n,Y_n) \Rightarrow (X,Y)$ if
 	\beqnn
 	(X_n, Y_n(f_1),\cdots,Y_n(f_m))\overset{\rm d}\to(X,Y(f_1),\cdots, Y(f_m))\quad \mbox{in }D([0,\infty),\mathbb{H}\times\mathbb{R}^m).
 	\eeqnn
 \end{enumerate} 
 	
 \end{definition}

 \begin{lemma}[Theorem 5.5 in \cite{KurtzProtter1996}] \label{kurztheorem}
  Let $\{Y_n\}_{n\geq 1}$ be a sequence of standard $\mathbb{H}^\#$-semimartingales and $\{X_n\}_{n\geq 1}$ a sequence of c\`adlà\`ag, $\mathbb{H}$-valued processes. 
  If $\{Y_n\}_{n\geq 1}$ is uniformly tight and 
 	$(X_n,Y_n)\Rightarrow (X,Y),$
  then there exists a filtration $\{\mathscr{G}_t\}$ such that $Y$ is an $\{\mathscr{G}_t\}$-adapted, standard, $\mathbb{H}^\#$-semimartingale, $X$ is $\{\mathscr{G}_t\}$-adapted, and 
  \beqnn
 	(X_n,Y_n,X_{n-}\cdot dY_n)\Rightarrow (X,Y,X_-\cdot dY). 
  \eeqnn	 
 \end{lemma}

 \bibliographystyle{plain}
 \bibliography{Reference}

 \end{document}